%% file: main.tex
\documentclass[11pt,twoside,reqno]{amsart}
\usepackage{import}
\usepackage{newclude}
\include{preamble}

\title{Eigenvalues of $p$-adic random matrices}
\author{Jiahe Shen and Roger Van Peski}

\date{\today}
\begin{document}

\maketitle

\begin{abstract} 
    We develop the basic theory of eigenvalues of $p$-adic random matrices, analogous to the classical theory for random matrices over $\R$ and $\C$. Such eigenvalue statistics were proposed as a model for the zeroes of $p$-adic $L$-functions by Ellenberg-Jain-Venkatesh, who computed the limiting distribution of the number of eigenvalues in a unit disc. We compute the full joint distribution of the $n$ eigenvalues of an $n \times n$ matrix with Haar distribution, obtaining Coulomb gas type formulas as in the archimedean case, with Vandermonde terms leading to eigenvalue repulsion. 

    From these Coulomb gas density functions we derive asymptotics of eigenvalue statistics as $n \to \infty$. These include exact computations, such as a closed form 
    $$\rho(x,y) = 1 - \theta_3(-\sqrt{p};||x-y||^2/p)$$
    for the limiting pair correlation of eigenvalues in $\Z_p$, and similar results in quadratic extensions. Such formulas yield concrete numerical predictions on zeroes of $p$-adic $L$-functions.
    
    For eigenvalues in arbitrary extensions of $\Q_p$ we also give precise estimates on their pair-repulsion and expected number of eigenvalues in each extension. Finally, we compute the asymptotic probability that all eigenvalues lie in $\Z_p$.
    
    Our proofs combine results from several distinct areas: $p$-adic orbital integrals, roots of random $p$-adic polynomials, the Sawin-Wood moment method for random modules, and Markov chains associated with measures on integer partitions.
\end{abstract}

\tableofcontents

\include*{Introduction}

\include*{Preliminaries}

\include*{joint_distribution}

\include*{joint_distribution_2}

\include*{Expected_determinants}

\include*{cor_function_asymptotic}

\include*{GL_r}

\include*{markov}

\include*{limiting_correlation_functions}

\include*{repulsion}

\include*{average_number_high_degree}

\include*{all_Zp}

\bibliographystyle{alpha_abbrvsort}
\bibliography{references.bib}

\end{document}

%% file: preamble.tex
\usepackage{anysize}
\usepackage{amsmath}
\usepackage{amsthm}
\usepackage{amsmath,amscd}
\usepackage[utf8]{inputenc}
\usepackage{amssymb}
\usepackage{scalerel}
\usepackage{stmaryrd}
\usepackage{wasysym}
\usepackage{mathrsfs}
\usepackage[pagebackref=true]{hyperref} 
\renewcommand*{\backref}[1]{}
\renewcommand*{\backrefalt}[4]{({\tiny%
   \ifcase #1 Not cited.%
         \or Cited on page~#2.%
         \else Cited on pages #2.%
   \fi%
   })}
\usepackage{cleveref}
\usepackage{enumitem}
\usepackage{dsfont}
\usepackage{soul}
\hypersetup{colorlinks=true,linkcolor=black,citecolor=black}
\usepackage{color}
\usepackage[all]{xy}
\usepackage{float}
\usepackage{bm}
\usepackage{mathtools}
\usepackage{makecell}
\usepackage{thmtools}
\usepackage{setspace}
\usepackage{comment}
\usepackage{tikz}
\usepackage{tikz-cd}
\usepackage{mathdots}
\usepackage{graphicx}
\usepackage[myheadings]{fullpage}
\usepackage{csquotes}

\DeclareFontFamily{U}{mathb}{}

\DeclareFontShape{U}{mathb}{m}{n}{
   <-5.5>  mathb5
   <5.5-6.5> mathb6
   <6.5-7.5> mathb7
   <7.5-8.5> mathb8
   <8.5-9.5> mathb9
   <9.5-11>  mathb10
   <11->     mathb12
}{}

\DeclareSymbolFont{mathb}{U}{mathb}{m}{n}

\DeclareMathSymbol{\lefttorightarrow}{3}{mathb}{"FC}
\DeclareMathSymbol{\righttoleftarrow}{3}{mathb}{"FD}

\makeatletter
\newcommand{\acts}{%
  \mathrel{\mathpalette\@acts\righttoleftarrow}}
\newcommand{\actedby}{%
  \mathrel{\mathpalette\@acts\lefttorightarrow}}
\newcommand{\@acts}[2]{\reflectbox{$\m@th#1#2$}}
\makeatother

\numberwithin{equation}{section}
\newcommand\mtop{.95in}
\newcommand\mbottom{.95in}
\newcommand\mleft{1in}
\newcommand\mright{1in}
\usepackage[top = \mtop, bottom = \mbottom, left = \mleft, right=\mright]{geometry}
\DeclareMathOperator{\val}{val}
\DeclareMathOperator{\Mat}{Mat}

\newtheorem{thm}{Theorem}[section]

\newtheorem{example}[thm]{Example}
\newtheorem{prop}[thm]{Proposition}

\newtheorem{lemma}[thm]{Lemma}

\newtheorem{cor}[thm]{Corollary}
\theoremstyle{definition}
\newtheorem{defi}[thm]{Definition}
\newtheorem{notation}[thm]{Notation}
\newtheorem{rmk}[thm]{Remark}

\newcommand\reallywidehat[1]{%
\savestack{\tmpbox}{\stretchto{%
  \scaleto{%
    \scalerel*[\widthof{\ensuremath{#1}}]{\kern-.6pt\bigwedge\kern-.6pt}%
    {\rule[-\textheight/2]{1ex}{\textheight}}
  }{\textheight}%
}{0.5ex}}%
\stackon[1pt]{#1}{\tmpbox}%
}
\DeclareSymbolFont{bbold}{U}{bbold}{m}{n}
\DeclareSymbolFontAlphabet{\mathbbold}{bbold}

\DeclarePairedDelimiter{\floor}{\lfloor}{\rfloor}

\makeatletter
\def\@tocline#1#2#3#4#5#6#7{\relax
  \ifnum #1>\c@tocdepth 
  \else
    \par \addpenalty\@secpenalty\addvspace{#2}%
    \begingroup \hyphenpenalty\@M
    \@ifempty{#4}{%
      \@tempdima\csname r@tocindent\number#1\endcsname\relax
    }{%
      \@tempdima#4\relax
    }%
    \parindent\z@ \leftskip#3\relax \advance\leftskip\@tempdima\relax
    \rightskip\@pnumwidth plus4em \parfillskip-\@pnumwidth
    #5\leavevmode\hskip-\@tempdima
      \ifcase #1
       \or\or \hskip 1em \or \hskip 2em \else \hskip 3em \fi%
      #6\nobreak\relax
    \hfill\hbox to\@pnumwidth{\@tocpagenum{#7}}\par
    \nobreak
    \endgroup
  \fi}
\makeatother

\renewcommand{\vec}{\mathbf}

\newcommand{\R}{\mathbb{R}}
\newcommand{\Z}{\mathbb{Z}}
\newcommand{\Q}{\mathbb{Q}}
\newcommand{\N}{\mathbb{N}}

\newcommand{\C}{\mathbb{C}}
\newcommand{\F}{\mathbb{F}}

\renewcommand{\H}{\mathbb{H}}

\newcommand{\E}{\mathbb{E}}

\newcommand{\mf}{\mathfrak}

\newcommand{\mc}{\mathcal}

\newcommand{\bbone}{\mathbbold{1}}
\renewcommand{\l}{\lambda}
\newcommand{\la}{\lambda}

\newcommand{\pfrac}[2]{\left(\frac{#1}{#2}\right)}

\newcommand{\ot}{\otimes}

\newcommand{\tth}{^{th}}

\renewcommand{\L}{\Lambda}
\newcommand{\Y}{\mathbb{Y}}

\newcommand{\tT}{\tilde{T}}
\newcommand{\tf}{\tilde{f}}

\DeclareMathOperator{\Hom}{Hom}
\DeclareMathOperator{\Sur}{Sur}
\DeclareMathOperator{\End}{End}
\DeclareMathOperator{\Res}{Res}
\DeclareMathOperator{\new}{new}
\DeclareMathOperator{\Disc}{Disc}
\DeclareMathOperator{\Nm}{Nm}
\DeclareMathOperator{\Mod}{Mod}
\DeclareMathOperator{\Cok}{Cok}
\DeclareMathOperator{\cok}{Cok}

\DeclareMathOperator{\Ext}{Ext}

\DeclareMathOperator{\poly}{poly}

\DeclareMathOperator{\Cl}{Cl}

\DeclareMathOperator{\Gal}{Gal}
\DeclareMathOperator{\Aut}{Aut}

\DeclareMathOperator{\Ker}{Ker}

\DeclareMathOperator{\GL}{GL}

\DeclareMathOperator{\Cov}{Cov}

\DeclareMathOperator{\diag}{diag}

\DeclareMathOperator{\rank}{rank}
\DeclareMathOperator{\id}{id}
\let\Im\relax
\DeclareMathOperator{\Im}{Im}
\newcommand{\tZ}{\tilde{Z}}
\newcommand{\tx}{\tilde{x}}
\newcommand{\mfx}{\mathfrak{x}}
\newcommand{\modzp}{\Mod_{\Z_p[\mathfrak{x}]}}
\newcommand{\emod}{E^\times \backslash \modzp}
\newcommand{\Den}{V}

\newcommand{\mcc}{\mathcal{C}}

\newcommand{\JC}[1]{{\color{purple} \sf $\clubsuit$ Jiahe's comment: [#1]}}

\newcommand{\RC}[1]{{\color{green} \sf $\clubsuit$ Roger's comment: [#1]}}

\renewcommand{\JC}[1]{}

\renewcommand{\RC}[1]{}


%% file: Introduction.tex
\section{Introduction}

\subsection{Preface} During teatime one day in April of 1972, Hugh Montgomery mentioned to Freeman Dyson his then-unpublished conjecture \cite{montgomery1973pair} on the distribution of the nontrivial zeroes of the Riemann zeta function on $\tfrac{1}{2}+i\R$. Suitably normalizing the imaginary parts, they appeared to behave like a random collection of points on the real line, with density $1$ and correlations given by
\begin{equation}\label{eq:GUE_pair}
    \E[\#\{\text{zeroes in $U$}\} \cdot \#\{\text{zeroes in $V$}\}] = \int_{U \times V}1 - \pfrac{\sin (\pi(x-y))}{\pi(x-y)}^2dx dy
\end{equation}
for any disjoint measurable sets $U,V \subset \R$. 

Dyson immediately recognized this integrand as the limiting pair correlation function of eigenvalues of large Hermitian or unitary random matrices, which he and others had been studying as models of the energy levels of heavy nuclei \cite{Thomas2013-IAS-PrimesRandomMatrices}. This now-famous episode sparked a blaze of activity devoted to extending these observations to other families of $L$-functions, and to finer statistics of the $L$-functions, such as the distribution of zeroes close to the real line averaged over a family of $L$-functions; see for instance the surveys of Katz-Sarnak \cite{katz2023random} or Keating-Snaith \cite[Chapter 24]{akemann2011oxford}. These and other random matrix statistics have now arisen in a surprising number of areas of pure and applied mathematics and physics.

Later, Ellenberg-Jain-Venkatesh \cite{ellenberg2011modeling} proposed a $p$-adic version of this story, namely that the eigenvalues of random matrices over $\Z_p$ also model the asymptotic distribution of zeroes of $p$-adic $L$-functions, which are points in $\bar{\Q}_p$. Computing the zeroes of large families of $L$-functions numerically and comparing with exact formulas on the random matrix side, they found very good agreement between the distributions of the number of zeroes and the number of eigenvalues in the open unit ball\footnote{This number, for the $L$-functions associated to number fields $K$ which were studied in \cite{ellenberg2011modeling}, is called the \emph{Iwasawa $\lambda$-invariant} and governs the growth of class numbers in cyclotomic towers over $K$.} of $\bar{\Q}_p$ centered at $1$. On the $L$-function side, several works such as Delbourgo-Chao \cite{delbourgo2015invariants}, Delbourgo-Knospe \cite{delbourgo2023iwasawa}, Knospe \cite{knospe2024special}, Kundu-Washington \cite{kundu2024first,kundu2024heuristics}, and Santato \cite{santato2021zeroes} have computed more data for other $L$-function families, finding good agreement with the $\lambda$-invariant conjecture of \cite{ellenberg2011modeling}. 

There has been little study of any statistics on the zeroes/eigenvalues other than the number which lie in the unit ball. This seems surprising, given both the success of the classical $L$-function/random matrix connection and the explicit push from the original work in this direction:

\begin{displayquote}
\emph{More ambitiously, one can speculate that the distribution of more refined
features of the $p$-adic $L$-function can also be predicted by random matrix
models; we discuss one example, concerning the spacing of zeroes. Here the
numerical evidence is less compelling but not wholly discouraging...} 

\quad \quad \quad -\cite[p3]{ellenberg2011modeling}
\end{displayquote}
The ``one example'' was the product of squared differences of all eigenvalues/zeroes in the open unit ball around $1$ (the discriminant of the Iwasawa polynomial). It is a highly nonlinear statistic, so the lack of complete agreement may or may not simply be slow convergence, but we are not aware of any works studying it following \cite{ellenberg2011modeling} to determine whether this is the case. 

From the perspective of Dyson and Montgomery's original interaction and its outgrowth, it seems obvious to ask for more refined eigenvalue statistics on the random matrix side. Surprisingly, it has remained open since \cite{ellenberg2011modeling} to compute explicit formulas or estimates even for basic statistics like the pair correlation functions analogous to \eqref{eq:GUE_pair}.

The goal of the present work is to develop the basic theory of $p$-adic random matrix eigenvalues along the same lines as was done classically for random matrices over $\R$ and $\C$, and in particular to provide such explicit eigenvalue statistics. These statistics may be compared with data coming from zeroes of $L$-functions, or from anywhere else where $p$-adic random matrix statistics may be suspected in the future. If the history of classical random matrix theory is any indication, these statistics may not appear in only one place.

To the question of why such statistics on the $p$-adic side were not computed earlier, one possible answer is that it is apparently more difficult to do so than in the archimedean case. Our proofs combine techniques from several threads of research which have developed separately---some later than \cite{ellenberg2011modeling}---and to our knowledge have not previously been combined:
\begin{itemize}[left=1em]
    \item \emph{$p$-adic orbital integrals}, specifically certain module-counting computations appearing previously in works on this subject such as Yun \cite{yun2013orbital}.
    \item \emph{Zeroes of random polynomials with independent coefficients over $\Z_p$}. We use $p$-adic Kac-Rice formula and related results introduced by Caruso \cite{caruso2022zeroes} (following Evans \cite{evans2002elementary}) to study these, and many proofs are inspired by those of \cite{caruso2022zeroes} in this setting.
    \item \emph{Cokernels of discrete random matrices and the moment method for random modules.} Cokernels are well-studied and amenable to the moment method of Wood \cite{wood2017distribution}, but are \emph{a priori} different statistics than eigenvalues. However, we reduce understanding eigenvalues to understanding cokernels of polynomials $Z(A)$ in a random matrix $A$, at which point we may use Sawin-Wood's generalized moment method \cite{sawin2022moment} as applied to this problem by Cheong-Yu \cite{cheong2023distribution}.
    \item \emph{Markov chains related to measures on integer partitions.} Specifically, we use Fulman's diagonalization \cite{fulman-RR} of a Markov chain found independently by Evans \cite{evans2002elementary} to sample the cokernel of $p$-adic random matrix.
\end{itemize}
We hope that in addition to the theorems, the present paper provides some value by giving a blueprint for how the above techniques may be used together to prove future results. Much remains to be explored in the $p$-adic case in comparison with the copious literature on random matrices over $\R$ and $\C$. For example, we have made no attempt to treat matrix distributions outside of the Haar distributions on $\Mat_n(\Z_p)$ and $\GL_n(\Z_p)$, but it would be natural to treat other `symmetry types' coming from different groups or homogeneous spaces, and we believe the same methods will be useful there.

One feature which makes the non-archimedean case difficult is the existence of infinitely many extensions of $\Q_p$ in which the eigenvalues can live. This is also one of the interesting features, but for simplicity we will first state those of our results which bypass it, either by considering only the eigenvalues in $\Q_p$ or by considering balls in $\bar{\Q}_p$ without differentiating by extension. These give a good sense of many features of the probabilistic picture with a minimum of complication, and are explained in \Cref{subsec:first_results}. In \Cref{subsec:joint_density_intro} we return to the issue of eigenvalues in different extensions, explain how to treat them through correlation functions, and characterize the eigenvalue distribution of an $n \times n$ matrix, which is the starting point of all of our asymptotic results. In \Cref{subsec:asymp_intro} we give some asymptotics of eigenvalues in extensions of $\Q_p$ which may be derived from it. \Cref{subsec: result in real ginibre} discusses related results in the archimedean case, and one we prove in the $p$-adic case which was inspired by these, and \Cref{subsec:random_polynomials} discusses relations to the random polynomials literature in more detail. \Cref{subsec:proof_methods} outlines the ideas of the proofs and structure of the paper.

\subsection{First results} \label{subsec:first_results} Recall that $\bar{\Q}_p$ is equipped with a norm, which we denote $||\cdot||$, which restricts to the usual norm on $\Q_p$ defined by $||u\cdot p^k||=p^{-k}$ for $u \in \Z_p^\times$. All eigenvalues of a matrix $A \in \Mat_n(\Z_p)$ lie in $\bar{\Z}_p$, the closed unit ball inside the algebraic closure $\bar{\Q}_p$ with respect to this norm, since the coefficients of the characteristic polynomial lie in the closed unit ball $\Z_p \subset \Q_p$. The different topology of this space leads to a very different limit than in classical archimedean random matrix theory:

\begin{thm}\label{thm:lifted_subspace_intro_revised}
    There exist sets $\mc{U}_1,\mc{U}_2,\ldots$, such that
    \begin{enumerate}[left=1em]
        \item $\bar{\Z}_p = \displaystyle \bigsqcup_{i \geq 1} \mc{U}_i$.\label{item:disjoint_union}
        \item Each set $\mc{U}_i$ is a finite union of some number $d_i$ of open balls of radius $1$ in $\bar{\Z}_p$, each of which is Galois-conjugate to the others.\label{item:balls}
        
        \item Let $A \in \Mat_n(\Z_p)$ be a random matrix with independent entries distributed by the additive Haar probability measure on $\Z_p$. Then for any finite collection $\mc{U}_{i_1},\ldots,\mc{U}_{i_k}$, the eigenvalues in each $\mc{U}_{i_j}$ are asymptotically independent of those in the other ones $\mc{U}_{i_{j'}}$ as $n \to \infty$, in the sense of factorization of correlation functions (see \Cref{thm: independent distribution over lifted subspaces}). Furthermore, the number of eigenvalues of $A$ lying in $\mc{U}_i$ converges in distribution as $n \to \infty$ to a $d_i \Z_{\geq 0}$-valued random variable $N_i$, with distribution defined by 
        \begin{equation}\label{eq:number_eigenvalues_on_island_intro}
            \mathbf{P}(N_i/d_i = j) = \left(\prod_{k \geq 1} (1-p^{-d_i k })\right) \cdot \frac{p^{-d_i j}}{\prod_{k=1}^j (1-p^{-d_i k})}
        \end{equation}
        for $j \geq 0$.\label{item:independence}

    \end{enumerate}
    The sets $\mc{U}_i$ have the following explicit description. Enumerate the nonzero monic irreducible polynomials in $\F_p[x]$ by $F_1,F_2,\ldots$, with degrees $d_1,d_2,\ldots$. Let $\mc{U}_i \subset \bar{\Z}_p$ be the set of all $x \in \bar{\Z}_p$ such that the minimal polynomial of $x$, reduced modulo $p$, is a power of $F_i$.
\end{thm}

This gives a picture of eigenvalues as living on an infinite archipelago of small islands $\mc{U}_i$, each island populated by a random small number of eigenvalues, which interact with each other but not with those on the other islands. One such island is the open unit ball around $0$, corresponding to the polynomial $F_1=x$, and the main random matrix result which \cite{ellenberg2011modeling} compared with zeroes of $L$-functions was a $\GL_n$ analogue of the formula \eqref{eq:number_eigenvalues_on_island_intro} in this case. \Cref{thm:lifted_subspace_intro_revised} may in fact be obtained by modifying the arguments of \cite{ellenberg2011modeling}, though we obtain it instead as a byproduct of machinery used for other results.

The question then becomes what the interactions between eigenvalues on a given island $\mc{U}_i$ look like asymptotically. To understand this, it becomes necessary to distinguish eigenvalues by which field extension $K/\Q_p$ they lie in. This distinction, and treating eigenvalue correlations rather than the number of eigenvalues lying in a ball, are core differences between our approach and that of \cite{ellenberg2011modeling}. 

The simpler results for those eigenvalues lying in $\Z_p$ already give a good sense of the general ones, so we state these first before setting up notation in general extensions. 

\begin{thm}\label{thm:Zp_2pt_intro}
    Let $U,V \subset \Z_p$ be two disjoint measurable sets. For each $n$, let $A \in \Mat_n(\Z_p)$ have independent entries distributed by the additive Haar probability measure on $\Z_p$, and let $Z_{U,n}, Z_{V,n}$ be the number of eigenvalues of $A$ in $U$ and $V$ respectively. Then 
    \begin{equation}\label{eq: eigenvalues in UXV_intro}
        \lim_{n \to \infty} \E[Z_{U,n} \cdot Z_{V,n}] = \int_{U \times V} \rho^{(\infty)}_{\Q_p,\Q_p}(x,y) dx dy,
    \end{equation}
    where integration is with respect to the Haar probability measure, and the limiting pair correlation function is given by
    \begin{equation}\label{eq: 2pt correlation_intro}
        \rho^{(\infty)}_{\Q_p,\Q_p}(x,y) = 1 - \theta_3(-\sqrt{p};||x-y||^2/p),
    \end{equation}
    where 
    $$
        \theta_3(z;t) = \prod_{i \geq 1} (1-t^i) (1+t^{i-1/2}z) (1+t^{i-1/2}/z) = \sum_{k \in \Z} t^{\frac{k^2}{2}}z^k
    $$
    is a Jacobi theta function. Furthermore, 
    \begin{equation} \label{eq:asymp_cov_Zp_intro}
        \lim_{n \to \infty} \Cov(Z_{U,n},Z_{V,n}) = - \int_{U \times V} \theta_3(-\sqrt{p};||x-y||^2/p) dx dy \leq 0.
    \end{equation}
\end{thm}

\Cref{thm:Zp_2pt_intro} may be viewed as an analogue of the sine kernel \eqref{eq:GUE_pair} in the archimedean case (though a closer archimedean analogue is the pair correlation of real eigenvalues of the real Ginibre ensemble discussed below). Taking $U,V$ to be very small balls around two points $x$ and $y$ in either \eqref{eq:GUE_pair} or \Cref{thm:Zp_2pt_intro}, the integrand can be viewed as the density of pairs of eigenvalues close to $x$ and $y$ respectively. In the classical case, the fact that the integrand in \eqref{eq:GUE_pair} converges to $0$ as $x-y \to 0$ witnesses the well-known repulsion of eigenvalues from one another. It is easily checked from the sum form of the theta function that the integrand in \Cref{thm:Zp_2pt_intro} converges to $0$ as $x - y \to 0$ as well, which shows pair-repulsion in the $p$-adic setting, as does the negative correlation in \eqref{eq:asymp_cov_Zp_intro}. While these properties of the pair correlation function tell a straightforward story, we do not have any conceptual reason for the mysterious appearance of the Jacobi theta function specifically, and this begs for explanation.

If the sets $U$ and $V$ are sufficiently far from one another so that $||x-y||=1$ for all $x \in U, y \in V$, then the numbers of eigenvalues in $U$ and $V$ are asymptotically uncorrelated as the equality case of \eqref{eq:asymp_cov_Zp_intro} holds. This corresponds to $U,V$ being in different sets $\mc{U}_i$ in \Cref{thm:lifted_subspace_intro_revised}. 

While \Cref{thm:Zp_2pt_intro} is just for the eigenvalues in $\Z_p$, we show pair-repulsion for eigenvalues in arbitrary extensions later in \Cref{thm:pair_correlation_estimate_intro}.

We obtain results such as \Cref{thm:Zp_2pt_intro}, and the independence statement in \Cref{thm:lifted_subspace_intro_revised}, by computing the joint eigenvalue distributions at fixed $n$ and taking the limit as $n \to \infty$. For eigenvalues in general extensions this requires some setup, and we explain it in the next subsection, but the result is of a similar form to the special case when all eigenvalues are in $\Z_p$ which we state now.

\begin{cor}[of {\Cref{thm: joint distribution}}]\label{thm:Zp_coulomb_intro}
Let $A \in \Mat_n(\Z_p)$ be a random matrix with independent entries distributed by the additive Haar probability measure on $\Z_p$, conditioned on the positive-probability event that all eigenvalues lie in $\Q_p$ (and hence in $\Z_p$). Then the joint density of its eigenvalues $x_1,\ldots,x_n $ is proportional to
\begin{equation}
    \prod_{1\le i<j\le n}||x_i-x_j||.
\end{equation}
Here, since the eigenvalues are not canonically ordered, this joint density is with respect to the pushforward to $\Z_p^{\vee n}$ of the Haar probability measure on $\Z_p^n$.
\end{cor}

\Cref{thm:Zp_coulomb_intro}, like \Cref{thm:Zp_2pt_intro}, has attractive similarities with previously studied real and complex random matrix ensembles. \Cref{thm:Zp_coulomb_intro} also shows a form of repulsion at finite $n$, since $||x_i-x_j||$ is small when $x_i$ is close to $x_j$. The same Vandermonde terms in eigenvalue densities over $\R$ and $\C$ are called \emph{Coulomb gas interactions} in the literature. To our knowledge this is their first appearance in $p$-adic random matrix theory, and we did not expect such close analogues with the archimedean case initially. However, such $p$-adic Coulomb gases (sometimes also called log-Coulomb gases) have been studied as objects in their own right by Webster \cite{webster2021log,webster2023log} and Sinclair \cite{sinclair2022non} in mathematics, and in physics by Z{\'u}{\~n}iga-Galindo with Torba \cite{zuniga2020non}, Zambrano-Luna and Le\'on-Cardenal \cite{zuniga2022graphs}, and Veys \cite{veys2025koba}. These works proceed without any relation to random matrices, but Sinclair writes in \cite{sinclair2022non} that ``\emph{Alternate titles for this paper include ``The p-adic Selberg Integral'' or ``p-adic Random Matrix Theory'' due to the appearance of a non-archimedean version of the Selberg Integral$\ldots$ The connection to random matrix theory is (currently) more tentative, since there are no random matrices
introduced in this paper}''. \Cref{thm:Zp_coulomb_intro} provides this connection.

\subsection{Extension-sorting and joint eigenvalue distribution for a fixed-size matrix}\label{subsec:joint_density_intro}

Eigenvalues of a matrix over $\Z_p$ can lie in any of the infinitely many extensions of $\Q_p$, and do so with positive probability asymptotically, making the general formulas more complicated than \Cref{thm:Zp_coulomb_intro} and \Cref{thm:Zp_2pt_intro}. 

We find the case of nonsymmetric real random matrices to be a helpful mental model. Here there are only two extensions, $\R$ and $\C$, which eigenvalues can lie in, but the same general behavior occurs. An $n \times n$ \emph{real Ginibre matrix} is one in $\Mat_n(\R)$ with independent real standard Gaussian entries. For a real matrix, the eigenvalues are either real or appear in conjugate pairs, and in general there are $\ell$ conjugate pairs and $n-2\ell$ real eigenvalues. Naively, one might hope that the joint distribution of all $n$ eigenvalues is absolutely continuous with respect to the Lebesgue measure on $n$ unordered copies of $\C$ and can be described by a single density function with respect to this measure. This turns out not to be the case, because the set $\R$ (which has Lebesgue measure $0$ in $\C$) typically contains a nonzero number of eigenvalues.

The complex conjugate pairs may each be described by a single representative in the open upper half-plane $\H_+$, which is just the quotient of $\C \setminus \R$ by complex conjugation. So one may identify the space on which the $n$-tuples of eigenvalues live with 
\begin{equation}\label{eq:ginibre_disjoint_union}
     \bigsqcup_{\ell=0}^{\floor{\frac{n}{2}}} \H_+^{\vee \ell} \times \R^{\vee n-2\ell},
\end{equation}
where we take symmetric products since the eigenvalues are not canonically ordered. Of course, $\R$ itself is canonically ordered, and the complex eigenvalues may be ordered by real part, so eigenvalue multisets may be identified with points on an appropriate ordered subset of $\H_+^{\ell} \times  \R^{n-2\ell}$. One may specify the full joint distribution on the space \eqref{eq:ginibre_disjoint_union} by a collection of density functions $\rho_{n,\ell}$ with respect to the Lebesgue measure on these subsets, for each $0 \leq \ell \leq \floor{\frac{n}{2}}$. 

In the $p$-adic case the situation is the same in broad strokes, but there are many more components in the disjoint union analogous to \eqref{eq:ginibre_disjoint_union} since there are many more extensions. A further issue is that it is not so simple to explicitly describe and order the quotient sets analogous to $\H_+$ in general. The eigenvalues of $A \in \Mat_n(\Q_p)$ live in the algebraic closure $\bar{\Q}_p$, which is infinite-dimensional over $\Q_p$. If $x \in \bar{\Q}_p$ is an eigenvalue of some matrix $A$, then all conjugates of $x$ by $\Gal(\bar{\Q}_p/\Q_p)$ are also eigenvalues, and the minimal polynomial of $x$ will appear as an irreducible factor of the characteristic polynomial of $A$. It is natural to view such an eigenvalue inside the smallest extension of $\Q_p$ in which it lies. For this purpose we define the set of \emph{new} elements in an extension $K/\Q_p$ by
\begin{equation}
    K^{\new} = K \setminus \bigcup_{\Q_p \subset L \subsetneq K} L
\end{equation}
where the union is over extensions $L$ strictly contained in $K$, so that 
\begin{equation}\label{eq:Qbar_disjoint_union}
    \bar{\Q}_p = \bigsqcup_{K/\Q_p} K^{\new}
\end{equation}
where the disjoint union is over all finite extensions $K$. 

The unit ball in $K$, equivalently $\bar{\Z}_p \cap K$, is the ring of integers $\mc{O}_K$, and the eigenvalues of a matrix $A \in \Mat_n(\Z_p)$ which lie in $K$ will always lie in $\mc{O}_K$. We denote $\mc{O}_K^{\new} = K^{\new} \cap \mc{O}_K$, and note that the set $\bar{\Z}_p$ in which the eigenvalues lie is a disjoint union of the sets $\mc{O}_K^{\new}$ by \eqref{eq:Qbar_disjoint_union}. Each $\mc{O}_K$ comes equipped with an additive Haar probability measure $\mu_K$, and the complement $\mc{O}_K \setminus \mc{O}_K^{\new}$ has measure $0$.

To put an explicit probability measure on the space analogous to \eqref{eq:ginibre_disjoint_union} in our setting would require describing the quotients of each $K^{\new}$ by the Galois action and putting an ordering on them. We find it much easier to consider the equivalent viewpoint of \emph{correlation functions}, which is standard in classical random matrix theory. The pair-correlation function of eigenvalues in $\Z_p$ already appeared in \Cref{thm:Zp_2pt_intro} above, and in general we define them as follows.

\begin{defi}
    \label{def:correlation_function}
    Let $X$ be a random subset of $\bar{\Z}_p$ such that $\sigma(X)=X$ with probability $1$ for any $\sigma \in \Gal(\bar{\Q}_p/\Q_p)$, and for any finite extension $K/\Q_p$, $X \cap \mathcal{O}_K$ is finite almost surely. Then for extensions $K_1,\ldots,K_m$ of $\Q_p$, we say that
    $$\rho_{K_1,\ldots,K_m}: \mc{O}_{K_1}^{\new} \times \cdots \times \mc{O}_{K_m}^{\new} \to \R_{\geq 0}$$ 
    is the \emph{$(K_1,\ldots,K_m)$-correlation function} of $X$ if for any bounded measurable function
    $$f: \mc{O}_{K_1}^{\new} \times \cdots \times \mc{O}_{K_m}^{\new} \to \C, $$
    we have
    \begin{equation*}
        \E\left[\sum_{\substack{\{x_1,\ldots,x_m\} \subset X}} f(x_1,\ldots,x_m)\right] = \displaystyle \int_{\mc{O}_{K_1}^{\new} \times \cdots \times \mc{O}_{K_m}^{\new}} f(x_1,\ldots,x_m) \rho_{K_1,\ldots,K_m}(x_1,\ldots,x_m) dx_1 \cdots dx_m
    \end{equation*}
    where the sum is over all tuples $(x_1,\ldots,x_m)  \in \mc{O}_{K_1}^{\new} \times \cdots \times \mc{O}_{K_m}^{\new}$ such that no two of the $x_i$ are equal or Galois-conjugate and $\{x_1,\ldots,x_m\} \in X$, and integration is with respect to the Haar probability measure on $\prod_{i=1}^m \mc{O}_{K_i}$.

    In the case where $X$ is the eigenvalue set of a random matrix $A \in \Mat_n(\Z_p)$ with i.i.d. additive Haar entries, we denote its $(K_1,\ldots,K_m)$-correlation function by $\rho^{(n)}_{K_1,\ldots,K_m}$.
    
\end{defi}

Correlation functions allow us to express expectations over subsets of the random point configuration, which do not come with a canonical ordering, using functions and integrals over ordered products. In general, taking $f$ in \Cref{def:correlation_function} to simply be an indicator function of a measurable set $U \subset  \mc{O}_{K_1}^{\new} \times \cdots \times \mc{O}_{K_m}^{\new}$, we have that 
 \begin{equation}\label{eq:number_roots_cor_fns_intro}
        \E[Z_{U,n}] = \int_U \rho^{(n)}_{K_1,\ldots,K_m}(x_1,\ldots,x_m) dx_1 \cdots dx_m
    \end{equation}
where $Z_{U,n}$ is the number of $m$-tuples $(x_1,\ldots,x_m)$ of distinct, non-conjugate eigenvalues lying in $U$. One should interpret $\rho^{(n)}_{K_1,\ldots,K_m}(x_1,\ldots,x_m)$ as the density of $m$-tuples of eigenvalues of $X$ near $(x_1,\ldots,x_m)$, in the sense that if $U_1,\ldots,U_m$ are sufficiently small balls around (distinct, non-conjugate) $x_1,\ldots,x_m$, then 

\begin{equation}
    \mathbf{P}(\text{Each set $U_i$ contains an eigenvalue}) \approx \rho_{K_1,\ldots,K_m}(x_1,\ldots,x_m) \mu_{K_1}(U_1) \cdots \mu_{K_m}(U_m).
\end{equation}

\noindent We will simply prove the equalities \eqref{eq:number_roots_cor_fns_intro}, which are equivalent to the general \Cref{def:correlation_function} by approximating $f$ by simple functions.

In particular, the correlation function $\rho^{(n)}_{K_1,\ldots,K_m}$ determines the joint distribution of eigenvalues on the event that they lie in $K_1,\ldots,K_m$. The set of these correlation functions, ranging over all tuples $(K_1,\ldots,K_m)$ with $\sum_{i=1}^m [K_i : \Q_p] = n$, hence gives the full description of the joint distribution of eigenvalues of a random Haar matrix $A \in \Mat_n(\Z_p)$. This set of correlation functions is also finite, since there are only finitely many extensions of $\Q_p$ in each degree by Krasner's theorem \cite{krasner1946nombre}.

\Cref{thm: joint distribution} below is an explicit formula for exactly these correlation functions. To state it, we introduce some notation.

\begin{defi}\label{def:delta_sigma}
    Let $x_1,\ldots,x_m \in \bar{\Q}_p$, and for each $i$ let $\sigma_{i,l_i}(x_i),1\le l_i\le r_i$ denote the full set of Galois conjugates of $x_i$. Then we write
    $$\Delta_\sigma(x_1,\ldots,x_m)=\Delta(\sigma_{1,1}(x_1),\ldots,\sigma_{1,r_1}(x_1),\ldots,\sigma_{m,1}(x_m),\ldots,\sigma_{m,r_m}(x_m))$$
    where $\Delta(z_1,\ldots,z_r)=\prod_{1\le i<j\le r}(z_i-z_j)$ is the usual Vandermonde determinant. 
\end{defi}

Note that the absolute value of the usual discriminant is related to $\Delta_\sigma$ by $$||\Disc(x_1,\ldots,x_m)||=||\Delta_\sigma(x_1,\ldots,x_m)||^2.$$

\begin{defi}\label{defi: Den and distance}
Let $K/\Q_p$ be an extension of degree $r$, $\mc{O}_K$ its ring of integers, and $x \in \mc{O}_K^{\new}$. Then we define
$$\Mod_{\Z_p[x]}=\{M\subset K\mid M \text{ is a $\Z_p[x]$-module, also a $\Z_p$-lattice of rank $r$}\}.$$ 
We have a canonical action of the group $K^\times$ on $\Mod_{\Z_p[x]}$ by multiplication, and similarly for its subgroup $\Lambda = \{\pi^n: n \in \Z\}$ where $\pi$ is a uniformizer of $K$. This latter action has finitely many orbits, and we denote the number by $\#(\Lambda \backslash \Mod_{\Z_p[x]})$. 

Finally, we define a function $\Den: \mc{O}_K \to \R$ (which depends implicitly on $K$) by
$$\Den(x)=\frac{1}{1-p^{-r/e}}\cdot ||\Delta_\sigma(x)||\cdot \#(\L\backslash \Mod_{\Z_p[x]})$$
where $e$ is the ramification index of $K/\Q_p$.
\end{defi}

\begin{thm}\label{thm: joint distribution}
Fix $n \geq 1$ and let $K_1,\ldots,K_m$ be extensions of $\Q_p$ with $\sum_{i=1}^m [K_i: \Q_p] = n$. Let $U \subset \mc{O}_{K_1}^{\new} \times \cdots \times \mc{O}_{K_m}^{\new}$ be any measurable set. Let $A \in \Mat_n(\Z_p)$ be random with independent Haar-distributed entries, and $Z_{U,n}$ be the number of $m$-tuples $(x_1,\ldots,x_m)$ of its distinct, non-conjugate eigenvalues which lie in $U$. Then
$$\E[Z_{U,n}]=\int_U\rho^{(n)}_{K_1,\ldots,K_m}(x_1,\ldots,x_m)dx_1\cdots dx_m,$$
where 
\begin{equation}
\rho^{(n)}_{K_1,\ldots,K_m}(x_1,\ldots,x_m)=(1-p^{-1})\cdots(1-p^{-n})||\Delta_\sigma(x_1,\ldots,x_m)|| \prod_{i=1}^m\Den(x_i)
\end{equation}
and integration (and definition of measurability of $U$) is with respect to the Haar probability measure on $\mc{O}_{K_1}\times \cdots \times \mc{O}_{K_m}$. 
\end{thm}

For explicit examples of $\#(\L\backslash \Mod_{\Z_p[x]})$ in the case of quadratic extensions, see \Cref{prop: orbital integral of quadratic extension}. For $x \in \Q_p$ it is just $1$, and so \Cref{thm:Zp_coulomb_intro} follows immediately from \Cref{thm: joint distribution}. \Cref{thm: joint distribution} shows that the Coulomb gas behavior seen over $\Z_p$ in \Cref{thm:Zp_coulomb_intro} is in fact a general phenomenon for all the eigenvalues, not just the ones in $\Z_p$. The Vandermonde leads directly to repulsion of eigenvalues in the general setting of \Cref{thm: joint distribution} as well, which we quantify in \Cref{thm:pair_correlation_estimate_intro}.

From the standpoint of seeing random matrix statistics in other phenomena such as zeroes of $L$-functions, it is not the finite-$n$ distributions of \Cref{thm: joint distribution} which matter, but the $n \to \infty$ limits. In these limits one expects to see universal distributions which are insensitive to the distribution of the matrix entries, and which model other objects outside of random matrix theory. 

Explicitly, we wish to compute the \emph{limiting $(K_1,\ldots,K_m)$-correlation functions}
    \begin{equation}\label{eq:limiting_cor_functions_intro}
        \rho_{K_1,\ldots,K_m}^{(\infty)}(x_1,\ldots,x_m) := \lim_{n \to \infty} \rho_{K_1,\ldots,K_m}^{(n)}(x_1,\ldots,x_m).
    \end{equation}
These capture the asymptotic probabilities of seeing an eigenvalue near each of $x_1,\ldots,x_m$ simultaneously, hence collectively govern all limiting local statistics of eigenvalues. 

\subsection{Asymptotic results for eigenvalues in extensions of $\Q_p$}\label{subsec:asymp_intro}

We have already given exact formulas for the pair correlation of eigenvalues in $\Z_p$ in \Cref{thm:Zp_2pt_intro}, and give a similar result for eigenvalues in quadratic extensions below. This result concerns the $K$-correlation function $\rho_K^{(n)}(x)$, often simply called a \emph{density function} of eigenvalues. This is not to be confused with probability density functions; the integral over $\mc{O}_K^{\new}$ of $\rho_K$ will typically not be $1$.

\begin{thm}
    \label{thm:quadratic_intro}
    Let $K/\Q_p$ be an unramified quadratic extension and $x \in K$ be a choice of generator with $||x||=1$, so that $\mc{O}_K = \Z_p[x]$ and every element of $\mc{O}_K$ may be written uniquely as $a_0+a_1 x$ for some $a_0,a_1 \in \Z_p$. Then the limiting density function depends only on $||a_1||$ and is given by
    \begin{equation}
        \rho_K^{(\infty)}(a_0+a_1x) = \frac{p^{-m}(1+p^{-1}-2p^{-m-1})}{1-p^{-1}}\sum_{k \geq 0} p^{-\frac{(2m+1)k^2+(4m+1)k}{2}}(1-p^{-k-1}),
    \end{equation}
    where $m$ is defined by $p^{-m}=||a_1||$. 

    In the case where $K$ is ramified, choosing $x \in K$ to be a generator with $||x||=p^{-1/2}$, we instead have 
    \begin{equation}
        \rho_K^{(\infty)}(a_0+a_1x) = ||\Disc_{K/\Q_p}||\frac{p^{-m}(1-p^{-m-1})}{1-p^{-1}}\sum_{k \geq 0} p^{-(m+1)k^2 - (2m+1)k}(1-p^{-2k-2})
    \end{equation}
    for any $a_0,a_1 \in \Z_p$ with $||a_1||=p^{-m}$. Note that when $p$ is odd, the term $||\Disc_{K/\Q_p}||$ is just $p^{-1}$.
    
\end{thm}

Because eigenvalues $x$ in a quadratic extension occur in conjugate pairs, the Vandermonde term in \Cref{thm: joint distribution} suggests that the two eigenvalues in the pair interact with one another by Coulomb interactions and so will be unlikely to be found nearby to one another. \Cref{thm:quadratic_intro} shows this to be the case: in the unramified case, the correlation function expands as 
\begin{equation}
\rho_K^{(\infty)}(a_0+a_1x) = p^{-m}(1+O_m(p^{-1}))    
\end{equation}
as a series in $p^{-1}$, and $p^{-m} = ||a_1|| = ||(a_0+a_1x) - \sigma(a_0+a_1x)||$ where $\sigma$ is the nontrivial Galois automorphism of $K$. Hence the main contribution to the density function comes from the Coulomb repulsion between the eigenvalue and its conjugate. The ramified case is the same apart from $m$-independent multiplicative constants.

We note also that these series in \Cref{thm:quadratic_intro} are similar to the series formula for the $\theta_3$ function in \Cref{thm:Zp_2pt_intro}. All three series are linear combinations of constants and `partial theta functions' $\theta_P(z;t) = \sum_{k \geq 0} t^{k^2/2}z^k$, which are also studied in the literature, see e.g. Warnaar \cite{warnaar2003partial}. Unlike the series in \Cref{thm:Zp_2pt_intro}, those in \Cref{thm:quadratic_intro} do not seem to have product forms.

Beyond the cases of $\Z_p$ and quadratic extensions, we do not presently have exact algebraic formulas, but we do have estimates. The following one says that the Vandermonde term governs repulsion between (Galois orbits of) eigenvalues in a straightforward fashion, in any field extension.

\begin{thm}[{\Cref{thm: estimate of two point correlation}} in text]
    \label{thm:pair_correlation_estimate_intro} 
    Let $\mc{U}$ be the lifted subspace corresponding to irreducible polynomial $F \in \F_p[x]$ of degree $d$. Let $x_1,x_2 \in \mc{U}$ be two distinct elements which are not Galois-conjugate, and let $\sigma_1(x_1),\ldots,\sigma_{\ell}(x_1)$ and $\sigma_1'(x_2),\ldots,\sigma_r'(x_2)$ be the full Galois orbits of $x_1$ and $x_2$. Then letting $K_1,K_2$ be the extensions (possibly the same) of $\Q_p$ for which $x_1 \in K_1^{\new},x_2 \in K_2^{\new}$, we have
    \begin{equation}\label{eq:2-pt_estimate_intro}
        \rho_{K_1,K_2}^{(\infty)}(x_1,x_2) = \rho_{K_1}^{(\infty)}(x_1) \cdot \rho_{K_2}^{(\infty)}(x_2) \cdot \left(\prod_{\substack{1 \leq i \leq \ell \\ 1 \leq j \leq r}} ||\sigma_i(x_1) - \sigma_j'(x_2) ||\right) \cdot I(x_1,x_2),
    \end{equation}
    with error term
    \begin{equation}
        \prod_{k=1}^\infty (1-p^{-kd}) < I(x_1,x_2) < \frac{1}{\prod_{k=1}^\infty (1-p^{-kd})}.
    \end{equation}
\end{thm}

The left-hand side of \eqref{eq:2-pt_estimate_intro} is the probability density of finding two eigenvalue orbits, one close to $x_1$ and one close to $x_2$, while the term $\rho_{K_1}^{(\infty)}(x_1)$ on the right-hand side gives the probability density of finding an eigenvalue orbit close to $x_1$ independent of what happens close to $x_2$, and similarly for $\rho_{K_2}^{(\infty)}(x_2)$.  The product term, which equivalently is the resultant of the minimal polynomials of $x_1$ and $x_2$, is always strictly less than $1$ because lifted subspaces are unions of Galois-conjugate open balls of radius $1$, and is in fact at most $p^{-1}$. Since the error term $I(x_1,x_2)$ is approximately $1$, finding eigenvalue orbits close to $x_1$ and $x_2$ simultaneously is less likely\footnote{A careful reader may notice that \Cref{thm:pair_correlation_estimate_intro} does not quite show the pair-repulsion inequality
\begin{equation}\label{eq:pair_repulsion_intro}
    \rho_{\Q_p[x_1],\Q_p[x_2]}^{(\infty)}(x_1,x_2) < \rho_{\Q_p[x_1]}^{(\infty)}(x_1) \cdot \rho_{\Q_p[x_1]}^{(\infty)}(x_2)
\end{equation}
when $p=2$ and the resultant is $p^{-1}$, but we prove an additional bound in \Cref{thm: refinement of two point correlation} to show \eqref{eq:pair_repulsion_intro} in this case.} than it would be if eigenvalues near both locations behaved independently, due to the Coulomb-type interactions witnessed by the product term. Together, Theorems~\ref{thm:lifted_subspace_intro_revised} and~\ref{thm:pair_correlation_estimate_intro} show that eigenvalues in different lifted subspaces behave independently, while those in the same lifted subspace repel according to Coulomb interactions.

It is also natural to ask how the eigenvalues sort into extensions of $\Q_p$, and the following theorem provides some answers.

\begin{thm}
    \label{thm:roots_in_extensions_intro}
    For any $n \geq 1$, the expected number $\E[Z_{\Z_p,n}]$ of eigenvalues in $\Z_p$ for $A \in \Mat_n(\Z_p)$ with independent Haar entries is $1$. For a quadratic extension $K/\Q_p$,
    \begin{equation*}
        \lim_{n \to \infty} \E[Z_{\mathcal{O}_K^{\new},n}] = \begin{cases}
            (1-p^{-1})\displaystyle \sum_{k\ge 0}\frac{(1-p^{-k-1})(1+p^{-k^2-2k-3})p^{-\frac{k^2+k}{2}}}{(1-p^{-k^2-2k-2})(1-p^{-k^2-2k-3})} & K \text{ unramified} \\ 
            ||\Disc_{K/\Q_p}||(1-p^{-1})\displaystyle\sum_{k\ge 0}\frac{(1-p^{-2k-2})p^{-k^2-k}}{(1-p^{-k^2-2k-2})(1-p^{-k^2-2k-3})} & K \text{ ramified.}
        \end{cases}
    \end{equation*}
    
    Furthermore, let $K/\Q_p$ be any arbitrary finite extension, and $p^f$ the order of its residue field. Then $\lim_{n \to \infty} \E[Z_{\mc{O}_K^{\new},n}]$ exists, is finite, and satisfies the estimate
    \begin{equation}
        \lim_{n \to \infty} \E[Z_{\mathcal{O}_K^{\new},n}] = ||\Disc_{K/\Q_p}|| \cdot \sum_{d|f} \mu\left(\frac{f}{d}\right)p^{d-f} + \mathcal{E}(K).
    \end{equation}
    Here $\mu$ is the M\"obius function, and the error term $\mathcal{E}(K)$ satisfies the bounds
    \begin{equation}
        -p^{-f}||\Disc_{K/\Q_p}|| < \mathcal{E}(K) < \frac{1+\tau(f)}{1-p^{-f}}p^{-f}||\Disc_{K/\Q_p}||
    \end{equation}
    where $\tau$ is the divisor function. 
\end{thm}

The expected number of eigenvalues in $\Z_p$ and quadratic extensions are natural statistics to test numerically against $L$-function zeroes. We also compute the variance of the number of eigenvalues in $\Z_p$ in \Cref{thm:variance_of_roots_in_Z_p}. If one is interested in a related statistic, e.g. the number of eigenvalues in a quadratic extension $K$ which lie in the open unit ball around $0$ rather than the closed one, this may be derived by the same method as we use for the formulas in \Cref{thm:roots_in_extensions_intro}, namely integrating the correlation function in \Cref{thm:quadratic_intro} over the appropriate set. We hope to see such tests of the Ellenberg-Jain-Venkatesh heuristics soon.

\begin{rmk}
    \label{rmk:gln_and_matn}
    One may also consider matrices distributed by the Haar measure on the (multiplicative) group $\GL_n(\Z_p)$, which are the ones originally treated in \cite{ellenberg2011modeling} as models for zeroes of $L$-functions. In fact, the limiting correlation functions of their eigenvalues are identical to those for the additive Haar measure on $\Mat_n(\Z_p)$, save for the restriction that eigenvalues cannot lie in the open unit ball around $0$ due to invertibility; this is shown in \Cref{thm:limit_cor_fns_exist_GL}. In particular, the eigenvalues in the open unit ball around $1$ considered in \cite{ellenberg2011modeling} have the same asymptotic distribution as the eigenvalues in the open unit ball around $0$ for $\Mat_n(\Z_p)$ after shifting by $1$.
    
    The finite-$n$ joint density formulas for $\GL_n(\Z_p)$, given in \Cref{sec:glr} as well, are also almost identical to those for $\Mat_n(\Z_p)$. We have not written down the analogue of \Cref{thm:roots_in_extensions_intro} for $\GL_n(\Z_p)$, but it is easy to do by the same argument, integrating the same density functions over a slightly restricted set.
\end{rmk}

\subsection{Corresponding results in the real Ginibre case} \label{subsec: result in real ginibre}

The history of the real Ginibre ensemble gives some indication of the added difficulty introduced by the sorting of eigenvalues into different extensions. For the simpler case of real symmetric matrices, in which all eigenvalues are automatically real, the computation of the joint eigenvalue distribution is very short and was done by Rosenzweig-Porter \cite{rosenzweig1960repulsion} around the same time such matrices were introduced as a physical model. For the real Ginibre ensemble, by contrast, Ginibre \cite{ginibre1965statistical} computed the joint distribution of eigenvalues under conditioning that all eigenvalues are real (the analogue of \Cref{thm:Zp_coulomb_intro}), but it was over $20$ years later that Lehmann-Sommers \cite{lehmann1991eigenvalue} computed the joint distributions of eigenvalues for any number of real eigenvalues (the analogue of \Cref{thm: joint distribution}).

Notably, the joint densities obtained in \cite{lehmann1991eigenvalue} are of exactly the same form as \Cref{thm: joint distribution}. They feature the Vandermonde of eigenvalues times a product over all eigenvalues of certain potential functions. Of course, the form of the potential functions is different than our $\Den(x)$, and perhaps more importantly, the different geometry of $\bar{\Z}_p$ versus $\C$ yields quite different probabilistic behavior.

The question of $n\to \infty$ limiting formulas for the correlation functions (such as the analogue of \Cref{thm:Zp_2pt_intro}) for the real Ginibre ensemble was still widely thought to be unsolvable after the joint densities were computed in \cite{lehmann1991eigenvalue}. After several more years, however, Borodin-Sinclair \cite{borodin2009ginibre} and Forrester-Nagao \cite{forrester2007eigenvalue} independently computed them using formulas of Sinclair \cite{sinclair2007averages}. Those works made use of a non-obvious Pfaffian point process structure, for which we do not know any $p$-adic analogue. It remains an exciting open question to find in $p$-adic random matrix theory an analogue, if it exists, of the determinantal and Pfaffian point process structures which have been so powerful in classical random matrix theory.

Another natural question is the probability that a real Ginibre matrix has $k$ real eigenvalues. This was treated by Edelman-Kostlan-Shub \cite{edelman1994many} and Edelman \cite{edelman1997probability}, the latter of which independently rederived the joint densities of \cite{lehmann1991eigenvalue}. In particular, the probability that all eigenvalues are real was shown in \cite{edelman1997probability} to be exactly $2^{-n(n-1)/4}$. In the $p$-adic setting, we give an asymptotic version of this result, which shows the same exponential-quadratic decay.

\begin{thm}\label{thm: all eigenvalues in Z_p intro}
Let $A \in \Mat_n(\Z_p)$ be distributed by the additive Haar measure. Then
$$\mathbf{P}(\text{all eigenvalues of $A$ lie in $\Z_p$})=p^{-\frac{n^2}{2(p-1)}-\frac{1}{2}n\log_p n+O(n)}$$
as $n \to \infty$, where the implicit constant depends on $p$.
\end{thm}

\Cref{thm: all eigenvalues in Z_p intro} is proven together with the analogous statement for $\GL_n(\Z_p)$ in \Cref{thm: all eigenvalues in Z_p} below. It would be interesting to find good exact expressions for these probabilities at finite $n$ in the $p$-adic setting, analogous to the elegant formulas for the real Ginibre case in \cite{edelman1997probability}. 

\subsection{Relation to random polynomials over $\Z_p$ with independent coefficients} \label{subsec:random_polynomials}
The series formulas in \Cref{thm:roots_in_extensions_intro} for the limiting expected number of eigenvalues in a quadratic extension are derived easily from \Cref{thm:quadratic_intro}. Similarly, \Cref{thm:Zp_2pt_intro} yields formulas for the variance of the number of eigenvalues in $\Z_p$ (\Cref{thm:variance_of_roots_in_Z_p}). Remarkably, in all of these cases, the formula is invariant under the substitution $p \leftrightarrow p^{-1}$. We are not aware of this phenomenon being observed in $p$-adic random matrix theory before. However, the same $p \leftrightarrow p^{-1}$ invariance was observed by Bhargava-Cremona-Fisher-Gajovi\'c \cite{bhargava2022densitypolynomialsdegreen} in the analogous formulas for the expected number of roots lying in given extensions for a random polynomial over $\Z_p$ with independent Haar-distributed coefficients. The invariance in this setting was conjectured in general in \cite[Conjecture 1.2]{bhargava2022densitypolynomialsdegreen}, proven for quadratic and cubic extensions and for $\Q_p$ itself in the same paper, and proven in a further special case by Shmueli \cite{shmueli2022probability}. It was finally settled in general by G-Wei-Yin \cite{g2023chebotarev}, by geometric machinery which eventually reduces it to Poincar\'e duality. As \cite{g2023chebotarev} notes, the methods used are quite general, and we hope they may apply in the random matrix setting as well.

Several other comparisons with the literature on $p$-adic random polynomials are worth remarking upon. The subject predates \cite{bhargava2022densitypolynomialsdegreen}: roots of a $p$-adic polynomial with independent random coefficients were considered by Evans \cite{evans2006expected}, who computed the expected number of roots in $\Q_p$ when the polynomial has Haar coefficients, and Buhler-Goldstein-Moews-Rosenberg \cite{buhler2006probability}, who computed asymptotics of the probability that all roots lie in $\Q_p$, analogous to \Cref{thm: all eigenvalues in Z_p intro} (the proof of which uses their result). This direction was continued by Kulkarni-Lerario \cite{kulkarni2021p}, Ait El Manssour-Lerario \cite{ait2022probabilistic}, the above \cite{bhargava2022densitypolynomialsdegreen}, Shmueli \cite{shmueli2023expected}, Caruso \cite{caruso2022zeroes} and others; the introduction of \cite{caruso2022zeroes} provides a good overview of these developments. 

The most similar of these works to the present one is \cite{caruso2022zeroes}, which inspired many parts of our proof of \Cref{thm: joint distribution} and from which we borrow various techniques detailed below. The same phenomenon of different correlation functions corresponding to different extensions of $\Q_p$ is present there, though the formulas are different. More specifically, \cite[Theorem C]{caruso2022zeroes} computes explicitly the analogue of $\rho_K^{(\infty)}(x)$ in \Cref{thm:quadratic_intro}, finding simpler rational functions in $p^{-1}$ rather than series. \cite[Theorem D]{caruso2022zeroes} computes the expected number of roots in a given extension, finding essentially the same asymptotic as \Cref{thm:roots_in_extensions_intro} when the residue field is large. 

In summary, while we obtain different limiting statistics, many of the same broad phenomena---repulsion by Coulomb interactions, asymptotically finite number of roots in each extension, independence of roots at distance $1$, $p \leftrightarrow p^{-1}$ invariance---appear for both eigenvalues of random matrices and roots of random polynomials.

\subsection{Proof methods and structure of the paper} \label{subsec:proof_methods}

Proving \Cref{thm: joint distribution} is one of the main technical components of the paper, and everything else is proven from it by taking an $n \to \infty$ limit. Most rely also on computations of expected determinants which use the moment method of Sawin-Wood \cite{sawin2022moment}, but with these inputs the proofs of most asymptotic results are largely independent. An exception is the exact computations of correlation functions in Theorems~\ref{thm:Zp_2pt_intro} and~\ref{thm:quadratic_intro}: they have quite similar proofs and rely on another technical input, namely Markov chains related to measures on partitions and cokernels of random matrices.

After giving some necessary preliminaries in \Cref{sec:preliminaries}, Sections \ref{sec:char_poly} and \ref{sec: joint distribution} are devoted almost entirely to the proof of \Cref{thm: joint distribution}. The main step in the proof of \Cref{thm: joint distribution}, carried out in \Cref{sec:char_poly}, is to compute a formula (\Cref{thm: points on variety_main theorem}) for the Haar measure of the set
\begin{equation}\label{eq:A_s_set_intro}
    \{A\in\Mat_n(\Z_p)\mid ||P_{A}(x_i)||\le p^{-s_i},\quad\forall 1\le i\le m\}
\end{equation}
of matrices whose characteristic polynomial $P_A$ takes small values on points $x_1,\ldots,x_m$ (and their Galois conjugates). These points should be viewed as proxies for the eigenvalues, so this set is a tractable notion of the set of matrices where the eigenvalues are close to $x_1,\ldots,x_m$ and their conjugates. This is inspired by consideration of the analogous set in Caruso's work \cite{caruso2022zeroes} on roots of polynomials with i.i.d. coefficients.

We find that the matrices in \eqref{eq:A_s_set_intro} may be expressed uniquely as polynomials $f(A)$ where $A$ has eigenvalues \emph{exactly} $x_1,\ldots,x_m$ and their conjugates, and $f$ satisfies certain properties. This leads to considering the map
\begin{align}
    \begin{split}
        \GL_n(\Q_p) &\to \Mat_n(\Q_p) \\ 
        U &\mapsto UC(Z)U^{-1}
    \end{split}
\end{align}
where $C(Z)$ is the companion matrix of the degree-$n$ monic polynomial $Z$ with roots $x_1,\ldots,x_m$. We care about the intersection of this map's image with the set $\Mat_n(\Z_p) \subset \Mat_n(\Q_p)$ on which our Haar measure is supported. The measure of this intersection is an orbital integral, which are treated in many works. We use various results related to orbital integrals which may be found in a note by Yun \cite{yun2013orbital}.

Using all this, we express the set \eqref{eq:A_s_set_intro} as a disjoint union of small pieces, breaking it up successively over several steps, and keeping track of the (finite) number of pieces at each step. In the end, we may compute the measures of the smallest pieces by linearizing and passing to Lie algebra computations. This part is analogous to the Jacobian computations by which eigenvalue joint densities are computed in the archimedean random matrix literature, though considerably more complicated.

Once this formula for the measure of \eqref{eq:A_s_set_intro} is computed, we must extract the correlation functions in \Cref{thm: joint distribution} by taking a limit as the $s_i$ become large. The strategy here follows the one used by Caruso \cite{caruso2022zeroes} for random polynomials, and the tool needed is a $p$-adic Kac-Rice formula introduced there. This is done in \Cref{sec: joint distribution}, proving \Cref{thm: joint distribution}. 

We remark that in all of these computations, it is very useful to endow the space of tuples $(x_1,\ldots,x_m)$ of eigenvalues with the algebra structure on $E = K_1 \times \cdots \times K_m$, as then the tuple itself is a root of the characteristic polynomial. This viewpoint of étale algebras was also used for instance in \cite{caruso2022zeroes} and in works on orbital integrals such as \cite{yun2013orbital}.

The remainder of the paper is concerned with deriving the above-stated asymptotic results from \Cref{thm: joint distribution}. Except for \Cref{thm: all eigenvalues in Z_p intro}, these results rely on correlation functions. All come from the following, which is essentially a corollary of \Cref{thm: joint distribution} but is important enough to state fully here.

\begin{lemma}
        \label{thm:general_cor_functions_intro}
    Let $K_1,\ldots,K_m$ be extensions of $\Q_p$ and $r := \sum_{i=1}^m [K_i:\Q_p]$. For any $n \geq r$, the $(K_1,\ldots,K_m)$-correlation function of eigenvalues of an $n \times n$ matrix $A$ with i.i.d. additive Haar entries is
    \begin{equation}
    \rho_{K_1,\ldots,K_m}^{(n)}(x_1,\ldots,x_m)=\prod_{i=0}^{r-1} (1-p^{-n+i}) \cdot ||\Delta_\sigma(x_1,\ldots,x_m)|| \cdot \left(\prod_{i=1}^m\Den(x_i) \right) \cdot \E||\det(Z(\tilde A))||,
    \end{equation}
    where $Z(x)=\prod_{i=1}^m Z_i(x)$ for $Z_i\in\Z_p[x]$ the minimal polynomial of $x_i$, and $\tilde A\in\Mat_{n-r}(\Z_p)$ is distributed by the additive Haar measure.
\end{lemma}

It is derived from \Cref{thm: joint distribution} in \Cref{sec: joint distribution} by integrating out the remaining $n-r$ eigenvalues which are not accounted for by $x_1,\ldots,x_m$ and their Galois conjugates, using the product structure of the correlation functions in \Cref{thm: joint distribution}. This integration is the source of the $\E||\det(Z(\tilde A))||$ term, which is an average over these other eigenvalues.

It might seem that \Cref{thm:general_cor_functions_intro} has not made the problem of computing limiting correlation functions of the eigenvalues of $A$ any easier, since it features a determinant of another random matrix $\tilde A$ with the same entry distribution! However, an expected determinant is an easier object than a joint eigenvalue distribution (though still not easy, for general $Z$). 

The terms $\Den(x_i)$ are independent of $n$ and are relatively easy to compute explicitly for small extensions. Hence the challenge for computing limiting eigenvalue statistics is to understand the $n \to \infty$ limit of $\E||\det(Z(\tilde A))||$. When $Z$ is a linear polynomial this is routine, but even for quadratic polynomials (needed for Theorems~\ref{thm:Zp_2pt_intro} and \ref{thm:quadratic_intro}) doing an exact computation is surprisingly difficult and we succeed only after serious effort. After the proof of \Cref{thm: joint distribution}, the bulk of the remaining technical work in this paper goes toward understanding such limits. 

The first key observation is that for any nonsingular $A \in \Mat_n(\Z_p)$,
\begin{equation}
    ||\det(A)|| = \frac{1}{\#\Cok(A)},
\end{equation}
where $\Cok(A) = \Z_p^n/A\Z_p^n$ is the cokernel. The cokernel is an abelian group (and, for $\Cok(Z(A))$, has an additional $\Z_p[T]/(Z(T))$-module structure), so \emph{a priori} it is more complicated than the determinant. However, this extra structure actually makes it easier to study. 

Passing to cokernels brings in an important input: general machinery has been developed by Wood and collaborators to handle limits of random groups and modules through the so-called moment method. The motivation for the original works of Wood \cite{wood2015random,wood2017distribution} was the study of the distribution of the cokernel itself, which appears in Cohen-Lenstra heuristics for class groups and Jacobians of random graphs. Many other papers on cokernels followed, including a body of literature on the cokernels of polynomials $Z(A)$ in a random matrix, such as Cheong-Huang \cite{cheong2021cohen,cheong2023cokernel}, Cheong-Kaplan \cite{cheong2022generalizations}, Cheong-Liang-Strand \cite{cheong2023polynomial}, and Lee \cite{lee2023joint}. General results of Sawin-Wood \cite{sawin2022moment} gave a moment method for random modules, which was applied to $\Cok(Z(A))$ by Cheong-Yu \cite{cheong2023cokernel}. This gave a complete answer to the limiting distribution of $\Cok(Z(A))$, though in terms of certain $\Ext$ groups which can be hard to compute in practice. We use these results of \cite{sawin2022moment} and \cite{cheong2023cokernel} in \Cref{sec:cokernels} to prove several results on expected determinants. The literature on cokernels of matrix polynomials is relatively new, and as far as we are aware it arose simply because it is a natural generalization of the cokernel of a single matrix with independent entries; the present work appears to be its first application to another problem.

\Cref{sec:limit_cor_func} uses these results to establish existence and basic properties of the limiting correlation functions $\rho^{(\infty)}_{K_1,\ldots,K_m}$. In particular, the independence stated in \Cref{thm:lifted_subspace_intro_revised} and given in detail in \Cref{thm: independent distribution over lifted subspaces} is deduced there from factorization properties of $\E||\det(Z(\tilde A))||$ which come from the moment method. It is worth saying clearly that we stated \Cref{thm:lifted_subspace_intro_revised} as the first theorem because it is an enlightening and simple-enough probabilistic statement which clarifies what comes later, but at a technical level it is not as new or difficult as some other results in the paper. We obtain it almost as a by-product of our attempts to understand correlation functions in general, but it should also be possible to establish it directly by a different approach slightly extending the arguments of \cite{ellenberg2011modeling}. This is not true of our other results, for which the methods are fairly orthogonal to those of \cite{ellenberg2011modeling}.

\Cref{sec:glr} makes precise the claim of \Cref{rmk:gln_and_matn} that Haar-distributed elements of $\GL_n(\Z_p)$ have the same asymptotic local statistics except for the lack of eigenvalues in the open unit ball around $0$. The finite-$n$ joint distribution formula for $\GL_n(\Z_p)$ which we prove in that section is also very similar to the one in \Cref{thm: joint distribution}.

For the finer understanding of $\Cok(Z(A))$ needed for the exact computations of Theorems~\ref{thm:Zp_2pt_intro} and~\ref{thm:quadratic_intro}, which together encompass the case where $Z$ has degree $2$, we must draw on another body of work. To explain, first note that a Haar-distributed $A \in \Mat_n(\Z_p)$ may be generated by first sampling a uniformly random matrix in $\Mat_n(\Z/p\Z)$, then a uniformly random lift of this matrix to $\Mat_n(\Z/p^2\Z)$, then a uniformly random lift of this to $\Mat_n(\Z/p^3\Z)$, and so on. The groups $\Cok(A)/p^d\Cok(A)$ depend only on the reduction $A \pmod{p^d} \in \Mat_n(\Z/p^d\Z)$, so if we sample $A$ as above, we also sample $\Cok(A)/p\Cok(A), \Cok(A)/p^2\Cok(A),\ldots$ step by step. At each step, we can parametrize $\Cok(A)/p^d\Cok(A)$ by a collection of $d$ integers $\la_1 \geq \ldots \geq \la_d \geq 0$, defined by 
\begin{equation}
    \la_i = \rank_{\F_p}(p^{i-1}\Cok(A)/p^i\Cok(A))
\end{equation}
(since $p^{i-1}\Cok(A)/p^i\Cok(A)$ is $p$-torsion, it is an $\F_p$-vector space). The full collection $\la_1 \geq \la_2 \geq \ldots$ is an integer partition, and for random $A$ it is a random partition. We mention that these measures on partitions are part of the broader framework of so-called \emph{Macdonald processes} developed by Borodin-Corwin \cite{borodin2014macdonald}, see e.g. \cite{van2020limits}, and also arose in the work of Fulman \cite{fulman_main} on conjugacy classes of random matrices over finite fields, though we use neither connection here.

A non-obvious fact which we do use is that the probability of a given $\la_{i+1}$ conditional on $\la_1,\ldots,\la_i$ depends only on $\la_i$, i.e. this process forms a discrete-time Markov chain. Additionally, this transition probability depends only on $\la_i$ but not on $i$, i.e. $\la_1,\la_2,\ldots$ is produced by a homogeneous descending Markov chain on the nonnegative integers. This was shown by Evans \cite{evans2002elementary}. From the explicit description of the transition probabilities of this Markov chain, the computation of $\E[||\det(A)||]$ is straightforward.

However, for a nonlinear polynomial $Z(A)$ we do not have such a simple Markovian structure of lifts modulo $p,p^2,\ldots$. The key to overcome this is a `linearization' result proven in \Cref{thm:compute_char_poly_det} using the moment method. It equates $\lim_{n \to \infty} \E[1/\#\Cok(Z(\tilde A))]$ to a mixed moment of sizes of cokernels of several matrices $B_1,\ldots,B_m$. Each matrix $B_i$ is marginally distributed according to the Haar measure on an appropriate extension of $\Q_p$, and hence its cokernel has the above Markov chain description. However, in general there are correlations for different $i$, and one has several correlated copies of this Markov chain which move together for a certain number of steps and then diverge. In the case where $Z$ has degree $2$, these reduce to expectations of certain functions of the $m\tth$ step of a single copy of this Markov chain, which are $m$-fold sums over integers $\la_1,\la_2,\ldots,\la_m$.

Fulman \cite{fulman-RR} had studied this same Markov chain earlier than \cite{evans2002elementary}, motivated by random matrix theory over finite fields. The work \cite{fulman-RR} gives an explicit diagonalization of this Markov chain, and uses this for a probabilistic proof of the Andrews-Gordon and Rogers-Ramanujan identities, the former of which are $m$-fold sums of a similar form. Using this diagonalization, we are able to reduce the $m$-fold sums to single sums where $m$ features as a parameter, seen in Theorems~\ref{thm:Zp_2pt_intro} and \ref{thm:quadratic_intro}. 

In \Cref{sec:markov_chains} we give background on the Markov chain and needed computational lemmas in our setting. In \Cref{sec: limiting correlation functions in low degree} we use this to prove Theorems~\ref{thm:Zp_2pt_intro} and~\ref{thm:quadratic_intro}, and use these to find series expression for the variance of the number of eigenvalues in $\Z_p$ (\Cref{thm:variance_of_roots_in_Z_p}) and prove the degree-$2$ part of \Cref{thm:roots_in_extensions_intro}.

The expectations we compute are in a sense deformed versions of the Andrews-Gordon identities, and we do not know whether the seeming algebraic miracles which allow their computations to simplify (see \Cref{rmk:miracle_in_diagonalization}) may hint at some deeper story. For $Z$ of degree $3$ and higher, we could compute many terms of the series expansions of $\lim_{n\rightarrow\infty}\E[1/\#\Cok(Z(\tilde A))]$ by this method, but understanding the structure of these series in general appears difficult and we have deferred it to future work. 

\Cref{sec:pair_repulsion} concerns the proof of the pair-repulsion estimate given in \Cref{thm:pair_correlation_estimate_intro}, and an additional estimate (\Cref{thm: refinement of two point correlation}) which is needed to confirm that repulsion still holds when $p=2$.

\Cref{sec:average_number_of_eigenvalues_in_high_degree_extensions} is devoted to proving the estimate of \Cref{thm:roots_in_extensions_intro} for extensions of higher degree. The method follows the one of \cite{caruso2022zeroes} for roots of random polynomials, using certain bounds of \cite{yun2013orbital} and results on expected determinants from \Cref{sec:cokernels} as input. The basic idea is that for a large portion of eigenvalue tuples (those which generate the ring of integers of the relevant field extension), the density function actually has a simple exact form, and we can estimate the relatively small contribution of the remaining eigenvalue tuples.

\Cref{sec:all_zp} concerns the proof of \Cref{thm: all eigenvalues in Z_p intro} and the analogous result for $\GL_n(\Z_p)$. These results concern a global rather than local property of the eigenvalues, namely whether they all lie in $\Z_p$, and hence cannot be accessed by the asymptotics of any given correlation function. Instead, we rely on an \emph{a priori} surprising connection to random polynomials with independent coefficients. It is essentially immediate from \Cref{thm: joint distribution} that the correlation function 
\begin{equation}
    \rho_{\Q_p^n}^{(n)}(x_1,\ldots,x_n) = \frac{(1-p^{-1}) \cdots (1-p^{-n})}{(1-p^{-1})^n} \prod_{1 \leq i < j \leq n} ||x_i-x_j||
\end{equation}
of the eigenvalues in $\Q_p^n$ has a simple Vandermonde form (this is how \Cref{thm:Zp_coulomb_intro} is proven). For a monic, degree $n$ polynomial with independent coefficients drawn from the additive Haar measure on $\Z_p$, the density function of roots on $\Q_p^n$ is the same Vandermonde without the constant in front (without the monic condition this was shown in \cite{caruso2022zeroes} and our proof is analogous). This match lets us port asymptotics from the setting of random polynomials, where the relevant ones were proven much earlier by Buhler-Goldstein-Moews-Rosenberg \cite{buhler2006probability}.

\textbf{Acknowledgments.} We thank Alexei Borodin, Jordan Ellenberg, Nathan Hayford, Chao Li, Will Sawin, Akshay Venkatesh, Joe Webster, and Zhiwei Yun for helpful conversations and comments. J.S. was supported by NSF grant DMS-2246576 and Simons Investigator grant 929852 (PI: Ivan Corwin).

%% file: Preliminaries.tex
\section{Preliminaries}\label{sec:preliminaries}

\subsection{Some general notation} The following notations will be used frequently. We always fix a prime $p$, and except in a few places indicated otherwise, $p$ does not have to be odd. We write $||\cdot||$ for the standard norm on $\bar{\Q}_p$ or its restriction to any subfield of $\Q_p$, and $\val(x) = -\log_p(||x||) \in \Q$ for the $p$-adic valuation. $K$ or $K_i$ will always denote a finite extension of $\Q_p$. $\bar{\Z}_p$ denotes the closed unit ball in $\bar{\Q}_p$. We write $\mathbf{P}(\cdots)$ for probability, which unless otherwise noted is with respect to the additive Haar probability measure on $\Mat_n(\Z_p)$. $P_A(x)$ denotes the characteristic polynomial of a matrix $A$. The $q$-Pochhammer symbol is
$$(a;q)_n := \prod_{i=0}^{n-1} (1-aq^i)$$
for $n \in \Z_{\geq 0} \cup \{\infty\}$. We often write RHS(x.yz) for the right-hand side of equation (x.yz) and similarly LHS(x.yz) for the left-hand side.

Many sections work in a specific setting detailed at the beginning of the section, so we recommend a reader who has skipped to a result in-text and finds unexplained notation to first glance at the beginning of the section.

\subsection{Étale algebras}

During our discussion in the Introduction, we frequently encountered a tuple of extensions $K_1,\ldots,K_m$ over $\Q_p$. As we will see, the forthcoming sections also involve the algebraic structure of such tuples, which makes it necessary to introduce the notion of an étale algebra. The readers may also refer to \cite[Chapter 5.6]{bourbaki1998commutative} or \cite[Chapter 2.1.2]{cohen2012advanced} for more information.

\begin{defi}
\cite[Definition 2.1.2]{cohen2012advanced} A finite étale algebra over $\Q_p$ is defined as a finite-degree algebra $E$ over $\Q_p$ without nonzero nilpotent elements.
\end{defi}

In the rest of the paper, we are only concerned with étale algebras that are finite-degree over $\Q_p$, and `étale algebra' will always implicitly include this condition.

\begin{prop}
\cite[Corollary 2.1.6]{cohen2012advanced} Every étale algebra $E/\Q_p$ is isomorphic to one of the form
$$E\cong K_1\times K_2\times\cdots\times K_m$$
where $K_1,\ldots,K_m$ are finite field extensions over $\Q_p$.
\end{prop}

Also, given $E=K_1\times\cdots\times K_m$, we denote
$$
E^{\new}:=\{\mathfrak{x}=(x_1,\ldots,x_m)\in E\mid\Q_p[\mathfrak{x}]=E\}
$$
as the \emph{new} elements in $E$, and $\mathcal{O}_E^{\new}:=\mathcal{O}_E\cap E^{\new}$. Given $\mfx\in E$, it is new if and only if its minimal polynomial (i.e., the monic polynomial $f\in\Q_p[x]$ of least degree such that $f(\mfx)=0$) has the same degree as the degree of $E/\Q_p$.

It is worth mentioning that we can also write the set $E^{\new}$ as
$$E^{\new}=\{\mathfrak{x}=(x_1,\ldots,x_m)\in K_1^{\new}\times\cdots\times K_m^{\new}\mid x_1,\ldots,x_m\text{ lie in distinct Galois orbits}\}.$$
In particular, $E^{\new}$ (resp. $\mc{O}_{E}^{\new}$) will be a proper subset of $K_1^{\new}\times\cdots\times K_m^{\new}$ (resp. $\mathcal{O}_{K_1}^{\new}\times\cdots\times \mathcal{O}_{K_m}^{\new}$) if and only if at least one pair of fields $K_i$ are isomorphic. $\mc{O}_{K_1} \times \cdots \times \mc{O}_{K_m}$ has a product Haar probability measure, with respect to which $\mc{O}_{E}^{\new}$ has full measure.

We endow $E\cong K_1\times K_2\times\cdots\times K_m$ with the absolute value $||\cdot||$ given by
$$||(x_1,\ldots,x_m)||=\max\{||x_1||,\ldots,||x_m||\}\quad(x_i\in K_i).$$
In this case, we let $\mathcal{O}_E=\mathcal{O}_{K_1}\times \mathcal{O}_{K_2}\times\cdots\times \mathcal{O}_{K_m}$ be the subring of $E$ consisting of elements of norm at most $1$. We define the norm map $\Nm_{E/\Q_p}:E\rightarrow \Q_p$ by
$$\Nm_{E/\Q_p}(x_1,\ldots,x_m)=\Nm_{K_1/\Q_p}(x_1)\cdots\Nm_{K_m/\Q_p}(x_m).$$
$\Nm_{E/\Q_p}(x_1,\ldots,x_m)$ is also equal to the determinant of the scalar multiplication $(x_1,\ldots,x_m)$, viewed as a linear map over the $\Q_p$-vector space $E$.

We also define the discriminant of a finite étale $\Q_p$-algebra $E=K_1\times\cdots\times K_m$ as $\Disc_{E/\Q_p}=\Disc_{K_1/\Q_p}\cdots \Disc_{K_m/\Q_p}$. For all $\mfx=(x_1,\ldots,x_m)\in E$, define $\Disc(\mfx):=0$ when $\mfx$ is not in $E^{\new}$, and 
$$\Disc(\mfx):=\prod_{(i,l_i)\ne(j,l_j)}(\sigma_{i,l_i}(x_i)-\sigma_{j,l_j}(x_j))=\prod_{i=1}^m\Nm_{K_i/\Q_p}Z'(x_i)\in\Q_p$$
when $\mfx\in E^{\new}$. Here, for all $1\le i\le m$, $\sigma_{i,1}(x_i)=x_i,\ldots,\sigma_{i,r_i}(x_i)$ are the Galois conjugates of $x_i$, and $Z\in\Q_p[x]$ is the minimal polynomial of $(x_1,\ldots,x_m)$.

For a finite étale algebra $E=K_1\times\cdots\times K_m$ over $\Q_p$, we denote by $\Aut_{\Q_p}(E)$ the group of $\Q_p$-algebra automorphisms of $E$, and $m_E:=m$ as the number of field components of $E$.

\begin{prop}\label{prop: aut_e and generators}
Let $E$ be a finite étale algebra over $\Q_p$, and $\mfx\in E^{\new}$, so that $E=\Q_p[\mfx]$. Then, every element $\tau\in\Aut_{\Q_p}(E)$ can be written in the form
\begin{align}\label{eq: automorphism as moving generator}
\begin{split}
\tau:E&\rightarrow E\\
f(\mfx)&\mapsto f(\mathfrak{y}),\quad f\in\Q_p[x].\\
\end{split}
\end{align}
where $\mathfrak{y}\in E^{\new}$ has the same minimal polynomial as $\mfx$. Consequently, we have
$$\#\Aut_{\Q_p}(E)=\#\{\mathfrak{y}\in E^{\new}: \mfx \text{ and }\mathfrak{y}\text{ have the same minimal polynomial in }\Q_p\}.$$
\end{prop}

\begin{proof}
On the one hand, for all $\tau\in\Aut_{\Q_p}(E)$, $\tau(\mfx)$ must have the same minimal polynomial as $\mfx$. On the other hand, every mapping $\tau$ given in \eqref{eq: automorphism as moving generator} is indeed an isomorphism.
\end{proof}

We will often simply write $\Aut(E)$ when the subscript is clear from the text. The definitions of 
$$||\cdot||,\Nm_{E/\Q_p},\Disc_{E/\Q_p},\#\Aut(E)$$
do not depend on the chosen identification $E\cong K_1\times\cdots\times K_m$. 

From now on, when $E=K_1\times\cdots\times K_m$, we will also use the notation $\rho_{E}^{(n)}(x_1,\ldots,x_m)$ or $\rho_{K_1\times\cdots\times K_m}^{(n)}(x_1,\ldots,x_m)$ for $\rho_{K_1,\ldots,K_m}^{(n)}(x_1,\ldots,x_m)$, and $\rho_E^{(\infty)}(x_1,\ldots,x_m)$ or $\rho_{K_1\times\cdots\times K_m}^{(\infty)}(x_1,\ldots,x_m)$ for $\rho_{K_1,\ldots, K_m}^{(\infty)}(x_1,\ldots,x_m)$.

\subsection{The ring {$\bar{\Z}_p$}, lifted subspaces, and the resultant}\label{subsec:zpbar}

Let $\bar\Z_p=\{y\in\bar\Q_p\mid ||y||\le 1\}$ be the ring of algebraic integers over $\Z_p$, which is the space that the eigenvalues of $A\in\Mat_n(\Z_p)$ fall in. Denote by $\{F_1=x,F_2,\ldots\}$ the set of monic irreducible polynomials over $\F_p$, and denote $d_i:=\deg F_i$. Then $\bar\Z_p$ is the disjoint union 
\begin{equation}\label{eq: decomposition of algebraic closure}
\bar\Z_p=\bigsqcup_{i=1}^\infty \mathcal{U}_i,
\end{equation}
where for all $i\ge 1$, $\mathcal{U}_i$, the \emph{lifted subspace} over $F_i$, is defined by
\begin{equation}\label{eq: lifted subspace}
\mathcal{U}_i:=\{y\in\bar\Z_p\mid\text{the residue of the minimal polynomial of $y$ is a power of $F_i$}\}.
\end{equation}
For instance, the irreducible polynomial $F_1=x\in\F_p[x]$ corresponds to the open unit disk
$\mathcal{U}_1=\{y\in\bar\Q_p\mid ||y||<1\}$. It is clear that $\mathcal{U}_1$ is the unique maximal ideal in $\bar\Z_p$.  

\begin{lemma}\label{lem: Hensel}
(Hensel's lemma, strong version) Let $Z(x)\in\Z_p[x]$ be monic and $\overline Z(x)\in\F_p[x]$ be the residue of $Z(x)$. Suppose that there exists a decomposition
$$\overline Z(x)=\overline Z_1(x)\overline Z_2(x),$$
where $\overline Z_1(x),\overline Z_2(x)\in\F_p[x]$ are relatively prime. Then $Z(x)$ admits a factorization
$$Z(x)=Z_1(x)Z_2(x),$$
such that $\overline Z_1(x)$ is the residue of $Z_1(x)$, and $\overline Z_2(x)$ is the residue of $Z_2(x)$.
\end{lemma}

\begin{proof}
See \cite[Theorem II.4.6]{neukirch2013algebraic}.
\end{proof}

Now, suppose 
$$Z_1(x)=a_{r_1}x^{r_1}+a_{r_1-1}x^{r_1-1}+\cdots+a_0,\quad Z_2(x)=b_{r_2}x^{r_2}+b_{r_2-1}x^{r_2-1}+\cdots+b_0\in\Z_p[x]$$
are nonzero polynomials of degree $r_1,r_2$ respectively. Denote by $x_{i,1},\ldots,x_{i,r_i}$ for the roots of $Z_i$, where $i=1,2$ (note that the roots may not be in $\Z_p$). Then the \emph{resultant} of the two polynomials $Z_1(x),Z_2(x)$, denoted by $\Res(Z_1,Z_2)$, is defined as
$$\Res(Z_1,Z_2)=a_{r_1}^{r_2}b_{r_2}^{r_1}\prod_{l_1=1}^{r_1}\prod_{l_2=1}^{r_2}(x_{1,l_1}-x_{2,l_2})\in\Z_p.$$
We furthermore define $\Res(Z_1,0)=\Res(0,Z_2)=0$. The resultant describes distances between the roots of two polynomials. It satisfies some basic properties which we record without proof.

\begin{prop}\label{prop: properties of resultant}
Let $Z_1,Z_2,Z_3\in\Z_p[x]$. Then we have the following properties for the resultant:
\begin{enumerate}
\item $\Res(Z_1,Z_2)=\Res(Z_2,Z_1)$.
\item $\Res(Z_1,Z_2)=0$ if and only if they share a root.
\item $\Res(Z_1Z_2,Z_3)=\Res(Z_1,Z_3)\Res(Z_2,Z_3)$.
\item If $Z_3$ is monic, then $\Res(Z_1,Z_3)=\Res(Z_1-Z_2Z_3,Z_3)$.
\end{enumerate}
\end{prop}

\begin{prop}\label{prop: resultant positive power}
Let $Z_1\in\Z_p[x]$ be monic, and $Z_2\in\Z_p[x]$ with non-zero residue. Furthermore, assume that the residue of $Z_1$ is a power of $F_i$, which has degree $d_i$. Then we have
$$||\Res(Z_1,Z_2)||=p^{-kd_i},$$
where $k\ge 0$ is an integer. In particular, $k=0$ if and only if $F_i \nmid Z_2 \pmod{p}$.
\end{prop}

\begin{proof}
We first clarify the case where the residue of $Z_2$ is not divisible by $F_i$. In this case, applying the Euclidean algorithm over the residue field, there exists $f_1(x),f_2(x)\in\Z_p[x]$ such that
$$f_1(x)Z_1(x)+f_2Z_2(x)\in1+p\Z_p[x].$$
Therefore, for any root $x_1$ of $Z_1$, we have $||f_2(x_1)Z_2(x_1)||=1$. Since $||f_2(x_1)||,||Z_2(x_1)||\le 1$, we must have $||f_2(x_1)||=||Z_2(x_1)||=1$. Therefore, we have
$$||\Res(Z_1,Z_2)||=\prod_{x_1}||Z_2(x_1)||=1,$$
where the product runs through all roots $x_1$ of $Z_1$.

Next, suppose that the residue of $Z_2$ is divisible by $F_i$. We only need to verify that the absolute value of the resultant $||\Res(Z_1,Z_2)||$ is a positive integer power of $p^{-d_i}$. The proof is to apply induction over $\deg Z_1$. If $\deg Z_1=d_i$, then we can write $Z_2=Z_3Z_1+p^kZ_4$, where $Z_3\in\Z_p[x]$, $k\ge 1$, and $Z_4\in\Z_p[x]$ is of degree $\le d_i-1$ with non-zero residue. In this case, applying \Cref{prop: properties of resultant} gives
\begin{align}
\begin{split}
||\Res(Z_1,Z_2)||&=||\Res(Z_1,Z_3Z_1+p^kZ_4)||\\
&=||\Res(Z_1,p^kZ_4)||\\
&=p^{-kd_i}||\Res(Z_1,Z_4)||\\
&=p^{-kd_i}.
\end{split}
\end{align}
Here, the last line holds because the residue of $Z_4$ is not divisible by $F_i$.

Now we suppose the claim is proven for $\deg Z_1 \leq (n-1)d_i$, and prove it for $\deg Z_1 = nd_i$. By \Cref{lem: Hensel}, we have a decomposition $Z_2=Z_3Z_4$, where $Z_3\in\Z_p[x]$ is monic with residue being a power of $F_i$, and the residue of $Z_4\in\Z_p[x]$ is not divisible by $F_i$. We have
$$||\Res(Z_1,Z_2)||=||\Res(Z_1,Z_3)||\cdot||\Res(Z_1,Z_4)||=||\Res(Z_1,Z_3)||.$$
If $\deg Z_3<\deg Z_1=nd_i$, then by relabeling $Z_1$ as $Z_3$ and vice versa, we are done by the induction hypothesis. If $\deg Z_3\ge\deg Z_1$, then we can write $Z_3=Z_4Z_1+p^kZ_5$, where $Z_4\in\Z_p[x]$, $k\ge 1$, and $Z_5\in\Z_p[x]$ of degree $\le nd_i-1$ with nonzero residue. Applying \Cref{prop: properties of resultant} and we have
\begin{align}
\begin{split}
||\Res(Z_1,Z_3)||&=||\Res(Z_1,Z_4Z_1+p^kZ_5)||\\
&=||\Res(Z_1,p^kZ_5)||\\
&=p^{-nkd_i}||\Res(Z_1,Z_5)||.\\
\end{split}
\end{align}
By \Cref{lem: Hensel} again, we have $Z_5=Z_6Z_7$, where $Z_6\in\Z_p[x]$ is monic with residue being a power of $F_i$, and the residue of $Z_7\in\Z_p[x]$ is not divisible by $F_i$. Therefore,
$$||\Res(Z_1,Z_3)||=p^{-nkd_i}||\Res(Z_1,Z_5)||=p^{-nkd_i}||\Res(Z_1,Z_6)||,$$
which must be a positive power of $p^{-d_i}$ by the induction hypothesis over $\deg Z_6$.
\end{proof}

The following proposition intuitively claims that two  elements in $\bar\Z_p$ are close if and only if they fall into the same lifted subspace. It is a direct corollary of \Cref{prop: resultant positive power}, so we omit the proof here.
\begin{prop}\label{prop: distance between lifted subspaces}
Let $x_1\in \mathcal{U}_i,x_2\in \mathcal{U}_j$, and $Z_1(x),Z_2(x)$ denote the minimal polynomial of $x_1,x_2$ respectively. Then we have the following:
\begin{enumerate}
\item If $i\ne j$, then $||\Res(Z_1,Z_2)||=1$.
\item If $i=j$, then $||\Res(Z_1,Z_2)||=p^{-kd_i}$, where $k\ge 1$ is an integer.
\end{enumerate}
\end{prop}

\begin{prop}\label{prop: distances between roots of Z_i}
Suppose we have a monic $Z(x)\in\Z_p[x]$ whose residue $F(x)\in\F_p[x]$ is irreducible, so that $Z(x)$ is also irreducible. Let $d=\deg F$, and let
$$\sigma_1(x),\ldots,\sigma_d(x)$$
be the roots of $Z$. Here $x$ is a root of $Z$, and $\sigma_1=\id,\ldots,\sigma_d\in\Gal(\bar\Q_p/\Q_p)$. Then, for all $1\le j<k\le d$, we have $||\sigma_j(x)-\sigma_k(x)||=1$. 
\end{prop}

\begin{proof}
We only need to show that the discriminant of $Z(x)$ is in $\Z_p^\times$. This is equivalent to the statement that the discriminant of $F(x)$ is nonzero. Since $F(x)$ does not have repeated roots, we are done.
\end{proof}

\begin{prop}\label{thm: lifted space as cosets of unit disk}
Let $F_i\in\F_p[x]$ be monic and irreducible, and  $\mathcal{U}_i$ be the lifted subspace that corresponds to $F_i$. Let $Z_i\in\Z_p[x]$ be monic, whose residue is $F_i$. Let $d_i=\deg F_i$, and let
$$\sigma_{i,1}(x_i),\ldots,\sigma_{i,d_i}(x_i)$$
be the roots of $Z_i$. Then, $\mathcal{U}_i$ can be written as the disjoint union
\begin{equation}\label{eq:lifted_space_coset_union}
    \mathcal{U}_i=\bigsqcup_{j=1}^{d_i} (\sigma_{i,j}(x_i)+\mathcal{U}_1).
\end{equation}
\end{prop}

\begin{proof} 
Recalling the strong triangle inequality and applying \Cref{prop: distances between roots of Z_i}, we know that union in \eqref{eq:lifted_space_coset_union} is indeed disjoint. 

$(\subset):$ Suppose $x_i'\in \mathcal{U}_i$. Denote $Z_{i'}$ as the minimal polynomial of $x_i'$. Applying \Cref{prop: distance between lifted subspaces}, we have $||\Res(Z_i,Z_{i'})||<1$, and $\Res(Z_i,Z_{i'})$ is a product of Galois conjugates of $Z_i(x_i')$, so this implies $||Z_i(x_i')||=\prod_{j=1}^{d_i}||\sigma_{i,j}(x_i)-x_i'||<1$. Therefore, there exists $1\le j\le d_i$ such that $||\sigma_{i,j}(x_i)-x_i'||<1$, and $x_i'\in\sigma_{i,j}(x_i)+\mathcal{U}_1$.

$(\supset):$ Take $y\in\mathcal{U}_1$, so that $||y||<1$. Denote by 
$$\sigma_1(y),\ldots,\sigma_r(y)$$
as all the Galois conjugates of $y$. Then the residue of the polynomial 
$$\prod_{1\le j\le d_i}\prod_{1\le l\le r}(x-\sigma_{i,j}(x_i)-\sigma_l(y))=\prod_{1\le l\le r}Z_i(x-\sigma_l(y))$$
is equal to $\overline Z_i(x)^r=F_i^r$. Therefore, we have $\sigma_{i,j}(x_i)+y\in\mathcal{U}_i$ for all $1\le j\le d_i$ and $y\in\mathcal{U}_1$.
\end{proof}

\begin{rmk}
The way to pick the polynomial $Z_i$ and its root $x_i$ is not unique, but it always leads to the same decomposition of disjoint unions as in the above. Therefore, from now on, we will write $\mathcal{U}_{i,j}$ for the set $\sigma_{i,j}(x_i)+\mathcal{U}_1$ in \eqref{eq:lifted_space_coset_union}, where $x_i$ will be chosen implicitly. Since the Galois group acts transitively on these pieces, we deduce that the number of eigenvalues that fall inside each $\mathcal{U}_{i,j},1\le j\le d_i$ must be equal. In fact, any other (Galois-invariant) eigenvalue statistic is the same on all discs $\mc{U}_{i,j}$, because the eigenvalues in any disc $\mc{U}_{i,j}$ are Galois conjugate to the ones in $\mc{U}_{i,1}$.
\end{rmk}

%% file: joint_distribution.tex
\section{Explicit joint CDF of characteristic polynomial values at multiple points}
\label{sec:char_poly}

In this section, we will work with the following setting. Let

\begin{enumerate}
\item $E=K_1\times K_2\times\cdots\times K_m$ be a finite étale algebra over $\Q_p$, which has dimension $r$ as a $\Q_p$-vector space. Here for all $1\le i\le m$, $K_i/\Q_p$ is a finite field extension of degree $r_i$, where $r_1+\cdots+r_m=r$. 
\item $\mathfrak{x}=(x_1,\ldots,x_m)\in \mathcal{O}_E^{\new}$, so that $x_i\in \mathcal{O}_{K_i}^{\new}$ for all $1\le i\le m$, and $E=\Q_p[\mfx]$.
\item $\pi_i$ be a uniformizer of $\mathcal{O}_{K_i}$, so that $\pi_i \mathcal{O}_{K_i}$ is the unique maximal ideal of $\mc{O}_{K_i}$.
\item $\L_E=\prod_{i=1}^m\pi_i^\Z\subset E^\times$ be the free $\Z$-module complementary to $\mathcal{O}_E^\times=\mathcal{O}_{K_1}^\times\times\cdots\times \mathcal{O}_{K_m}^\times$.
\item $\L_{K_i}=\pi_i^\Z\subset K_i^\times$ be the free $\Z$-module complementary to $\mc{O}_{K_i}^\times$, so that $\L_E=\prod_{i=1}^m\L_{K_i}$.
\item $L_i$ be the normal closure of $K_i/\Q_p$.
\item $\sigma_{i,1}(x_i)=x_i,\ldots,\sigma_{i,r_i}(x_i)$ be conjugates of $x_i$ in $L_i$, where $\sigma_{i,1}=\id,\ldots,\sigma_{i,r_i}\in\Gal(\bar\Q_p/\Q_p)$.
\item $Z_i(x)$ be the minimal polynomial of $x_i$, which is of degree $r_i$. Since the elements $x_1,\ldots,x_m$ lie in different orbits under Galois automorphisms, the polynomials $Z_1,\ldots,Z_m$ are distinct. 
\item $Z(x)=\prod_{i=1}^m Z_i(x)$, which is monic of degree $r=r_1+\cdots+r_m$.
\item $e_i$ be the ramification index of the extension $K_i/\Q_p$, so $\#(\mathcal{O}_{K_i}:\pi_i \mathcal{O}_{K_i})=p^{r_i/e_i}$.
\item $\mathbf{P}(\cdots)$ be the (additive) Haar probability measure on $\Mat_r(\Z_p)$.
\end{enumerate}

This section is devoted to the proof of the following result. Informally, it computes the probability that the characteristic polynomial of an $r \times r$ matrix takes very small values at a given collection of $r$ points. It turns out that this packages most of the work needed to prove \Cref{thm: joint distribution}, because the eigenvalues must lie close to these $r$ points.

\begin{thm}\label{thm: points on variety_main theorem}
Suppose $A\in\Mat_r(\Z_p)$ is Haar-distributed. Then for all tuples $s=(s_1,\ldots,s_m)\in(\frac{1}{e_1}\N,\ldots,\frac{1}{e_m}\N)$ such that $s_i>r\val\Disc(\mathfrak{x})$ for all $1\le i\le m$, we have
\begin{equation}\label{eq: points on the variety}
\prod_{i=1}^m p^{s_ir_i}\cdot\mathbf{P}(\val P_{A}(x_i)\ge s_i,\forall 1\le i\le m)=\frac{(1-p^{-1})\cdots(1-p^{-r})}{||\Delta_\sigma(x_1,\ldots,x_m)||}\prod_{i=1}^m \Den(x_i),
\end{equation}
where $\Delta_\sigma$ and $\Den$ are as in Definitions~\ref{def:delta_sigma} and~\ref{defi: Den and distance}.
\end{thm}

\subsection{High-level summary of the proof}

In \Cref{thm: points on variety_main theorem}, we wish to compute the measure of the set
\begin{equation}\label{eq:A_s_set}
    \{A\in\Mat_r(\Z_p)\mid\val(P_{A}(x_i))\ge s_i,\quad\forall 1\le i\le m\}.
\end{equation}
Recall the overview explained in \Cref{subsec:proof_methods}; let us now be more specific. Letting
\begin{equation}
\mathcal{X}=\{A\in\Mat_r(\Z_p)\mid P_{A}=Z\},
\end{equation}
one may view \eqref{eq:A_s_set} as the set of matrices `close' to $\mathcal{X}$, where `close' is defined purely in terms of the characteristic polynomial. First, we prove in \Cref{thm: polynomial expression of matrix near variety} that the matrices in \eqref{eq:A_s_set} can be uniquely expressed by $A=f(A_0)$ with $A_0\in \mathcal{X}$, and $f\in\Z_p[x]$ is a polynomial which `only moves the matrix a little bit'. More precisely, $f$ belongs to the following set.

\begin{defi}
    \label{def:f_set}
    Let $s=(s_1,\ldots,s_m)\in(\frac{1}{e_1}\N,\ldots,\frac{1}{e_m}\N)$ such that $s_i>r\val\Disc(\mathfrak{x})$ for all $1\le i\le m$. Denote
    \begin{equation}
        F_{Z,s} := \{f \in \Z_p[x] \mid \deg f \leq r-1, \val(f(x_i) - x_i) + \val Z'(x_i) \geq s_i \text{ for all }1 \leq i \leq m\}.
    \end{equation}
    We will abbreviate $F_{Z,s}$ as $F_Z$ when $s_1,\ldots,s_m$ are clear from context.
\end{defi}

\begin{prop}\label{thm: polynomial expression of matrix near variety}
Let $s=(s_1,\ldots,s_m)\in(\frac{1}{e_1}\N,\ldots,\frac{1}{e_m}\N)$ such that $s_i>r\val\Disc(\mathfrak{x})$ for all $1\le i\le m$, and let $F_Z = F_{Z,s}$ be as in \Cref{def:f_set}. Then the map
\begin{align}\label{eq:f_forward_map}
    \begin{split}
        \{(f,A_0)\mid A_0\in \mc{X},f\in F_Z\} & \to \{A\in\Mat_r(\Z_p)\mid\val(P_{A}(x_i))\ge s_i,\quad\forall 1\le i\le m\} \\ 
        (f,A_0) & \mapsto f(A_0)
    \end{split}
\end{align}
is a bijection.
\end{prop}

Because $f$ is polynomial and hence 
\begin{equation}
    f(TA_0T^{-1}) = Tf(A_0)T^{-1},
\end{equation}
\Cref{thm: polynomial expression of matrix near variety} makes it natural to study the orbits in $\mathcal{X}$ under conjugation. If one allows conjugation by all of $\GL_r(\Q_p)$, by rational canonical form there is only one orbit, with a representative given by the companion matrix $C(Z)$ of $Z$ (see \Cref{def:companion_matrix}). However, this orbit $\{T C(Z) T^{-1}: T \in \GL_r(\Q_p)\}$ is not contained in $\Mat_r(\Z_p)$, so it is bigger than our set $\mathcal{X}$. 

Hence it is more natural to consider orbits under conjugation by the smaller group $\GL_r(\Z_p)$, which will naturally lie inside $\Mat_r(\Z_p)$. The structure of these orbits is more complicated. However, it turns out (\Cref{thm:matrix_disjoint_union}) that the natural way to parametrize these orbits is by the set of modules
\begin{equation}\label{eq:def_ModZp}
    \Mod_{\Z_p[\mathfrak{x}]}=\{M\subset E\mid M \text{ is a $\Z_p[\mathfrak{x}]$-module, also a $\Z_p$-lattice of rank $r$}\},
\end{equation}
or rather by the quotient $E^\times \backslash \Mod_{\Z_p[\mfx]}$ under scalar multiplication, which is analogous to the class group. This parametrization has already been understood in the literature on orbital integrals
$$\int_{\GL_r(\Q_p)} \bbone_{gC(Z)g^{-1} \in \Mat_r(\Z_p)} dg$$
such as \cite{yun2013orbital}, but we give a self-contained presentation of what we need (though after proving the parametrization, we use without proof certain facts about the aforementioned class group analogues shown in \cite{yun2013orbital}).

This makes the set \eqref{eq:A_s_set} a disjoint union of sets
\begin{equation}\label{eq:T_f_set}
    \{Tf(A_{\mcc})T^{-1}\mid T \in \GL_r(\Z_p),f\in F_Z\}
\end{equation}
where $A_{\mcc}$ are certain representative matrices, and we are reduced to computing the measure of these sets.

The sets in \eqref{eq:T_f_set} may be further decomposed into a finite number of translates (by conjugation) of sets
\begin{equation}\label{eq:diff_element_set}
    \{(I_r+\Delta T)f(A_{\mcc})(I_r+\Delta T)^{-1} \mid \Delta T \in p^{\mc{M}}\Mat_r(\Z_p), f \in F_Z\},
\end{equation}
all of which have the same measure. This is done in \Cref{thm:decompose_TT_sets} and \Cref{thm:how_many_sets_M}. 

The matrices $\Delta T$ in \eqref{eq:diff_element_set} are small, and the set in \eqref{eq:diff_element_set} may be described explicitly in a linear way using the image of the Lie bracket, see \Cref{lem: same image of Lie} and \Cref{lem: local conjugate and lie bracket}. From this description, we compute its measure in \Cref{prop: volume of differential element}.

\subsection{Proof of \Cref{thm: polynomial expression of matrix near variety}}

We begin proving the results needed for \Cref{thm: polynomial expression of matrix near variety}. The following lemma generalizes Krasner's lemma \cite{krasner1946theorie} to the étale algebra case, and essentially says that any sequence $(\tx_1,\ldots,\tx_m)$ of elements close to $(x_1,\ldots,x_m)$ has metric and algebraic properties similar to $(x_1,\ldots,x_m)$.

\begin{lemma}\label{lem: close elements}
Suppose we have a sequence of elements
$$\tilde{\mathfrak{x}}=(\tilde x_1,\ldots,\tilde x_m),\quad \tilde x_1,\ldots,\tilde x_m\in\bar\Q_p,$$ 
such that $||x_i-\tilde x_i||<||\Disc(\mathfrak{x})||$ for all $1\le i\le m$. Then 
\begin{enumerate}
\item\label{item: same distance} For all pairs of distinct automorphisms $\sigma_{i,l_i}, \sigma_{j,l_j}$ with $(i,l_i) \neq (j,l_j)$, we have
$$||\sigma_{i,l_i}(\tilde x_i)-\sigma_{j,l_j}(\tilde x_j)||=||\sigma_{i,l_i}(x_i)-\sigma_{j,l_j}(\tilde x_j)||=||\sigma_{i,l_i}(x_i)-\sigma_{j,l_j}(x_j)||.$$
\item\label{item: distinct conjugates} The elements 
$$\sigma_{1,1}(\tilde x_1),\ldots,\sigma_{1,r_1}(\tilde x_1),\ldots,\sigma_{m,1}(\tilde x_m),\ldots,\sigma_{m,r_m}(\tilde x_m)$$
are all distinct.
\item \label{item: etale algebra contain} For all $1\le i\le m$, we have $\Q_p[\tilde x_i]\supset\Q_p[x_i]$. Also, we have $\Q_p[\tilde{\mathfrak{x}}]\supset E=\Q_p[\mathfrak{x}]$. 
\end{enumerate} 
\end{lemma}

\begin{proof}
First, notice that for every pair $(i,l_i) \neq (j,l_j)$,
\begin{equation}\label{eq:sigma_disc}
    ||x_j - \tilde x_j|| < ||\Disc(\mfx)|| \leq ||\sigma_{i,l_i}(x_i)-\sigma_{j,l_j}(x_j)||,
\end{equation}
where the first inequality is by hypothesis and the second is because $||\Disc(\mathfrak{x})||$ is a product of factors $\leq 1$ of which one is $||\sigma_{i,l_i}(x_i)-\sigma_{j,l_j}(x_j)||$. In other words, $x_j$ is closer to $\tilde x_j$ than it is to any of its Galois conjugates or the Galois conjugates of the other $x_i, i \neq j$. Because a Galois automorphism $\sigma_{j,l_j}$ does not affect the distance on the left-hand side of \eqref{eq:sigma_disc}, we have 
\begin{equation}
    ||\sigma_{j,l_j}(x_j)-\sigma_{j,l_j}(\tilde x_j)|| < ||\sigma_{i,l_i}(x_i)-\sigma_{j,l_j}(x_j)||.
\end{equation}
Because $||a+b|| = ||a||$ when $||a|| > ||b||$, this implies 
\begin{align}
    \begin{split} \label{eq:tilde_swap_1}
     ||\sigma_{i,l_i}(x_i)-\sigma_{j,l_j}(\tilde x_j)||&=||\sigma_{i,l_i}(x_i)-\sigma_{j,l_j}(x_j)+\sigma_{j,l_j}(x_j)-\sigma_{j,l_j}(\tilde x_j)||\\
     &=||\sigma_{i,l_i}(x_i)-\sigma_{j,l_j}(x_j)||,   
    \end{split}
\end{align}
which is strictly greater than $||x_i-\tilde x_i||=||\sigma_{i,l_i}(x_i)-\sigma_{i,l_i}(\tilde x_i)||$ by \eqref{eq:sigma_disc}. Hence we also have
\begin{align}
    \begin{split}\label{eq:tilde_swap_2}
        ||\sigma_{i,l_i}(\tilde x_i)-\sigma_{j,l_j}(\tilde x_j)||&=||\sigma_{i,l_i}(\tilde x_i)-\sigma_{i,l_i}(x_i)+\sigma_{i,l_i}(x_i)-\sigma_{j,l_j}(\tilde x_j)||\\ 
        &=||\sigma_{i,l_i}(x_i)-\sigma_{j,l_j}(\tilde x_j)||.
    \end{split}
\end{align} 
Combining \eqref{eq:tilde_swap_1} with \eqref{eq:tilde_swap_2} proves \eqref{item: same distance}.

Next, for all $(i,l_i)\ne (j,l_j)$, we have $||\sigma_{i,l_i}(\tilde x_i)-\sigma_{j,l_j}(\tilde x_j)||=||\sigma_{i,l_i}(x_i)-\sigma_{j,l_j}(x_j)||>0$. This proves \eqref{item: distinct conjugates}.

In the end, by \eqref{item: distinct conjugates}, for all $\sigma\in\Gal(\bar\Q_p/\Q_p)$ that $\sigma(x_i)\ne x_i$, we have $\sigma(\tilde x_i)\ne \tilde x_i$. Taking the contrapositive, for every $\sigma\in\Gal(\bar{\Q}_p/\Q_p)$, $\sigma(\tilde x_i)=\tilde x_i$ implies $\sigma(x_i)=x_i$. Hence $\Gal(\bar{\Q}_p/\Q_p[\tilde x_i])\subseteq\Gal(\bar{\Q}_p/\Q_p[x_i])$, which implies $\Q_p[\tilde x_i]\supset\Q_p[x_i]$.

Furthermore, for every $\sigma\in\Gal(\bar{\Q}_p/\Q_p)$ and $1\le i\le m$, we have $\sigma(x_i)=\sigma_{i,l_i}(x_i)$ for some $1\le l_i\le r_i$. Therefore, for all $j\ne i$,
\begin{align}\label{eq:some_more_sigmas}
\begin{split}
||\sigma(\tilde x_i)-\sigma_{i,l_i}(x_i)||&=||\sigma(\tilde x_i)-\sigma(x_i)||\\
&=||\tilde x_i-x_i||.
\end{split}
\end{align}
By \eqref{eq:sigma_disc} with $i,j$ reversed and $\sigma_{j,l_j}=\id$, we have
\begin{align}
    \begin{split}
       \text{RHS\eqref{eq:some_more_sigmas}} &<||x_j-\sigma_{i,l_i}(x_i)||\\
&=||x_j-\sigma_{i,l_i}(x_i)+\tilde x_j- x_j||\\
&=||\tilde x_j-\sigma_{i,l_i}(x_i)||,
    \end{split}
\end{align}
where for the last line we use the strong triangle inequality and the fact that $||\tilde x_j - x_j|| < ||x_j - \sigma_{i,l_i}(x_i)||$, which follows again by \eqref{eq:sigma_disc}.

Hence, $\sigma(\tilde x_i)$ and $\tilde x_j$ have different distances to $\sigma_{i,l_i}(x_i)$, so in particular $\sigma(\tilde x_i)\ne\tilde x_j$. Since such assertion holds for all $i\ne j$ and $\sigma\in\Gal(\bar\Q_p/\Q_p)$, the minimal polynomials of $\tilde x_1,\ldots,\tilde x_m$ are distinct, and
$$\Q_p[\tilde{\mathfrak{x}}]=\Q_p[\tilde x_1]\times\cdots\times\Q_p[\tilde x_m]\supset\Q_p[ x_1]\times\cdots\times\Q_p[x_m]=E=\Q_p[\mathfrak{x}],$$
proving \eqref{item: etale algebra contain}.
\end{proof}

The goal of the following lemma is to show that if a polynomial $f$ of degree less than $r-1$ takes small values on every $x_i,1\le i\le m$, then all the coefficients of $f$ are small.

\begin{cor}\label{cor: bound of coefficient}
Let $f(x)\in\Q_p[x]$ be of degree $\le r-1$, $u\in\R$, and for all $1\le i\le m$, $\val f(x_i)-\val Z'(x_i)\ge u$. Then the coefficients of $f(x)$ all have valuations $\ge u$.
\end{cor}

\begin{proof}
Due to the appearance of $u$ in both the hypothesis and the conclusion, if we show the claim for $p^k f$ for some $k$ then it holds for $f$ as well. Hence we may assume without loss of generality that $u>r\val\Disc(\mathfrak{x})$. Consider the roots of $\tilde Z(x)=Z(x)+f(x)$, which is a monic polynomial of degree $r$. For any fixed $1\le i\le m$, the $r$ roots of $\tilde Z$ give $r$ distances to $x_i$. The product of these distances is exactly 
$$||\tilde Z(x_i)||=||f(x_i)||\le ||Z'(x_i)||p^{-u}\le p^{-u}<||\Disc(\mathfrak{x})||^r.$$ Therefore, we can find one root $\tilde x_i$ of $\tilde Z$ that satisfies $||\tilde x_i-x_i||<||\Disc(\mathfrak{x})||$. Since such $\tilde x_i$ exists for all $1\le i\le m$, we can apply \eqref{item: distinct conjugates} of \Cref{lem: close elements} to deduce that the roots $$\sigma_{1,1}(\tilde x_1),\ldots,\sigma_{1,r_1}(\tilde x_1),\ldots,\sigma_{m,1}(\tilde x_m),\ldots,\sigma_{m,r_m}(\tilde x_m)$$
are all distinct, thus must be all roots of the polynomial $\tilde Z$. Also, for all $1\le i\le m$, we have 
\begin{align}
\begin{split}
\val(x_i-\tilde x_i)&=\val \tilde Z(x_i)-\sum_{(j,l_j)\ne (i,1)}\val(x_i-\sigma_{j,l_j}(\tilde x_j))\\
&=\val \tilde Z(x_i)-\sum_{(j,l_j)\ne (i,1)}\val(x_i-\sigma_{j,l_j}(x_j))\\
&=\val f(x_i)-\val Z'(x_i)\\
&\ge u.
\end{split}
\end{align}
Here the second equality comes from \eqref{item: same distance} of \Cref{lem: close elements}, and the third since $\tZ(x_i) = f(x_i)$ because $x_i$ is a root of $Z$. Hence $\val(\sigma_{i,l_i}(x_i-\tilde x_i))\ge u$ for all $1\le l_i\le r_i$, and so when we expand the expression
$$f(x)=\tilde Z(x)-Z(x)=\prod_{i=1}^m\prod_{l_i=1}^{r_i}(x-\sigma_{i,l_i}(x_i)+\sigma_{i,l_i}(x_i-\tilde x_i))-\prod_{i=1}^m\prod_{l_i=1}^{r_i}(x-\sigma_{i,l_i}(x_i)),$$
all the coefficients of this polynomial must have valuations $\ge u$.
\end{proof}

The below result says intuitively that if $Z_0$ is a polynomial that takes small values on the roots of $Z$, then the roots of $Z$ and $Z_0$ may naturally be paired up into close together pairs, and the elements of these pairs can be interchanged by polynomials.

\begin{prop}
    \label{prop: f_and_g}
    Let $s=(s_1,\ldots,s_m)\in(\frac{1}{e_1}\N,\ldots,\frac{1}{e_m}\N)$ such that $s_i>r\val\Disc(\mathfrak{x})$ for all $1\le i\le m$. Let $Z_0 \in \Z_p[x]$ be any monic polynomial of degree $r$ for which $\val Z_0(x_i) = s_i$ for all $1 \leq i \leq m$. Then there exists a unique set of roots $x_1',\ldots,x_m'$ of $Z_0$ satisfying
    \begin{equation}
        \label{eq:roots_close_xi'}
        ||x_i' - x_i|| < ||\Disc(\mfx)|| \quad \quad \quad \quad \text{ for all }1 \leq i \leq m.
    \end{equation}
    Furthermore, there exist unique polynomials $f,g \in \Z_p[x]$ of degree $\leq r-1$ such that
    \begin{equation}\label{eq:f_and_g}
        x_i' = f(x_i) \quad \quad \quad \text{and }\quad \quad \quad x_i = g(x_i') \quad \quad \quad \text{ for all }1 \leq i \leq m,
    \end{equation}
    where $f$ satisfies
    $$\val(f(x_i)-x_i)=s_i-\val Z'(x_i),\quad\text{ for all }1\le i\le m,$$
    in particular $f \in F_Z$ (recall \Cref{def:f_set}).
\end{prop}

\begin{proof}
For any fixed $1\le i\le m$, the $r$ roots of $Z_0$ give $r$ distances to $x_i$. The product of these distances is exactly 
$$||Z_0(x_i)||=p^{-s_i}<||\Disc(\mathfrak{x})||^r.$$ Therefore, we can find one root $x_i'$ of $Z_0$ that satisfies 
\begin{equation}\label{eq:roots_close_disc}
||x_i'-x_i||<||\Disc(\mathfrak{x})||.  
\end{equation}
Since such $x_i'$ exists for all $1\le i\le m$, we can apply \eqref{item: distinct conjugates} of \Cref{lem: close elements} to deduce that the roots $$\sigma_{1,1}(x_1'),\ldots,\sigma_{1,r_1}(x_1'),\ldots,\sigma_{m,1}(x_m'),\ldots,\sigma_{m,r_m}(x_m')$$
are all distinct, thus must be all the roots of the polynomial $Z_0$. Denote $\mathfrak{x}'=(x_1',\ldots,x_m')$. On the one hand, by \eqref{item: etale algebra contain} of \Cref{lem: close elements} we have $\Q_p[\mathfrak{x'}]\supset E=\Q_p[\mathfrak{x}]$. On the other hand, $x_1',\ldots,x_m'$ are the roots of a polynomial $Z_0$ of degree $r$, so $\Q_p[\mathfrak{x}']$ is of dimension $\le r$ viewed as a $\Q_p$-vector space. These two observations together imply $\Q_p[\mathfrak{x'}]=E=\Q_p[\mathfrak{x}]$. Hence we have 
\begin{equation}
x_i'=f(x_i),x_i=g(x_i'),\quad\forall 1\le i\le m
\end{equation}
for some polynomials $f,g\in\Q_p[x]$, which may be taken to be of degree $\le r-1$ by subtracting off multiples of the minimal polynomials of $x_i$ and $x_i'$ respectively if necessary. Also, for all $1\le i\le m$, we have 
\begin{align}\label{eq:val_s_i_align}
\begin{split}
\val(f(x_i)-x_i)&=\val(x_i'-x_i)\\
&=\val Z_0(x_i)-\sum_{(j,l_j)\ne (i,1)}\val(x_i-\sigma_{j,l_j}(x_j'))\\
&=s_i-\sum_{(j,l_j)\ne (i,1)}\val(x_i-\sigma_{j,l_j}(x_j))\\
&=s_i-\val Z'(x_i).
\end{split}
\end{align} 
Here the third equality comes from \eqref{item: same distance} of \Cref{lem: close elements}. Applying \Cref{cor: bound of coefficient}, the coefficients $c_0,\ldots,c_{r-1}$ of $f(x)-x$ all satisfy
\begin{align}
    \begin{split}
        \val(c_i) & \ge\min_i(s_i-2\val Z'(x_i)) \\ 
        &>\min_i(r\val\Disc(\mathfrak{x})-2\val Z'(x_i))\\ 
        &\ge 0,
    \end{split}
\end{align}
where the last line comes as usual from the fact that the discriminant is a product of terms of norm $\leq 1$ and $Z'(x_i)$ is a product of a subset of these terms. Thus, we have $f\in\Z_p[x]$. Since
\begin{equation}
    \val(g(x_i') - x_i') = \val(x_i'-x_i) = s_i-\val Z'(x_i)
\end{equation}
by \eqref{eq:val_s_i_align}, the rest of the proof applies to $g$ as well and shows $g\in\Z_p[x]$.

In the end, we show uniqueness of our constructions. By \eqref{item: same distance} of \Cref{lem: close elements}, we have
$$||\sigma_{j,l_j}(x_j')-x_i||=||\sigma_{j,l_j}(x_j)-x_i||\ge||\Disc(\mathfrak{x})||,\quad\forall (j,l_j)\ne (i,1),$$
thus the way to pick the roots $x_1',\ldots,x_m'$ of $Z_0$ satisfying \eqref{eq:roots_close_xi'} is unique. Furthermore, given any polynomial $f_1\in\Q_p[x]$ with degree $\le r-1$ such that 
$$f_1(x_i)=x_i',\quad\forall 1\le i\le m,$$ 
we have that $f_1-f$ is a polynomial of degree $\le r-1$ such that $(f_1-f)(x_i)=0$ for all $1\le i\le m$. Hence $(f_1-f)(\mfx)=0$, and since $\mfx$ has minimal polynomial $Z$ of degree $r$, we must have $f_1-f=0$. The same argument shows uniqueness of $g$.
\end{proof}

\begin{proof}[Proof of {\Cref{thm: polynomial expression of matrix near variety}}]
We first check that the map \eqref{eq:f_forward_map} is well-defined. Since $A_0$ has roots $x_1,\ldots,x_m$ and their Galois conjugates, the matrix $A=f(A_0)$ has roots $f(x_1),\ldots,f(x_m)$ and their conjugates. We wish to show that the matrix $A$ lies in the desired set in \eqref{eq:f_forward_map}. In order to apply \Cref{lem: close elements} to the sequence $f(x_1),\ldots,f(x_m)$, we must check that $||f(x_i) - x_i|| < ||\Disc(\mfx)||$. For this, notice that
\begin{align}
\begin{split}
\val(f(x_i)-x_i)&\ge s_i-\val Z'(x_i)\\
&>r\val\Disc(\mathfrak{x})-\val Z'(x_i)\\
&\ge\val\Disc(\mathfrak{x}),
\end{split}
\end{align}
where the last inequality holds because when $r=1$, we have $\val(Z'(x_i))=0$, and when $r\ge 2$, we have $\val\Disc(\mathfrak{x})\ge\val Z'(x_i)$. Now, applying \eqref{item: same distance} of \Cref{lem: close elements} in the second equality, we have for all $1\le i\le m$ that 
\begin{align}
\begin{split}
\val(P_A(x_i))&=\sum_{j=1}^m\sum_{l_j=1}^{r_j}\val(f(\sigma_{j,l_j}(x_j))-x_i)\\
&=\val(f(x_i)-x_i)+\sum_{(j,l_j)\ne (i,1)}\val(\sigma_{j,l_j}(x_j)-x_i)\\
&=\val(f(x_i)-x_i)+\val(Z'(x_i))\\
&\ge s_i.
\end{split}
\end{align}
Hence the map \eqref{eq:f_forward_map} is well-defined.

We now define the inverse map, so let us start with an arbitrary $A\in\Mat_r(\Z_p)$ such that $\val(P_A(x_i))\ge s_i,\forall 1\le i\le m$. Applying \Cref{prop: f_and_g}, we obtain two polynomials $f_A,g_A\in\Z_p[x]$ that switch the roots of $P_A$ and $Z$, so the map
\begin{align}\label{eq:reverse_map}
    \begin{split}
         \{A\in\Mat_r(\Z_p)\mid\val(P_{A}(x_i))\ge s_i,\quad\forall 1\le i\le m\} & \to \{(f,A_0)\mid A_0\in \mathcal{X},f\in F_Z\} \\ 
        A & \mapsto (f_A,g_A(A))
    \end{split}
\end{align}
is well-defined.

Finally, we prove that the maps we constructed are inverses. Let $f_B,g_B$ be the polynomials associated to the eigenvalues of a matrix $B$ by \Cref{prop: f_and_g}. Then we must show that
\begin{enumerate}
    \item\label{item: polynomial of matrix_injection} Let $A\in\Mat_r(\Z_p)$ with $\val(P_{A}(x_i))\ge s_i$ for all $ 1\le i\le m$. Then $f_A(g_A(A))=A$, where $f_A$ and $g_A$ are as above.
    \item\label{item: polynomial of matrix_surjection} For any $f \in F_Z, A_0 \in \mathcal{X}$, we have $f_{f(A_0)} = f$ and $g_{f(A_0)}(f(A_0)) = A_0$.
\end{enumerate}

We first prove \eqref{item: polynomial of matrix_injection}. Denote $x_i'=f_A(x_i)$ for all $1\le i\le m$. In fact, we have $f_A(g_A(x_i'))=f_A(x_i)=x_i',\forall 1\le i\le m$. This implies $P_A(x)\mid f_A(g_A(x))-x$, and therefore $f_A(g_A(A))=A$ by the Cayley–Hamilton theorem.

We next prove \eqref{item: polynomial of matrix_surjection}. The polynomial $f$ has the property that $f(x_i)$ is an eigenvalue of $f(A_0)$ for each $i$, hence the uniqueness part of \Cref{prop: f_and_g} implies that $f_{f(A_0)}=f$. Furthermore, since the polynomials $f$ and $g_{f(A_0)}$ switch the corresponding eigenvalues of $A_0$ and $f(A_0)$, we have $g_{f(A_0)}(f(x_i))=x_i$ for all $1\le i\le m$. This implies $Z(x)\mid g_{f(A_0)}(f(x))-x$, and therefore $g_{f(A_0)}(f(A_0)) = A_0$ by the Cayley–Hamilton theorem. This ends the proof.
\end{proof}

\begin{cor}\label{cor: norm of differential}
Let $s=(s_1,\ldots,s_m)\in(\frac{1}{e_1}\N,\ldots,\frac{1}{e_m}\N)$ such that $s_i>r\val\Disc(\mathfrak{x})$ for all $1\le i\le m$. Suppose $A\in\Mat_r(\Z_p)$ satisfies $\val(P_A(x_i))\ge s_i$ for all $1\le i\le m$. Let $f$ be the same as the bijection in \Cref{thm: polynomial expression of matrix near variety}, so that
$$f(\sigma_{1,1}(x_1)),\ldots,f(\sigma_{1,r_1}(x_1)),\ldots,f(\sigma_{m,1}(x_m)),\ldots,f(\sigma_{m,r_m}(x_m))$$
gives all the eigenvalues of $A$. Denote $x_i'=f(x_i)$ for all $1\le i\le m$. Then we have
$$||P_A'(x_i)||=||P_A'(x_i')||=||Z'(x_i)||.$$
\end{cor}

\begin{proof}
First, we have
\begin{align}
\begin{split}
||P_A'(x_i')||&=\prod_{(j,l_j)\ne(i,1)}||\sigma_{j,l_j}(x_j')-x_i'||\\
&=\prod_{(j,l_j)\ne(i,1)}||\sigma_{j,l_j}(x_j)-x_i||\\
&=||Z'(x_i)||,
\end{split}
\end{align}
where the second line comes from \Cref{lem: close elements}. Next, notice that $P_A'\in\Z_p[x]$, we have $$||P_A'(x_i')-P_A'(x_i)||\le||x_i'-x_i||<||\Disc(\mathfrak{x})||\le||Z'(x_i)||,$$
which implies
$$||P_A(x_i)||=||P_A'(x_i)-P_A'(x_i')+P_A'(x_i')||=||P_A'(x_i')||=||Z'(x_i)||.$$
\end{proof}

\subsection{Matrices and $\Z_p[\mfx]$-modules}

Recall the set $\Mod_{\Z_p[\mfx]}$ of full rank $\Z_p$-sublattices of $E$ which are $\Z_p[\mfx]$-modules, defined in \eqref{eq:def_ModZp}. We have a canonical group action $E^\times \acts \Mod_{\Z_p[\mathfrak{x}]}$ by scalar multiplication. As noted on page 4 of \cite{yun2013orbital}, the action of its subgroup $\L_E \subset E^\times$ on $\Mod_{\Z_p[\mathfrak{x}]}$ has finitely many orbits, denoted as $\#(\L_E\backslash\Mod_{\Z_p[\mathfrak{x}]})$. We write $E^\times \backslash \Mod_{\Z_p[\mathfrak{x}]}$ for the set of equivalence classes in $\Mod_{\Z_p[\mathfrak{x}]}$ under constant multiplication, and for every $M\in\Mod_{\Z_p[\mathfrak{x}]}$, we denote by $[M]=\{xM\mid x\in E^\times\}$ the equivalence class that contains $M$.

Denote by $\GL_r(\Z_p)\backslash \mathcal{X}$ the orbits of elements in $\mathcal{X}$ under the conjugation action of $\GL_r(\Z_p)$. For $\Q_p$-conjugacy there is only one orbit: recall that the rational canonical form states that every $A\in\Mat_r(\Q_p)$ with $P_A=Z$ can be written as $A=TC(Z)T^{-1}$, where $T\in\GL_r(\Q_p)$ and we recall the definition of $C(Z)$:

\begin{defi}\label{def:companion_matrix}
For every monic $Z(x)=x^r+c_{n-1}x^{r-1}+\cdots+c_0\in\Z_p[x]$ of degree $r$, we denote by $$C(Z)=\begin{pmatrix}0 & 0 & \ldots & 0 & -c_0\\
1 & 0 & \ldots & 0 & -c_1\\
0 & 1 & \ldots & 0 & -c_2\\
\vdots & \vdots & \ddots & \vdots & \vdots\\
0 & 0 & \ldots & 1 & -c_{n-1}\\
\end{pmatrix}\in\Mat_r(\Z_p)$$ 
the companion matrix of $Z$. 
\end{defi}

The following identifies these orbits in $\GL_r(\Z_p)\backslash \mc{X}$ with lattice equivalence classes.

\begin{prop}\label{prop: etale matrix variety}
    Given any lattice $M \in \Mod_{\Z_p[\mfx]}$, fix a choice of $\Z_p$-basis for $M$, and let $A_M$ be the matrix by which $\mfx$ acts in this basis. Then the map
    \begin{align}\label{eq:map_before_quotient}
        \begin{split}
            \Mod_{\Z_p[\mfx]} &\to \GL_r(\Z_p)\backslash \mathcal{X} \\ 
            M & \mapsto [A_M]
        \end{split}
    \end{align}
    is surjective. Furthermore, it factors through the quotient $\emod$, and the resulting map 
        \begin{align}\label{eq:map_thru_quotient}
        \begin{split}
            E^\times \backslash \Mod_{\Z_p[\mfx]} &\to \GL_r(\Z_p)\backslash \mathcal{X} \\ 
            [M] & \mapsto [A_M]
        \end{split}
    \end{align}
    is bijective. 
\end{prop}
\begin{proof}
    The map \eqref{eq:map_before_quotient} is clearly well-defined, since changing the choice of $\Z_p$-basis simply corresponds to conjugating $A_M$ by a matrix in $\GL_r(\Z_p)$. 
    
    For surjectivity, let $A \in \mathcal{X}$ be arbitrary. By rational canonical form, $A = TC(Z)T^{-1}$ for some $T \in \GL_r(\Q_p)$. Let $M \subset E$ be the $\Z_p$-span of the elements $T^{-1}(1),\ldots,T^{-1}(\mfx^{r-1})$, where we identify $T$ with a linear transformation of $E$ by identifying the standard basis vectors of $\Q_p^r$ with $1,\mfx,\ldots,\mfx^{r-1} \in E$. This $M$ is clearly full rank. Furthermore, because $\mfx$ acts on the basis $1,\mfx,\ldots,\mfx^{r-1}$ by the matrix $C(Z)$, it acts on the basis $T^{-1}(1),\ldots,T^{-1}(\mfx^{r-1})$ by $TC(Z)T^{-1} = A$. In particular, since $A  \in \Mat_{r \times r}(\Z_p)$, it takes the $\Z_p$-span of this basis to itself, i.e. $AM \subset M$. Hence, because $A$ represents the action of $\mfx$ on $M$ in this basis, $M$ is not just a $\Z_p$-lattice but is also a $\Z_p[\mfx]$-module, $M \in \Mod_{\Z_p[\mfx]}$. Since $\mfx$ acts on $M$ by the matrix $TC(Z)T^{-1} = A$, $M \mapsto [A]$ under the map \eqref{eq:map_before_quotient}, showing surjectivity.

    We now show the map factors through the quotient $\emod$. Let $\mcc\in \emod$ be any equivalence class. If 
    $$M_1,M_2 \in \mcc,$$ 
    then $M_2 = uM_1$ where $u \in E^\times$. In particular $u \in E$ so $u=f(\mfx)$ for some polynomial $f \in \Q_p[x]$. Hence $u$ commutes with the action of $\mfx$, so fixing a basis $\vec{v}_1,\ldots,\vec{v}_r$ of $M_1$, the action of $\mfx$ on $\vec{v}_1,\ldots,\vec{v}_r$ and on the basis $u\vec{v}_1,\ldots,u\vec{v}_r$ of $M_2$ is described by the same matrix. Thus $[A_{M_1}] = [A_{M_2}]$.

    Finally, we show the map \eqref{eq:map_thru_quotient} is injective. Suppose for some $M_1,M_2 \in \modzp$ that $[A_{M_1}] = [A_{M_2}]$. Then there exists a $\Z_p$-basis $\vec{v}_1,\ldots,\vec{v}_r$ of $M_1$ and a $\Z_p$-basis $\vec{w}_1,\ldots,\vec{w}_r$ of $M_2$ such that $\mfx$ acts by the same matrix $A$ on both bases. Letting $U \in \End_{\Q_p}(E)$ be the linear map with $U\vec{v}_i = \vec{w}_i$ for each $i$, we have that $UM_1 = M_2$ and $U$ commutes with $A$. Because $U$ commutes with $A$, they are simultaneously diagonalizable. Since $A$ is conjugate over $\Q_p$ to $C(Z)$, its eigenvalues are given by the roots of $Z$, namely the coordinates of $\mfx$ and their Galois conjugates. Because $\mfx \in E^{\new}$, polynomials in $\mfx$ generate $E$. Hence $U = f(A)$ for some polynomial $f$, i.e. $U$ acts by a polynomial $f(\mfx)$. Since $U$ is invertible, $f(\mfx) \in E^\times$. Hence $M_2 = u M_1$ for some $u = f(\mfx) \in E^\times$, so $M_1$ and $M_2$ lie in the same equivalence class $\mcc \in \emod$. This shows injectivity.
\end{proof}

The remainder of this subsection consists of building on and establishing further properties of the dictionary given in \Cref{prop: etale matrix variety}.

\begin{defi}
    \label{def:end_of_module}
    Given $M \in \modzp$, we denote by $\End(M) \subset \End_{\Q_p}(E)$ the set of $\Q_p$-linear maps $T$ on $E$ which satisfy $TM \subset M$ and commute with the action of $\mfx$ by multiplication. We similarly denote by $\Aut(M)  \subset \End(M)$ the maps with $TM=M$.
\end{defi}

\begin{lemma}\label{thm:end_of_class}
    Given a lattice class $\mcc \in E^\times \backslash \Mod_{\Z_p[\mathfrak{x}]}$ and $M_1, M_2 \in \mcc$, 
    \begin{equation}
        \End(M_1) = \End(M_2)
    \end{equation}
    and 
    \begin{equation}
        \Aut(M_1) = \Aut(M_2).
    \end{equation}
\end{lemma}
\begin{proof}
    $M_2 = uM_1$ for some $u \in E^\times$, and given any $T \in \End(M_1)$, trivially $uTu^{-1} \in \End(M_2)$. But $u$ polynomial in $\mfx$ and $T$ commutes with the action of $\mfx$, so $uTu^{-1} = T$. Hence $T \in \End(M_1) \Rightarrow T \in \End(M_2)$, and the result follows.
\end{proof}

\begin{defi}\label{def:end_of_class}
    For any class $\mcc \in \emod$, let $\End(\mcc) = \End(M)$ where $M \in \mcc$ is any choice of representative; by \Cref{thm:end_of_class} this does not depend on the choice of $M$. Define $\Aut(\mcc)$ similarly. 
\end{defi}

\begin{rmk}\label{rmk:one_basis_per_class}
    To reiterate a useful point in the proof of \Cref{prop: etale matrix variety}: Given any lattice $M \in \modzp$ with a fixed choice of basis $\vec{v}_1,\ldots,\vec{v}_r$, any other lattice $uM$ for $u \in E^\times$ has a corresponding basis $u\vec{v}_1,\ldots,u\vec{v}_r$. The matrix $A_M$ by which $\mfx$ acts on this basis of $M$ is the same matrix by which $\mfx$ acts on the basis $u\vec{v}_1,\ldots,u\vec{v}_r$ of $uM$, because $u$ is a polynomial in $\mfx$ and hence commutes with $\mfx$. Hence while the choice of basis for $M$ and matrix $A_M$ is non-canonical, one only needs to make it once for the whole class $[M]$, rather than once for each module in the class.
\end{rmk}

\begin{notation}\label{notation:basis_choice}
    From now on, whenever we are speaking of a lattice $M \in \Mod_{\Z_p[\mfx]}$, we will implicitly (and non-canonically) fix a $\Z_p$-basis $\vec{v}_1,\ldots,\vec{v}_r$ of $M$. Similarly, if we are speaking of an equivalence class $\mcc \in \emod$, we will fix a $\Z_p$-basis of one of the lattices $M \in \mcc$, and thereby fix a (different) basis of all other lattices in $\mcc$ by \Cref{rmk:one_basis_per_class}. We denote by $A_{\mcc}$ the matrix by which $\mfx$ acts on this basis. 
    
    Using this basis, we will freely identify $E$ with $\Q_p^r$ and the set of linear maps $\End_{\Q_p}(E)$ (and thereby its various subsets such as $\End(\mcc),\Aut(\mcc)$) with $\Mat_r(\Q_p)$ similarly to \Cref{prop: etale matrix variety}. We will also often identify $E$ and subsets such as $\mathcal{O}_E$ with linear maps in $\End_{\Q_p}(E)$ by considering their action by multiplication, and thereby write inclusions such as $\End(\mcc) \subset \mathcal{O}_E$ (see \Cref{thm: the deadly inclusion}). In summary\footnote{The identifications of $\End(\mcc)$ and $\Aut(\mcc)$ in the table are not entirely trivial, see \Cref{thm:end_and_aut}.}:

    \begin{center}
\begin{tabular}{|c|c|}
\hline
\text{Linear maps over $E$ and $M \subset E$} & \text{Matrices} \\ 
\hline
$\mathfrak{x}$ & $A_{\mcc}\in\Mat_r(\Z_p)$\\
\hline
$E$ & $\Q_p[A_{\mcc}]\subset\Mat_r(\Q_p)$\\ 
\hline
$\Z_p[\mathfrak{x}]$ & $\Z_p[A_{\mcc}]$\\
\hline 
$\End(\mcc)$ & $\Q_p[A_{\mcc}] \cap \Mat_{r \times r}(\Z_p)$ \\ 
\hline 
 $\Aut(\mcc)$ & $\Q_p[A_{\mcc}] \cap \GL_r(\Z_p)$ \\ 
\hline
\makecell{$E^\times=\{f(\mathfrak{x})\mid$\\$f\in\Q_p[x],f(x_i)\ne 0,\forall 1\le i\le m\}$} & \makecell{$\Q_p[A_{\mcc}]\cap \GL_r(\Q_p)=\{f(A_{\mcc})\mid$ \\ $f\in\Q_p[x],f(x_i)\ne 0,\forall 1\le i\le m\}$}\\
\hline
$\mc{O}_E$ & $\{B \in \Q_p[A_{\mcc}]: P_B \in \Z_p[x]\}$ \\ 
\hline 
\end{tabular}
\end{center}

Note that we are sometimes identifying $E$ with the vector space $\Q_p^r$ and sometimes with linear maps on this vector space given by multiplication.
\end{notation}

\begin{lemma}
    \label{thm:end_and_aut}
    Under a choice of basis as in \Cref{notation:basis_choice}, we have
    \begin{equation}\label{eq:end_and_matrix}
        \End(\mcc) = \Q_p[A_{\mcc}] \cap \Mat_{r \times r}(\Z_p)
    \end{equation}
    and 
    \begin{equation}\label{eq:aut_and_matrix}
        \Aut(\mcc) = \Q_p[A_{\mcc}] \cap \GL_r(\Z_p).
    \end{equation}
 
\end{lemma}
\begin{proof}
    Pick $M \in \mcc$. The set $\Mat_{r \times r}(\Z_p)$, in our chosen basis of $M$, is exactly the set of matrices taking $M$ to itself. Elements of $\End(\mcc) = \End(M)$ (identified via \Cref{def:end_of_class}) furthermore commute with the action of $\mfx$. This forces them to be polynomials in the matrix $A_{\mcc}$ by which $\mfx$ acts, as in the proof of \Cref{prop: etale matrix variety}. The identification \eqref{eq:aut_and_matrix} follows immediately from \eqref{eq:end_and_matrix}.
\end{proof}

In light of \Cref{prop: etale matrix variety} it is natural to name the components of $\mathcal{X}$ and choose a representative of each.

\begin{defi}
    \label{def:X_M}
    For each class $\mcc \in E^\times \backslash \Mod_{\Z_p[\mathfrak{x}]}$, denote by $\mathcal{X}_{\mcc}$ the orbit in $\GL_r(\Z_p)\backslash \mathcal{X}$ that corresponds to $\mcc$ by \Cref{prop: etale matrix variety}. 
\end{defi}

Given all of the $\Z_p$-lattices just defined, the following basic terminology is useful. 

\begin{defi}
Let $L_0\subset\Z_p^k$ be a $\Z_p$-submodule of rank $k_0$. Then there exists a unique maximal submodule $L_0\subset L\subset \Z_p^k$ of the same rank $k_0$, given explicitly by
$$L=(L_0\otimes_{\Z_p}\Q_p)\cap\Z_p^k.$$
We say $L$ is the \emph{maximal extension} of $L_0$. If $L_0=L$, we say that $L_0$ is \emph{maximal} in $\Z_p^k$.
\end{defi}

\begin{prop}\label{prop: properties of maximal submodule}
Let $L_0\subset\Z_p^k$ be a submodule of rank $k_0$, and $e_1,\ldots,e_{k_0}\in\Z_p^k$ be a $\Z_p$-basis of $L_0$. Then the following are equivalent:
\begin{enumerate}
\item $L_0$ is maximal. \label{item: maximal}
\item There exists $e_{k_0+1},\ldots,e_k\in\Z_p^k$ such that the sequence $e_1,\ldots,e_k$ gives a $\Z_p$-basis of $\Z_p^k$. \label{item: basis expansion}
\item $\Z_p^k/L_0$ is a free $\Z_p$-module of rank $k-k_0$. \label{item: free quotient}
\end{enumerate}
\end{prop}

\begin{proof}
We will prove equivalence of the above assertions in the cyclic way.

$\eqref{item: maximal}\rightarrow\eqref{item: free quotient}$: We prove this by contradiction. Suppose $\Z_p^k/L_0$ has a non-trivial torsion. By the second isomorphism theorem of modules, there exists a module $L\subset\Z_p^{k}$ of rank $k_0$, such that $L_0$ is a proper submodule of $L$. This contradicts our assumption that $L_0$ is maximal.

$\eqref{item: free quotient}\rightarrow\eqref{item: basis expansion}$: Suppose $\Z_p^k/L_0$ is a free $\Z_p$-module of rank $k-k_0$. Pick its $\Z_p$-basis $e_{k_0+1}+L_0,\ldots,e_k+L_0$, where $e_{k_0+1},\ldots,e_k\in\Z_p^k$. Then the sequence $e_1,\ldots,e_k$ gives a $\Z_p$-basis of $\Z_p^k$.

$\eqref{item: basis expansion}\rightarrow\eqref{item: maximal}$: Suppose there exists $e_{k_0+1},\ldots,e_k\in\Z_p^k$ such that the sequence $e_1,\ldots,e_k$ gives a $\Z_p$-basis of $\Z_p^k$. Let $L\supset L_0$ be the maximal extension of $L_0$, and set $L_1=\bigoplus_{i=k_0+1}^k\Z_p e_i$. Then we have $L_0+L_1=\Z_p^k\supset L$. Notice that $L_1\cap L=0$, we have $L_0=L$, and therefore $L_0$ is maximal.
\end{proof}

\begin{prop}
    \label{thm: the deadly inclusion}
    We have inclusion relations 
    \begin{equation}
        \label{eq:deadly_inclusion}
        \Z_p[A_{\mcc}]\subset \End(\mcc)\subset \mathcal{O}_E,
    \end{equation}
which are all free $\Z_p$-modules of rank $r$. Furthermore, the maximal extension of $\Z_p[A_{\mcc}]$ in $\Mat_r(\Z_p)$ is $\End(\mcc)$.  
\end{prop}
\begin{proof}
    By \Cref{thm:end_and_aut}, 
    $$\End(\mcc) = \Q_p[A_{\mcc}] \cap \Mat_r(\Z_p).$$
    Because $A_{\mcc} \in \Mat_r(\Z_p)$ itself, the first inclusion in \eqref{eq:deadly_inclusion} is clear. 
    
    For the second inclusion, since $\End(\mcc) \subset \Mat_r(\Z_p)$, any element $A \in \End(\mcc)$ must have a characteristic polynomial with $\Z_p$ coefficients. Hence by the Cayley-Hamilton theorem, $A$ satisfies a monic polynomial in $\Z_p[x]$. Furthermore, $\End(\mcc) \subset \Q_p[A_{\mcc}] = E$ (under our identification of $E$ with linear maps given by multiplication), so $A$ is identified with an element of $E$ that satisfies a monic polynomial with coefficients in $\Z_p$, i.e. an element of $\mc{O}_E$. This shows the second inclusion.

    Finally, $\End(\mcc)$ is maximal in $\Mat_r(\Z_p)$ since it is the intersection of a $\Q_p$-vector space with $\Mat_r(\Z_p)$. Together with the first inclusion of \eqref{eq:deadly_inclusion}, this shows that the maximal extension of $\Z_p[A_{\mcc}]$ in $\Mat_r(\Z_p)$ is $\End(\mcc)$.
\end{proof}

\begin{example}
Note that both containments in \Cref{thm: the deadly inclusion} may be strict or non-strict. When $\mcc=[\Z_p[\mfx]]$, we have $\End(\mcc)=\Z_p[\mfx]$, which in general may be a proper subset of $\mc{O}_E$. When $\mcc=[\mathcal{O}_E]$, we have $\End(\mcc)=\mathcal{O}_E$, so the second containment is non-strict. However, one obtains strict containment of $\Z_p[\mfx]$ inside $\mathcal{O}_E$ in the following example. Let $E=\Q_p[\sqrt p]$ be a field and $\mfx=p\sqrt p$. Then $\Z_p[\mfx]=\Z_p[p\sqrt p]$ becomes a proper subset of $\mathcal{O}_E=\Z_p[\sqrt p]$.
\end{example}

\subsection{Decomposing and computing volumes}

We may now decompose the set of matrices in \Cref{thm: points on variety_main theorem} into subsets indexed by these module classes:

\begin{cor}
    \label{thm:matrix_disjoint_union}
    Let $(s_1,\ldots,s_m)\in(\frac{1}{e_1}\N,\ldots,\frac{1}{e_m}\N)$ be such that $s_i > r \val \Disc(\mfx)$ for all $1\le i\le m$. Then
\begin{multline}\label{eq: matrix disjoint union}
\{A \in \Mat_{r \times r}(\Z_p)\mid\val (P_{A}(x_i))\ge s_i,\quad\forall 1\le i\le m\}=\\
\bigsqcup_{\mcc \in E^\times \backslash \Mod_{\Z_p[\mathfrak{x}]}}\{Tf(A_{\mcc})T^{-1}\mid f\in F_Z, T \in \GL_r(\Z_p)\},
\end{multline}
where we recall $F_Z$ from \Cref{def:f_set} and $A_{\mcc}$ from \Cref{notation:basis_choice}.
\end{cor}
\begin{proof}
From \Cref{prop: etale matrix variety}, we have
\begin{align}
\begin{split}
\text{LHS}\eqref{eq: matrix disjoint union}&=\{f(A_0)\mid f\in F_Z,A_0\in \mathcal{X}\}\\
&=\bigcup_{\mcc\in E^\times \backslash \Mod_{\Z_p[\mathfrak{x}]}}\{f(A_0)\mid f\in F_Z,A_0\in \mathcal{X}_\mcc\}\\
&=\bigcup_{\mcc\in E^\times \backslash \Mod_{\Z_p[\mathfrak{x}]}}\{Tf(A_{\mcc})T^{-1}\mid f\in F_Z, T \in \GL_r(\Z_p)\}.
\end{split}
\end{align}
Now we only need to show that the above expression is indeed a disjoint union. Suppose we have $\mcc_1,\mcc_2\in\Mod_{\Z_p[\mfx]}$, $T_1,T_2\in\GL_r(\Z_p)$ and $f_1,f_2\in F_Z$ such that $$T_1f_1(A_{\mcc_1})T_1^{-1}=f_1(T_1A_{\mcc_1}T_1^{-1})=f_2(T_2A_{\mcc_2}T_2^{-1})=T_2f_2(A_{\mcc_2})T_2^{-1}.$$ Then \Cref{thm: polynomial expression of matrix near variety} implies that $f_1=f_2$, and $T_1A_{\mcc_1}T_1^{-1}=T_2A_{\mcc_2}T_2^{-1}$. Thus we have $\mc{X}_{\mcc_1}=\mc{X}_{\mcc_2}$, and $\mcc_1=\mcc_2$, which ends the proof.
\end{proof}

Now we wish to compute the measure of the sets
\begin{equation}\label{eq:conj_f_set}
    \{Tf(A_{\mcc})T^{-1}\mid f\in F_Z, T \in \GL_r(\Z_p)\}
\end{equation}
which appear in \eqref{eq: matrix disjoint union}. We will do this by further breaking them up into smaller sets as in the following lemma, and then computing the number of such sets and their measure. How to naturally break them up? Because the centralizer of $f(A_{\mcc})$ in $\Mat_r(\Z_p)$ is $\End(\mcc)$, and its intersection with $\GL_r(\Z_p)$ is $\Aut(\mcc)$, it is natural to quotient out by this group for the $T$ in \Cref{thm:matrix_disjoint_union}. Because this quotient is still infinite ($\Aut(\mcc)$ has rank $r$ while $\GL_r(\Z_p)$ has rank $r^2$), it is natural to also take all matrices modulo a large power of $p$ in order to deal with a finite number of sets. 

\begin{lemma}
    \label{thm:decompose_TT_sets}
    Let $\mcc\in \emod$ and let $\mc{M} \in \N$ be sufficiently large. Let $\tT_1,\ldots,\tT_k \in \GL_r(\Z/p^{\mc{M}}\Z)$ be left coset representatives of $\Aut(\mcc) \pmod{p^{\mc{M}}}$ in $\GL_r(\Z_p/p^{\mc{M}}\Z_p)$, so $$\GL_r(\Z_p/p^{\mc{M}}\Z_p)=\bigsqcup_{i=1}^k\tT_i\Aut(\mcc) \pmod{p^{\mc{M}}}.$$
    Let $T_1,\ldots,T_k$ be any lifts of $\tT_1,\ldots,\tT_k$ to $\GL_r(\Z_p)$, so that 
    $$\GL_r(\Z_p)=\bigsqcup_{i=1}^k \{(T_i + \Delta T)\Aut(\mcc)\mid\Delta T\in p^{\mc{M}}\Mat_r(\Z_p)\}.$$
    Then
    \begin{multline}\label{eq:decompose_TT_sets}
        \{Tf(A_{\mcc})T^{-1}\mid f\in F_Z, T \in \GL_r(\Z_p)\} \\ 
        = \bigsqcup_{i=1}^k \{(T_i + \Delta T) f(A_{\mcc}) (T_i + \Delta T)^{-1} \mid \Delta T \in p^{\mc{M}}\Mat_r(\Z_p), f \in F_Z \}.
    \end{multline}
    Furthermore, each set on the right-hand side of \eqref{eq:decompose_TT_sets} has the same Haar measure as the set $\{(I_r + p^{\mc{M}} \Delta T) f(A_{\mcc}) (I_r + p^{\mc{M}}\Delta T)^{-1} \mid \Delta T \in p^{\mc{M}}\Mat_r(\Z_p), f \in F_Z \}$.
\end{lemma}
\begin{proof}
We first prove that the components in RHS\eqref{eq:decompose_TT_sets} are disjoint. Suppose we have $\Delta T_i,\Delta T_j\in p^{\mc{M}}\Mat_r(\Z_p)$ and $f_i,f_j\in F_Z$ such that 
\begin{equation*}
    (T_i+\Delta T_i)f_i(A_{\mcc})(T_i+\Delta T_i)^{-1}=(T_j+\Delta T_j)f_j(A_{\mcc})(T_j+\Delta T_j)^{-1},
\end{equation*} 
which we rewrite as 
\begin{equation*}
    \label{eq:f_and_T_different}
    f_i((T_i+\Delta T_i)A_{\mcc}(T_i+\Delta T_i)^{-1})=f_j((T_j+\Delta T_j)A_{\mcc}(T_j+\Delta T_j)^{-1}).
\end{equation*} 
Then \Cref{thm: polynomial expression of matrix near variety} implies that $f_i=f_j$, and 
\begin{equation}
    \label{eq:in_aut}
    (T_i+\Delta T_i)A_{\mcc}(T_i+\Delta T_i)^{-1}=(T_j+\Delta T_j)A_{\mcc}(T_j+\Delta T_j)^{-1}.
\end{equation}
By \eqref{eq:in_aut} we have $(T_j+\Delta T_j)^{-1}(T_i+\Delta T_i)\in\Aut(\mcc)$, hence $i=j$, proving disjointness.

Now to prove \eqref{eq:decompose_TT_sets}. Consider an arbitrary element on the left-hand side of \eqref{eq:decompose_TT_sets}. Every element $T \in \GL_r(\Z_p)$ can be expressed in the form $(T_i+\Delta T)B$, where $1\le i\le k,\Delta T\in p^{\mc{M}}\Mat_r(\Z_p)$, and $B\in\Aut(\mcc)$. Since $B$ commutes with the action of $A_{\mcc}$, $Bf(A_{\mcc})=f(A_{\mcc})B$, and therefore
$$(T_i+\Delta T)Bf(A_{\mcc})B^{-1}(T_i+\Delta T)^{-1}=(T_i+\Delta T)f(A_{\mcc})(T_i+\Delta T)^{-1}.$$
This shows the inclusion $\text{LHS\eqref{eq:decompose_TT_sets}$\subset$RHS\eqref{eq:decompose_TT_sets}}$, and the reverse inclusion is clear.

Finally, the fact that all such sets have the same measure as the one where $T_i = I_r$ follows by conjugating by $T_i$ and the defining property of the Haar measure.
\end{proof}

Now it remains to compute the number of these sets and the measure of each one. The first task is easier. 

\begin{prop}\label{thm:how_many_sets_M}
    Let $\mcc \in E^\times\backslash\Mod_{\Z_p[\mfx]}$. Then
    \begin{equation*}
        \#(\GL_r(\Z_p/p^{\mc{M}}\Z_p) / \Aut(\mcc) \pmod{p^{\mc{M}}}) = p^{\mathcal{M}(r^2-r)}\frac{(1-p^{-1})\cdots(1-p^{-r})\#(\mathcal{O}_E^\times/\Aut(\mcc))}{\prod_{i=1}^m(1-p^{-r_i/e_i})\#(\mathcal{O}_E/\End(\mcc))} 
    \end{equation*}
     for all positive integers $\mc{M}$. For the quotients on the right-hand side, note that $\mathcal{O}_E^\times/\Aut(\mcc)$ is a quotient of groups while $\mathcal{O}_E/\End(\mcc)$ is a quotient of rings, and to make sense of these we identify $E$ and $\End(\mcc)$ with rings of matrices as in \Cref{notation:basis_choice}.
\end{prop}
\begin{proof}
On the one hand, we have
\begin{equation}\label{eq: order of GL}
\#\GL_r(\Z_p/p^{\mc{M}}\Z_p)=(1-p^{-1})\cdots(1-p^{-r})p^{\mathcal{M}r^2}.
\end{equation}
On the other hand, $\End(\mcc)\subset\Mat_r(\Z_p)$ is a maximal lattice of rank $r$, and therefore 
\begin{equation}\label{eq: order of endomorphisms mod p^M}
\#\End(\mcc) \pmod{p^{\mc{M}}}=\#(\Z/p^{\mc{M}}\Z)^r=p^{\mathcal{M}r}.
\end{equation}

Recall that $\End(\mcc) \subset \mc{O}_E$ is a full-rank submodule. Also, for all $T\in\Aut(M)$, we have $T+p^\mathcal{M}\End(\mcc)\subset\Aut(M)$. Therefore, $
\Aut(\mcc)$ is the disjoint union of subsets of $\End(\mcc)$ of the form $T+p^\mathcal{\mathcal{M}}\End(\mcc)$, where  $T\in\Aut(\mcc)$. Then
\begin{align}\label{eq: proportion of Aut in End}
\begin{split}
\frac{\#\Aut(\mcc)\pmod{p^{\mc{M}}}}{\#\End(\mcc)\pmod{p^{\mc{M}}}}&=\frac{\#\Aut(\mcc)\pmod{p^{\mc{M}}}}{\#\mathcal{O}_E^\times\pmod{p^{\mc{M}}}}\frac{\#\mathcal{O}_E^\times\pmod{p^{\mc{M}}}}{\#\mathcal{O}_E\pmod{p^{\mc{M}}}}\frac{\#\mathcal{O}_E\pmod{p^{\mc{M}}}}{\#\End(\mcc)\pmod{p^{\mc{M}}}}\\
&=\frac{\#\mathcal{O}_E^\times\pmod{p^{\mc{M}}}}{\#\mathcal{O}_E\pmod{p^{\mc{M}}}}\cdot\frac{\#(\mathcal{O}_E/\End(\mcc))}{\#(\mathcal{O}_E^\times/\Aut(\mcc))}\\
&=\prod_{i=1}^m(1-p^{-r_i/e_i})\cdot \frac{\#(\mathcal{O}_E/\End(\mcc))}{\#(\mathcal{O}_E^\times/\Aut(\mcc))},
\end{split}
\end{align}
where in the last step we use that an element of $\mc{O}_E$ lies in $\mc{O}_E^\times$ if and only if each $K_i$ coordinate projects to a nonzero element in the residue field $\mc{O}_{K_i}/\pi_i \mc{O}_{K_i}$, which has size $p^{r_i/e_i}$.
The results in \eqref{eq: order of endomorphisms mod p^M} and \eqref{eq: proportion of Aut in End} together yields 
\begin{equation}\label{eq: order of automorphisms mod p^M}
\#(\Aut(\mcc) \pmod{p^{\mc{M}}})=\frac{\#(\mathcal{O}_E/\End(\mcc))}{\#(\mathcal{O}_E^\times/\Aut(\mcc))}\prod_{i=1}^m(1-p^{-r_i/e_i})p^{\mathcal{M}r}.
\end{equation}
The ratio of the orders in \eqref{eq: order of GL} and \eqref{eq: order of automorphisms mod p^M} becomes the order of the quotient we want.

\end{proof}

Finally, computing the measure of the set \eqref{eq:conj_f_set} is reduced to computing the Haar measure of the single set 
$$\{(I_r + p^{\mc{M}} \Delta T) f(A_{\mcc}) (I_r + p^{\mc{M}}\Delta T)^{-1} \mid \Delta T \in \Mat_r(\Z_p), f \in F_Z \}.$$ 
We are conjugating a matrix $f(A_{\mcc})$ by matrices very close to the identity, so it is natural to try to pass to the Lie algebra and linearize the problem. 

Recall that $\Mat_r(\Z_p)$ as a $\Z_p$-module of rank $r^2$, and for $A \in \mathcal{X}$, $\Z_p[\mathfrak{x}]\cong\Z_p[A]\subset\Mat_r(\Z_p)$ is a $\Z_p$-submodule of rank $r$. Also, the maximal extension of $\Z_p[A_{\mcc}]$ is $\End(\mcc)$ by \Cref{thm: the deadly inclusion}. 

\begin{defi}\label{defi: End^perp}
    For $\mcc \in E^\times\backslash\Mod_{\Z_p[\mfx]}$, we define
$$\End(\mcc)^\perp=\Im([\cdot,A_{\mcc}]\mid_{\Mat_r(\Z_p)})=\{TA_{\mcc}-A_{\mcc}T\mid T\in\Mat_r(\Z_p)\},$$ 
the image of the Lie bracket $[\cdot,A_{\mcc}]$ on $\Mat_{r \times r}(\Z_p)$.
\end{defi}

The below lemma justifies the notation $\End(\mcc)^\perp$.

\begin{lemma}\label{lem: order of lattice decomposition}
Let $\mcc \in E^\times\backslash\Mod_{\Z_p[\mfx]}$. Then $$\End(\mcc)\cap\End(\mcc)^\perp=0$$
and $\End(\mcc)+\End(\mcc)^\perp\subset\Mat_r(\Z_p)$ is a $\Z_p$-submodule of full rank $r^2$. Furthermore, it has index
$$\#(\Mat_r(\Z_p)/\End(\mcc)+\End(\mcc)^\perp)=p^{\val\Disc(\mathfrak{x})}=||\Disc(\mathfrak{x})||^{-1}.$$
\end{lemma}

\begin{proof}
Identifying $\Mat_{r \times r}(\Z_p) \cong \Z_p^r \ot \Z_p^r$, the linear operator $[\cdot,A]:\Mat_r(\Z_p)\rightarrow\Mat_r(\Z_p)$ is given in the basis $e_1 \ot e_1,\ldots,e_r \ot e_1,\ldots,e_1 \ot e_r,\ldots,e_r \ot e_r$ by the matrix
$$A \ot I_r - I_r \ot A = \begin{pmatrix}a_{1,1}I_r-A & a_{1,2}I_r & \cdots & a_{1,r}I_r \\
a_{2,1}I_r & a_{2,2}I_r-A & \cdots & a_{2,r}I_r \\
\vdots & \vdots & \ddots & \vdots \\
a_{r,1}I_r & a_{r,2}I_r & \cdots & a_{r,r}I_r-A \\
\end{pmatrix},$$
where $\otimes$ is the Kronecker product. We claim that this matrix has characteristic polynomial $x^{r^2}+\cdots+\Disc(\mathfrak{x})x^r$.

To see this, note that the above matrix is conjugate by an $r$-block matrix in $\GL_{r^2\times r^2}(\bar\Q_p)$ to the block diagonal matrix
\begin{multline}\label{eq: block form of Lie bracket}
\diag_{r^2\times r^2}(\sigma_{1,1}(x_1)I_r-A,\ldots,\sigma_{1,r_1}(x_1)I_r-A,\sigma_{2,1}(x_2)I_r-A,\ldots,\sigma_{2,r_2}(x_2)I_r-A,\ldots,\\
\sigma_{m,1}(x_m)I_r-A,\ldots,\sigma_{m,r_m}(x_m)I_r-A),
\end{multline}
by diagonalizing the $A \ot I_r$ part and recalling that $\sigma_{i,j}(x_i)$ are the eigenvalues of any $A \in \mathcal{X}$. The block form in \eqref{eq: block form of Lie bracket} has a characteristic polynomial equal to 
\begin{align}
    \begin{split}
        \prod_{i=1}^m\prod_{l_i=1}^{r_i}\det((x-\sigma_{i,l_i}(x_i))I_r+A)&=\prod_{(i,l_i)}\prod_{(j,l_j)}(x-\sigma_{i,l_i}(x_i)+\sigma_{j,l_j}(x_j)) \\ 
        &=x^{r^2}+\cdots+\Disc(\mathfrak{x})x^r.
    \end{split}
\end{align}
Now let us turn back to the matrix discussion. Denote by $e_1,\ldots,e_r$ a $\Z_p$-basis of $\End(\mcc)$. Since $\End(\mcc)\subset\Mat_r(\Z_p)$ is maximal, by \Cref{prop: properties of maximal submodule} we can extend this sequence to $e_1,\ldots,e_{r^2}$, which gives a $\Z_p$-basis of $\Mat_r(\Z_p)$. Because $A$ commutes with elements of $\End(\mcc)$, the matrix of the Lie bracket $[\cdot,A]$ over the basis $e_1,\ldots,e_{r^2}$ has the form
$$\begin{pmatrix} 0 & H \\ 0 & D\end{pmatrix}\in\Mat_{r^2\times r^2}(\Z_p).$$
From our above discussion of the characteristic polynomial, $D\in\Mat_{(r^2-r)\times(r^2-r)}(\Z_p)$ has determinant $\Disc(\mathfrak{x})$. In this case, $\End(\mcc)^\perp$ is the free $\Z_p$-module generated by the column vectors of $\begin{pmatrix} 0 & H \\ 0 & D\end{pmatrix}$. Thus $\End(\mcc)\cap\End(\mcc)^\perp=0$, and since $\End(\mcc)$ is the kernel of $[\cdot,A]$ it follows that $\End(\mcc)+\End(\mcc)^\perp\subset\Mat_r(\Z_p)$ is a $\Z_p$-submodule of full rank $r^2$. Finally, since $\det(D) = \Disc(\mathfrak{x})$ we have
$$\#(\Mat_r(\Z_p)/(\End(\mcc)+\End(\mcc)^\perp))=\left|\left|\det\begin{pmatrix} I_r & H \\ 0 & D\end{pmatrix}\right|\right|^{-1}=||\Disc(\mathfrak{x})||^{-1}.$$
\end{proof}

\begin{lemma}\label{lem: same image of Lie}
Let $\mcc \in E^\times\backslash\Mod_{\Z_p[\mfx]}$ and $f\in F_Z$. Then
$$\{Tf(A_{\mcc})-f(A_{\mcc})T\mid T\in\Mat_r(\Z_p)\}=\{TA_{\mcc}-A_{\mcc}T\mid T\in\Mat_r(\Z_p)\},$$
where we recall that the right-hand side is just $\End(\mcc)^\perp$.
\end{lemma}

\begin{proof}
Let $A = A_{\mcc}$. The inclusion 
$$\{f(A)T-Tf(A)\mid T\in\Mat_r(\Z_p)\}\subset\{AT-TA\mid T\in\Mat_r(\Z_p)\}$$
actually holds for any polynomial $f \in \Z_p[x]$; to prove this, it suffices to prove for $f(x)=x^k$ for all $k \geq 1$. This must be true since
\begin{equation*}
A^kT-TA^k=A(A^{k-1}T+A^{k-2}TA+\cdots +TA^{k-1})-(A^{k-1}T+A^{k-2}TA+\cdots +TA^{k-1})A,
\end{equation*}
which is in $\End(\mcc)^\perp$. 

To prove the reverse inclusion, denote $A'=f(A)$. Let $x_i'=f(x_i)$ for all $1\le i\le m$. By \Cref{prop: f_and_g}, there exists $g\in\Z_p[x]$ of degree $\le r-1$ such that $x_i=g(x_i')$. In this case, we have $g(f(x_i))=g(x_i')=x_i$ for all $1\le i\le m$, and therefore $Z(x)\mid g(f(x))-x$. By the Cayley-Hamilton theorem, $A=g(f(A))=g(A')$ can be expressed as a polynomial of $A'$. The same argument we used for $f$ yields
$$\{g(A')T-Tg(A')\mid T\in\Mat_r(\Z_p)\}\subset\{A'T-TA'\mid T\in\Mat_r(\Z_p)\},$$
and recalling that $g(A')=A, A' = f(A)$ this is exactly what we wanted to show.
\end{proof}

The following proposition intuitively measures the volume of a differential element.

\begin{prop}\label{prop: volume of differential element}
Fix $s=(s_1,\ldots,s_m)\in(\frac{1}{e_1}\N,\ldots,\frac{1}{e_m}\N)$ with $s_i>r\val\Disc(\mathfrak{x})$ for all $1\le i\le m$. Let $\mcc \in E^\times\backslash\Mod_{\Z_p[\mfx]}$. Then for all $\mathcal{M}\in\Z$ sufficiently large, we have
$$\mathbf{P}(\{(I_r+\Delta T)f(A_{\mcc})(I_r+\Delta T)^{-1}\mid\Delta T \in p^{\mc{M}}\Mat_r(\Z_p), f \in F_Z\})=\\
\frac{\#(\mathcal{O}_E/\End(\mcc))}{p^{\sum_is_ir_i+\mathcal{M}(r^2-r)}},$$
where for the quotient $\mathcal{O}_E$ and $\End(\mcc)$ are identified with matrix rings as in \Cref{thm:how_many_sets_M}.
\end{prop}
We first state a crucial lemma, prove \Cref{prop: volume of differential element} conditional on it, and then return to the proof of the lemma.

\begin{lemma}\label{lem: local conjugate and lie bracket}
In the same setup as \Cref{prop: volume of differential element}, for all $\mathcal{M}\in\Z$ sufficiently large we have 
\begin{multline}\label{eq: etale local conjugate and lie bracket}
\{(I_r+\Delta T)f(A_{\mcc})(I_r+\Delta T)^{-1}\mid\Delta T \in p^{\mc{M}}\Mat_r(\Z_p), f \in F_Z\}=\\
A_{\mcc} + \prod_{i=1}^m\frac{\pi_i^{e_is_i}}{||Z'(x_i)||}\mathcal{O}_{K_i}+p^{\mathcal{M}}\End(\mcc)^\perp,
\end{multline}
where we identify $\mathcal{O}_{K_i} \subset E$ with a subset of $\Mat_{r \times r}(\Q_p)$ via \Cref{notation:basis_choice} as usual.
\end{lemma}

\begin{proof}[Proof of \Cref{prop: volume of differential element}, based on \Cref{lem: local conjugate and lie bracket}]
Translating the set in \Cref{lem: local conjugate and lie bracket} does not change its Haar measure, so we must show 
$$\mathbf{P}\left(\prod_{i=1}^m\frac{\pi_i^{e_is_i}}{||Z'(x_i)||}\mathcal{O}_{K_i}+p^{\mathcal{M}}\End(\mcc)^\perp\right)=\frac{\#(\mathcal{O}_E/\End(\mcc))}{p^{\sum_is_ir_i+\mathcal{M}(r^2-r)}}.$$
Since 
$$\prod_{i=1}^m\frac{\pi_i^{e_is_i}}{||Z'(x_i)||}\mathcal{O}_{K_i}+p^{\mathcal{M}}\End(\mcc)^\perp$$ 
is a $\Z_p$-submodule of $\Mat_r(\Z_p)$ and the latter has measure $1$, it is equivalent to show 
\begin{equation}\label{eq: order of quotient of differential element}
\#\left(\Mat_r(\Z_p)/\left(\prod_{i=1}^m\frac{\pi_i^{e_is_i}}{||Z'(x_i)||}\mathcal{O}_{K_i}+p^{\mathcal{M}}\End(\mcc)^\perp\right)\right)=\frac{p^{\sum_is_ir_i+\mathcal{M}(r^2-r)}}{\#(\mathcal{O}_E/\End(\mcc))}.
\end{equation}
In fact, by \Cref{lem: order of lattice decomposition}, we have
\begin{align}
\begin{split}
\text{LHS}\eqref{eq: order of quotient of differential element}&=\#\left(\left(\End(\mcc)+\End(\mcc)^\perp\right)/\left(\prod_{i=1}^m\frac{\pi_i^{e_is_i}}{||Z'(x_i)||}\mathcal{O}_{K_i}+p^{\mathcal{M}}\End(\mcc)^\perp\right)\right)\\
&\times\#(\Mat_r(\Z_p)/(\End(\mcc)+\End(\mcc)^\perp))\\
&=\#\left((\End(\mcc)+\End(\mcc)^\perp)/\left(\prod_{i=1}^m\frac{\pi_i^{e_is_i}}{||Z'(x_i)||}\mathcal{O}_{K_i}+p^{\mathcal{M}}\End(\mcc)^\perp\right)\right)||\cdot \Disc(\mathfrak{x})||^{-1}\\
&=\frac{p^{\mathcal{M}(r^2-r)}}{||\Disc(\mathfrak{x})||}\#\left(\End(\mcc)/\prod_{i=1}^m\frac{\pi_i^{e_is_i}}{||Z'(x_i)||}\mathcal{O}_{K_i}\right),
\end{split}
\end{align}
where in the last step we use that $\End(\mcc) \cap \End(\mcc)^\perp = 0$ and $\End(\mcc)^\perp$ has rank $r^2-r$. Furthermore, recall that $\mc{O}_E = \prod_{i=1}^m \mc{O}_{K_i}$. This implies
\begin{align}
    \begin{split}           \#\left(\End(\mcc)/\prod_{i=1}^m\frac{\pi_i^{e_is_i}}{||Z'(x_i)||}\mathcal{O}_{K_i}\right)&=\frac{\#(\mathcal{O}_E/\prod_{i=1}^m\frac{\pi_i^{e_is_i}}{||Z'(x_i)||}\mathcal{O}_{K_i})}{\#(\mathcal{O}_E/\End(\mcc))} \\ 
        &= \frac{\prod_{i=1}^m \#(\mc{O}_{K_i} / \frac{\pi_i^{e_is_i}}{||Z'(x_i)||}\mathcal{O}_{K_i})}{\#(\mathcal{O}_E/\End(\mcc))}\\
        &=\frac{\prod_{i=1}^mp^{s_i r_i}||Z'(x_i)||^{r_i}}{\#(\mathcal{O}_E/\End(\mcc))}.
    \end{split}
\end{align}
Thus we have
$$\text{LHS}\eqref{eq: order of quotient of differential element}=\frac{p^{\mathcal{M}(r^2-r)}}{||\Disc(\mathfrak{x})||}\frac{p^{\sum_is_ir_i}\prod_{i=1}^m||Z'(x_i)||^{r_i}}{\#(\mathcal{O}_E/\End(\mcc))}=\text{RHS}\eqref{eq: order of quotient of differential element},$$
where the last equality holds because $\prod_{i=1}^m||Z'(x_i)||^{r_i}=||\Disc(\mathfrak{x})||$. This ends the proof.
\end{proof}

We now turn back to the proof of \Cref{lem: local conjugate and lie bracket}.

\begin{proof}[Proof of \Cref{lem: local conjugate and lie bracket}]

Let $A = A_{\mcc}$. We prove ($\subset$) and ($\supset$) for \eqref{eq: etale local conjugate and lie bracket}.

($\subset$) For all such $f$ and $\Delta T$,
$$(I_r+\Delta T)f(A)(I_r+\Delta T)^{-1}-A=(I_r+\Delta T)f(A)(I_r-\Delta T+(\Delta T)^2-\cdots)-A$$
lies in the set
$$f(A)-A+\Delta Tf(A)-f(A)\Delta T+p^{2\mathcal{M}}\Mat_r(\Z_p).$$
Therefore, we only need to prove that for sufficiently large $\mc{M}$,
$$f(A)-A+\Delta Tf(A)-f(A)\Delta T+p^{2\mathcal{M}}\Mat_r(\Z_p)\subset\prod_{i=1}^m\frac{\pi_i^{e_is_i}}{||Z'(x_i)||}\mathcal{O}_{K_i}+p^{\mathcal{M}}\End(\mcc)^\perp.$$
We already know $f(A)-A\subset\prod_{i=1}^m\frac{\pi_i^{e_is_i}}{||Z'(x_i)||}\mathcal{O}_{K_i}$. By \Cref{lem: same image of Lie}, we have $\Delta Tf(A)-f(A)\Delta T\in p^{\mathcal{M}}\End(\mcc)^\perp$. So we are left to show $p^{2\mathcal{M}}\Mat_r(\Z_p)\subset\prod_{i=1}^m\frac{\pi_i^{e_is_i}}{||Z'(x_i)||}\mathcal{O}_{K_i}+p^{\mathcal{M}}\End(\mcc)^\perp$. When $\mc{M}$ is sufficiently large, we have $p^{\mc{M}}\Mat_r(\Z_p)\subset\End(\mcc)+\End(\mcc)^\perp$ since the latter is a full-rank sublattice, and for the same reason $p^{\mc{M}}\End(\mcc)\subset\prod_{i=1}^m\frac{\pi_i^{e_is_i}}{||Z'(x_i)||}\mathcal{O}_{K_i}$. In this case, 
$$p^{2\mathcal{M}}\Mat_r(\Z_p)\subset p^{\mc{M}}\End(\mcc)+p^{\mc{M}}\End(\mcc)^\perp\subset\prod_{i=1}^m\frac{\pi_i^{e_is_i}}{||Z'(x_i)||}\mathcal{O}_{K_i}+p^{\mathcal{M}}\End(\mcc)^\perp,$$
hence ($\subset$) is proved.

($\supset$) We first claim that the right-hand side of \eqref{eq: etale local conjugate and lie bracket} may be rewritten as 
\begin{multline}
    A + \prod_{i=1}^m\frac{\pi_i^{e_is_i}}{||Z'(x_i)||}\mathcal{O}_{K_i}+p^{\mathcal{M}}\End(\mcc)^\perp \\ 
    = \{\tf(A) + \Delta \tT A - A \Delta \tT \mid \tf \in F_Z, \Delta \tT \in p^{\mathcal{M}} \Mat_{r \times r}(\Z_p)\}.
\end{multline}
This holds by renaming $\Delta \tT = p^{\mathcal{M}}\Delta T$ in \Cref{defi: End^perp}, and noting that a slight reformulation of \Cref{def:f_set} gives
\begin{equation}\label{eq:prod_to_f-x}
    \prod_{i=1}^m\frac{\pi_i^{e_is_i}}{||Z'(x_i)||}\mathcal{O}_{K_i}=\{\tf(\mathfrak{x})-\mathfrak{x}\mid\tf \in F_Z\}
\end{equation}
and the right-hand side corresponds to the set of matrices $\{\tf(A)-A\mid\tf \in F_Z\}$. 

Hence we only need to prove that given any such $\tf\in F_Z$ and $\Delta \tT\in p^{\mc{M}}\Mat_r(\Z_p)$, there exist $\Delta T\in p^{\mathcal{M}}\Mat_r(\Z_p)$ and $f \in F_Z$ such that 
\begin{equation}
    \label{eq:supset_wts}
    (I_r+\Delta T)^{-1}(\tf(A)+\Delta\tT A-A\Delta\tT)(I_r+\Delta T)=f(A).
\end{equation}
The idea is to mimic the proof of Hensel's lemma and construct $\Delta T$ and $f$ step by step, by conjugating $\tf(A)+\Delta\tT A-A\Delta\tT$ by a sequence of matrices $I_r + \Delta T_j$ in succession to bring it into the desired form.

We claim that there exists a sequence $(f_j)_{j\ge 1}$ of polynomials in $F_Z$ and sequences of matrices $(\Delta T_j)_{j\ge 1}, (\Delta \tT_j)_{j\ge 1}$, with initial terms $f_1=\tf,\Delta\tT_1=\Delta \tT$, such that 
\begin{equation}\label{eq: new f and Delta}
(I_r+\Delta T_j)^{-1}(f_j(A)+\Delta \tT_jA-A\Delta \tT_j)(I_r+\Delta T_j)=f_{j+1}(A)+\Delta \tT_{j+1}A-A\Delta \tT_{j+1}, 
\end{equation}
\begin{equation}\label{eq:where_matrices_live}
    \Delta T_j,\Delta \tT_j\in p^{\lfloor (j+1)\mathcal{M}/2\rfloor}\Mat_r(\Z_p),
\end{equation}
(i.e. both matrices have smaller norms at each step), and 
\begin{equation}\label{eq:f_j_cauchy}
    f_{j+1}(A)-f_j(A)\in p^{\lfloor (j+1)\mathcal{M}/2\rfloor}\Mat_r(\Z_p).
\end{equation}
It suffices to prove \eqref{eq: new f and Delta}, \eqref{eq:where_matrices_live}, and \eqref{eq:f_j_cauchy}, because then the polynomials $f_j$ converge to a limit $f_\infty$ because the nonarchimedean norm satisfies the strong triangle inequality, and due to \eqref{eq:where_matrices_live} the $\Delta \tT_jA-A\Delta \tT_j$ terms converge to $0$ as $j \to \infty$, so we have
\begin{equation}
    \prod_{j = 1}^\infty (I_r + \Delta T_j)^{-1}(\tf(A)+\Delta\tT A-A\Delta\tT) \prod_{j = 1}^\infty (I_r + \Delta T_j)^{-1} = f_\infty(A).
\end{equation}
Taking $f = f_\infty$ and $\Delta T=\prod_{j=1}^\infty(I_r+\Delta T_j)-I_r$ then yields \eqref{eq:supset_wts}. Now let us give the construction.

Let $j \geq 1$. At the $j\tth$ step, we will take as input $f_j$ and $\Delta \tT_j$ (note that we already have these when $j=1$, since we set $f_1=\tf$ and $\Delta \tT_1 = \Delta \tT$), and find the polynomial $f_{j+1}$ and the matrices $\Delta T_j,\Delta \tT_{j+1}$ satisfying \eqref{eq: new f and Delta}, \eqref{eq:where_matrices_live}, and \eqref{eq:f_j_cauchy}. There are two steps.

\begin{enumerate}
\item We find $\Delta T_j\in p^{\lfloor (j+1)\mathcal{M}/2\rfloor}\Mat_r(\Z_p)$ such that $\Delta T_j f_j(A)-f_j(A)\Delta T_j=\Delta \tT_jA-A\Delta \tT_j$. This is always achievable due to \Cref{lem: same image of Lie}.
\item Consider the conjugation action of $(I_r+\Delta T_j)$ over $f_j(A)+\Delta \tT_jA-A\Delta \tT_j$:
\begin{multline*}
(I_r+\Delta T_j)^{-1}(f_j(A)+\Delta\tT_jA-A\Delta\tT_j)(I_r+\Delta T_j)=(I_r-\Delta T_j+\Delta T_j^2-\cdots)\cdot\\(f_j(A)+\Delta \tT_jA-A\Delta \tT_j)(I_r+\Delta T_j)\subset f_j(A)+p^{(j+1)\mathcal{M}-1}\Mat_r(\Z_p).
\end{multline*}
When $\mathcal{M}$ is sufficiently large (uniformly and does not rely on the specific value of $j\ge 1$), we must have  $$p^{(j+1)\mathcal{M}-1}\Mat_r(\Z_p)\subset\prod_{i=1}^m\frac{\pi_i^{e_is_i}}{||Z'(x_i)||}\mathcal{O}_{K_i}+p^{\lfloor (j+2)\mathcal{M}/2\rfloor}\End(\mcc)^\perp.$$
By \eqref{eq:prod_to_f-x},
\begin{equation}
    f_j(A) + \prod_{i=1}^m\frac{\pi_i^{e_is_i}}{||Z'(x_i)||}\mathcal{O}_{K_i} = \{g(A): g \in F_Z\},
\end{equation}
hence there exists $f_{j+1}\in F_Z$ and $\Delta \tT_{j+1}\in p^{\lfloor (j+2)\mathcal{M}/2\rfloor}\Mat_r(\Z_p)$ such that \eqref{eq: new f and Delta} holds.
\end{enumerate}
This constructs the sequences $(f_j)_{j\ge 1},(\Delta T_j)_{j\ge 1}, (\Delta \tT_j)_{j\ge 1}$ satisfying \eqref{eq: new f and Delta}, \eqref{eq:where_matrices_live}, and \eqref{eq:f_j_cauchy}, so we are done.
\end{proof}

\begin{cor}\label{cor: limit of one orbit}
Fix $s=(s_1,\ldots,s_m)\in(\frac{1}{e_1}\N,\ldots,\frac{1}{e_m}\N)$ with $s_i>r\val\Disc(\mathfrak{x})$ for all $1\le i\le m$. Let $\mcc \in E^\times\backslash\Mod_{\Z_p[\mfx]}$. Then we have 
\begin{equation}\label{eq: measure of one orbit}
\mathbf{P}(\{Tf(A_{\mcc})T^{-1}\mid T\in\GL_r(\Z_p),f\in F_Z\})=\frac{(1-p^{-1})\cdots(1-p^{-r})}{\prod_{i=1}^m(1-p^{-r_i/e_i})p^{s_ir_i}}\#(\mathcal{O}_E^\times/\Aut(\mcc)).
\end{equation}
\end{cor}

\begin{proof}
Recall that every set in the disjoint union in \eqref{eq:decompose_TT_sets} in \Cref{thm:decompose_TT_sets} has the same measure. By \Cref{thm:how_many_sets_M}, there are 
$$\#(\GL_r(\Z_p/p^{\mc{M}}\Z_p) / \Aut(\mcc) \pmod{p^{\mc{M}}}) = p^{\mathcal{M}(r^2-r)}\frac{(1-p^{-1})\cdots(1-p^{-r})\#(\mathcal{O}_E^\times/\Aut(\mcc))}{\prod_{i=1}^m(1-p^{-r_i/e_i})\#(\mathcal{O}_E/\End(\mcc))}$$
such elements in total. Therefore, we have 
\begin{align}
\begin{split}
\text{LHS}\eqref{eq: measure of one orbit}&=p^{\mathcal{M}(r^2-r)}\frac{(1-p^{-1})\cdots(1-p^{-r})\#(\mathcal{O}_E^\times/\Aut(\mcc))}{\prod_{i=1}^m(1-p^{-r_i/e_i})\#(\mathcal{O}_E/\End(\mcc))}\\
&\times\mathbf{P}(\{(I_r+\Delta T)f(A_{\mcc})(I_r+\Delta T)^{-1}\mid\Delta T \in p^{\mc{M}}\Mat_r(\Z_p), f \in F_Z\})\\
&=\frac{(1-p^{-1})\cdots(1-p^{-r})}{\prod_{i=1}^m(1-p^{-r_i/e_i})p^{s_ir_i}}\frac{\#(\mathcal{O}_E^\times/\Aut(\mcc))}{\#(\mathcal{O}_E/\End(\mcc))}\#(\mathcal{O}_E/\End(\mcc))\\
&=\frac{(1-p^{-1})\cdots(1-p^{-r})}{\prod_{i=1}^m(1-p^{-r_i/e_i})p^{s_ir_i}}\#(\mathcal{O}_E^\times/\Aut(\mcc)),
\end{split}
\end{align}
where the second equality comes from \Cref{prop: volume of differential element}. This gives the proof.
\end{proof}

The following results from \cite{yun2013orbital} allow us to simplify the sum over $\mcc \in \emod$ appearing in our formulas.

\begin{prop}
The following equalities hold for the quotients $\#(\L_E\backslash \Mod_{\Z_p[\mathfrak{x}]})$:
\begin{equation}\label{eq: orbital integral as sum of orbit} \sum_{\mcc\in E^\times \backslash \Mod_{\Z_p[\mathfrak{x}]}}\#(\mathcal{O}_E^\times/\Aut(\mcc))=\#(\L_E\backslash \Mod_{\Z_p[\mathfrak{x}]}),
\end{equation} 
and
\begin{equation}\label{eq: orbital integral as elliptic components} \#(\L_E\backslash \Mod_{\Z_p[\mathfrak{x}]})||\Delta_\sigma(x_1,\ldots,x_m)||=\prod_{i=1}^m\#(\L_{K_i}\backslash \Mod_{\Z_p[x_i]})||\Delta_\sigma(x_i)||.
\end{equation}
\end{prop}

\begin{proof}
\eqref{eq: orbital integral as sum of orbit} comes from \cite[Section 2.14]{yun2013orbital}, specifically the sentence `Now we claim...' after \cite[(2.17)]{yun2013orbital}; $\bar{\Cl}(R)$ and $X_R$ there correspond to $\emod$ and $\Mod_{\Z_p[\mfx]}$ in our notation. For \eqref{eq: orbital integral as elliptic components}, notice that 
\begin{align}
\begin{split}
\frac{\prod_{i=1}^m||\Delta_\sigma(x_i)||}{||\Delta_\sigma(x_1,\ldots,x_m)||}&=\frac{\prod_{i=1}^m||\Delta(\sigma_{i,1}x_i,\ldots,\sigma_{i,r_i}x_i)||}{||\Delta(\sigma_{1,1}x_1,\ldots,\sigma_{1,r_1}x_1,\ldots,\sigma_{m,1}x_m,\ldots,\sigma_{m,r_m}x_m)||}\\
&=\prod_{1\le i<j\le m}\prod_{l_i=1}^{r_i}\prod_{l_j=1}^{r_j}||\sigma_{i,l_i}(x_i)-\sigma_{j,l_j}(x_j)||^{-1}\\
&=\prod_{1\le i<j\le m}||\Res(Z_i,Z_j)||^{-1}\\
&=p^{\sum_{1\le i<j\le m}\val(\Res(Z_i,Z_j))}.
\end{split}
\end{align}
Here for all $1\le i<j\le m$, $\Res(Z_i,Z_j)$ denotes the resultant of the polynomials $Z_i,Z_j$ from \Cref{subsec:zpbar}. Therefore, in order to prove \eqref{eq: orbital integral as elliptic components}, we only have to show
$$\#(\L_E\backslash \Mod_{\Z_p[\mathfrak{x}]})=p^{\sum_{1\le i<j\le m}\val(\Res(Z_i,Z_j))}\cdot\prod_{i=1}^m\#(\L_{K_i}\backslash \Mod_{\Z_p[x_i]}),$$
which is exactly the same as \cite[Corollary 4.10]{yun2013orbital}\footnote{In the notation of \cite[Corollary 4.10]{yun2013orbital}, we are taking $q=p$, $O_\gamma=\#(\L_E\backslash \Mod_{\Z_p[\mathfrak{x}]})$, $B(\gamma)=\{1,2,\ldots,m\}$, and $O_{\gamma_i}^{L_i}=\#(\L_{K_i}\backslash \Mod_{\Z_p[x_i]})$ for all $1\le i\le m$. According to the definition given in \cite[Subsection 4.1]{yun2013orbital}, the notation $\rho(\gamma)$ is equal to $\sum_{1\le i<j\le m}\val(\Res(Z_i,Z_j))$, which is exactly the expression we need here.}.
\end{proof}

\begin{rmk}
Notice that we study the group action of $E^\times$ and $\Lambda$ over the same set of modules $\Mod_{\Z_p[\mfx]}$. The group $\Lambda$ a subgroup of $E^\times$, and therefore generates more orbits. In fact, the order of the quotient $\#(\mathcal{O}_E^\times/\Aut(\mcc))$ in \eqref{eq: orbital integral as sum of orbit} is the number of orbits in class $\mcc$, under the group action of $\Lambda$. 
\end{rmk}

Now, we are ready for the proof of \Cref{thm: points on variety_main theorem}.

\begin{proof}[Proof of \Cref{thm: points on variety_main theorem}] When $s=(s_1,\ldots,s_m)\in(\frac{1}{e_1}\N,\ldots,\frac{1}{e_m}\N)$ such that $s_i>r\val\Disc(\mathfrak{x})$ for all $1\le i\le m$, we have (recalling the notation of \Cref{def:X_M}) that 
\begin{align}\label{eq: sum up all lattice classes}
\begin{split}
\text{LHS}\eqref{eq: points on the variety}&=\prod_{i=1}^m p^{s_ir_i}\cdot\sum_{\mcc\in E^\times \backslash \Mod_{\Z_p[\mathfrak{x}]}}\mathbf{P}(\{Tf(A_{\mcc})T^{-1}\mid T\in\GL_r(\Z_p),f\in F_Z\})\\
&=\sum_{\mcc\in E^\times \backslash \Mod_{\Z_p[\mathfrak{x}]}}\frac{(1-p^{-1})\cdots(1-p^{-r})}{\prod_{i=1}^m(1-p^{-r_i/e_i})}\#(\mathcal{O}_E^\times/\Aut(\mcc))\\
&=\frac{(1-p^{-1})\cdots(1-p^{-r})}{\prod_{i=1}^m(1-p^{-r_i/e_i})}\#(\L_E\backslash \Mod_{\Z_p[\mathfrak{x}]})\\
&=\frac{(1-p^{-1})\cdots(1-p^{-r})}{||\Delta_\sigma(x_1,\ldots,x_m)||}\prod_{i=1}^m\frac{\#(\L_{K_i}\backslash\Mod_{\Z_p[x_i]})||\Delta_\sigma(x_i)||}{(1-p^{-r_i/e_i})}\\
&=\text{RHS}\eqref{eq: points on the variety}.
\end{split}
\end{align}
where the first equality is by \Cref{thm:matrix_disjoint_union}, the second equality comes from \Cref{cor: limit of one orbit}, the third equality comes from \eqref{eq: orbital integral as sum of orbit}, and the fourth equality comes from \eqref{eq: orbital integral as elliptic components}. Therefore, we are done.
\end{proof}

\begin{rmk}
We remark that similar formulas and techniques appeared in the unpublished thesis of Boreico \cite[Section 3]{boreico2016statistics} (we are grateful to Akshay Venkatesh for bringing this to our attention after we found the proofs here). This work studied the related question of the joint distribution of coefficients of the characteristic polynomial of a random matrix. An expression for the joint density $\rho(f)$ of these coefficients with respect to the Haar measure is given in \cite[Theorem 3.3.26(b)]{boreico2016statistics}. The formula there is similar to the second line of \eqref{eq: sum up all lattice classes}, featuring an infinite sum. It appears that the additional simplification of \eqref{eq: sum up all lattice classes}, in which the term $\#(\L_E\backslash\Mod_{\Z_p[\mfx]})$ emerges when we sum up the probability measure corresponding to each lattice class $\mcc$, could also lead to a simpler expression of $\rho(f)$ similar to the right hand side of \eqref{eq: points on the variety}, but this was not done in \cite{boreico2016statistics}.
\end{rmk}

%% file: joint_distribution_2.tex
\section{From characteristic polynomial to joint distribution of eigenvalues}\label{sec: joint distribution}

Unless otherwise stated, we assume the same notations in this section as given in the beginning of the previous one. 

\subsection{Proof of {\Cref{thm: joint distribution}}} The following proposition allows us to focus on new elements in the étale algebras.

\begin{prop}\label{prop: almost surely distinct eigenvalues}
Let $A\in\Mat_n(\Z_p)$ be Haar-distributed. Then we have 
$$\mathbf{P}(A\text{ has repeated eigenvalues})=0.$$
\end{prop}

\Cref{prop: almost surely distinct eigenvalues} is a direct corollary of the following lemma, which is a $p$-adic analog of \cite{caron2005zero}.

\begin{lemma}\label{lem: zero set of polynomial has measure zero}
Let $n\ge 1$ be an integer, and $f\in\Z_p[\mathbf{x}_1,\ldots,\mathbf{x}_n]:\Z_p^n\rightarrow\Z_p$ be a nonzero polynomial with coefficients in $\Z_p$. Let $(a_1,\ldots,a_n)\in\Z_p^n$ be random with respect to the Haar probability measure. Then we have
$$\mathbf{P}(f(a_1,\ldots,a_n)=0)=0.$$
\end{lemma}

\begin{proof}[Proof of \Cref{prop: almost surely distinct eigenvalues}, based on \Cref{lem: zero set of polynomial has measure zero}]
Let us consider the characteristic polynomial $P_A$. Its discriminant can be viewed as a polynomial over the $n^2$ entries of $A$, and is given by 
\begin{align}
\begin{split}
f:\Mat_n(\Z_p)&\rightarrow\Z_p\\
A&\mapsto\Disc(P_A).
\end{split}
\end{align} 
The zero set of $f$ is exactly the matrices $A\in\Mat_n(\Z_p)$ that have repeated roots. It is clear that $f$ is not always zero, because there exists a matrix in $\Mat_n(\Z_p)$ without repeated roots. This leads to the proof by applying \Cref{lem: zero set of polynomial has measure zero}.
\end{proof}

\begin{proof}[Proof of \Cref{lem: zero set of polynomial has measure zero}]
We give the proof by induction over $n$. The statement is trivial when $n=1$. Now, suppose the assertion holds for $n-1$. Given $f\in\Z_p[\mathbf{x}_1,\ldots,\mathbf{x}_n]$, we write
$$f=\sum_{j\in\N}f_j(\mathbf{x}_1,\ldots,\mathbf{x}_{n-1})\mathbf{x}_n^j,$$
where the polynomials $f_j$ have coefficients in $\Z_p$. Since $f$ is nonzero, there exists some $j\in\N$ such that $f_j$ is nonzero. When $f(a_1,\ldots,a_n)=0$, one of the following cases must be true for the sequence $(a_1,\ldots,a_{n-1})$:
\begin{enumerate}
\item $f_j(a_1,\ldots,a_{n-1})=0$. By the induction hypothesis, this case has probability zero.
\item $f_j(a_1,\ldots,a_{n-1})\ne 0$. In this case, there are only finitely many $a_n$ that satisfy $f(a_1,\ldots,a_n)=0$, which also has probability zero.
\end{enumerate}
And therefore we are done.
\end{proof}

To prepare for the proof of \Cref{thm: joint distribution}, we state the following $p$-adic Kac-Rice formula, which is a generalized version of \cite[Theorem 2.1]{caruso2022zeroes}:

\begin{lemma}\label{lem: Kac-Rice formula}
For fixed $A\in\Mat_r(\Z_p)$, and a compact open subset $U\subset \mathcal{O}_E^{\new}$, we have for all integers $s\ge 0$,
$$\# \{\mathfrak{x}\in U\mid P_A(\mathfrak{x})=0\}=\lim_{s\rightarrow\infty}\int _UI_s(\mathfrak{x},A)d\mathfrak{x},$$
where the limit is taken along nonnegative integers $s\in\N$, and for $\mathfrak{x}=(x_1,\ldots,x_m)$,
$$I_s(\mathfrak{x},A)=p^{sr}\prod_{i=1}^m ||P_A'(x_i)||^{r_i}\cdot\bbone_{\val P_A(x_i)\ge s,\forall 1\le i\le m}.$$
\end{lemma}

\begin{proof}
Suppose we have $\mathfrak{y}=(y_1,\ldots,y_m)\in U$ with $P_A(\mathfrak{y})=0$. Applying \cite[Lemma 3.4]{caruso2014tracking}\footnote{The statement of \cite[Lemma 3.4]{caruso2014tracking} is somewhat notation-heavy, but exactly the same special case we use here is used later in the proof of \cite[Theorem 2.1]{caruso2022zeroes}. Both results are stated for strictly differentiable functions, a class which in particular includes all polynomials.}, for each $1 \leq i \leq m$ there exists a positive integer $S_{y_i}$ such that the following holds:

\begin{enumerate}
\item $||P_A'(x_i)||=||P_A'(y_i)||$ and $x_i\in y_i+\frac{\pi_i^{e_iS_{y_i}}}{||P_A'(y_i)||}\mathcal{O}_{K_i}$.
\item For any integer $s>S_{y_i}$, the polynomial $P_A$ gives a bijection
\begin{equation}
    \label{eq:caruso_bij_coordinates}
    y_i+\frac{\pi_i^{e_is}}{||P_A'(y_i)||}\mathcal{O}_{K_i}\leftrightarrow\pi_i^{e_is}\mathcal{O}_{K_i}.
\end{equation}
\end{enumerate}

Taking $S_{\mf{y}} = \max_i S_{y_i}$, it follows immediately from \eqref{eq:caruso_bij_coordinates} that for any integer $s > S_{\mf{y}}$, $P_A$ gives a bijection 
\begin{equation}
    \label{eq:P_A_bijective}
    \mathfrak{y}+\prod_{i=1}^m\frac{\pi_i^{e_is}}{||P_A'(y_i)||}\mathcal{O}_{K_i}\leftrightarrow\prod_{i=1}^m\pi_i^{e_is}\mathcal{O}_{K_i}.
\end{equation}

In this case, for all $s>S_{\mathfrak{y}}$, we have
\begin{align}
\begin{split}
\int_{\mathfrak{y}+\prod_{i=1}^m\frac{\pi_i^{e_iS_{\mathfrak{y}}}}{||P_A'(y_i)||}\mathcal{O}_{K_i}}I_s(\mathfrak{x},A)d\mathfrak{x}&=\int_{\mathfrak{y}+\prod_{i=1}^m\frac{\pi_i^{e_iS_{\mathfrak{y}}}}{||P_A'(y_i)||}\mathcal{O}_{K_i}}p^{sr}\prod_{i=1}^m ||P_A'(y_i)||^{r_i}\bbone_{P_A(\mfx)\in\prod_{i=1}^m\pi_i^{e_is}\mathcal{O}_{K_i}}d\mfx\\
&=p^{sr}\prod_{i=1}^m ||P_A'(y_i)||^{r_i}\int_{\mathfrak{y}+\prod_{i=1}^m\frac{\pi_i^{e_is}}{||P_A'(y_i)||}\mathcal{O}_{K_i}}d\mfx\\
&=1.
\end{split}
\end{align}
Also, because $P_A$ is a bijection between the sets in \eqref{eq:P_A_bijective} and the latter set includes $0$, $\mathfrak{y}$ is the unique zero of $P_A$ in the set $\mathfrak{y}+\prod_{i=1}^m\frac{\pi_i^{e_iS_{\mathfrak{y}}}}{||P_A'(y_i)||}\mathcal{O}_{K_i}$. Therefore, the zeros of $P_A$ form a discrete subset of $U$. Since $U$ is compact, there are thus only finitely many zeros in $U$, say $$\mathfrak{y}_1=(y_{1,1},\ldots,y_{1,m}),\ldots,\mathfrak{y}_n=(y_{n,1},\ldots,y_{n,m}).$$ 
Take an integer $S$ such that $S>\max\{S_{\mathfrak{y}_1},\ldots,S_{\mathfrak{y}_n}\}$, and the sets
$$U_k=\mathfrak{y}_k+\prod_{i=1}^m\frac{\pi_i^{e_iS}}{||P_A'(y_{k,i})||}\mathcal{O}_{K_i},\quad 1\le k\le n$$
are disjoint. Denote by $V$ the complement in $U$ of $U_1\cup\cdots\cup U_n$. Since $V$ is compact and $P_A$ does not vanish on it, we deduce that when $s$ is sufficiently large, $||P_A(\mathfrak{x})||\ge p^{-s}$ for all $\mathfrak{x}\in V$. The indicator function in $I_s(\mfx,A)$ is then always $0$ when $\mfx \in V$. Hence
$$\int _UI_s(\mathfrak{x},A)d\mathfrak{x}=\int_VI_s(\mathfrak{x},A)d\mathfrak{x}+\sum_{i=1}^n\int_{U_i}I_s(\mathfrak{x},A)d\mathfrak{x}=0+\sum_{i=1}^n 1=n,$$
which ends the proof.
\end{proof}

\begin{rmk}
    While we do not assume that $A$ has distinct eigenvalues yet in \Cref{lem: Kac-Rice formula}, tuples of eigenvalues with repeated entries cannot occur on $\mathcal{O}_{E}^{\new}$, so when restricted to $\mathcal{O}_{E}^{\new}$ all roots of $P_A$ are distinct.
\end{rmk}

The following lemma is essential for our proof of \Cref{thm: joint distribution}, which is an application of \cite[Proposition 2.3]{evans2006expected}.

\begin{lemma}\label{lem: integral and number of solutions}
For fixed $A\in\Mat_r(\Z_p)$, and a open subset $U\subset \mathcal{O}_E^{\new}$, we have for all integers $s\ge 0$, 
\begin{equation}\label{eq: integral and number of solutions}
\int_U I_s(\mfx,A)d\mfx=p^{sr}\int_{\prod_{i=1}^m\pi_i^{e_is}\mathcal{O}_{K_i}}\#\{\mathfrak{x}\in U\mid P_A(\mathfrak{x})=\mathfrak{y}\}d\mathfrak{y}.
\end{equation}
\end{lemma}

\begin{proof}
We identify $E=K_1\times\cdots\times K_m$ with the $\Q_p$-vector space $\Q_p^r$. Similarly, each component $K_i$ is identified with $\Q_p^{r_i}$. We recall the definition of Jacobian of a map of $\Q_p$-vector spaces from e.g. \cite[p282]{evans2006expected}. For any polynomial $f \in \Q_p[x]$ and any element $\mfx = (x_1,\ldots,x_m) \in E$, we note that $f(\mfx) = (f(x_1),\ldots,f(x_m))$, so the Jacobian $Jf$ is a block-diagonal matrix with an $r_i \times r_i$ block corresponding to each field component $K_i$. We denote this block by $J_{K_i}f$, which is just the Jacobian of the map $f: K_i \to K_i$ where we again identify $K_i$ with $\Q_p^{r_i}$. 

By \cite[Proposition 2.3]{evans2006expected}\footnote{In the notation of \cite{evans2006expected}, we are taking $f=P_A$, $g=\bbone_{\prod_{i=1}^m\pi_i^{e_is}\mathcal{O}_{K_i}}$ and $X=U$.}, we have
\begin{equation}\label{eq: Evan's formula}
\text{RHS}\eqref{eq: integral and number of solutions}=p^{sr}\int_U ||\det(JP_A(\mfx))||\bbone_{\val P_A(x_i)\ge s,\forall 1\le i\le m}d\mfx.
\end{equation}
Here, $JP_A(\mfx)$ is the scalar multiplication $P_A'(\mfx)$ over $E$, and $\det(JP_A(\mfx))$ is the determinant of $JP_A(\mfx)$, viewed as a linear map. Likewise, for each $1\le i\le m$, we denote by $J_{K_i}P_A(x_i)$ the scalar multiplication of $P_A'(x_i)$ over $K_i$, and by $\det(J_{K_i}P_A(x_i))$ its determinant. If we pick a basis of $K_i$ for each $1\le i\le m$, then these vectors together form a natural basis of $E$. This implies that
\begin{align}
\begin{split}
\label{eq: decomposition of Jacobian}
||\det(JP_A(\mfx))||&=\prod_{i=1}^m||\det(J_{K_i}P_A(x_i))||\\
&=\prod_{i=1}^m||\Nm_{K_i/\Q_p}P_A'(x_i)||\\
&=\prod_{i=1}^m||P_A'(x_i)||^{r_i}.
\end{split}
\end{align}
Combining \eqref{eq: Evan's formula} and \eqref{eq: decomposition of Jacobian}, we have
$$\text{RHS}\eqref{eq: integral and number of solutions}=\int_Up^{sr}\prod_{i=1}^m ||P_A'(x_i)||^{r_i}\cdot\bbone_{\val P_A(x_i)\ge s,\forall 1\le i\le m}dx_1\cdots dx_m=\text{LHS}\eqref{eq: integral and number of solutions},$$
which gives the proof.
\end{proof}

\begin{rmk}
The reader may note that the statements of \Cref{lem: Kac-Rice formula} and \Cref{lem: integral and number of solutions} are closely related. In fact, after justifying that
\begin{equation}
    \label{eq:justify_lemma_rmk}
    \#\{\mathfrak{x}\in U\mid P_A(\mathfrak{x})=\mathfrak{y}\}=\#\{\mathfrak{x}\in U\mid P_A(\mathfrak{x})=0\}
\end{equation}
for all $\mathfrak{y}\in\prod_{i=1}^m\pi_i^{e_is}\mathcal{O}_{K_i}$ when $s$ is sufficiently large, one may deduce \Cref{lem: Kac-Rice formula} from \Cref{lem: integral and number of solutions}. Specifically, for such $s$,
\begin{align}
\begin{split}
\int_UI_s(\mfx,A)d\mfx&=p^{sr}\int_{\prod_{i=1}^m\pi_i^{e_is}\mathcal{O}_{K_i}}\#\{\mathfrak{x}\in U\mid P_A(\mathfrak{x})=\mathfrak{y}\}d\mathfrak{y}\\
&=p^{sr}\int_{\prod_{i=1}^m\pi_i^{e_is}\mathcal{O}_{K_i}}\#\{\mathfrak{x}\in U\mid P_A(\mathfrak{x})=0\}d\mathfrak{y}\\
&=\#\{\mathfrak{x}\in U\mid P_A(\mathfrak{x})=0\},
\end{split}
\end{align}
and \Cref{lem: Kac-Rice formula} follows. However, justifying \eqref{eq:justify_lemma_rmk} encompasses most of the work in the proof of \Cref{lem: Kac-Rice formula} already.

\end{rmk}

We now turn back to our proof of \Cref{thm: joint distribution}. Our argument follows the proof of \cite[Theorem 2.6]{caruso2022zeroes}. 

\begin{proof}[Proof of \Cref{thm: joint distribution}]
To keep the notation consistent with the preceding lemmas, we continue to denote the size of the matrix by $r$, and $\sum_{i=1}^m[K_i:\Q_p]=r$, although the theorem itself uses $n$ in place of $r$. By \Cref{prop: almost surely distinct eigenvalues}, it suffices to prove the case where $U$ is a measurable set in $\mathcal{O}_E^{\new}$.

To start with, we assume that $U\subset \mathcal{O}_E^{\new}$ is compact and open. In this case, we deduce from \Cref{lem: Kac-Rice formula} that
$$\E[Z_{U,r}]=\int_{\Mat_r(\Z_p)}\lim_{s\rightarrow\infty}\int_U I_s(\mathfrak{x},A)d\mathfrak{x}dA.$$
It is clear that $I_s(\mathfrak{x},A)$ is nonnegative everywhere, and
\begin{align}
\begin{split}
\int_UI_s(\mathfrak{x},A)d\mathfrak{x}
&=p^{sr}\int_{\prod_{i=1}^m\pi_i^{e_is}\mathcal{O}_{K_i}}\#\{\mathfrak{x}\in U\mid P_A(\mathfrak{x})=\mathfrak{y}\}d\mathfrak{y}\\
&\le p^{sr}\int_{\prod_{i=1}^m\pi_i^{e_is}\mathcal{O}_{K_i}}r^r d\mathfrak{y}\\
&=r^r,
\end{split}
\end{align}
where the first line comes from \Cref{lem: integral and number of solutions}, and the second line is a naive bound since $P_A$ has $r$ roots in $\bar{\Q}_p$ and hence $r^r$ $r$-tuples of roots. Since this bound is independent of $s$, Lebesgue's dominated convergence theorem yields
\begin{equation}
    \label{eq:dct_swap_Z}
\E[Z_{U,r}]=\lim_{s\rightarrow\infty}\int_{\Mat_r(\Z_p)}\int_U I_s(\mathfrak{x},A)d\mathfrak{x}dA=\lim_{s\rightarrow\infty}\int_U\int_{\Mat_r(\Z_p)} I_s(\mathfrak{x},A)dAd\mathfrak{x},
\end{equation}
where the second equality comes from Tonelli's theorem. 

We now want to apply \Cref{thm: points on variety_main theorem} to compute the above integral when $s$ is sufficiently large. Since $\val\Disc(\cdot)$ is a continuous function over the compact set $U$, there exists $\mathfrak{y}\in U$ such that $\val\Disc(\mathfrak{y})$ takes the maximum. In this case, when $s>r\val\Disc(\mathfrak{y})$, we have
\begin{align}
\begin{split}
\int_{\Mat_r(\Z_p)} I_s(\mathfrak{x},A)dA&=\int_{\Mat_r(\Z_p)}p^{sr}\prod_{i=1}^m ||P_A'(x_i)||^{r_i}\cdot\bbone_{\val P_A(x_i)\ge s,\forall 1\le i\le m}dA\\
&=p^{sr}\prod_{i=1}^m ||Z'(x_i)||^{r_i}\int_{\Mat_r(\Z_p)}\bbone_{\val P_A(x_i)\ge s,\forall 1\le i\le m}dA\\
&=p^{sr}||\Disc(\mathfrak{x})||\mathbf{P}(\val P_A(x_i)\ge s,\forall 1\le i\le m)\\
&=\rho_E^{(r)}(\mathfrak{x}).
\end{split}
\end{align}
Here the second equality comes from \Cref{cor: norm of differential}, and the last equality comes from \Cref{thm: points on variety_main theorem}. Therefore, when $s>r\val\Disc(\mathfrak{y})$, we have
\begin{equation}\label{eq:I_to_rho}
    \int_U\int_{\Mat_r(\Z_p)} I_s(\mathfrak{x},A)dAd\mathfrak{x}=\int_U\rho_E^{(r)}(\mathfrak{x})d\mathfrak{x}.
\end{equation}
Combining \eqref{eq:dct_swap_Z} and \eqref{eq:I_to_rho} proves
\Cref{thm: joint distribution} in the case when $U$ is compact and included in $\mathcal{O}_E^{\new}$.

We now extend the result to any open subset of $\mathcal{O}_E^{\new}$. Specifically, if $U$ is open in $\mathcal{O}_E^{\new}$, we can construct an increasing sequence $(U_m)_{m\ge 1}$ of compact open subsets of $\mathcal{O}_E^{\new}$ such that $\bigcup_{m\ge 1}U_m=U$. By the case we have just shown, we have 
\begin{equation}\label{eq: expectation of number of roots in U_m}
\E[Z_{U_m,r}]=\int_{U_m}\rho_E^{(r)}(\mathfrak{x})d\mathfrak{x}.
\end{equation}
Moreover, the sequence $(Z_{U_m,r})_{m \geq 1}$ of functions on $A \in \Mat_r(\Z_p)$ is nonnegative, non-decreasing, and converges pointwise to $Z_{U,r}$. Hence we have
\begin{align}
\begin{split}
\E[Z_{U,r}]&=\lim_{m \to \infty} \E[Z_{U_m,r}]\\
&=\lim_{m \to \infty}\int_{U_m}\rho_E^{(r)}(\mathfrak{x})d\mathfrak{x}.\\
&=\int_U\rho_E^{(r)}(\mathfrak{x})d\mathfrak{x},
\end{split}
\end{align} 
where the first and last line are deduced by the monotone convergence theorem.

Following a similar approach as the above, by considering a descending sequence of open sets, we deduce that the result holds for any $G_\delta$-set (that is, a countable intersection of open sets). Recall also that any measurable set $U$ has measure $\mu_E(U) = \inf_V \mu_E(V)$, where the infimum is over all open sets $V \supset U$. Therefore, for any zero measure set $U_0$, we have a descending sequence of open sets $U_1 \supset U_2 \supset \ldots \supset U_0$ with $U_0' := \bigcap_{m \geq 1} U_m$ having measure $0$. Since $U_0' \supset U_0$ is a $G_\delta$ set, this gives $\E[Z_{U_0,r}]\le\E[Z_{U_0',r}]=0$, hence $\E[Z_{U_0,r}]=0$. 

Now, let $U\subset\mathcal{O}_E^{\new}$ be any measurable set. As before, there is a descending sequence $(U_m)_{m\ge 1}$ of open subsets of $\mathcal{O}_E^{\new}$ such that each contain $U$ and the measure of $U_m\backslash U$ converges to zero when $m$ goes to infinity. In this case, $\bigcap_{m\ge 1}U_m$ is a $G_\delta$-set that contains $U$, and $(\bigcap_{m\ge 1}U_m)\backslash U$ has measure zero. Therefore,
$$\E[Z_{U,r}]=\E[Z_{\bigcap_{m\ge 1}U_m,r}]=\int_U\rho_E^{(r)}(\mathfrak{x})d\mathfrak{x}+\int_{(\bigcap_{m\ge 1}U_m)\backslash U}\rho_E^{(r)}(\mathfrak{x})d\mathfrak{x}=\int_U\rho_E^{(r)}(\mathfrak{x})d\mathfrak{x}.$$
This ends the proof. 
\end{proof}

\subsection{Averages of correlation functions and proof of {\Cref{thm:general_cor_functions_intro}}}
Based on the form of the joint distribution $\rho_E^{(n)}$, we can write down the expectation of a function over $\Mat_n(\Z_p)$ when it depends only on the eigenvalues, which is the purpose of the following corollary. Recall from the beginning of \Cref{sec:char_poly} that $n_E(K)$ is the number of factors isomorphic to $K$ in the product decomposition $E \cong K_1 \times \cdots \times K_{m_E}$.

\begin{prop}\label{cor: expectation of function depends on eigenvalues}
Let $A\in\Mat_n(\Z_p)$ be Haar-distributed. Then for any bounded measurable function $h:\Mat_n(\Z_p)\rightarrow\R$ that only depends on the eigenvalues, we have
$$\E[h(A)]=\sum_E\int_{\mathcal{O}_E}\frac{1}{\#\Aut(E)}h(x_1,\ldots,x_{m_E})\rho_E^{(n)}(x_1,\ldots,x_{m_E})dx_1\cdots dx_{m_E}.$$
Here the sum $\sum_E$ ranges over all isomorphism classes of étale algebras $E=K_1\times\cdots\times K_{m_E}$ of degree $r$ over $\Q_p$. The symbol $h(x_1,\ldots,x_{m_E})$ denotes the value of $h$ on any matrix $A\in\Mat_n(\Z_p)$ that has $x_1,\ldots,x_{m_E}$ and their Galois conjugates as its eigenvalues.
\end{prop}

The following lemma will be useful for our proof of \Cref{cor: expectation of function depends on eigenvalues}.

\begin{lemma}\label{lem: zero as sum of aut}
Let $E:=K_1\times\cdots\times K_m$, where $K_1,\ldots,K_m$ are the same as in \Cref{thm:general_cor_functions_intro}. Let $E'=K_{m+1}\times\cdots\times K_{m+m_{E'}}$ be an étale algebra of degree $n-r$. Moreover, let
$U\subset\mathcal{O}_{K_1}^{\new}\times\cdots\times\mathcal{O}_{K_m}^{\new}$ be a measurable set, and write $Z_{U,n}(x_1,\ldots,x_{m+m_{E'}})$ for the value of $Z_{U,n}$ (defined in \Cref{thm:general_cor_functions_intro}) on any matrix $A\in\Mat_n(\Z_p)$ that has $x_1,\ldots,x_{m+m_{E'}}$ and their Galois conjugates as its eigenvalues. 

Then, for all $(x_1,\ldots,x_{m+m_{E'}})\in\mathcal{O}_{E\times E'}^{\new}$, we have
$$Z_{U,n}(x_1,\ldots,x_{m+m_{E'}})=\frac{1}{\#\Aut(E')}\sum_{\tau\in\Aut(E\times E')}\bbone_{\tau(x_1,\ldots,x_{m+m_{E'}})\in U\times\mathcal{O}_{E'}}.$$
In particular, when $n=r$, we have 
$$Z_{U,n}(x_1,\ldots,x_m)=\sum_{\tau\in\Aut(E)}\bbone_{\tau(x_1,\ldots,x_m)\in U}.$$
\end{lemma}

\begin{proof}
Let $A \in \Mat_n(\Z_p)$ be a matrix with $x_1,\ldots,x_{m+m_{E'}}$ and their Galois conjugates as its eigenvalues. Let us consider the set 
$$\mathcal{S}:=\{\tau(x_1,\ldots,x_{m+m_{E'}}):\tau\in\Aut(E\times E')\},$$
which is exactly the set of elements in $E\times E'$ that have minimal polynomial equal to $P_A$, the characteristic polynomial of $A$. 

For all $\mathfrak{y}:=(y_1,\ldots,y_m)\in U$ such that $y_1,\ldots,y_m$ are distinct, non-conjugate eigenvalues of $A$, denote by $P_{\mathfrak{y}}\in\Z_p[x]$ the minimal polynomial of $\mathfrak{y}$. Then, we have 
\begin{align}
\begin{split}
\#\{\mathfrak{y'}\in\mathcal{O}_{E'}:(\mathfrak{y},\mathfrak{y}')\in\mathcal{S}\}&=\{\mathfrak{y'}\in\mathcal{O}_{E'}:P_{\mathfrak{y}}\cdot P_{\mathfrak{y'}}=P_A\}\\
&=\{\mathfrak{y'}\in\mathcal{O}_{E'}:P_{\mathfrak{y'}}=P_A/P_{\mathfrak{y}}\}\\
&=\#\Aut(E').
\end{split}
\end{align}
Here, $P_{\mathfrak{y'}}\in\Z_p[x]$ is the minimal polynomial of $\mathfrak{y'}$, and the last line comes from \Cref{prop: aut_e and generators}. Therefore, we have
$$\#(\mathcal{S}\cap (U\times\mathcal{O}_{E'}))=\#\Aut(E')\cdot Z_{U,n}.$$
The $n=r$ case follows immediately as a special case.
\end{proof}

\begin{proof}[Proof of \Cref{cor: expectation of function depends on eigenvalues}]
Since the equality is linear in $h$ on both sides, it suffices to prove the case when $h$ (viewed as a function of the eigenvalues) is supported in some given isomorphism class $E=K_1\times\cdots\times K_m$. First, suppose $h=Z_{U,n}$, where $U\subset \mathcal{O}_{K_1}^{\new}\times\cdots\times\mathcal{O}_{K_m}^{\new}$ is a measurable set. Applying \Cref{thm: joint distribution}, we have
\begin{align}
\begin{split}
\E[Z_{U,n}]&=\int_{U}\rho_E^{(n)}(x_1,\ldots,x_{m_E})dx_1\cdots dx_{m_E}\\
&=\int_{\mathcal{O}_E}\bbone_{(x_1,\ldots,x_{m_E})\in U}\rho_E^{(n)}(x_1,\ldots,x_{m_E})dx_1\cdots dx_{m_E}\\
&=\int_{\mathcal{O}_E}\frac{1}{\#\Aut(E)}\sum_{\tau\in\Aut(E)}\bbone_{\tau(x_1,\ldots,x_{m_E})\in U}\rho_E^{(n)}(x_1,\ldots,x_{m_E})dx_1\cdots dx_{m_E}\\
&=\int_{\mathcal{O}_E}\frac{1}{\#\Aut(E)}Z_{U,n}(x_1,\ldots,x_{m_E})\rho_E^{(n)}(x_1,\ldots,x_{m_E})dx_1\cdots dx_{m_E}.
\end{split}
\end{align}
Here, the third lines holds because the integral of the summand corresponding to any $\tau\in\Aut(E)$ is the same, and the last line comes from \Cref{lem: zero as sum of aut}. 

After we verify the case $h=Z_{U,n}$ in the above, we know that the equality holds when $h$ is a simple function, i.e., a finite linear combination of indicator functions of measurable sets. Finally, since a bounded measurable function can be uniformly approximated by a sequence of simple functions, we are done.
\end{proof}

\begin{rmk}
Later in \Cref{thm: upper bound of correlation functions}, we will see that the correlation function $\rho_E^{(n)}(x_1,\ldots,x_{m_E})$ is bounded in $\mathcal{O}_E$. Therefore, we can extend \Cref{cor: expectation of function depends on eigenvalues} to the case when $h$ is integrable.
\end{rmk}

Based on the expression of the correlation function in \Cref{thm: joint distribution}, we now turn to \Cref{thm:general_cor_functions_intro}. 

\begin{proof}[Proof of \Cref{thm:general_cor_functions_intro}] 
Following the notation of \Cref{thm: joint distribution}, for all measurable set $U\subset\mathcal{O}_{K_1}^{\new}\times\cdots\times\mathcal{O}_{K_m}^{\new}$, we still denote by $Z_{U,n}$ the number of $m$-tuples $(x_1,\ldots,x_m)$ of distinct, non-conjugate eigenvalues of $A$ which lie in $U$. Since a bounded measurable function can be uniformly approximated by a sequence of simple functions,
it suffices to verify the lemma for $Z_{U,n}$ associated to an arbitrary measurable set $U\subset\mathcal{O}_{K_1}^{\new}\times\cdots\times\mathcal{O}_{K_m}^{\new}$.

Denote $E:=K_1\times\cdots \times K_m$. Then, the étale algebras over $\Q_p$ of degree $n$ that contain $E$ are isomorphic to an \'etale algebra $E\times E'$, where $E'=K_{m+1}\times\cdots\times K_{m+m_{E'}}$ is of degree $n-r$. $Z_{U,n}$ is a bounded measurable function on $\Mat_{n\times n}(\Z_p)$ that only depends on the eigenvalues, so by \Cref{cor: expectation of function depends on eigenvalues} we have
\begin{equation}\label{eq:Z_and_aut}
    \E[Z_{U,n}] = \sum_{E'}\int_{\mc{O}_E \times \mathcal{O}_{E'}}\frac{Z_{U,n}(x_1,\ldots,x_{m+m_E})}{\#\Aut(E \times E')}\rho_{E \times E'}^{(n)}(x_1,\ldots,x_{m+m_E})dx_1\cdots dx_{m+m_E}
\end{equation}
where we write $Z_{U,n}(x_1,\ldots,x_{m+m_{E'}})$ for the value of $Z_{U,n}$ on any matrix with $x_1,\ldots,x_{m+m_{E'}}$ and their Galois conjugates as its eigenvalues, and $\sum_{E'}$ ranges over all isomorphism classes of étale algebras $E'=K_{m+1}\times\cdots\times K_{m+m_{E'}}$ of degree $n-r$. By \Cref{prop: almost surely distinct eigenvalues}, with probability $1$, $A$ has no repeated eigenvalues, and so we may restrict the domain of integration to the subset where $(x_1,\ldots,x_{m+m_E}) \in \mc{O}_{E \times E'}^{\new}$ without changing the value of the integral. Hence we may use the expression in \Cref{lem: zero as sum of aut} for $Z_{U,n}$ to obtain
\begin{align}
    \begin{split}
        \text{RHS\eqref{eq:Z_and_aut}} &=\sum_{E'}\frac{1}{\#\Aut(E')}\int_{\mathcal{O}_E\times \mathcal{O}_{E'}}\sum_{\tau\in\Aut(E\times E')}\frac{\bbone_{\tau(x_1,\ldots,x_{m+m_{E'}})\in U\times\mathcal{O}_{E'}}}{\#\Aut(E\times E')}\\
&\times\rho_{E\times E'}^{(n)}(x_1,\ldots,x_{m+m_{E'}})dx_1\cdots dx_{m+m_{E'}}\\
&=\sum_{E'}\int_{\mathcal{O}_E\times \mathcal{O}_{E'}}\sum_{\tau\in\Aut(E\times E')}\frac{\bbone_{\tau(x_1,\ldots,x_{m+m_{E'}})\in U\times\mathcal{O}_{E'}}}{\#\Aut(E')\#\Aut(E\times E')}(1-p^{-1})\cdots(1-p^{-n})\\
&\times \prod_{i=1}^{m+m_{E'}}\Den(x_i)\cdot||\Delta_\sigma(x_1,\ldots,x_{m+m_{E'}})||dx_1\cdots dx_mdx_{m+1}\cdots dx_{m+m_{E'}}.
    \end{split}
\end{align}

Notice that for all $E'$ in the above sum, every étale algebra isomorphism $\tau\in\Aut(E\times E')$ contributes equally to the whole sum. Therefore, we only need to handle the case that $\tau=\id$ is the identity map, so the above implies
\begin{align}
    \begin{split}\label{eq: tau is the identity map}
   \E[Z_{U,n}]&=\sum_{E'}\int_{\mathcal{O}_E\times \mathcal{O}_{E'}}\frac{\bbone_{(x_1,\ldots,x_m)\in U}}{\#\Aut(E')}(1-p^{-1})\cdots(1-p^{-n})\\
&\times \prod_{i=1}^{m+m_{E'}}\Den(x_i)\cdot||\Delta_\sigma(x_1,\ldots,x_{m+m_{E'}})||dx_1\cdots dx_mdx_{m+1}\cdots dx_{m+m_{E'}} 
\end{split}
\end{align}
For all $E'$ in the above sum, the integral over $\mathcal{O}_E\times \mathcal{O}_{E'}$ is finite because $\E[Z_{U,n}]$ is finite, and the integrand is non-negative, thus it must be integrable. Therefore, by Fubini's theorem, we can switch the integral in \eqref{eq: tau is the identity map} into iterated integrals, i.e., we first integrate over $x_{m+1},\ldots,x_{m+m_{E'}}$, then integrate over $x_1,\ldots,x_m$. By Krasner's theorem \cite{krasner1946nombre} there are only finitely many extensions of $\Q_p$ of each degree, hence there are only finitely many isomorphism classes of étale algebras $E'$ of degree $n-r$. Therefore we can move the sum $\sum_{E'}$ past the integral over $\mc{O}_E$, so \eqref{eq: tau is the identity map} yields 
\begin{multline}\label{eq:factorized_integral}
    \E[Z_{U,n}]= \int_{U} (1-p^{r-n-1}) \cdots (1-p^{-n}) ||\Delta_\sigma(x_1,\ldots,x_m)|| \left(\prod_{i=1}^{m} V(x_i)\right) (1-p^{-1}) \cdots (1-p^{r-n}) \\ 
\times \left(\sum_{E'} \int_{\mc{O}_{E'}} \frac{1}{\#\Aut(E')} \frac{||\Delta_\sigma(x_1,\ldots,x_{m+m_{E'}})||}{||\Delta_\sigma(x_1,\ldots,x_m)||} \prod_{i=m+1}^{m_{E'}} V(x_i)  dx_{m+1} \cdots dx_{m+m_{E'}}\right) dx_1 \cdots dx_m. 
\end{multline}
We will now write the inner sum in terms of the expected determinant. Notice that when $E'=K_{m+1}\times\cdots\times K_{m+m_{E'}}$, and a fixed matrix $\tilde A\in\Mat_{n-r}(\Z_p)$ has $x_{m+1}\in\mathcal{O}_{K_{m+1}}^{\new},\ldots,x_{m+m_{E'}}\in\mathcal{O}_{K_{m+m_{E'}}}^{\new}$ and their Galois conjugates as its eigenvalues, we have
$$||\det(Z(\tilde A))||=\prod_{i=m+1}^{m+m_E}||\Nm_{K_i/\Q_p}Z(x_i)||=\frac{||\Delta_{\sigma}(x_1,\ldots,x_{m+m_{E'}})||}{||\Delta_{\sigma}(x_1,\ldots,x_m)||\cdot||\Delta_{\sigma}(x_{m+1},\ldots,x_{m+m_{E'}})||}.$$
Therefore, applying \Cref{cor: expectation of function depends on eigenvalues} by considering $h(\tilde A)=||\det(Z(\tilde A))||$ where $\tilde A\in\Mat_{n-r}(\Z_p)$ is Haar-distributed, we have
\begin{align}\label{eq:write_out_E_det}
\begin{split}
\E||\det(Z(\tilde A))||&=\sum_{E'}\int_{\mathcal{O}_{E'}}\frac{1}{\#\Aut(E')}\cdot\frac{||\Delta_{\sigma}(x_1,\ldots,x_{m+m_{E'}})||}{||\Delta_{\sigma}(x_1,\ldots,x_m)||\cdot||\Delta_{\sigma}(x_{m+1},\ldots,x_{m+m_{E'}})||}\\
&\times\rho_{E'}^{(n-r)}(x_{m+1},\ldots,x_{m+m_{E'}})dx_{m+1}\cdots dx_{m+m_{E'}}\\
&=\sum_{E'}\int_{\mathcal{O}_{E'}}\frac{(1-p^{-1})\cdots(1-p^{r-n})}{\#\Aut(E')}\prod_{i=m+1}^{m+m_E'}\Den(x_i)\\
&\times \frac{||\Delta_{\sigma}(x_1,\ldots,x_{m+m_{E'}})||}{||\Delta_{\sigma}(x_1,\ldots,x_m)||}dx_{m+1}\cdots dx_{m+m_{E'}}.
\end{split}
\end{align}
Combining \eqref{eq:factorized_integral} with \eqref{eq:write_out_E_det} completes the proof.
\end{proof}

%% file: Expected_determinants.tex
\section{Expected determinants, cokernels, and the moment method}\label{sec:cokernels}

The main results of the section, Theorems~\ref{thm:compute_char_poly_det} and~\ref{thm: decomposition of limit expectation}, concern the limit of the expected determinants $\lim_{n\rightarrow\infty}\E||\det(Z(A))||$. \Cref{thm: decomposition of limit expectation} states that under some restricted hypotheses these expectations factorize, and is behind the independence of eigenvalues in different lifted subspaces in \Cref{thm:lifted_subspace_intro_revised}. \Cref{thm:compute_char_poly_det} is true more generally, but simply reduces the desired expectations to expectations of a product of correlated random variables. This `linearization' is very useful, but even after linearization the expectations will still be quite difficult to compute. The proofs of both results rely on the moment method machinery for random modules due to Sawin-Wood \cite{sawin2022moment}, as applied by Cheong-Yu \cite{cheong2023distribution}. 

\begin{thm}\label{thm:compute_char_poly_det}
Let $n\ge 1$, and $A$ be distributed by the additive Haar measure on $\Mat_n(\Z_p)$. Suppose $Z=\prod_{i=1}^m Z_i(x)\in\Z_p[x]$ has degree $r$, where the polynomials $Z_i(x)\in\Z_p[x],1\le i\le m$ are monic, distinct, and irreducible. Let $r_i$ be the degree of $Z_i(x)$, so that $r=r_1+\cdots+r_m$, and $x_i$ be a root of $Z_i(x)$. Then the limit $\lim_{n\rightarrow\infty}\E||\det(Z(A))||^k$ exists for any $k\geq 0$, and is given by
\begin{equation}
    \label{eq:compute_char_poly_det}
    \lim_{n\rightarrow\infty}\E[||\det(Z(A))||^k]=\lim_{n\rightarrow\infty}\E\left[\prod_{i=1}^m||\det(A_0+x_iA_1+\cdots+x_i^{r-1}A_{r-1})||^{r_ik}\right],
\end{equation}
where $A_0,\ldots,A_{r-1}$ are i.i.d. Haar-distributed over $\Mat_n(\Z_p)$. 
\end{thm}

\begin{thm}\label{thm: decomposition of limit expectation}
Let $n\ge 1$, and $A$ be distributed by the additive Haar measure on $\Mat_n(\Z_p)$. Let $Z=\prod_{i=1}^m Z_i(x)\in\Z_p[x]$, where the polynomials $Z_i(x)\in\Z_p[x],1\le i\le m$ are monic and do not have repeated roots. Furthermore, suppose that the residues $\bar Z_i(x)\in\F_p[x],1\le i\le m$ are powers of distinct irreducible polynomials in $\F_p[x]$. Then, when $n$ tends to infinity, the joint random variable 
$$(||\det(Z_1(A))||,\ldots,||\det(Z_m(A))||)$$ 
converges in distribution to a tuple of $m$ independent random variables, thus
\begin{equation}\label{eq: decomposition of limit expectation}
\lim_{n\rightarrow\infty}\E[||\det (Z(A))||^k]=\prod_{i=1}^m\lim_{n\rightarrow\infty}\E[||\det(Z_i(A))||^k]
\end{equation}
for all $k\ge 0$.
\end{thm}

We remark that the limit is necessary and the equalities in the above theorems are not true for fixed $n$. 

\subsection{Algebraic preparations}

In this subsection, we continue to work with the same notations as in \Cref{sec:char_poly} and \Cref{sec: joint distribution}, and we also abbreviate $R:=\Z_p[x]/(Z(x))\cong\Z_p[\mfx]$. Recall that for a ring $L$, the \emph{cokernel} of a matrix $A \in \Mat_n(L)$ is 
\begin{equation}
    \Cok_L(A) := L^n/AL^n.
\end{equation}
Sometimes we wish to make $A$ act on $M^n$ where $M$ is an $L$-module, and write $\Cok_M(A) := M^n/AM^n$. We will often simply write $\Cok$ when the subscript is clear from context.

For $A \in \Mat_n(\Z_p)$, abelian $p$-group
\begin{equation}
    \Cok_{\Z_p}(Z(A)) = \Z_p^n/Z(A)\Z_p^n
\end{equation}
also carries the structure of an $R$-module for $R$ as defined above, where $x \in R$ acts by $A$. 

Moreover, as pointed out by Lee \cite{lee2023joint}, we have
\begin{equation}\label{eq:Z_vs_linear_iso}
     \Cok_{\Z_p}(Z(A)) \cong \Cok_{\Z_p[\mfx]}(A-\mfx I_n),
 \end{equation}
where the isomorphism is as $R$-modules. We will use both expressions in \eqref{eq:Z_vs_linear_iso} interchangeably; the former naturally connects to the desired determinants, while the latter is expressed in a linear way and easier to work with. 

\begin{lemma}\label{lem: cok and det}
Let $K$ be an extension of $\Q_p$ of degree $r$, $\mathcal{O}_K$ be its ring of integers, and $A\in\Mat_n(\mathcal{O}_K)$ be nonsingular. Then, we have $\#\Cok(A)=||\det(A)||^{-r}$.
\end{lemma}

\begin{proof}
Let $\pi$ be a uniformizer of $\mathcal{O}_K$, and
$e$ be the ramification index of $K/\Q_p$, so that $\val(\pi)=\frac{1}{e}$, and $\#(\mathcal{O}_K:\pi\mathcal{O}_K)=p^{r/e}$. By the Smith normal form (i.e., after multiplying by elements of $\GL_n(\mathcal{O}_K)$ on both sides), there is no loss of generality in assuming that
$$A=\diag(\pi^{\l_1},\ldots,\pi^{\l_n}),\quad\l_1\ge\ldots\ge\l_n\ge0.$$
In this case, we have 
$$\#\Cok(A)=\#(\mathcal{O}_K^n:A\mathcal{O}_K^n)=\prod_{i=1}^n\#(\mathcal{O}_K:\pi^{\l_i}\mathcal{O}_K)=\prod_{i=1}^n p^{r\l_i/e}=p^{r(\l_1+\cdots+\l_n)/e}.$$
Also, we have $\det(A)=\pi^{\l_1+\cdots+\l_n}$, and thus $||\det(A)||^{-r}=p^{r(\l_1+\cdots+\l_n)/e}$. This completes the proof.
\end{proof}

The following lemma is a useful intermediate step. 

\begin{lemma}\label{cor: Cok as product of determinant over fields}
Let $A_\mfx=A_0+A_1\mfx+\cdots+A_{r-1}\mfx^{r-1}\in\Mat_n(\Z_p[\mfx])$ be nonsingular, where $\mfx=(x_1,\ldots,x_m)\in \mathcal{O}_{K_1}^{\new}\times\cdots\times \mathcal{O}_{K_m}^{\new}$, and $A_0,\ldots,A_{r-1}\in\Mat_n(\Z_p)$. Then we have
$$\#\Cok(A_\mfx)=\prod_{i=1}^m||\det(A_0+A_1x_i+\cdots+A_{r-1}x_i^{r-1})||^{-r_i}.$$
\end{lemma}

\begin{proof}
We first pass to cokernels over $\mc{O}_E$, which have better factorization properties than cokernels over $\Z_p[\mfx]$. We know $\#(\mathcal{O}_E:\Z_p[\mfx])<\infty$ is finite. If we write
$$a_1+\Z_p[\mfx]^n,\ldots,a_m+\Z_p[\mfx]^n$$
as the cosets of $\Z_p[\mfx]^n$ in $\mathcal{O}_E^n$, then
$$A_\mfx(a_1+\Z_p[\mfx]^n),\ldots,A_\mfx(a_m+\Z_p[\mfx]^n)$$
gives all the cosets of $A_\mfx\Z_p[\mfx]^n$ in $A_\mfx\mathcal{O}_E^n$. Therefore, we have $\#(A_\mfx\mathcal{O}_E^n:A_\mfx\Z_p[\mfx]^n)=\#(\mathcal{O}_E^n:\Z_p[\mfx]^n)$, and
\begin{align}\label{eq: cok to field factors}
\begin{split}
\#\Cok(A_\mfx)&=\#(\Z_p[\mfx]^n:A_\mfx\Z_p[\mfx]^n)\\
&=\#(\mathcal{O}_E^n:A_\mfx\mathcal{O}_E^n)\frac{\#(A_\mfx\mathcal{O}_E^n:A_\mfx\Z_p[\mfx]^n)}{\#(\mathcal{O}_E^n:\Z_p[\mfx]^n)}\\
&=\#(\mathcal{O}_E^n:A_\mfx\mathcal{O}_E^n)\\
&=\prod_{i=1}^m\#(\mathcal{O}_{K_i}^n:(A_0+A_1x_i+\cdots+A_{r-1}x_i^{r-1})\mathcal{O}_{K_i}^n).
\end{split}
\end{align}
Here, the last line follows because $\mc{O}_E \cong \prod_{i=1}^m \mc{O}_{K_i}$ and $A_{\mfx}$ acts on the $i\tth$ factor by $A_0+A_1x_i+\cdots+A_{r-1}x_i^{r-1}$. Furthermore, by \Cref{lem: cok and det}, we have
\begin{align}\label{eq: field factors to det}
\begin{split}
\#(\mathcal{O}_{K_i}^n/(A_0+A_1x_i+\cdots+A_{r-1}x_i^{r-1})\mathcal{O}_{K_i}^n)&=\#\Cok_{\mathcal{O}_{K_i}}(A_0+A_1x_i+\cdots+A_{r-1}x_i^{r-1})\\
&=||\det(A_0+A_1x_i+\cdots+A_{r-1}x_i^{r-1})||^{-r_i}.\\
\end{split}
\end{align}
\eqref{eq: cok to field factors} and \eqref{eq: field factors to det} together give the proof.
\end{proof}

\subsection{Proof of \Cref{thm:compute_char_poly_det}}
The goal of this subsection is to prove \Cref{thm:compute_char_poly_det}. We assume the same notations as given in the previous preparations subsection. The following lemma lets us replace $A-\mfx I_n$ by a matrix with i.i.d. entries that are Haar-distributed.

\begin{lemma}\label{lem: Cok of Haar matrix with mfx}
For all finite $\Z_p[\mfx]$-modules $M$, we have
$$\lim_{n\rightarrow\infty}\mathbf{P}(\Cok(A-\mathfrak{x}I_n)\cong M)=\lim_{n\rightarrow\infty}\mathbf{P}(\Cok(A_\mathfrak{x})\cong M),$$
where $A$ is Haar-distributed in $\Mat_n(\Z_p)$, and $A_\mathfrak{x}$ is Haar-distributed in $\Mat_n(\Z_p[\mfx])$. In other words, the random cokernels $\Cok(A-\mathfrak{x}I),\Cok(A_\mathfrak{x})$ converge in distribution to the same limit.
\end{lemma}

\begin{proof}
Our proof is based on the moment method over random $\Z_p[\mfx]$-modules. For any finite $\Z_p[\mfx]$-module $M$, by \cite[Theorem 1.12]{cheong2023distribution}\footnote{This theorem concerns $(\Z/p^k\Z)[\mfx]$-modules, but taking $k$ large enough so that $p^k M = 0$ yields the form that we use.} 
$$\lim_{n\rightarrow\infty}\E[\#\Sur_{\Z_p[\mfx]}(\Cok(A-\mfx I_n),M)]=1.$$
Here $\#\Sur_{\Z_p[\mfx]}(\cdot,\cdot)$ denotes the number of $\Z_p[\mfx]$-module epimorphisms. By \cite[Lemma 6.3]{sawin2022moment}, the moments $1$ are well-behaved (a growth condition defined in the same paper), and hence by \cite[Theorem 1.8]{sawin2022moment} we have that convergence of moments to $1$ implies convergence in distribution of the random module to the limit distribution characterized by these moments. Hence it suffices to prove $$\lim_{n\rightarrow\infty}\E[\#\Sur_{\Z_p[\mfx]}(\Cok(A_\mfx),M)]=1$$
for all finite $\Z_p[\mfx]$-modules $M$.

To prove this, notice that
\begin{align}\label{eq: expection of surjection of Haar matrix}
\begin{split}
\E[\#\Sur_{\Z_p[\mfx]}(\Cok(A_\mfx),M)]&=\int_{\Mat_n(\Z_p[\mfx])}\#\Sur_{\Z_p[\mfx]}(\Cok(A_\mfx),M)dA_\mfx\\
&=\int_{\Mat_n(\Z_p[\mfx])}\sum_{F\in\Sur_{\Z_p[\mfx]}(\Cok(A_\mfx),M)}1dA_\mfx\\
&=\sum_{F\in\Sur_{\Z_p[\mfx]}(\Z_p[\mfx]^n,M)}\mathbf{P}(F(A_\mfx \Z_p[\mfx]^n)=0).
\end{split}
\end{align}
Since the columns of $A_\mfx$ are independent, we have
\begin{align}
\begin{split}
\mathbf{P}(F(A_\mfx \Z_p[\mfx]^n)=0)&=\mathbf{P}(F(\omega)=0)^n\\
&=\mathbf{P}(\omega\in\Ker F)^n\\
&=\frac{1}{\#(\Z_p[\mfx]^n/\Ker F)^n}\\
&=\frac{1}{(\#M)^n},
\end{split}
\end{align}
where $\omega\in\Z_p[\mfx]^n$ is Haar-distributed. Since $\#\Hom_{\Z_p[\mfx]}(\Z_p[\mfx]^n,M) = (\#M)^n$, we have 
\begin{align}
\begin{split}
\text{RHS}\eqref{eq: expection of surjection of Haar matrix}&=\sum_{F\in\Sur_{\Z_p[\mfx]}(\Z_p[\mfx]^n,M)}\frac{1}{(\#M)^n}\\
&=1-\sum_{F\in\Hom_{\Z_p[\mfx]}(\Z_p[\mfx]^n,M)\backslash\Sur_{\Z_p[\mfx]}(\Z_p[\mfx]^n,M)}\frac{1}{(\#M)^n}.\\
\end{split}
\end{align}
The above expression is less than $1$. On the other hand, since any non-surjective homomorphism has some image $N \subsetneqq M$, we have
\begin{align}
\begin{split}
\sum_{F\in\Hom_{\Z_p[\mfx]}(\Z_p[\mfx]^n,M)\backslash\Sur_{\Z_p[\mfx]}(\Z_p[\mfx]^n,M)}\frac{1}{(\#M)^n}&\le\sum_{N\subsetneqq M}\sum_{F\in\Hom_{\Z_p[\mfx]}(\Z_p[\mfx]^n,N)}\frac{1}{(\#M)^n}.\\
&=\sum_{N\subsetneqq M}\frac{(\#N)^n}{(\#M)^n}\\
&\le\#\{N\subsetneqq M\mid N \text{ is a $\Z_p[\mfx]$-module}\}\frac{1}{2^n}.
\end{split}
\end{align}
The limit as $n \to \infty$ of this quantity is $0$, and thus $\lim_{n\rightarrow\infty}\E[\#\Sur_{\Z_p[\mfx]}(\Cok(A_\mfx),M)]=1$. This completes the proof.
\end{proof}

\begin{proof}[Proof of \Cref{thm:compute_char_poly_det}]
By \Cref{lem: cok and det}, we have
$$||\det(Z(A))||^k=\#\Cok(Z(A))^{-k}=\#\Cok(A-\mfx I_n)^{-k}.$$
When $n$ goes to infinity, by \Cref{lem: Cok of Haar matrix with mfx}, the random variables
$$\#\Cok(A-\mfx I_n)^{-k},\#\Cok(A_0+A_1\mfx+\cdots+A_{r-1}\mfx^{r-1})^{-k}$$
weakly converge to the same distribution. Since these random variables are always less or equal to $1$, we have
\begin{align}
\begin{split}
\lim_{n\rightarrow\infty}\E||\det(Z(A))||^k&=\lim_{n\rightarrow\infty}\E\left[\#\Cok(A-\mfx I_n)^{-k}\right]\\
&=\lim_{n\rightarrow\infty}\E\left[\#\Cok(A_0+A_1\mfx+\cdots+A_{r-1}\mfx^{r-1})^{-k}\right]\\
&=\lim_{n\rightarrow\infty}\E\prod_{i=1}^m||\det(A_0+x_iA_1+\cdots+x_i^{r-1}A_{r-1})||^{r_ik},
\end{split}
\end{align}
where the last line comes from \Cref{cor: Cok as product of determinant over fields}. This completes the proof.
\end{proof}

\subsection{Proof of \Cref{thm: decomposition of limit expectation}}
In this subsection, we prove \Cref{thm: decomposition of limit expectation}. We will work with the following setting, which is slightly different from the version in the previous subsections:
\begin{enumerate}
\item $Y_1(x),\ldots,Y_m(x)\in\Z_p[x]$ be monic and the residues $\overline Y_i(x)\in\F_p[x],1\le i\le m$ are distinct and irreducible.
\item $d_i\ge 1$ be degree of $Y_i$. 
\item $Z_1(x),\ldots,Z_m(x)\in\Z_p[x]$ be monic, such that for all $1\le i\le m$, the residue $\overline Z_i(x)\in\F_p[x]$ is a power of $\overline Y_i$. Here, we require the polynomials $Z_i(x),1\le i\le m$ have no repeated roots, but we do not require them to be irreducible.
\item $r_i$ be the degree of $Z_i$, which is divisible by $d_i$.
\item $Z(x)=\prod_{i=1}^mZ_i(x)$, which is monic of degree $r:=r_1+\cdots+r_m$.
\item $R=\Z_p[x]/(Z(x))$ be the quotient ring.
\item $R_i:=\Z_p[x]/(Z_i(x))$, so that we can write $R$ as the product $R=R_1\times\cdots\cdots R_m$.
\item $\mathfrak{m}_i=(p,Y_i(x))/(Z_i(x))$ be the unique maximal ideal of the ring $R_i$. 
\item $\F_{p^{d_i}}\cong\Z_p[x]/(p,Y_i(x))=R_i/\mathfrak{m}_i$ be the corresponding quotient field.
\end{enumerate}

The following proposition is essentially contained in \cite[Theorem 1.3]{cheong2023distribution}, which is furthermore a concrete example of \cite[Lemma 6.3]{sawin2022moment}. It describes the limit distribution of $\Cok(A-\mfx I_{n})$.

\begin{prop}\label{thm: limit distribution of cok}
For a finite $R$-module $M$, we have
$$\lim_{n\rightarrow\infty}\mathbf{P}(\Cok(Z(A))\cong M)=\frac{1}{\#\Aut_{R}(M)}\prod_{j=1}^m\prod_{i=1}^\infty(1-p^{-id_j})$$
if there exists a nonnegative integer $s\ge 0$ and a short exact sequence $$0\rightarrow R^s\rightarrow R^s\rightarrow M\rightarrow 0$$ (in other words, $M$ is the quotient of two free $R$-modules of the same rank); otherwise, we have $\lim_{n\rightarrow\infty}\mathbf{P}(\Cok(Z(A))\cong M)=0$.
\end{prop}

The following lemma is a necessary preparation for the proof of \Cref{thm: limit distribution of cok}.

\begin{lemma}\label{lem: ext equals hom when free quotient}
Suppose there exists a nonnegative integer $s\ge 0$ and a short exact sequence $$0\rightarrow R^s\rightarrow R^s\rightarrow M\rightarrow 0.$$
Then we have $\#\Ext^1_R(M,\F_{p^{d_j}})=\#\Hom_R(M,\F_{p^{d_j}})$.
\end{lemma}

\begin{proof}
Since $R^s$ is free, we always have
$$\Ext^1_R(R^s,\F_{p^{d_j}})=0.$$
Applying the $\Ext$ functor to the short exact sequence $0\rightarrow R^s\rightarrow R^s\rightarrow M\rightarrow 0$, we derive the following exact sequence of $R$-modules:
$$0\rightarrow\Hom_R(M,\F_{p^{d_j}})\rightarrow\Hom_R(R^s,\F_{p^{d_j}})\rightarrow\Hom_R(R^s,\F_{p^{d_j}})\rightarrow\Ext_R^1(M,\F_{p^{d_j}})\rightarrow 0.$$
Therefore,
$$\#\Ext^1_R(M,\F_{p^{d_j}})=\frac{\#\Hom_R(R^s,\F_{p^{d_j}})}{\#\Hom_R(R^s,\F_{p^{d_j}})}\#\Hom_R(M,\F_{p^{d_j}})=\#\Hom_R(M,\F_{p^{d_j}}),$$
which ends the proof.
\end{proof}

\begin{proof}[Proof of \Cref{thm: limit distribution of cok}]
On the one hand, by definition, for every fixed $n\ge 1$ and $A\in\Mat_n(\Z_p)$, $\Cok(Z(A))$ is the quotient of two free $R$-modules of rank $n$, so for those $M$ that cannot be expressed as the quotient of two free $R$-modules, we must have $\lim_{n\rightarrow\infty}\mathbf{P}(\Cok(A-\mfx I_{n})\cong M)=0$. On the other hand, suppose $M$ comes from a short exact sequence $$0\rightarrow R^s\rightarrow R^s\rightarrow M\rightarrow 0$$ 
with $s\ge 0$. By \cite[Theorem 1.3]{cheong2023distribution}, we have
\begin{equation}
\lim_{n\rightarrow\infty}\mathbf{P}(\Cok(Z(A))\cong M)=\frac{1}{\#\Aut_R(M)}\prod_{j=1}^m\prod_{i=1}^\infty\left(1-\frac{\#\Ext^1_R(M,\F_{p^{d_j}})}{\#\Hom_R(M,\F_{p^{d_j}})}p^{-id_j}\right).
\end{equation}
Applying \Cref{lem: ext equals hom when free quotient}, we complete the proof.
\end{proof}

\begin{cor}\label{cor: limit expectation}
For all integers $k\ge 0$, we have
$$\lim_{n\rightarrow\infty}\E||\det(Z( A))||^k=\prod_{j=1}^m\prod_{i=1}^\infty(1-p^{-id_j}) \cdot \sum_M\frac{1}{\#\Aut_R(M)(\#M)^k},$$
where the sum ranges over all isomorphism classes of finite $R$-modules $M$ that could be expressed as the quotient of two free $R$-modules of the same rank.
\end{cor}

\begin{proof}
By \Cref{lem: cok and det}, we always have $||\det(Z(A))||^k=\frac{1}{\#\Cok(Z(A))^k}$. Then, the proof follows from \Cref{thm: limit distribution of cok}.
\end{proof}

\begin{proof}[Proof of \Cref{thm: decomposition of limit expectation}]
Notice that the number of module automorphisms satisfies the relation
\begin{equation}\label{eq: decomposition of number of automorphisms}
\#\Aut_{R_1\times\cdots\times R_m}(M_1\times\cdots\times M_m)=\prod_{i=1}^m\#\Aut_{R_i}(M_i).
\end{equation}
Here for all $1\le i\le m$, $M_i$ is a finite $R_i$-module. Thus, by applying \Cref{thm: limit distribution of cok} twice, we have
\begin{align}
\begin{split}
\lim_{n\rightarrow\infty}\mathbf{P}(\Cok(Z_i(A))\cong M_i,\forall 1\le i\le m)&=\lim_{n\rightarrow\infty}\mathbf{P}(\Cok(Z(A))\cong M_1\times\cdots\times M_m)\\
&=\frac{1}{\#\Aut_{R_1\times\cdots\times R_m}(M_1\times\cdots\times M_m)}\prod_{j=1}^m\prod_{i=1}^\infty(1-p^{-id_j})\\
&=\prod_{j=1}^m\frac{1}{\#\Aut_{R_j}(M_j)}\prod_{i=1}^\infty(1-p^{-id_j})\\
&=\prod_{i=1}^m\lim_{n\rightarrow\infty}\mathbf{P}(\Cok(Z_i(A))\cong M_i).
\end{split}
\end{align}
By \Cref{lem: cok and det}, we have
$||\det(Z_i(A))||=\frac{1}{\#\Cok(Z_i(A))}$ for all $1\le i\le m$. Therefore, the proof follows.
\end{proof}

\begin{rmk}\label{rmk: pairwise coprime independence}
The independence statement in \Cref{thm: decomposition of limit expectation} also extends to the case when $\bar Z_1,\ldots,\bar Z_m\in\F_p[x]$ are not powers of distinct irreducible polynomials anymore, but still pairwise coprime. Indeed, applying Hensel's lemma in \Cref{lem: Hensel}, we can further decompose $Z_1,\ldots,Z_m$ as
$$Z_i=\prod_{l_i=1}^{j_i}Z_{i,l_i},\forall 1\le i\le m,$$
where the residues $\bar Z_{1,1},\ldots,\bar Z_{1,j_1},\ldots,\bar Z_{m,1},\ldots,\bar Z_{m,j_m}\in\F_p[x]$ are powers of distinct irreducible polynomials. In this case, the independence assertion for the case $Z_1,\ldots,Z_m$ immediately follows from the assertion over $Z_{1,1},\ldots, Z_{1,j_1},\ldots,Z_{m,1},\ldots, Z_{m,j_m}$.
\end{rmk}

In particular, let us consider the case where every $Z_i(x)$ is the minimal polynomial of some $x_i\in\bar\Z_p$ such that $\Z_p[x_i]=\mc{O}_{K_i}$, where $K_i:=\Q_p[x_i]$. In this case, the elements $x_1,\ldots,x_m$ lie in different lifted subspaces, and the polynomials $Z_1,\ldots,Z_m$ are irreducible. Under such assumptions, the limit of the expected determinant has a simple form:

\begin{prop}\label{prop: limit expectation over generator}
Suppose that for all $1\le i\le m$, $Z_i(x)$ is the minimal polynomial of some $x_i\in\bar\Z_p$ such that $\Z_p[x_i]=\mc{O}_{K_i}$, where $K_i:=\Q_p[x_i]$. Then for all integers $k\ge 0$, we have $$\lim_{n\rightarrow\infty}\E||\det(Z(A))||^k=\prod_{i=1}^m\prod_{j=1}^k(1-p^{-jd_i}).$$
\end{prop}

\begin{proof} 
Let $E:=K_1\times\cdots\times K_m$. In this case,  we have $R_i\cong\mc{O}_{K_i}$ for all $1\le i\le m$, and $R=R_1\times\cdots\times R_m\cong \mc{O}_E$. Note that every finite $\mathcal{O}_E$-module $M$ is isomorphic to one of the form
$$M\cong M_1\times\cdots\times M_m$$
where for all $1\le i\le m$, $M_i\cong\prod_{j=1}^\infty (\mathcal{O}_{K_i}/\pi_i^{\lambda^i_j}\mathcal{O}_{K_i})$ is a finite $\mathcal{O}_{K_i}$-module determined by some $\la^i_1 \geq \la^i_2 \geq \ldots \geq 0$ (with only finitely many nonzero). This implies that every finite $\mathcal{O}_E$-module is the quotient of two free $\mathcal{O}_E$-module of the same rank. Applying \cite[Lemma 6.3]{sawin2022moment} by taking $u=k$ and recalling \Cref{lem: ext equals hom when free quotient}, we deduce that
$$\nu(M)=\frac{\prod_{j=1}^m\prod_{i=1}^\infty(1-p^{-(i+k)d_j})}{\#\Aut_{\mathcal{O}_E}(M)(\#M)^k}$$
gives a probability measure over all finite $\mathcal{O}_E$-modules. Therefore, we have $$\sum_M\frac{1}{\#\Aut_{\mathcal{O}_E}(M)(\#M)^k}=\prod_{j=1}^m\prod_{i=1}^\infty(1-p^{-(i+k)d_j})^{-1}.$$ 
Combining the above equality with \Cref{cor: limit expectation} gives the proof.
\end{proof}

%% file: cor_function_asymptotic.tex
\section{Limiting correlation functions and their properties} \label{sec:limit_cor_func}

This section is mainly devoted to proofs of the following basic results on the limiting correlation functions $\rho_E^{(\infty)}$ using results of the previous section. In proving these, we will also establish some bounds on correlation functions which are useful later. We recall that for an étale algebra $E = K_1 \times \cdots \times K_m$, the notation $\rho^{(\infty)}_E$ is interchangeable with the notation $\rho^{(\infty)}_{K_1,\ldots,K_m}$ used in the Introduction.

\begin{thm}\label{thm:limit_cor_fns_exist}
Let $K_1,\ldots,K_m$ be finite extensions of $\Q_p$, and $E := K_1\times\cdots\times K_m$. Moreover, let $(x_1,\ldots,x_m)\in\mathcal{O}_{K_1}^{\new}\times\cdots\times\mathcal{O}_{K_m}^{\new}$. Then the limits
\begin{equation}
    \rho_{E}^{(\infty)}(x_1,\ldots,x_m) := \lim_{n\rightarrow\infty}\rho_E^{(n)}(x_1,\ldots,x_m)
\end{equation}
exist, where $\rho_E^{(n)}(x_1,\ldots,x_m)=\rho_{K_1,\ldots,K_m}^{(n)}(x_1,\ldots,x_m)$ is as defined in \Cref{thm:general_cor_functions_intro}, and are given explicitly by
\begin{align}
    \begin{split}
       \rho_{E}^{(\infty)}(x_1,\ldots,x_m) =||\Delta_\sigma(x_1,\ldots,x_m)|| \cdot \left(\prod_{i=1}^m\Den(x_i) \right) \cdot \lim_{n\rightarrow\infty}\E||\det(Z(A))||.
    \end{split}
\end{align}
Here, $A\in\Mat_n(\Z_p)$ is Haar-distributed, and $Z(x)=\prod_{i=1}^m Z_i(x)$, where $Z_i(x)\in\Z_p[x]$ is the minimal polynomial of $x_i$.
\end{thm}

\begin{thm}\label{thm: limit of expectation of zeros}
Let $E:=K_1\times\cdots\times K_m$, where $K_1,\ldots,K_m$ be finite extensions of $\Q_p$. Then, for all measurable $U\subset\mathcal{O}_{K_1}^{\new}\times\cdots\times\mathcal{O}_{K_m}^{\new}$, 
$$\lim_{n\rightarrow\infty}\E[Z_{U,n}]=\int_U\rho_{E}^{(\infty)}(x_1,\ldots,x_m) dx_1 \cdots dx_m.$$
\end{thm}

\begin{thm}\label{thm: independent distribution over lifted subspaces}
Let $k_1<k_2<\ldots<k_m$ be positive integers. For each $i$ let $x_{i,1},\ldots,x_{i,j_i}$ be a sequence of elements in $\mathcal{U}_{k_i}$, which gives the corresponding étale algebra $E_i=\Q_p[(x_{i,1},\ldots,x_{i,j_i})]$. Let $E = E_1 \times \cdots \times E_m$. Then the limiting correlation functions of \Cref{thm:limit_cor_fns_exist} satisfy
\begin{equation}\label{eq: decomposition of correlation function}
\rho_E^{(\infty)}(x_{1,1},\ldots,x_{1,j_1},\ldots,x_{m,1},\ldots,x_{m,j_m})=\prod_{i=1}^m\rho_{E_i}^{(\infty)}(x_{i,1},\ldots,x_{i,j_i}).
\end{equation}
\end{thm}

\Cref{thm: independent distribution over lifted subspaces} shows that local statistics of eigenvalues in different lifted subspaces are asymptotically independent.

\begin{proof}[Proof of {\Cref{thm:limit_cor_fns_exist}}]
We only need to prove the case where the elements $x_1,\ldots,x_m\in\bar\Z_p$ lie in distinct Galois orbits, since otherwise the limit must be zero. Recall from \Cref{thm:compute_char_poly_det} that the limit $\lim_{n\rightarrow\infty}\E||\det(Z(A))||$ exists. Therefore, this theorem can be directly deduced from the expression of the finite size correlation functions $\rho_E^{(n)}(x_1,\ldots,x_m)$ in \Cref {thm:general_cor_functions_intro}. 
\end{proof}

\Cref{thm: limit of expectation of zeros} may appear to be evident from \Cref{thm:limit_cor_fns_exist} at first glance, but we need to verify that we can indeed interchange the order of integration and taking limits, which is the purpose of the following proposition.

\begin{prop}\label{thm: upper bound of correlation functions}
Let $E:=K_1\times\cdots\times K_m$, where $K_1,\ldots,K_m$ be finite extensions of $\Q_p$. Moreover, let $(x_1,\ldots,x_m)\in\mathcal{O}_{K_1}^{\new}\times\cdots\times\mathcal{O}_{K_m}^{\new}$. 
For all $1\le i\le m$, let $r_i:=[K_i:\Q_p]$,  $e_i$ be the ramification index of $K_i/\Q_p$, and $Z_i\in\Z_p[x]$ be the minimal polynomial of $x_i$. Then, we have 
\begin{equation}
\rho_E^{(n)}(x_1,\ldots,x_m)\le||\Disc_{E/\Q_p}||\cdot\prod_{1\le i\le m}\frac{1}{1-p^{-r_i/e_i}}\cdot\prod_{1\le i<j\le m}||\Res(Z_i,Z_j)||.
\end{equation}
\end{prop}

\begin{proof}[Proof of {\Cref{thm: limit of expectation of zeros}}, assuming \Cref{thm: upper bound of correlation functions}]
For all $1\le i\le m$, let $r_i:=[K_i:\Q_p]$, and $e_i$ be the ramification index of $K_i/\Q_p$. By \Cref{thm: upper bound of correlation functions}, for all $(x_1,\ldots,x_m)\in\mathcal{O}_{K_1}^{\new}\times\cdots\times\mathcal{O}_{K_m}^{\new},$ we have
$$\rho_{K_1,\ldots,K_m}^{(n)}(x_1,\ldots,x_m)\le\prod_{i=1}^m\frac{1}{1-p^{-r_i/e_i}}.$$
The result is now immediate from the bounded convergence theorem.
\end{proof}

Now, let us return to \Cref{thm: upper bound of correlation functions}. The estimate of orbital integrals in the following lemma is a weaker version of \cite[Theorem 1.5]{yun2013orbital}.

\begin{lemma}\label{lem: estimate of orbital integral}
Let $K/\Q_p$ be a finite extension. Suppose $x\in\mathcal{O}_K^{\new}$ satisfies $\F_p[x\pmod{p}]=\F_{p^d}$, where $x\pmod{p}$ is the image of $x$ in the residue field of $K$. When $\Z_p[x]$ is a proper subset of $\mathcal{O}_K$, we have
$$\frac{\#(\Lambda\backslash\Mod_{\Z_p[x]})}{\#(\mathcal{O}_K/\Z_p[x])^2}\le p^{-d}+2p^{-2d}.$$
\end{lemma}

\begin{proof}
Let $\delta$ be the Serre invariant of $\Z_p[x]$ (see \cite[Subsection 4.1]{yun2013orbital}), so $\#(\mathcal{O}_K/\Z_p[x])=p^{d\delta}$. Since  $\Z_p[x]$ is a proper subset of $\mathcal{O}_K$, we have $\delta>0$.

We give a few combinatorial calculations as preparation. We denote by $\Y$ the set of partitions $\lambda=(\lambda_1,\lambda_2,\ldots)$, which are (finite or infinite) sequences of nonnegative integers $\lambda_1\ge\lambda_2\ge\cdots$ that are eventually zero. We do not distinguish between two such sequences that differ only by a string of zeros at the end. Let $|\l| := \sum_{i\ge 1}\l_i$, and $l(\l):= \#\{i\mid \l_i >0\}$. Observe that for all $n,k\ge 1$,
\begin{align}
\begin{split}
\#\{\l\in\Y\mid |\l|=n,l(\l)=k\}&=\#\{\l_1\ge\ldots\ge\l_k\ge 1\mid \l_1+\cdots+\l_k=n\}\\
&\le\#\{\lambda_1,\ldots,\lambda_k\ge 1\mid\l_1+\cdots+\l_k=n\}\\
&=\begin{pmatrix}n-1 \\ k-1\end{pmatrix}.
\end{split}
\end{align}
As a consequence, 
\begin{align}\label{eq: estimate of partitions of given length}
\begin{split}
\#\{\l\in\Y\mid |\l|\le \delta,l(\l)=k\}&=\sum_{n\le\delta}\#\{\l\in\Y\mid |\l|=n,l(\l)=k\}\\
&\le \sum_{n\le\delta}\begin{pmatrix}n-1 \\ k-1\end{pmatrix}\\
&=\begin{pmatrix}\delta \\ k\end{pmatrix},
\end{split}
\end{align}
where the last step follows by iterating Pascal's identity. By the upper bound in \cite[Theorem 1.5]{yun2013orbital}\footnote{In the notation of \cite[Theorem 1.5]{yun2013orbital}, $q$ corresponds to our $p$, and $O_\gamma$ is our $\#(\Lambda\backslash\Mod_{\Z_p[x]})$. The product $\prod_{i\in B(\gamma)}$ only has one term in our case and $d_i=d,\delta_i=\delta$, so $\rho(\gamma)=0$.} we have
\begin{align}
\begin{split}
\#(\Lambda\backslash\Mod_{\Z_p[x]})&\le\sum_{|\l|\le\delta}p^{d(\delta-l(\l))}+\sum_{|\l|<\delta}p^{d(|\l|-l(\l))}\\
&\le\sum_{|\l|\le\delta}p^{d(\delta-l(\l))}+\sum_{|\l|\le\delta-1}p^{d(\delta-1-l(\l))}\\
&\le\sum_{k=0}^\delta p^{d(\delta-k)}\begin{pmatrix}\delta \\ k\end{pmatrix}+\sum_{k=0}^{\delta-1} p^{d(\delta-1-k)}\begin{pmatrix}\delta-1 \\ k\end{pmatrix}\\
&=(1+p^d)^{\delta}+(1+p^d)^{\delta-1}.
\end{split}
\end{align}
Here the third line comes from \eqref{eq: estimate of partitions of given length}, and the last line is the binomial theorem. Therefore, we obtain the inequality
$$\frac{\#(\Lambda\backslash\Mod_{\Z_p[x]})}{\#(\mathcal{O}_K/\Z_p[x])^2}\le\frac{(1+p^d)^{\delta}+(1+p^d)^{\delta-1}}{p^{2d\delta}} \le p^{-d}+2p^{-2d}$$
by plugging in $\delta=1$ since the middle expression is decreasing in $\delta$.
\end{proof}

\begin{proof}[Proof of \Cref{thm: upper bound of correlation functions}]
For all $1\le i\le m$, if $\Z_p[x_i]=\mathcal{O}_{K_i}$, then we have $\#(\Lambda_i\backslash\Mod_{\Z_p[x_i]})=1$. Otherwise, by \Cref{lem: estimate of orbital integral}, we have that 
$$\frac{\#(\Lambda_i\backslash\Mod_{\Z_p[x_i]})}{\#(\mathcal{O}_{K_i}/\Z_p[x_i])^2}\le p^{-d_i}+2p^{-2d_i}\le 1,$$
where $d_i:=[\F_p[x_i\pmod{p}]:\F_p]\ge 1$. Thus, regardless of which of the above scenarios is true, we always have
$$\frac{\#(\Lambda_i\backslash\Mod_{\Z_p[x_i]})}{\#(\mathcal{O}_{K_i}/\Z_p[x_i])^2}\le 1.$$
Following the expression given in \Cref{thm:general_cor_functions_intro}, we have
\begin{align}
\begin{split}
\rho_E^{(n)}(x_1,\ldots,x_m)&\le||\Delta_\sigma(x_1,\ldots,x_m)|| \cdot \prod_{i=1}^m\Den(x_i) \\
&=\prod_{1\le i<j\le m}||\Res(Z_i,Z_j)||\cdot\prod_{i=1}^m\frac{\#(\Lambda_i\backslash\Mod_{\Z_p[x_i]})||\Delta_\sigma(x_i)||^2}{1-p^{-r_i/e_i}}\\
&=\prod_{1\le i<j\le m}||\Res(Z_i,Z_j)||\cdot\prod_{i=1}^m\frac{\#(\Lambda_i\backslash\Mod_{\Z_p[x_i]})||\Disc_{K_i/\Q_p}||}{(1-p^{-r_i/e_i})\#(\mathcal{O}_{K_i}/\Z_p[x_i])^2}\\
&\le\prod_{1\le i<j\le m}||\Res(Z_i,Z_j)||\cdot\prod_{i=1}^m\frac{||\Disc_{K_i/\Q_p}||}{1-p^{-r_i/e_i}}\\
&=||\Disc_{E/\Q_p}||\cdot\prod_{1\le i\le m}\frac{1}{1-p^{-r_i/e_i}}\cdot\prod_{1\le i<j\le m}||\Res(Z_i,Z_j)||.
\end{split}
\end{align}
Here, the third line is obtained from the result on page 10 of \cite{caruso2022zeroes}.
\end{proof}

\begin{proof}[Proof of \Cref{thm: independent distribution over lifted subspaces}]
We only need to prove the case where the elements $$x_{1,1},\ldots,x_{1,j_1},\ldots,x_{m,1},\ldots,x_{m,j_m}\in\bar\Z_p$$ 
lie in distinct Galois orbits, since otherwise both sides must be zero. For $1\le i\le m$, let $Z_i$ denote the minimal polynomial of $\{x_{i,l}:1\le l\le j_i\}$. Then, the polynomials $Z_1,\ldots,Z_m$ satisfy the condition in \Cref{thm: decomposition of limit expectation}, and thus
$$\lim_{n\rightarrow\infty}\E||\det(Z(A))||=\prod_{i=1}^m\lim_{n\rightarrow\infty}\E||\det(Z_i(A))||.$$
Here, $Z(x):=\prod_{i=1}^mZ_i(x)$, and $A\in\Mat_n(\Z_p)$ is Haar-distributed. Moreover, as shown in \Cref{prop: distance between lifted subspaces}, the distances between the eigenvalues in different lifted subspaces are always $1$, which implies 
$$||\Delta_\sigma(x_{1,1},\ldots,x_{1,j_1},\ldots,x_{m,1},\ldots,x_{m,j_m})||=\prod_{i=1}^m||\Delta_\sigma(x_{i,1},\ldots,x_{i,j_i})||.$$
Therefore, by \Cref{thm:limit_cor_fns_exist}, we have 
\begin{align}
\begin{split}
\text{LHS}\eqref{eq: decomposition of correlation function}&=\left(\prod_{i=1}^m\prod_{l_i=1}^{j_i}\Den(x_{i,l_i})\right)\cdot||\Delta_\sigma(x_{1,1},\ldots,x_{1,j_1},\ldots,x_{m,1},\ldots,x_{m,j_m})||\cdot\lim_{n\rightarrow\infty}\E||\det(Z(A))||\\
&=\left(\prod_{i=1}^m\prod_{l_i=1}^{j_i}\Den(x_{i,l_i})\cdot||\Delta_\sigma(x_{i,1},\ldots,x_{i,j_i})||\right)\cdot\lim_{n\rightarrow\infty}\E||\det(Z(A))||\\
&=\prod_{i=1}^m\left(\prod_{l_i=1}^{j_i}\Den(x_{i,l_i})\cdot||\Delta_\sigma(x_{i,1},\ldots,x_{i,j_i})||\cdot\lim_{n\rightarrow\infty}\E||\det(Z_i(A))||\right)\\
&=\text{RHS}\eqref{eq: decomposition of correlation function}.
\end{split}
\end{align}
This ends the proof.
\end{proof}

\begin{proof}[Proof of \Cref{thm:lifted_subspace_intro_revised}]
The sets $\mathcal{U}_1,\mathcal{U}_2,\ldots$ appearing in \Cref{thm:lifted_subspace_intro_revised} are exactly the lifted subspaces introduced in \Cref{subsec:zpbar}. The disjoint union \eqref{item:disjoint_union} is clear, and \eqref{item:balls} is proven in \Cref{thm: lifted space as cosets of unit disk}. The asymptotic independence of eigenvalues lying in different lifted subspaces is proven in \Cref{thm: independent distribution over lifted subspaces}, so all that remains is the limit of the number of eigenvalues in each lifted subspace. 

In the end, by passing $A\in\Mat_n(\Z_p)$ to the residue field version with entries in $\F_p$, we can apply \cite[Corollary 1.2]{shen2025universality}\footnote{In the notation of \cite[Corollary 1.2]{shen2025universality}, we are taking $f_i$ be $F_i\in\F_p[x]$, which is the irreducible polynomial that correspond to $\mathcal{U}_i$ as in \eqref{eq: lifted subspace}. Therefore, this polynomial has degree $d_i$. Furthermore, we are taking $\nu_{f_i}$ be the law of $N_i/d_i$.} to conclude that the number of eigenvalues of $A$ lying in $\mathcal{U}_i$ converges in distribution, and the law of the limiting random variable $N_i$ is given there. 
\end{proof}

We note that the results of \cite{shen2025universality} used above were proven for general matrix entry distribution. In the case of uniformly random entries in $\F_p$ (which is all we need here) they were also proven in Corollary 3.1.13 of the second author's unpublished undergraduate thesis \cite{van2018random}, using the same argument as was used for the $\GL_n(\F_q)$ case in Fulman’s work \cite{fulman1997probability}. We do not currently know how to prove universality for the finer eigenvalue statistics we consider beyond the numbers of eigenvalues in lifted subspaces, but this is a very interesting question.

Below is a simple exact computation of the limiting correlation functions in the (not very explicitly defined) special case where $\mathcal{O}_E=\Z_p[\mathfrak{x}]$. This will be useful for later estimates in \Cref{sec:average_number_of_eigenvalues_in_high_degree_extensions}.

\begin{thm}\label{thm: limit density over generator}
Let $E=K_1\times\cdots\times K_m$ be a finite étale algebra over $\Q_p$, and $\mathfrak{x}=(x_1,\ldots,x_m)\in \mathcal{O}_E$ such that $\mathcal{O}_E=\Z_p[\mathfrak{x}]$. Then we have $\rho_E^{(\infty)}(\mfx)=||\Disc_{E/\Q_p}||$.
\end{thm}

\begin{proof}
In this case, for all $1\le i\le m$, we have $\Z_p[x_i] = \mc{O}_{K_i}$ and hence
$$\#(\L_i\backslash \Mod_{\Z_p[x_i]})=\#(\L_i\backslash \Mod_{\mathcal{O}_{K_i}})=1.$$ 
The result on page 10 of \cite{caruso2022zeroes} implies that when $\Z_p[x_i] = \mc{O}_{K_i}$, $||\Delta_\sigma(x_i)||^2=||\Disc_{K_i/\Q_p}||$. Hence $||\Delta_\sigma(x_1,\ldots,x_m)||=||\Disc_{E/\Q_p}||^{1/2}$. In addition, let $p^{f_i}$ denote the order of the residue field of $\mathcal{O}_{K_i}$, so that $\Den(x_i)=\frac{||\Disc_{K_i/\Q_p}||^{1/2}}{1-p^{-f_i}}$ where $\Den$ is as in \Cref{defi: Den and distance}.

By \Cref{thm:limit_cor_fns_exist}, we have 
\begin{align}
\begin{split}
\rho_E^{(\infty)}(\mfx)&=\prod_{i=1}^m\Den(x_i)\cdot ||\Delta_\sigma(x_1,\ldots,x_m)||\cdot\lim_{n\rightarrow\infty}\E||\det(Z(A))||\\
&=\prod_{i=1}^m\frac{||\Disc_{K_i/\Q_p}||^{1/2}}{1-p^{-f_i}}\cdot ||\Disc_{E/\Q_p}||^{1/2}\cdot\prod_{i=1}^m(1-p^{-f_i})\\
&=||\Disc_{E/\Q_p}||.
\end{split}
\end{align}
Here, the second line comes from \Cref{prop: limit expectation over generator}. This completes the proof. 
\end{proof}

%% file: GL_r.tex
\section{Eigenvalues of Haar matrices in {$\GL_n(\Z_p)$}}\label{sec:glr}

The main result of this section, \Cref{thm:limit_cor_fns_exist_GL}, shows that the limiting correlation functions of $\GL_n(\Z_p)$ are the same as for $\Mat_n(\Z_p)$ after taking the obvious restrictions into account. We will also prove results about prelimit correlation functions. We assume the same notation as in \Cref{sec:char_poly}, except that in this section we furthermore require that for all $1\le i\le m$, we have $||x_i||=1$, or equivalently, the constant term of $Z_i(x)$ belongs to $\Z_p^\times$.

\begin{thm}\label{thm: points on variety_GL}
Suppose $A$ is random and distributed according to the Haar probability measure on $\GL_r(\Z_p)$. Then for all tuples $s=(s_1,\ldots,s_m)\in(\frac{1}{e_1}\N,\ldots,\frac{1}{e_m}\N)$ such that $s_i>r\val\Disc(\mathfrak{x})$ for all $1\le i\le m$, we have
\begin{equation}\label{eq: points on the variety_GL}
\prod_{i=1}^mp^{s_ir_i}\cdot\mathbf{P}(\val P_{A}(x_i)\ge s_i,\forall 1\le i\le m)=\frac{1}{||\Delta_\sigma(x_1,\ldots,x_m)||}\prod_{i=1}^m \Den(x_i).
\end{equation}
\end{thm}

\begin{lemma}\label{lem: the same set between Mat and GL}
Suppose $A\in\Mat_r(\Z_p)$ satisfies $\val P_A(x_i)> r\val\Disc(\mathfrak{x})$ for all $1\le i\le m$. Then $A\in\GL_r(\Z_p)$.
\end{lemma}

\begin{proof}[Proof of \Cref{thm: points on variety_GL}, assuming \Cref{lem: the same set between Mat and GL}]
\Cref{lem: the same set between Mat and GL} implies that the set that we want to compute the probability of in \eqref{eq: points on the variety_GL} is exactly the same as in \eqref{eq: points on the variety}. The difference is that in \eqref{eq: points on the variety_GL} it is regarded as a subset in $\GL_r(\Z_p)$, while in \eqref{eq: points on the variety} it is regarded as a subset in $\Mat_r(\Z_p)$. Therefore, we have 
\begin{align}
\begin{split}
\text{LHS}\eqref{eq: points on the variety_GL}&=\frac{\text{RHS}\eqref{eq: points on the variety}}{\mathbf{P}(A\in\GL_r(\Z_p)\mid A\text{ Haar-distributed in }\Mat_r(\Z_p))}\\
&=\frac{\frac{(1-p^{-1})\cdots(1-p^{-r})}{||\Delta_\sigma(x_1,\ldots,x_m)||}\prod_{i=1}^m \Den(x_i)}{(1-p^{-1})\cdots(1-p^{-r})}\\
&=\text{RHS}\eqref{eq: points on the variety_GL},
\end{split}
\end{align}
which ends the proof.
\end{proof}

Now we turn back to the proof of \Cref{lem: the same set between Mat and GL}.

\begin{proof}[Proof of \Cref{lem: the same set between Mat and GL}]
Applying \Cref{thm: polynomial expression of matrix near variety}, the roots of $P_A$ can be written as
$$f(\sigma_{i,l_i}(x_i)),\quad 1\le i\le m,1\le l_i\le r_i,$$
where $f\in F_Z$ (recall \Cref{def:f_set}). Due to the strong triangle inequality, $||\sigma_{i,l_i}(x_i)||=1$ and $f\in F_Z$ implies $$||f(\sigma_{i,l_i}(x_i))||=||f(\sigma_{i,l_i}(x_i))-\sigma_{i,l_i}(x_i)+\sigma_{i,l_i}(x_i)||=||\sigma_{i,l_i}(x_i)||=1.$$ 
Thus the constant term of $P_A$ belongs to $\Z_p^\times$, and $A\in\GL_r(\Z_p)$.
\end{proof}

The next two results give the expressions of the correlation functions.

\begin{thm}\label{thm: joint distribution_GL}
Fix $n \geq 1$ and let $K_1,\ldots,K_m$ be extensions of $\Q_p$ with $\sum_{i=1}^m [K_i: \Q_p] = n$. Let $U \subset \mc{O}_{K_1}^{\new} \times \cdots \times \mc{O}_{K_m}^{\new}$ be any measurable set. Let $A \in \GL_n(\Z_p)$ be multiplicatively Haar-distributed, and $Z_{U,n}$ be the number of $m$-tuples $(x_1,\ldots,x_m)$ of its eigenvalues which lie in $U$. Then
$$\E[Z_{U,n}]=\int_U\rho^{(n),\GL}_{K_1,\ldots,K_m}(x_1,\ldots,x_m)dx_1\cdots dx_m,$$
where the correlation function is explicitly given by
\begin{equation}
\rho_{K_1,\ldots,K_m}^{(n),\GL}(x_1,\ldots,x_m)=||\Delta_\sigma(x_1,\ldots,x_m)||\cdot \left(\prod_{i=1}^m\Den(x_i) \right)\cdot\bbone_{x_i\in\mc{O}_{K_i}^\times,\forall 1\le i\le m},
\end{equation}
and integration (and definition of measurability of $U$) is with respect to the Haar probability measure on $\mc{O}_{K_1}\times \cdots \times \mc{O}_{K_m}$. 
\end{thm}

\begin{thm}\label{thm: local distribution_GL}
Let $K_1,\ldots,K_m$ be extensions of $\Q_p$ and $r := \sum_{i=1}^m [K_i:\Q_p]$. Let $A\in\GL_n(\Z_p)$ be multiplicatively Haar-distributed of size $n\ge r$. Then, the $(K_1,\ldots,K_m)$-correlation function $\rho_{K_1,\ldots,K_m}^{(n),\GL}(x_1,\ldots,x_m)$ of the eigenvalues of $A$ has the form
\begin{equation}
\rho_{K_1,\ldots,K_m}^{(n),\GL}(x_1,\ldots,x_m)=||\Delta_\sigma(x_1,\ldots,x_m)|| \cdot\left(\prod_{i=1}^m\Den(x_i) \right) \cdot \E||\det(Z(\tilde A))||.
\end{equation}
Here $Z(x)=\prod_{i=1}^m Z_i(x)$ has degree $r$, $Z_i\in\Z_p[x]$ is the minimal polynomial of $x_i$ of degree $r_i$, and $\tilde A\in\GL_{n-r}(\Z_p)$ is distributed by the multiplicative Haar measure.
\end{thm}

\begin{proof}[Proof of \Cref{thm: joint distribution_GL} and \Cref{thm: local distribution_GL}]
They can be deduced from \Cref{thm: points on variety_GL} by exactly the same proofs by which \Cref{thm: joint distribution} and \Cref{thm:general_cor_functions_intro} are deduced from \Cref{thm: points on variety_main theorem} in \Cref{sec: joint distribution}.
\end{proof}

\begin{thm}\label{thm:limit_cor_fns_exist_GL}
Let $K_1,\ldots,K_m$ be finite extensions of $\Q_p$, and $E := K_1\times\cdots\times K_m$. For all $1\le i\le m$, let $x_i\in\mc{O}_{K_i}^{\new}$, such that $x_1,\ldots,x_m$ lie in distinct Galois orbits. Then, the limit of the correlation functions
$$\rho_{E}^{(\infty),\GL}(x_1,\ldots,x_m):=\lim_{n\rightarrow\infty}\rho_{E}^{(n),\GL}(x_1,\ldots,x_m)$$
exists, and is given by
$$\rho_{E}^{(\infty),\GL}(x_1,\ldots,x_m)=\rho_{E}^{(\infty)}(x_1,\ldots,x_m)\cdot\bbone_{x_i\in\mc{O}_{K_i}^\times,\forall 1\le i\le m}.$$
\end{thm}

The following lemma will be useful for our proof of \Cref{thm:limit_cor_fns_exist_GL}.

\begin{lemma}\label{cor: expected determinant GL equals Mat}
Let $\tilde A\in\GL_n(\Z_p)$ be multiplicatively Haar-distributed, and $A\in\Mat_n(\Z_p)$ be additively Haar-distributed. Furthermore, let $Z\in\Z_p[x]$ be the same as in \Cref{thm: local distribution_GL}, so that its constant term is in $\Z_p^\times$. Then, we have
$$\lim_{n\rightarrow\infty}\E||\det(Z(\tilde A))||=\lim_{n\rightarrow\infty}\E||\det(Z(A))||.$$
\end{lemma}

\begin{proof}
Notice that the distribution of $\tilde A$ is the same as the distribution of $A$ conditioned on $||\det(A)||=1$. Therefore, we have
$$
\E||\det(Z(\tilde A))||=\E\left[||\det(Z(A))||\bigg|||\det(A)||=1\right].
$$
Furthermore, the residue $\overline Z\in\F_p[x]$ is coprime to $x\in\F_p[x]$. Therefore, by \Cref{rmk: pairwise coprime independence}, we deduce that when $n$ tends to infinity, the joint random variable $(||\det(Z(A))||,||\det(A)||)$ weakly converges to a tuple of two independent variables. Hence
$$
\lim_{n\rightarrow\infty}\E||\det(Z(\tilde A))||=\lim_{n\rightarrow\infty}\E\left[||\det(Z(A))||\bigg|||\det(A)||=1\right]=\lim_{n\rightarrow\infty}\E||\det(Z(A))||,
$$
completing the proof.
\end{proof}

\begin{proof}[Proof of \Cref{thm:limit_cor_fns_exist_GL}]
We only need to prove the case where the elements $x_1,\ldots,x_m\in\bar\Z_p$ lie in distinct Galois orbits, since otherwise the limit must be zero. If one of the $x_i$ does not belong to $\mathcal{O}_{K_i}^\times$, then it is clear that $\rho_{E}^{(\infty),\GL}(x_1,\ldots,x_m)=\lim_{n\rightarrow\infty}\rho_{E}^{(n),\GL}(x_1,\ldots,x_m)=0$. Otherwise, by \Cref{thm: local distribution_GL} and \Cref{cor: expected determinant GL equals Mat}, we have
$$\rho_{E}^{(\infty),\GL}(x_1,\ldots,x_m)=||\Delta_\sigma(x_1,\ldots,x_m)|| \cdot\left(\prod_{i=1}^m\Den(x_i) \right) \cdot \lim_{n\rightarrow\infty}\E||\det(Z(A))||.$$
Here, the matrix $A\in\Mat_n(\Z_p)$ is Haar-distributed. The right hand side of the above equality is exactly the expression of $\rho_{E}^{(\infty)}(x_1,\ldots,x_m)$ in \Cref{thm:limit_cor_fns_exist}.
\end{proof}

\begin{rmk}
Another viewpoint on the results for $\GL_n(\Z_p)$ is as follows. It is a simple fact that for the $\GL_n(\Z_p)$ case, the eigenvalues cannot lie in the open unit disk around $0$. Therefore, we can also regard the correlation functions as functions defined on $\mathcal{O}_{K_1}^{\times,\new}\times\cdots\times\mathcal{O}_{K_m}^{\times,\new}$, where $\mathcal{O}_{K_i}^{\times,\new}:=\mathcal{O}_{K_i}^{\times}\cap\mathcal{O}_{K_i}^{\new}$. Moreover, we may also equip $\mathcal{O}_{K_i}^{\times}$ with the multiplicative Haar probability measure. In this case, the expressions of correlation functions with respect to the product of these multiplicative Haar measures are the same as in \Cref{thm:limit_cor_fns_exist_GL}, except that we remove the indicator function and replace $V: \mathcal{O}_{K_i}^{\new}\rightarrow\R$ by 
\begin{align}
\begin{split}
V^{\GL}:\mathcal{O}_{K_i}^{\times,\new}&\rightarrow\R\\
x_i&\mapsto||\Delta_\sigma(x_i)||\cdot\#(\Lambda\backslash\Mod_{\Z_p[x_i]}).
\end{split}
\end{align}
One can see that the value of $V^{\GL}$ in $\mathcal{O}_{K_i}^{\times,\new}$ differs from $V$ only by a constant factor $1-p^{-r_i/e_i}$, which is exactly the additive Haar measure of $\mathcal{O}_{K_i}^\times$ in $\mathcal{O}_{K_i}$.
\end{rmk} 

%% file: markov.tex
\section{Markov chain description of the cokernel}\label{sec:markov_chains}

In \Cref{thm:compute_char_poly_det} we showed that the terms $\lim_{n \to \infty} \E[||\det(Z(A))||]$, which we need to compute in order to find the limiting correlation functions, are given by certain mixed moments of the $p$-adic valuation of determinants of related matrices $A_0 + x_i A_1 + \ldots + x_i^{r-1} A_{r-1}$ as $x_i$ ranges over the roots of $Z$. Our strategy to show \Cref{thm:compute_char_poly_det} was to consider cokernels rather than determinants, and we will use the same perspective to compute these mixed moments more explicitly. The basic idea is that for a single Haar matrix $A$ over $\Z_p$ or some more general $\mc{O}_K$, one may sample $\cok(A)/(\pi), \cok(A)/(\pi^2),\ldots$ successively, and (not so obviously) this may be stated in terms of a homogeneous Markov chain. Each cokernel of a matrix $A_0 + x_i A_1 + \ldots + x_i^{r-1} A_{r-1}$ may be individually sampled in this way, but for different $i$ one will in general have dependent copies of the underlying Markov chain. This section begins with known results of Fulman and Evans on this Markov chain, and uses them to prove a general formula for the types of expectations we will encounter.

Recall the $q$-Pochhammer notation 
\begin{equation}
    (a;t)_n := \prod_{i=0}^{n-1} (1-at^i)
\end{equation}
for $n \in \Z_{\geq 0} \cup \{\infty\}$.

\begin{defi}[{\cite[Theorem 2]{fulman-RR}}]
    \label{def:fulman_markov_chain}
    For fixed parameters $t \in (0,1)$ and $u \in (0,t^{-1})$, define a Markov kernel $\mc{K}(\cdot,\cdot)$ on $\Z_{\geq 0}$ by 
    \begin{equation}
        \mc{K}(a,b) = \mc{K}_{t,u}(a,b) = \begin{cases}
            \frac{t^{b^2}u^b(t;t)_a (ut;t)_a}{(t;t)_{a-b} (t;t)_b (ut;t)_b} & 0 \leq b \leq a \\ 
            0 & \text{else}
        \end{cases}.
    \end{equation}
    We will sometimes allow initial state $\infty$, and define transition probabilities
    \begin{equation}
        \mc{K}(\infty,b) = \lim_{a \to \infty} \mc{K}(a,b) = (ut;t)_\infty \frac{t^{b^2}u^b}{(t;t)_b (ut;t)_b}
    \end{equation}
    to finite states.
\end{defi}

It is not immediate that $\mc{K}$ defines a Markov kernel, but this is proven in {\cite[Theorem 2]{fulman-RR}}. Evans \cite{evans2002elementary} found the same Markov chain independently in his study of elementary divisors (which are also sometimes called singular numbers, and are equivalent to cokernels) of $p$-adic random matrices.

\begin{defi}
     Let $K$ be a non-archimedean local field with uniformizer $\pi$ and residue field isomorphic to $\F_q$, and let $M$ be an $\mc{O}_K$-module. Then for any $i \in \Z_{\geq 1}$ we define
     \begin{equation}
         M_i' := \rank_{\F_q}(\pi^{i-1}M/\pi^iM),
     \end{equation}
     noting that $\pi^{i-1}M/\pi^iM$ has the structure of an $\mc{O}_K/\pi\mc{O}_K \cong \F_q$-vector space.
\end{defi}

Alternatively, by the structure theorem for modules over PIDs, any module $M$ as above is isomorphic to $\bigoplus_{i \geq 0} \mc{O}_K/\pi^{\la_i}\mc{O}_K$ for some integers $\la_1 \geq \la_2 \geq \ldots \geq 0$, finitely many of which are nonzero. One may equivalently define
\begin{equation}
    M_i' = \#\{j: \la_j \geq i\}.
\end{equation}

\begin{prop}[{\cite[Theorem 3.5]{evans2002elementary}}]\label{thm:from_evans}
    Let $K$ be a non-archimedean local field with residue field of size $q$, and let $A$ be an $n \times n$ matrix with iid Haar-distributed entries in $\mc{O}_K$. Then the random sequence $\Cok(A)_1',\Cok(A)_2',\ldots$ is distributed according to the Markov chain $\mc{K}_{q^{-1},1}$ of \Cref{def:fulman_markov_chain} started at state $n$. In other words,
    \begin{equation}
        \Pr(\Cok(A)_1' = b) = \mc{K}_{q^{-1},1}(n,b), 
    \end{equation}
    and
    \begin{equation}
        \Pr(\Cok(A)_{i+1}' = b | \Cok(A)_i' = a) = \mc{K}_{q^{-1},1}(a,b)
    \end{equation}
    for all $i \geq 1$. 
\end{prop}

The simple fact that it is a Markov kernel for every $u$, verified in \cite[Theorem 2]{fulman-RR}, already allows computation of $\E[||\det(A)||]$ for any Haar matrix $A$:

\begin{prop}
    \label{thm:compute_t_moment_single_markov}
    Let $\la_1,\la_2,\ldots$ be distributed as the steps of the Markov chain $\mc{K}_{t,u}$ of \Cref{def:fulman_markov_chain}, with fixed initial state $\la_0 = n \in \N \cup \{\infty\}$. Then for $k \in \Z_{\geq 0}$,
    \begin{equation}
        \E[t^{k(\la_1+\la_2+\ldots)}] = \frac{(ut;t)_k}{(ut^{n+1};t)_k},
    \end{equation}
    where when $n=\infty$ we interpret the denominator as $1$. 
\end{prop}
\begin{proof}
    The factors $\frac{(t;t)_a(ut;t)_a}{(t;t)_b(ut;t)_b}$ in the explicit formula for $\mc{K}$ cancel in any infinite product 
    $$\mc{K}(n,\la_1)\mc{K}(\la_1,\la_2)\cdots$$ 
    except for the first $(t;t)_n(ut;t)_n$. This yields
    \begin{align}
        \begin{split}
            \E[t^{k(\la_1+\la_2+\ldots)}] &= \sum_{\substack{\la_1 \geq \la_2 \geq \ldots \geq 0 \\ \la_1 \leq n}} \left((t;t)_n (ut;t)_n \prod_{i \geq 0} \frac{t^{\la_{i+1}^2} u^{\la_{i+1}}}{(t;t)_{\la_i - \la_{i+1}}}\right) \cdot t^{k (\la_1+\ldots)} \\ 
            &= \frac{(ut;t)_n}{(ut^{k+1};t)_n} \sum_{\substack{\la_1 \geq \la_2 \geq \ldots \geq 0 \\ \la_1 \leq n}}(t;t)_n (ut^{k+1};t)_n \prod_{i \geq 0} \frac{t^{\la_{i+1}^2} (ut^k)^{\la_{i+1}}}{(t;t)_{\la_i - \la_{i+1}}} \\ 
            &=  \frac{(ut;t)_n}{(ut^{k+1};t)_n} \sum_{\substack{\la_1 \geq \la_2 \geq \ldots \geq 0 \\ \la_1 \leq n}} \prod_{i \geq 0} \mc{K}_{t,ut^k}(\la_i,\la_{i+1}).
        \end{split}
    \end{align}
    The sum in the last line is just $1$, and the quotient of $q$-Pochhammers is equal to the desired one after canceling factors.
\end{proof}

\begin{cor}
    \label{thm:moments_of_determinant}
    Let $K/\Q_p$ be an extension of degree $r$ and residue field of size $q$. Let $A \in \Mat_n(\mc{O}_K)$ be Haar-distributed. Then for any $k \geq 0$,
    \begin{equation}\label{eq:expected_determinant_single_matrix}
        \E[||\det(A)||^{rk}] = \frac{(q^{-1};q^{-1})_k}{(q^{-n-1};q^{-1})_k}.
    \end{equation}
\end{cor}
\begin{proof}
    By \Cref{lem: cok and det} we have
    \begin{equation}
        \E[||\det(A)||^{rk}] = \E[\#\Cok(A)^{-k}].
    \end{equation}
    By \Cref{thm:from_evans} this is equal to $\E[(q^{-1})^{k(\la_1+\la_2+\ldots)}]$ where $\la_0=n$ and $\la_1,\la_2,\ldots$ are the successive steps of the Markov chain $\mc{K}_{q^{-1},1}$. The result now follows from \Cref{thm:compute_t_moment_single_markov} with $u=1,t=q^{-1}$.
\end{proof}

For more general expectations of the type in \Cref{thm:compute_char_poly_det}, we need finer properties of the Markov kernel $\mc{K}$. As a $\Z_{\geq 0} \times \Z_{\geq 0}$ matrix, it has an explicit diagonalization, which was used to give a new proof of the Rogers-Ramanujan identities in \cite{fulman-RR}. To state, it is easier to work with the conjugated matrix $M$ given by
    \begin{equation}
        M(a,b) = M_{t,u}(a,b) = \begin{cases}
            \frac{t^{b^2}u^b}{(t;t)_{a-b}} &  b \leq a \\ 
            0 & \text{else}
            \end{cases},
    \end{equation}
which is no longer a Markov kernel but is algebraically simpler, and already essentially appeared in the proof of \Cref{thm:compute_t_moment_single_markov}.

\begin{prop}[{\cite[Theorem 3]{fulman-RR}}]\label{thm:fulman_diagonalization}
    Define the matrix $U = (U(i,j))_{i,j \in \Z_{\geq 0}}$ by 
    \begin{equation}
        U(i,j) = U_{t,u}(i,j) = \begin{cases}
            \frac{1}{(t;t)_{i-j}(ut;t)_{i+j}} & i \geq j \\ 
            0 & i < j
        \end{cases}.
    \end{equation}
    Then the columns $(U(i,j))_{i \in \Z_{\geq 0}}$ of $U$ are eigenvectors of $M$, and $M = UEU^{-1}$ where $E = E_{t,u}$ is the $\Z_{\geq 0} \times \Z_{\geq 0}$ diagonal matrix with $(i,i)$ entry $u^i t^{i^2}$ for $i=0,1,\ldots$.

    Furthermore the inverse matrix in the diagonalization has explicit entries
    \begin{equation}
        U^{-1}(i,j) = \begin{cases}
            \frac{(-1)^{i-j}t^{\binom{i-j}{2}}(1-u t^{2i}) (ut;t)_{i+j-1}}{(t;t)_{i-j}} & i \geq j \\ 
            0 & i < j
        \end{cases}.
    \end{equation}
    Here in the case $i=j=0$ we must interpret $(ut;t)_{-1} = 1/(1-u)$ according to the standard convention for $q$-Pochhammer symbols of negative index.
\end{prop}

Note that our $U$ is called $A$ in \cite{fulman-RR}, but we already use $A$ for Haar matrices. 

\begin{cor}
    \label{thm:markov_expectations}
    Let $m \geq 1$, and let $\xi \in \C$. Let $\la_1,\ldots,\la_m$ be successive states of the Markov chain $\mc{K} = \mc{K}_{t,u}$ started from $\la_0 = n$ where $n \in \Z_{\geq 0} \cup \{\infty\}$. Then
    \begin{equation}
             \E[\xi^{\la_1+\cdots+\la_m} f(\la_m)]  = (t;t)_n(ut;t)_n\sum_{\substack{d \geq \ell \geq 0 \\ d \leq n}} U_{t,u\xi}(n,d) E_{t,u\xi}(d,d)^m U_{t,u\xi}^{-1}(d,\ell) \frac{f(\ell)}{(t;t)_\ell (ut;t)_\ell},
    \end{equation}
    for any function $f: \Z_{\geq 0} \to \R$. Here if $n = \infty$ we assume $f$ has sufficient decay so that the left-hand side is defined, and set
    \begin{equation}\label{eq:u_infty_def}
        U_{t,u\xi}(\infty,d) = \lim_{n \to \infty} U_{t,u\xi}(n,d) = \frac{1}{(t;t)_\infty (u\xi t;t)_\infty}.
    \end{equation}
\end{cor}
\begin{proof}
The desired expectation is
        \begin{align}\label{eq:absorb_xi}
            \begin{split}
           &\sum_{\la_1 \geq \ldots \geq \la_m \geq 0} \prod_{i=0}^{m-1} \mc{K}(\la_i,\la_{i+1}) \xi^{\la_1+\cdots+\la_m} f(\la_m)\\ 
           &= \sum_{\la_1 \geq \ldots \geq \la_m \geq 0}(t;t)_n (ut;t)_n \prod_{i=0}^{m-1} \frac{t^{\la_{i+1}^2} u^{\la_{i+1}}}{(t;t)_{\la_i - \la_{i+1}}} \cdot \frac{1}{(t;t)_{\la_m} (ut;t)_{\la_m}} \xi^{\la_1+\cdots+\la_m} f(\la_m) \\ 
           &= (t;t)_n(ut;t)_n \sum_{\la_1 \geq \la_1 \geq \ldots \geq \la_m\geq 0} \prod_{i=0}^{m-1} M_{t,u\xi}(\la_i,\la_{i+1}) \cdot \frac{f(\la_m)}{(t;t)_{\la_m} (ut;t)_{\la_m}}.
            \end{split}
        \end{align}
    Diagonalizing the matrices $M_{t,u\xi}$ using \Cref{thm:fulman_diagonalization} yields
        \begin{equation}\label{eq:use_diagonalization}
        \text{RHS\eqref{eq:absorb_xi}} = (t;t)_n(ut;t)_n \sum_{\la_m \geq 0} (U_{t,u\xi} E_{t,u\xi}^m U_{t,u\xi}^{-1})(n,\la_m)  \frac{f(\la_m)}{(t;t)_{\la_m} (ut;t)_{\la_m}},
    \end{equation}
    which becomes the desired expression after writing out the matrix product and replacing $\la_m$ by $\ell$.    
\end{proof}

To simplify such expectations still further, \Cref{thm:markov_expectations} suggests a natural route: write the function $f$ in terms of the eigenvectors of $M_{t,u\xi}$ (i.e. columns of $U_{t,u\xi}$). The next result does this for the functions which arise for our correlation functions in the next section, which are cases of the $u=1,\xi=t^2$ case of \Cref{thm:markov_expectations}.

\begin{lemma}
    \label{thm:expand_in_eigenbasis_for_2pt}
    Set
    \begin{equation}
        F_j(i) = \begin{cases}
            \frac{1}{(t;t)_{i-j}(t^3;t)_{i+j}} & i \geq j \\ 
            0 & i < j
        \end{cases},
    \end{equation}
    the function on $\Z_{\geq 0}$ given by the $j\tth$ column of $U_{t,t^2}$. Let $g_1,g_2,g_3$ be functions on $\Z_{\geq 0}$ given by 
    \begin{align}
        g_1(\ell) &= \frac{1}{(t;t)_\ell (t;t)_\ell} \pfrac{1-t}{1-t^{\ell+1}}^2 \\ 
        g_2(\ell) &= \frac{1}{(t;t)_\ell (t;t)_\ell} \frac{1-t}{1-t^{\ell+1}} \\ 
        g_3(\ell) &= \frac{1}{(t;t)_\ell (t;t)_\ell} \frac{1-t^2}{1-t^{2\ell+2}}.
    \end{align}
    Then these functions have eigenbasis expansions
    \begin{align}\label{eq:expand_eigenbasis_1}
         g_1 &= \sum_{j \geq 0} \frac{(-1)^j t^{\binom{j+1}{2}}(1+t^{j+1})}{1+t} F_j \\ \label{eq:expand_eigenbasis_2}
        g_2 &= \sum_{j \geq 0} \frac{t^{j^2+j}(1-t^{2j+2})}{1-t^2}F_j \\ \label{eq:expand_eigenbasis_3}
        g_3 &= \sum_{j \geq 0} \frac{t^{\binom{j+1}{2}}(1-t^{j+1})}{1-t} F_j.
    \end{align}
\end{lemma}
\begin{proof}
    Each of \eqref{eq:expand_eigenbasis_1}, \eqref{eq:expand_eigenbasis_2} and \eqref{eq:expand_eigenbasis_3} is a similar computation, and we give the first and sketch the others. In each case we just have to check the identity when the functions are evaluated every $\ell \in \Z_{\geq 0}$. It will be useful to slightly rewrite $F_j(i)$ as 
    \begin{equation}
        F_j(i) = \begin{cases}
            \frac{(t;t)_2}{(t;t)_{i-j}(t;t)_{i+j+2}} & i \geq j \\ 
            0 & i < j
        \end{cases}.
    \end{equation}
    Now let us prove \eqref{eq:expand_eigenbasis_1} evaluated at a fixed $\ell$. Setting 
    \begin{equation}
        E_j^{(1)} = \frac{(-1)^j t^{\binom{j+1}{2}}(1-t)^2}{(1-t^{\ell+1})(t;t)_{\ell-j}(t;t)_{\ell+j+1}},
    \end{equation}
    it is an elementary check that
    \begin{equation}
        E_j^{(1)} - E_{j+1}^{(1)} = (1-t)^2\frac{(-1)^j t^{\binom{j+1}{2}}(1+t^{j+1})}{(t;t)_{\ell-j}(t;t)_{\ell+j+2}} = \frac{(-1)^j t^{\binom{j+1}{2}}(1+t^{j+1})}{1+t} F_j(\ell)
    \end{equation}
    for $j=0,\ldots,\ell-1$. 
    Hence the right-hand side of \eqref{eq:expand_eigenbasis_1} evaluated at $\ell$ is just 
    \begin{equation}\label{eq:with_E_l}
        \sum_{j=0}^{\ell} \frac{(-1)^j t^{\binom{j+1}{2}}(1+t^{j+1})}{1+t} F_j(\ell)= E_0^{(1)} - E_{\ell}^{(1)} + \frac{(-1)^\ell t^{\binom{\ell+1}{2}}(1+t^{\ell+1})}{1+t}F_\ell(\ell).
    \end{equation}
    Using the identity $(1-t^{\ell+1})(t;t)_{2\ell+1} = (t;t)_{2\ell+2}/(1+t^{\ell+1})$, we obtain that the last two terms of \eqref{eq:with_E_l} are equal, so \eqref{eq:with_E_l} is just $E_0^{(1)} = g_1(\ell)$. This proves \eqref{eq:expand_eigenbasis_1}.

    To prove \eqref{eq:expand_eigenbasis_2} one uses the same telescoping argument, with a different formula 
    \begin{equation}
        E_j^{(2)} = \frac{(1-t)t^{j^2+j}}{(t;t)_{\ell-j}(t;t)_{\ell+j+1}}
    \end{equation}
    instead of $E_j^{(1)}$. For \eqref{eq:expand_eigenbasis_3} one uses 
    \begin{equation}
        E_j^{(3)} = \frac{(1-t^2)t^{\binom{j+1}{2}}}{(1+t^{\ell+1})(t;t)_{\ell-j}(t;t)_{\ell+j+1}}.
    \end{equation}
    In each case $E_0^{(k)} = g_k(\ell)$ and $E_j^{(k)}-E_{j+1}^{(k)}$ gives the terms in the sum on the right-hand side of \eqref{eq:expand_eigenbasis_2} or \eqref{eq:expand_eigenbasis_3} evaluated at $\ell$, which telescope to leave only the $E_0^{(k)}$ term.
\end{proof}

 \begin{rmk}
     \label{rmk:miracle_in_diagonalization} 
     We remark that the eigenbasis expansions of $g_1,g_2,g_3$ in \Cref{thm:expand_in_eigenbasis_for_2pt} take a particularly simple form, and the eigenbasis expansions of an arbitrary rational function $g$ are typically not nearly so explicit. This algebraic miracle is responsible for the simplicity of our series formulas in \Cref{thm:compute_2pt_expectations}, and consequently of the formulas in \Cref{thm:Zp_2pt_intro} and \Cref{thm:quadratic_intro}, and we do not have a good understanding of its source. The three formulas in \Cref{thm:compute_2pt_expectations}, with the left-hand sides written as explicit sums over $\la_1 \geq \ldots \geq \la_m \geq 0$, may be viewed as deformed versions of the Andrews-Gordon identities \cite{andrews1974analytic}.
 \end{rmk}

\begin{cor}
    \label{thm:compute_2pt_expectations}
    Fix $m \geq 1$. For any $n\ge 1$, and $\la_1,\ldots,\la_m$ be successive states of the Markov chain $\mc{K}=\mc{K}_{t,1}$ started from $\la_0 = n$. Then 
    \begin{align}
        \lim_{n\rightarrow\infty}\E\left[t^{2\la_1+\cdots+2\la_m} \pfrac{1-t}{1-t^{\la_m+1}}^2\right] &= (1-t)^2 \sum_{k \geq 0} (-1)^k t^{\frac{(2m+1)k^2 + (4m+1)k}{2}}(1+t^{k+1}) \\ 
        \lim_{n\rightarrow\infty}\E\left[t^{2\la_1+\cdots+2\la_m} \frac{1-t}{1-t^{\la_m+1}}\right] &=(1-t)\sum_{k \geq 0} t^{(m+1)k^2 + (2m+1)k}(1-t^{2k+2}) \\ 
        \lim_{n\rightarrow\infty}\E\left[t^{2\la_1+\cdots+2\la_m} \frac{1-t^2}{1-t^{2\la_m+2}}\right] &=  
        (1-t^2)\sum_{k \geq 0} t^{\frac{(2m+1)k^2+(4m+1)k}{2}}(1-t^{k+1}).
    \end{align}
\end{cor}
\begin{proof}
Let $\nu_1,\ldots,\nu_m$ be successive states of the Markov chain $\mathcal{K}=\mathcal{K}_{t,1}$ started from $\infty$. When $n$ goes to infinity, the joint variable $(\l_1,\ldots,\l_m)$ weakly converges to $(\nu_1,\ldots,\nu_m)$ since the transition probabilities of the first step of the Markov chain converge. This weak convergence implies
$$\lim_{n\rightarrow\infty}\E[t^{2\la_1+\cdots+2\la_m}f(\la_m)]=\E[t^{2\nu_1+\cdots+2\nu_m}f(\nu_m)]$$
for each of the three choices of $f$ appearing in \Cref{thm:compute_2pt_expectations}, since $t^{2\nu_1+\cdots+2\nu_m}f(\nu_m)$ is a bounded function in these cases. Now, we use \Cref{thm:markov_expectations} with $u=1$ and $\xi = t^2$. Note that 
\begin{equation}\label{eq:F_U_inverse}
\sum_{\ell \geq 0} U_{t,t^2}^{-1}(d,\ell) F_j(\ell) = \bbone_{j=d}
\end{equation}
because the functions $F_j$ (viewed as vectors $(F_j(\ell))_{\ell \in \Z_{\geq 0}}$) are columns of $U_{t,t^2}$. Hence if the functions $f$ in \Cref{thm:markov_expectations} expand in the basis $F_j$ of \Cref{thm:expand_in_eigenbasis_for_2pt} as 
\begin{equation}
    \frac{f(\ell)}{(t;t)_\ell(t;t)_\ell} = \sum_{j \geq 0} c_j F_j(\ell),
\end{equation}
then \Cref{thm:markov_expectations} implies
\begin{align}
\begin{split}
\E[t^{2\nu_1+\cdots+2\nu_m}f(\nu_m)] & = (t;t)_\infty^2 \sum_{d,\ell \geq 0} U_{t,t^2}(\infty,d) E_{t,t^2}(d,d)^m U_{t,t^2}^{-1}(d,\ell) \sum_{j \geq 0} c_j F_j(\ell) \\ 
&= (t;t)_\infty^2 \sum_{d \geq 0} U_{t,t^2}(\infty,d) E_{t,t^2}(d,d)^m c_d
\end{split}
\end{align}
where the last line is by \eqref{eq:F_U_inverse}. Each of the cases we must show corresponds to a different $f$, for which the expansion coefficients $c_j$ were computed in \Cref{thm:expand_in_eigenbasis_for_2pt}, and the eigenvalues $E_{t,t^2}(d,d)^m$ are explicit powers of $t$ given in \Cref{thm:fulman_diagonalization}. Elementary calculation concludes the proof.
\end{proof}

%% file: limiting_correlation_functions.tex
\section{Limiting correlation functions in low-degree extensions}\label{sec: limiting correlation functions in low degree}

The goal of this section is to compute some specific correlation functions of low degree extensions. Specifically, we prove the formulas for correlation functions in \Cref{thm:Zp_2pt_intro} and \Cref{thm:quadratic_intro}. Using these, we show the exact formulas for expected number of roots in quadratic extensions in \Cref{thm:roots_in_extensions_intro}, and also compute a formula for the variance of the number of eigenvalues in $\Z_p$ (\Cref{thm:variance_of_roots_in_Z_p}) which was not stated in the Introduction.

\subsection{Eigenvalues in $\Z_p$}

In this subsection, we deal with the eigenvalues in $\Z_p$. We start from the density of eigenvalues in $\Q_p$.

\begin{thm}\label{thm: one point correlation function over Zp}
For all positive integers $n$ and all $x_1\in\Z_p$, we have $\rho_{\Q_p}^{(n)}(x_1)=1$.
\end{thm}

\begin{proof} 
In this case, we have $r=e=1$, $\Disc_{K/\Q_p}=\#(\mathcal{O}_K:\Z_p[x_1])=\#(\L_E\backslash \Mod_{\Z_p[x_1]})=1$, and $\Den(x_1)=\frac{1}{1-p^{-1}}$. By \Cref{thm:general_cor_functions_intro}, we have 
$$
\rho_{\Q_p}^{(n)}(x_1)=\frac{1-p^{-n}}{1-p^{-1}}\E||\det(\tilde A)||=1.
$$
Here, $\tilde A\in\Mat_{n-1}(\Z_p)$ is Haar-distributed, and the second equality comes from \Cref{thm:moments_of_determinant}. This completes the proof.
\end{proof}

\Cref{thm: one point correlation function over Zp} and \Cref{thm:general_cor_functions_intro} lead to the following corollary.

\begin{cor}\label{cor: expected eigenvalues in Zp}
We have $\E[Z_{\Z_p,n}]=1$ for all $n\ge 1$. 
\end{cor}

Let us turn to pair correlation functions over $\Z_p$.

\begin{thm}\label{thm: two point correlation function in Z_p}
Let $x_1,x_2\in\Z_p$, and $m:=\val(x_1-x_2)$. Then we have
$$\rho_{\Q_p\times\Q_p}^{(\infty)}(x_1,x_2)=p^{-m}\sum_{k \geq 0} (-1)^k p^{-\frac{(2m+1)k^2 + (4m+1)k}{2}}(1+p^{-k-1}).$$
\end{thm}

\Cref{thm: two point correlation function in Z_p} can be directly deduced from the following lemma.

\begin{prop}\label{lem: expectation of two point determinant}
Let $x_1,x_2\in\Z_p$, and $m:=\val(x_1-x_2)$. When $A\in\Mat_n(\Z_p)$ is Haar-distributed, we have
$$\lim_{n\rightarrow\infty}\E||(x_1I_n-A)(x_2I_n-A)||=(1-p^{-1})^2\sum_{k \geq 0} (-1)^k p^{-\frac{(2m+1)k^2 + (4m+1)k}{2}}(1+p^{-k-1}).$$
\end{prop}

\begin{proof}[Proof of \Cref{thm: two point correlation function in Z_p}, assuming \Cref{lem: expectation of two point determinant}]
Following the formula in \Cref{thm:limit_cor_fns_exist}, we have
\begin{align}
\begin{split}
\rho_{\Q_p\times\Q_p}^{(\infty)}(x_1,x_2)&=||x_1-x_2||\cdot\frac{1}{(1-p^{-1})^2}\cdot\lim_{n\rightarrow\infty}\E||(x_1I_n-A)(x_2I_n-A)||\\
&=p^{-m}\cdot\frac{1}{(1-p^{-1})^2}\cdot(1-p^{-1})^2\sum_{k \geq 0} (-1)^k p^{-\frac{(2m+1)k^2 + (4m+1)k}{2}}(1+p^{-k-1})\\
&=p^{-m}\sum_{k \geq 0} (-1)^k p^{-\frac{(2m+1)k^2 + (4m+1)k}{2}}(1+p^{-k-1}).
\end{split}
\end{align}
This gives the proof.
\end{proof}

The following lemma is needed for the proof of \Cref{lem: expectation of two point determinant}.

\begin{lemma}\label{lem: joint distribution of chains_2pt}
Let $n\ge 1$, and $A,B\in\Mat_n(\Z_p)$ be independent and Haar-distributed. Let $x_1,x_2$ be fixed elements in $\Z_p$, and $m:=\val(x_1-x_2)$. Then, the joint distribution of $\Z_p$-modules 
$$(\Cok(A-x_1B),\Cok(A-x_2B))$$
is given by the following way. First, the random sequence
\begin{multline}
\Cok(A-x_1B)_1'=\l_1,\Cok(A-x_1B)_2'=\l_2,\ldots,\Cok(A-x_1B)_m'=\l_m,\\
\Cok(A-x_1B)_{m+1}'=\l_{m+1}^{(1)},\Cok(A-x_1B)_{m+2}'=\l_{m+2}^{(1)},\ldots
\end{multline}
is distributed according to the Markov chain $\mathcal{K}_{p^{-1},1}$ of \Cref{def:fulman_markov_chain} started at state $n$. After we determine $\Cok(A-x_1B)$ in the above way, we have $\Cok(A-x_2B)_i=\Cok(A-x_1B)_i$ for all $1\le i\le m$, and the random sequence 
$$\Cok(A-x_2B)_{m+1}'=\l_{m+1}^{(2)},\Cok(A-x_2B)_{m+2}'=\l_{m+2}^{(2)},\ldots$$
is distributed according to the Markov chain $\mathcal{K}_{p^{-1},1}$ started at state $\l_m$.
\end{lemma}

\begin{proof}
Since $A\in\Mat_n(\Z_p)$ is Haar-distributed and independent of $B$, there is no loss of generality to assume $x_1=0$. Then, the distribution of $\Cok(A-x_1B)=\Cok(A)$ is already given by $\mc{K}_{p^{-1},1}$ by \Cref{thm:from_evans}. Now, given a fixed sequence of integers
\begin{equation}
    n\ge\l_1\ge\ldots\ge\l_m\ge\l_{m+1}^{(1)}\ge\ldots\ge 0,\label{eq:lambda_and_parens}
\end{equation}
it remains to understand the distribution of the random cokernel $\Cok(A-x_2B)$, conditioned on
\begin{equation}\label{eq: condition sequence_2pt in Z_p}
\Cok(A)_1'=\l_1,\ldots,\Cok(A)_m'=\l_m,\Cok(A)_{m+1}'=\l_{m+1}^{(1)},\ldots
\end{equation}
By the Smith normal form, there exists a sequence of integers $\nu_1\ge\ldots\ge\nu_n\ge 0$, such that every $A\in\Mat_n(\Z_p)$ that satisfies \eqref{eq: condition sequence_2pt in Z_p} can be expressed in the form
$$A=U_1\diag(p^{\nu_1},\ldots,p^{\nu_n})U_2$$
for some $U_1,U_2\in\GL_n(\Z_p)$; note that these are entirely determined by the partition \eqref{eq:lambda_and_parens} since $\nu_j=\#\{i\ge 1:\lambda_i\ge j\}$. Since the Haar distribution over $\Mat_n(\Z_p)$ is invariant under the left- and right-multiplication of $\GL_n(\Z_p)$, we have that conditioned on $A$ satisfying \eqref{eq: condition sequence_2pt in Z_p},
\begin{equation}
    \Cok(A-x_2 B) = \Cok\left(\diag(p^{\nu_1},\ldots,p^{\nu_n})-x_2B\right) \quad \quad \quad \quad \text{in distribution.}
\end{equation}
We will analyze the distribution of the cokernel on the right-hand side. Since $x_2 \in p^m \Z_p$, 
\begin{equation}
    \Cok\left(\diag(p^{\nu_1},\ldots,p^{\nu_n})-x_2B\right)_i' = \la_i \quad \quad \quad \quad \text{ for $1 \leq i \leq m$,}
\end{equation}
and it remains to find the distribution of $\Cok\left(\diag(p^{\nu_1},\ldots,p^{\nu_n})-x_2B\right)_i'$ for $i=m+1,m+2,\ldots$. Denote by $A_1:=\diag(p^{\nu_1},\ldots,p^{\nu_{\l_m}})$, and $A_4:=\diag(p^{\nu_{\l_m+1}},\ldots,p^{\nu_n})$. $A_1$ has size $\l_m$, and its entries have valuation $\ge m$. $A_4$ has size $n-\l_m$, and the entries on its diagonal have valuation $<m$. Denote 
$$B=\begin{pmatrix} B_1 & B_2\\ B_3 & B_4\end{pmatrix},$$
where 
$$B_1\in\Mat_{\lambda_m}(\Z_p),B_2\in\Mat_{\lambda_m\times(n-\lambda_m)}(\Z_p),B_3\in\Mat_{(n-\lambda_m)\times\lambda_m}(\Z_p),B_4\in\Mat_{n-\lambda_m}(\Z_p)$$ 
are all independently Haar-distributed. We have
\begin{align}
\begin{split}
\Cok\left(\diag(p^{\nu_1},\ldots,p^{\nu_n})-x_2B\right)&=\Cok\begin{pmatrix}A_1-x_2B_1 & -x_2B_2\\ -x_2B_3 & A_4-x_2B_4\end{pmatrix}\\
&\cong\Cok\begin{pmatrix} A_1-B_0-x_2B_1 & -x_2B_2\\ 0 & A_4-x_2B_4\end{pmatrix}\\
&\cong\Cok\begin{pmatrix} A_1-B_0-x_2B_1 & 0\\ 0 & A_4-x_2B_4\end{pmatrix},\\
\end{split}
\end{align}
where $B_0=x_2^2B_2(A_4-x_2B_4)^{-1}B_3\in p^m\Mat_{\l_m}(\Z_p)$ is independent of $B_1$. Therefore, the matrix $A_0:=A_1-B_0-x_2B_1$ is Haar-distributed in $\Mat_{\l_m}(p^m\Z_p)$. By \Cref{thm:from_evans}, the sequence 
$$\Cok(A_0)_{m+1}',\Cok(A_0)_{m+2}',\ldots$$
is distributed according to the Markov chain $\mathcal{K}_{p^{-1},1}$ started at state $\l_m$.
Moreover, because the entries on the diagonal of $A_4$ have valuations strictly less than $m$, we deduce that $\Cok(A_4-x_2B_4)\cong\Cok(A_4)$. Therefore, we have $\Cok\begin{pmatrix} A_0 & 0\\ 0 & A_4-x_2B_4\end{pmatrix}_i'=\l_i$ for all $1\le i\le m$, and $\Cok\begin{pmatrix} A_0 & 0\\ 0 & A_4-x_2B_4\end{pmatrix}_{m+i}'=\Cok(A_0)_{m+i}'$ for all $i\ge 1$. This completes the proof.
\end{proof}

\begin{proof}[Proof of \Cref{lem: expectation of two point determinant}]
Applying \Cref{thm:compute_char_poly_det} and \Cref{lem: cok and det}, we have
\begin{align}\label{eq:two_cokernels}
\begin{split}
\lim_{n\rightarrow\infty}\E||\det((A-x_1I_n)(A-x_2I_n))||&=\lim_{n\rightarrow\infty}\E||\det((A-x_1B)(A-x_2B))||\\
&=\lim_{n\rightarrow\infty}\E\left[\frac{1}{\#\Cok(A-x_1B)}\cdot\frac{1}{\#\Cok(A-x_2B)}\right],
\end{split}
\end{align}
where $B\in\Mat_{n\times n}(\Z_p)$ is independent of $A$ and also Haar-distributed over $\Mat_n(\Z_p)$. Therefore, we only need to compute the limit on the last line.

Now, let $n\ge 1$ be a fixed integer. Denote by
$\l_i=\Cok(A-x_1B)_i'$ for all $1\le i\le m$, and $\l_{m+i}^{(j)}=\Cok(A-x_jB)_{m+i}'$ for all $i\ge 1$ and $j\in\{1,2\}$. Following the distribution given in \Cref{lem: joint distribution of chains_2pt}, we have 
\begin{align}
\begin{split}
\E\left[\frac{1}{\#\Cok(A-x_1B)}\cdot\frac{1}{\#\Cok(A-x_2B)}\right]&=\E\left[p^{-2(\l_1+\cdots+\l_m)-(\l_{m+1}^{(1)}+\l_{m+2}^{(1)}+\cdots)-(\l_{m+1}^{(2)}+\l_{m+2}^{(2)}+\cdots)}\right]\\
&=\E \bigg[p^{-2(\l_1+\cdots+\l_m)}\\
&\cdot\E[p^{-(\l_{m+1}^{(1)}+\l_{m+2}^{(1)}+\cdots)-(\l_{m+1}^{(2)}+\l_{m+2}^{(2)}+\cdots)}\mid \l_m]\bigg]\\
&=\E \left[p^{-2(\l_1+\cdots+\l_m)}\cdot(\E[p^{-(\l_{m+1}^{(1)}+\l_{m+2}^{(1)}+\cdots)}\mid \l_m])^2\right]\\
&=\E \left[p^{-2(\l_1+\cdots+\l_m)}\cdot\frac{(1-p^{-1})^2}{(1-p^{-\l_m-1})^2}\right].\\
\end{split}
\end{align}
Here, the last line comes from \Cref{thm:compute_t_moment_single_markov}, and the fourth line holds because when we condition on $\l_m$, the two sequences $(\l_{m+1}^{(1)},\l_{m+2}^{(1)},\ldots)$ and $(\l_{m+1}^{(2)},\l_{m+2}^{(2)},\ldots)$ are independent and have the same law. Applying \Cref{thm:compute_2pt_expectations}, we have
\begin{align}
\begin{split}
\lim_{n\rightarrow\infty}\E||\det((A-x_1I_n)(A-x_2I_n))||&=\lim_{n\rightarrow\infty}\E\left[\frac{1}{\#\Cok(A-x_1B)}\cdot\frac{1}{\#\Cok(A-x_2B)}\right]\\
&=\lim_{n\rightarrow\infty}\E\left[p^{-2(\l_1+\cdots+\l_m)}\cdot\frac{(1-p^{-1})^2}{(1-p^{-\l_m-1})^2}\right]\\
&=(1-p^{-1})^2 \sum_{k \geq 0} (-1)^k p^{-\frac{(2m+1)k^2 + (4m+1)k}{2}}(1+p^{-k-1}).
\end{split}
\end{align}
This completes the proof.
\end{proof}

\begin{cor}\label{thm:variance_of_roots_in_Z_p}
The number of eigenvalues in $\Z_p$ has limiting variance
$$\lim_{n \to \infty} \E[(Z_{\Z_p,n} - \E[Z_{\Z_p,n}])^2] = \lim_{n\rightarrow\infty}\E[Z_{\Z_p\times\Z_p,n}]=\sum_{k\ge 0}\frac{(-1)^k(1-p^{-1})(1+p^{-k-1})p^{-\frac{k^2+k}{2}}}{1-p^{-k^2-2k-2}}.$$ 
\end{cor}

\begin{proof}
First, recall from \Cref{cor: expected eigenvalues in Zp} that we have $\E[Z_{\Z_p,n}]=1$ for all $n\ge 1$. Therefore, we have
\begin{align}
\begin{split}
\lim_{n \to \infty} \E[(Z_{\Z_p,n} - \E[Z_{\Z_p,n}])^2]&=\lim_{n\rightarrow\infty}\E[Z_{\Z_p,n}^2]-1\\
&=\lim_{n\rightarrow\infty}\E[Z_{\Z_p,n}(Z_{\Z_p,n}-1)]\\
&=\lim_{n\rightarrow\infty}\E[Z_{\Z_p\times\Z_p,n}].
\end{split}
\end{align}
Here, the last line holds because we always have $Z_{\Z_p\times\Z_p,n}=Z_{\Z_p,n}(Z_{\Z_p,n}-1)$, since both count ordered pairs of eigenvalues in $\Z_p$. Next, for all $m\ge 0$, the subset $$\mathcal{U}_m:=\{(x_1,x_2):x_1,x_2\in\Z_p,\val(x_1-x_2)=m\}\subset\Z_p\times\Z_p$$
has Haar probability measure $\mu_{\Q_p \times \Q_p}(\mathcal{U}_m)=(1-p^{-1})p^{-m}$.
Moreover, by \Cref{thm: two point correlation function in Z_p}, the limit correlation function $\rho_{\Q_p\times\Q_p}^{(\infty)}$ takes value $\rho_{\Q_p\times\Q_p}^{(\infty)}(p^m,0)$ in $\mathcal{U}_m$. Therefore, by \Cref{thm: limit of expectation of zeros}, we have 
\begin{align}
\begin{split}
\lim_{n\rightarrow\infty}\E[Z_{\Z_p\times\Z_p,n}]&=\int_{\Z_p\times\Z_p}\rho_{\Q_p\times\Q_p}^{(\infty)}(x_1,x_2)dx_1dx_2\\
&=\sum_{m\ge 0}\rho_{\Q_p\times\Q_p}^{(\infty)}(p^m,0)\cdot \mu_{\Q_p \times \Q_p}(\mathcal{U}_m)\\
&=\sum_{m\ge 0}(1-p^{-1})p^{-m}\rho_{\Q_p\times\Q_p}^{(\infty)}(p^m,0).
\end{split}
\end{align}
In the end, applying \Cref{thm: two point correlation function in Z_p}, we have
\begin{align}
\begin{split}
\sum_{m\ge 0}(1-p^{-1})p^{-m}\rho_{\Q_p\times\Q_p}^{(\infty)}(p^m,0)&=\sum_{m,k\ge 0}(-1)^k(1-p^{-1})(1+p^{-k-1})p^{-\frac{(2m+1)k^2+(4m+1)k}{2}-2m}\\
&=\sum_{k\ge 0}(-1)^k(1-p^{-1})(1+p^{-k-1})p^{-\frac{k^2+k}{2}}\sum_{m\ge 0}p^{-(k^2+2k+2)m}\\
&=\sum_{k\ge 0}\frac{(-1)^k(1-p^{-1})(1+p^{-k-1})p^{-\frac{k^2+k}{2}}}{1-p^{-k^2-2k-2}}.
\end{split}
\end{align}
This completes the proof.
\end{proof}

\begin{proof}[Proof of \Cref{thm:Zp_2pt_intro}]
First of all, given $x,y\in\Z_p$, let $m:=\val(x-y)$, so that $||x-y||=p^{-m}$. By \Cref{thm: two point correlation function in Z_p}, we have 
\begin{align}
\begin{split}
\rho_{\Q_p\times\Q_p}^{(\infty)}(x,y)&=p^{-m}\sum_{k \geq 0} (-1)^k p^{-\frac{(2m+1)k^2 + (4m+1)k}{2}}(1+p^{-k-1})\\
&=-\sum_{k \geq 1} (-1)^k p^{-\frac{(2m+1)k^2}{2}}(p^{k/2}+p^{-k/2})\\
&=1-\sum_{k\in\Z}\left(\frac{1}{p^{2m+1}}\right)^{\frac{k^2}{2}}\cdot(-\sqrt p)^k\\
&=1-\theta_3(-\sqrt{p}; ||x-y||^2/p).
\end{split}
\end{align}
Hence, we have verified \eqref{eq: 2pt correlation_intro}. Next, when $U,V\subset\Z_p$ are two disjoint measurable sets, we have $Z_{U,n}Z_{V,n}=Z_{U\times V,n}$, thus \eqref{eq: eigenvalues in UXV_intro} is a specific case of \Cref{thm: limit of expectation of zeros}. In the end, we recall from \Cref{thm: one point correlation function over Zp} that $\rho_{\Q_p}^{(n)}(x)=1$ for all $x\in\Z_p$ and $n\ge 1$. Therefore, we have
\begin{align}
\begin{split}
\lim_{n\rightarrow\infty}\Cov(Z_{U,n},Z_{V,n})
&=\lim_{n\rightarrow\infty}(\E[Z_{U,n}\cdot Z_{V,n}]-\E[Z_{U,n}]\E[Z_{V,n}])\\
&=\lim_{n \to \infty} \int_{U \times V} (\rho^{(\infty)}_{\Q_p,\Q_p}(x,y)-1) dx dy\\
&=-\lim_{n \to \infty} \int_{U \times V} \theta_3(-\sqrt{p}; ||x-y||^2/p) dx dy.
\end{split}
\end{align}
By the Jacobi triple product identity, the integrand in the last line (with $p^{-m}=||x-y||$, $m \geq 0$) is 
\begin{equation}
    \prod_{i\ge 1}\left(1-\left(\frac{1}{p^{2m+1}}\right)^i\right)\left(1-\left(\frac{1}{p^{2m+1}}\right)^{i-1/2}\cdot \sqrt p\right)\left(1-\left(\frac{1}{p^{2m+1}}\right)^{i-1/2}\cdot \frac{1}{\sqrt p}\right) \geq 0,
\end{equation}
so
\begin{equation}
    \lim_{n\rightarrow\infty}\Cov(Z_{U,n},Z_{V,n}) \leq 0.
\end{equation}
This verifies \eqref{eq:asymp_cov_Zp_intro} and  finishes the proof.
\end{proof}

\subsection{Quadratic extensions}

The main goal of this subsection is to prove \Cref{thm:quadratic_intro}. We will always work with the following setting. Let $K/\Q_p$ be a quadratic extension, and $\mathcal{O}_K=\Z_p[x_1]$ be the subring of integers. Here $\val(x_1)=\frac{1}{2}$ if $K/\Q_p$ is ramified, and $\val(x_1)=0$ if $K/\Q_p$ is unramified. Let $a_0,a_1 \in \Z_p$ with $\val(a_1)=m\ge 0$, and denote by $Z(x)$ the minimal polynomial of $a_0+a_1x_1$, which is monic of degree two. 

Before giving the proof of \Cref{thm:quadratic_intro}, let us start with some necessary related calculations.

\begin{prop}\label{thm: expection of quadratic extension}
Let $A\in\Mat_{n\times n}(\Z_p)$ be Haar-distributed. In the setting of this subsection, 
\begin{enumerate}
\item When $K/\Q_p$ is ramified, we have
$$\lim_{n\rightarrow\infty}\E||\det(Z(A))||=(1-p^{-1})\sum_{k \geq 0} p^{-(m+1)k^2 - (2m+1)k}(1-p^{-2k-2}).$$
\item When $K/\Q_p$ is unramified, we have
$$\lim_{n\rightarrow\infty}\E||\det(Z(A))||=(1-p^{-2})\sum_{k \geq 0} p^{-\frac{(2m+1)k^2+(4m+1)k}{2}}(1-p^{-k-1}).$$
\end{enumerate}
\end{prop}

\begin{lemma}\label{lem: distribution of chains_quadratic}
Let $n\ge 1$, and $A,B\in\Mat_n(\Z_p)$ be independent and Haar-distributed. Then, the distribution of the $\mathcal{O}_K$-module $\Cok(A-(a_0+a_1x_1)B)$ is given by the following way. 
\begin{enumerate}
\item When $K/\Q_p$ is ramified, $\Cok(A-(a_0+a_1x_1)B)_{2j-1}'=\Cok(A-(a_0+a_1x_1)B)_{2j}'$ for $1 \leq j \leq m$, and the random sequence
\begin{multline}
\Cok(A-(a_0+a_1x_1)B)_1'=\Cok(A-(a_0+a_1x_1)B)_2'=\l_1,\ldots,
\\\Cok(A-(a_0+a_1x_1)B)_{2m-1}'=\Cok(A-(a_0+a_1x_1)B)_{2m}'=\l_m
\end{multline}
is distributed according to the Markov chain $\mathcal{K}_{p^{-1},1}$ of \Cref{def:fulman_markov_chain} started at state $n$. Then, the sequence of subsequent terms 
$$\Cok(A-(a_0+a_1x_1)B)_{2m+1}'=\l_{m+1},\Cok(A-(a_0+a_1x_1)B)_{2m+2}'=\l_{m+2},\ldots$$
is distributed according to the Markov chain $\mathcal{K}_{p^{-1},1}$ started at state $\l_m$.
\item When $K/\Q_p$ is unramified, the random sequence
$$
\Cok(A-(a_0+a_1x_1)B)_1'=\l_1,\Cok(A-(a_0+a_1x_1)B)_2'=\l_2,\ldots,\Cok(A-(a_0+a_1x_1)B)_m'=\l_m
$$
is distributed according to the Markov chain $\mathcal{K}_{p^{-1},1}$ of \Cref{def:fulman_markov_chain} started at state $n$. Then, the sequence of subsequent terms 
$$\Cok(A-(a_0+a_1x_1)B)_{m+1}'=\l_{m+1},\Cok(A-(a_0+a_1x_1)B)_{m+2}'=\l_{m+2},\ldots$$
is distributed according to the Markov chain $\mathcal{K}_{p^{-2},1}$ started at state $\l_m$.
\end{enumerate}
\end{lemma}

\begin{proof}
Since $A\in\Mat_n(\Z_p)$ is Haar-distributed, there is no loss of generality to assume $a_0=0$. We first deal with the unramified case. Since $\val(a_1)=m$, we have $\Cok(A-a_1x_1B)_i=\Cok(A)_i$ for all $1\le i\le m$. Therefore, the distribution of the random sequence
$$\Cok(A-a_1x_1B)_1'=\Cok(A)_1',\Cok(A-a_1x_1B)_2'=\Cok(A)_2',\ldots,\Cok(A-a_1x_1B)_m'=\Cok(A)_m'$$
is already given in \Cref{thm:from_evans}. Now, given a fixed sequence of integers
$$n\ge\l_1\ge\ldots\ge\l_m,$$
we wish to understand the law of the subsequent terms $\Cok(A-a_1x_1B)_{m+i}'$ for all $i\ge 1$, conditioned on
\begin{equation}\label{eq: condition sequence_quadratic}
\Cok(A)_1'=\l_1,\ldots,\Cok(A)_m'=\l_m.
\end{equation}
For all $\l_m+1\le j\le n$, let
$\nu_j:=\#\{1\le i\le m:\l_i\ge j\}$. Let $A_1\in\Mat_{\lambda_m}(p^m\Z_p)$ be Haar-distributed, and $A_4:=\diag(p^{\nu_{\l_m+1}},\ldots,p^{\nu_n})$ be a fixed matrix. Then we have $\Cok(A_4)_i'=\l_i$ for all $1\le i\le m$. Moreover, by \Cref{thm:from_evans}, the sequence 
$$\Cok(A_1)_{m+1}',\Cok(A_1)_{m+2}',\ldots$$
is distributed according to the Markov chain $\mathcal{K}_{p^{-1},1}$ started at state $\l_m$. This implies that conditioned on $A$ satisfying \eqref{eq: condition sequence_quadratic},
$$\Cok(A)=\Cok\begin{pmatrix}A_1 & 0 \\ 0 & A_4\end{pmatrix},\quad \quad \quad \quad \text{in distribution.}$$
Since the Haar distribution of $B\in\Mat_n(\Z_p)$ is invariant under the left and right-multiplication of $\GL_n(\Z_p)$, we deduce that conditioned on $A$ satisfying \eqref{eq: condition sequence_quadratic},
$$\Cok(A-a_1x_1B)=\Cok\left(\begin{pmatrix}A_1 & 0 \\ 0 & A_4\end{pmatrix}-a_1x_1B\right),\quad \quad \quad \quad \text{in distribution.}$$
Denote 
$$B=\begin{pmatrix} B_1 & B_2\\ B_3 & B_4\end{pmatrix},$$
where 
$$B_1\in\Mat_{\lambda_m}(\Z_p),B_2\in\Mat_{\lambda_m\times(n-\lambda_m)}(\Z_p),B_3\in\Mat_{(n-\lambda_m)\times\lambda_m}(\Z_p),B_4\in\Mat_{n-\lambda_m}(\Z_p)$$ 
are all independently Haar-distributed. We have
\begin{align}
\begin{split}
\Cok\left(\begin{pmatrix}A_1 & 0 \\ 0 & A_4\end{pmatrix}-a_1x_1B\right)&=\Cok\begin{pmatrix}A_1-a_1x_1B_1 & -a_1x_1B_2\\ -a_1x_1B_3 & A_4-a_1x_1B_4\end{pmatrix}\\
&\cong\Cok\begin{pmatrix} A_1-a_1x_1B_1-B_0 & -a_1x_1B_2\\ 0 & A_4-a_1x_1B_4\end{pmatrix}\\
&\cong\Cok\begin{pmatrix} A_1-a_1x_1B_1-B_0 & 0\\ 0 & A_4-a_1x_1B_4\end{pmatrix},\\
\end{split}
\end{align}
where $B_0=a_1^2x_1^2B_2(A_4-a_1x_1B_4)^{-1}B_3\in p^m\Mat_{\l_m}(\mathcal{O}_K)$ is independent of the matrix $A_1-a_1x_1B_1$, which is Haar-distributed in $\Mat_{\l_m}(p^m\mathcal{O}_K)$. Therefore, the matrix $A_0:=A_1-a_1x_1B_1-B_0$ is Haar-distributed in $\Mat_{\l_m}(p^m\mathcal{O}_K)$. By \Cref{thm:from_evans}, the sequence 
$$\Cok(A_0)_{m+1}',\Cok(A_0)_{m+2}',\ldots$$
is distributed according to the Markov chain $\mathcal{K}_{p^{-2},1}$ started at state $\l_m$.
Moreover, because the entries on the diagonal of $A_4$ have valuations strictly less than $m$, we deduce that $\Cok(A_4-a_1x_1B_4)\cong\Cok(A_4)$. Therefore, we have $\Cok\begin{pmatrix} A_0 & 0\\ 0 & A_4-a_1x_1B_4\end{pmatrix}_i'=\l_i$ for all $1\le i\le m$, and $\Cok\begin{pmatrix} A_0 & 0\\ 0 & A_4-a_1x_1B_4\end{pmatrix}_{m+i}'=\Cok(A_0)_{m+i}'$ for all $i\ge 1$. 

For the ramified case, we can regard $A$ either as a matrix in $\Mat_n(\Z_p)$ or a matrix in $\Mat_n(\mathcal{O}_K)$. Moreover, we have $\Cok_{\mathcal{O}_K}(A)'_{2i-1}=\Cok_{\mathcal{O}_K}(A)_{2i}'=\Cok_{\Z_p}(A)_{i}'$ for all $i\ge 1$, so that the distribution of the random sequence
\begin{multline}
\Cok(A-a_1x_1B)_1'=\Cok(A-a_1x_1B)_2'=\Cok_{\Z_p}(A)'_1,\ldots,
\\\Cok(A-a_1x_1B)_{2m-1}'=\Cok(A-a_1x_1B)_{2m}'=\Cok_{\Z_p}(A)'_m
\end{multline}
is already given in \Cref{thm:from_evans}. The rest of the proof follows similarly.
\end{proof}

\begin{proof}[Proof of \Cref{thm: expection of quadratic extension}]
Applying \Cref{thm:compute_char_poly_det} and \Cref{lem: cok and det}, we have
\begin{align}
\begin{split}
\lim_{n\rightarrow\infty}\E||\det(Z(A))||&=\lim_{n\rightarrow\infty}\E||\det(A-(a_0+a
_1x_1)B)||^2\\
&=\lim_{n\rightarrow\infty}\E\left[\frac{1}{\#\Cok(A-(a_0+a
_1x_1)B)}\right].
\end{split}
\end{align}
where $B\in\Mat_{n\times n}(\Z_p)$ is independent of $A$ and also Haar-distributed. Therefore, we only need to compute the limit on the last line.

Now, let $n\ge 1$ be a fixed integer. When $K/\Q_p$ is ramified, denote by $\l_i:=\Cok(A-(a_0+a_1x_1)B)_{2i}'$ for all $1\le i\le m$, and $\l_{m+i}:=\Cok(A-(a_0+a_1x_1)B)_{2m+i}'$ for all $i\ge 1$. Following the distribution given in \Cref{lem: distribution of chains_quadratic}, we have 
\begin{align}
\begin{split}
\E\left[\frac{1}{\#\Cok(A-(a_0+a
_1x_1)B)}\right]&=\E \left[p^{-2\l_1-\cdots-2\l_m-\l_{m+1}-\l_{m+2}-\cdots}\right]\\
&=\E \left[p^{-2\l_1-\cdots-2\l_m}\cdot\E[p^{-(\l_{m+1}+\l_{m+2}+\cdots)}\mid \l_m]\right]\\
&=\E \left[p^{-2\l_1-\cdots-2\l_m}\cdot\frac{1-p^{-1}}{1-p^{-\l_m-1}}\right].\\
\end{split}
\end{align}
Here, the last line comes from \Cref{thm:compute_t_moment_single_markov}. Applying \Cref{thm:compute_2pt_expectations}, we have 
\begin{align}
\begin{split}
\lim_{n\rightarrow\infty}\E\left[\frac{1}{\#\Cok(A-(a_0+a
_1x_1)B)}\right]&=\lim_{n\rightarrow\infty}\E \left[p^{-2\l_1-\cdots-2\l_m}\cdot\frac{1-p^{-1}}{1-p^{-\l_m-1}}\right]\\
&=(1-p^{-1})\sum_{k \geq 0} p^{-(m+1)k^2 - (2m+1)k}(1-p^{-2k-2}).
\end{split}
\end{align}
When $K/\Q_p$ is unramified, denote by $\l_i:=\Cok(A-(a_0+a
_1x_1)B)_i'$ for all $i\ge 1$. Following the distribution given in \Cref{lem: distribution of chains_quadratic}, we deduce that for fixed $n\ge 1$,
\begin{align}
\begin{split}
\E\left[\frac{1}{\#\Cok(A-(a_0+a
_1x_1)B)}\right]&=\E \left[p^{-\l_1-\cdots-\l_m-\l_{m+1}-\l_{m+2}-\cdots}\right]\\
&=\E \left[p^{-\l_1-\cdots-\l_m}\cdot\E[p^{-(\l_{m+1}+\l_{m+2}+\cdots)}\mid \l_m]\right]\\
&=\E \left[p^{-\l_1-\cdots-\l_m}\cdot\frac{1-p^{-2}}{1-p^{-2\l_m-2}}\right].\\
\end{split}
\end{align}
Here, the last line comes from \Cref{thm:compute_t_moment_single_markov}. Applying \Cref{thm:compute_2pt_expectations}, we have 
\begin{align}
\begin{split}
\lim_{n\rightarrow\infty}\E\left[\frac{1}{\#\Cok(A-(a_0+a
_1x_1)B)}\right]&=\lim_{n\rightarrow\infty}\E \left[p^{-\l_1-\cdots-\l_m}\cdot\frac{1-p^{-2}}{1-p^{-2\l_m-2}}\right]\\
&=(1-p^{-2})\sum_{k \geq 0} p^{-\frac{(2m+1)k^2+(4m+1)k}{2}}(1-p^{-k-1}).
\end{split}
\end{align}
This completes the proof.
\end{proof}

\begin{rmk}
We consider two degenerate cases of the calculations in \Cref{thm: expection of quadratic extension}. First, when $m=0$, we have
$$\lim_{n\rightarrow\infty}\E||\det(Z(A))||=\begin{cases}
1-p^{-1} & K/\Q_p\text{ ramified}\\
1-p^{-2} & K/\Q_p\text{ unramified},
\end{cases}$$
which coincides with \Cref{prop: limit expectation over generator}. Also, when $m$ goes to infinity, we have
$$\lim_{m\rightarrow\infty}\lim_{n\rightarrow\infty}\E||\det(Z(A))||=(1-p^{-1})(1-p^{-2}).$$
This coincides with \Cref{thm:moments_of_determinant}, which provides a formula for $\lim_{n\rightarrow\infty}\E||\det(A-a_0I)||^2=\lim_{n\rightarrow\infty}\E||\det(A)||^2$ that yields the same result.
\end{rmk}

\begin{prop}\label{prop: orbital integral of quadratic extension}
We have 
$$\#(\Lambda\backslash\Mod_{\Z_p[a_0+a_1x_1]})=\frac{1-p^{m+1}}{1-p}
$$
when $K/\Q_p$ is ramified, and
$$\#(\Lambda\backslash\Mod_{\Z_p[a_0+a_1x_1]})=p^m+2p^{m-1}+2p^{m-2}+\cdots+2=\frac{1-p^m}{1-p}+\frac{1-p^{m+1}}{1-p}
$$
when $K/\Q_p$ is unramified. 
\end{prop}

\begin{proof}
Notice that $\Z_p[a_0+a_1x_1]=\Z_p+p^m\Z_px_1$. Our calculation follows the equation \eqref{eq: orbital integral as sum of orbit}. Under the group action of $K^\times$, every orbit in $\Mod_{\Z_p+p^m\Z_px_1}$ has a unique representative that contains $1$, but does not contain an element with absolute value strictly greater than $1$. Such a representative by definition lies inside $\mathcal{O}_K=\Z_p+\Z_px_1$, and it contains the lattice $\Z_p+p^m\Z_px_1$ because it contains the element $1$. Notice that the $\Z_p$-submodules of $\mathcal{O}_K/(\Z_p+p^m\Z_px_1)\cong \Z_p/p^m\Z_p$ are exactly $p^k\Z_p/p^m\Z_p$, where $0\le k\le m$. Lifting it back to elements in $\Mod_{\Z_p+p^m\Z_px_1}$, we deduce that a set of representatives of $K^\times\backslash\Mod_{\Z_p[a_0+a_1x_1]}$ has the form 
$$\mathcal{O}_K=\Z_p+\Z_px_1,\Z_p+p\Z_px_1,\ldots,\Z_p+p^m\Z_px_1.$$
For all $1\le k\le m$, we have $\Aut(\Z_p+p^k\Z_px_1)=\Z_p^\times+p^k\Z_px_1$ with the right-hand side acting by multiplication. When $k=0$ we have $\Aut(\Z_p+p^k\Z_px_1) = \mc{O}_K^\times$, which is given explicitly by $$\mathcal{O}_K^\times=\Z_p^\times+\Z_px_1$$ for the ramified case, and 
$$\mathcal{O}_K^\times=\{b_0+b_1x_1\mid b_0,b_1\text { are not both in }p\Z_p\}$$
for the unramified case. 

When $1\le k\le m$, both $\Aut(\Z_p+p^k\Z_px_1)=\Z_p^\times+p^k\Z_px_1$ and $\mathcal{O}_K^\times$ can be written (as sets) as union of cosets of the $\Z_p$-submodule $p\Z_p+p^k\Z_px_1 \subset \mc{O}_K$:
\begin{equation}
    \Z_p^\times+p^k\Z_px_1 = \bigsqcup_{b=1}^{p-1} b + p\Z_p+p^k\Z_px_1,
\end{equation}
and 
\begin{equation}
\mc{O}_K^\times =\begin{cases}
\bigsqcup_{1\le b_0\le p-1, 0\le b_1\le p^k-1}b_0+b_1x_1+p\Z_p+p^k\Z_px_1 & K/\Q_p \text{ ramified}\\
\bigsqcup_{\substack{0\le b_0\le p-1, 0\le b_1\le p^k-1\\ b_0\ne 0 \text{ or }p\nmid b_1}}b_0+b_1x_1+p\Z_p+p^k\Z_px_1 & K/\Q_p \text{ unramified}
\end{cases}.
\end{equation}

Therefore, the number of cosets that make up $\Aut(\Z_p+p^k\Z_px_1)$ is $p-1$, while the number of cosets that make up $\mathcal{O}_K^\times$ is $(p-1)p^k$ for the ramified case, or $(p^2-1)p^{k-1}$ for the unramified case. Hence, for $0\le k\le m$,
$$\#(\mathcal{O}_K^\times/\Aut(\Z_p+p^k\Z_px_1))=\begin{cases}
1 & k=0\\
p^k & k \geq 1, K/\Q_p\text{ ramified}\\
p^k+p^{k-1} & k \geq 1, K/\Q_p\text{ unramified}
\end{cases}.$$
Now, applying \eqref{eq: orbital integral as sum of orbit}, we have
$$\#(\Lambda\backslash\Mod_{\Z_p[a_0+a_1x_1]})=1+\sum_{k=1}^m p^k=\frac{1-p^{m+1}}{1-p}$$
for the ramified case, and 
$$\#(\Lambda\backslash\Mod_{\Z_p[a_0+a_1x_1]})=1+\sum_{k=1}^m (p^k+p^{k-1})=\frac{1-p^m}{1-p}+\frac{1-p^{m+1}}{1-p}$$
for the unramified case.
\end{proof}

\begin{proof}[Proof of \Cref{thm:quadratic_intro}]
We have 
\begin{equation}\label{eq:delta_squared_m}
    ||\Delta_\sigma(a_0+a_1x_1)||^2=p^{-2m}||\Delta_\sigma(x_1)||^2=p^{-2m}||\Disc_{K/\Q_p}||.
\end{equation}
Therefore, by \Cref{thm:limit_cor_fns_exist}, when $K/\Q_p$ is ramified,
\begin{align}
\begin{split}
\rho_K^{(\infty)}(a_0+a_1x_1)&=||\Delta_\sigma(a_0+a_1x_1)|| \cdot \Den(a_0+a_1x_1)  \cdot \lim_{n\rightarrow\infty}\E||\det(Z(A))||\\
&=\frac{p^{-2m}||\Disc_{K/\Q_p}||}{1-p^{-1}}\cdot\frac{1-p^{m+1}}{1-p}\cdot \lim_{n\rightarrow\infty}\E||\det(Z(A))||\\
&=||\Disc_{K/\Q_p}||\frac{p^{-m}(1-p^{-m-1})}{1-p^{-1}}\sum_{k \geq 0} p^{-(m+1)k^2 - (2m+1)k}(1-p^{-2k-2}).
\end{split}
\end{align}
Here, the second line comes from \Cref{prop: orbital integral of quadratic extension} (and we recall from \Cref{defi: Den and distance} that $V(a_0+a_1x_1)$ contributes an additional $\Delta_\sigma$ term to give us \eqref{eq:delta_squared_m}), and the third line comes from \Cref{thm: expection of quadratic extension}. The unramified case follows similarly, noting that $||\Disc_{K/\Q_p}||=1$ in this case.
\end{proof}

By integrating the correlation function in \Cref{thm:quadratic_intro}, we obtain the following.

\begin{cor}\label{thm:number_of_roots_in_quadratic_extension}
Let $K/\Q_p$ be a quadratic extension. Then we have 
$$\lim_{n\rightarrow\infty}\E[Z_{\mathcal{O}_K^{\new},n}]=||\Disc_{K/\Q_p}||(1-p^{-1})\sum_{k\ge 0}\frac{(1-p^{-2k-2})p^{-k^2-k}}{(1-p^{-k^2-2k-2})(1-p^{-k^2-2k-3})}$$
when $K/\Q_p$ is ramified, and
$$\lim_{n\rightarrow\infty}\E[Z_{\mathcal{O}_K^{\new},n}]=(1-p^{-1})\sum_{k\ge 0}\frac{(1-p^{-k-1})(1+p^{-k^2-2k-3})p^{-\frac{k^2+k}{2}}}{(1-p^{-k^2-2k-2})(1-p^{-k^2-2k-3})}$$
when $K/\Q_p$ is unramified.
\end{cor}

\begin{proof}
Following the setting stated at the beginning of this subsection, denote by $x_1$ a generator of $\mathcal{O}_K$. For all $m\ge 0$, the subset $$\mathcal{U}_m:=\{b_0+b_1x_1:b_0,b_1\in\Z_p,\val(b_1)=m\}=(\Z_p+p^m\Z_px_1)\backslash(\Z_p+p^{m+1}\Z_px_1)\subset\mathcal{O}_K$$
has Haar probability measure $\mu_{K}(\mathcal{U}_m)=(1-p^{-1})p^{-m}$.
Moreover, by \Cref{thm:quadratic_intro}, the limit correlation function $\rho_K^{(\infty)}$ takes value $\rho_K^{(\infty)}(p^mx_1)$ in $\mathcal{U}_m$. Therefore, by \Cref{thm: limit of expectation of zeros}, we have 
\begin{align}
\begin{split}
\lim_{n\rightarrow\infty}\E[Z_{\mathcal{O}_K^{\new},n}]&=\int_{\mathcal{O}_{K}}\rho_K^{(\infty)}(x)dx\\
&=\sum_{m\ge 0}\rho_K^{(\infty)}(p^mx_1)\cdot \mu_K(\mathcal{U}_m)\\
&=\sum_{m\ge 0}\frac{1-p^{-1}}{p^m}\rho_K^{(\infty)}(p^mx_1).
\end{split}
\end{align}
By \Cref{thm:quadratic_intro}, when $K/\Q_p$ is ramified, we have 
\begin{align}
\begin{split}
\sum_{m\ge 0}\frac{1-p^{-1}}{p^m}\rho_K^{(\infty)}(p^mx_1)&=\sum_{m,k\ge 0}||\Disc_{K/\Q_p}||\frac{(1-p^{-m-1})(1-p^{-2k-2})}{p^{(m+1)k^2+(2m+1)k+2m}}\\
&=\sum_{k\ge 0}||\Disc_{K/\Q_p}||\frac{1-p^{-2k-2}}{p^{k^2+k}}\sum_{m\ge 0}\frac{1-p^{-m-1}}{p^{(k^2+2k+2)m}}\\
&=\sum_{k\ge 0}||\Disc_{K/\Q_p}||\frac{1-p^{-2k-2}}{p^{k^2+k}}\left(\frac{1}{1-p^{-k^2-2k-2}}-\frac{p^{-1}}{1-p^{-k^2-2k-3}}\right)\\
&=||\Disc_{K/\Q_p}||(1-p^{-1})\sum_{k\ge 0}\frac{(1-p^{-2k-2})p^{-k^2-k}}{(1-p^{-k^2-2k-2})(1-p^{-k^2-2k-3})}.
\end{split}
\end{align}
When $K/\Q_p$ is unramified, we have 
\begin{align}
\begin{split}
\sum_{m\ge 0}\frac{1-p^{-1}}{p^m}\rho_K^{(\infty)}(p^mx_1)&=\sum_{m,k\ge 0}(1+p^{-1}-2p^{-m-1})(1-p^{-k-1})p^{-\frac{(2m+1)k^2+(4m+1)k}{2}-2m}\\
&=\sum_{k\ge 0}(1-p^{-k-1})p^{-\frac{k^2+k}{2}}\sum_{m\ge 0}p^{-(k^2+2k+2)m}(1+p^{-1}-2p^{-m-1})\\
&=\sum_{k\ge 0}(1-p^{-k-1})p^{-\frac{k^2+k}{2}}\left(\frac{1+p^{-1}}{1-p^{-k^2-2k-2}}-\frac{2p^{-1}}{1-p^{-k^2-2k-3}}\right)\\
&=(1-p^{-1})\sum_{k\ge 0}\frac{(1-p^{-k-1})(1+p^{-k^2-2k-3})p^{-\frac{k^2+k}{2}}}{(1-p^{-k^2-2k-2})(1-p^{-k^2-2k-3})}.
\end{split}
\end{align}
This completes the proof.
\end{proof}

%% file: repulsion.tex
\section{Estimates on repulsion of eigenvalues}\label{sec:pair_repulsion}

It is natural to ask for the two point correlation of eigenvalues. That is to say, given two eigenvalues $x_1,x_2\in\bar\Z_p$, are they attracted, independent, or repelled? In \Cref{thm: independent distribution over lifted subspaces}, we proved that 
two eigenvalues $x_1,x_2$ are independent when they belong to different lifted subspaces, i.e.,
$$\rho_{\Q_p[x_1]\times\Q_p[x_2]}^{(\infty)}(x_1,x_2)=\rho_{\Q_p[x_1]}^{(\infty)}(x_1)\rho_{\Q_p[x_2]}^{(\infty)}(x_2).$$
In this section, we prove that two eigenvalues repel when they belong to the same lifted subspace. We let 
\begin{enumerate}
\item $\mathcal{U}_i$ be the lifted subspace corresponding to $F_i\in\F_p[x]$ with $\deg F_i=d_i$.
\item $x_1,x_2\in\mathcal{U}_i$ be two elements in $\mathcal{U}_i$.
\item $K_1=\Q_p[x_1],K_2=\Q_p[x_2]$.
\item $Z_1,Z_2$ be the monic minimal polynomials of $x_1,x_2$, respectively. We furthermore require $Z_1\ne Z_2$, i.e., $x_1$ and $x_2$ do not lie in the same Galois orbit.
\end{enumerate}

The following result asserts that when $p^{d_i}$ is large, the ratio of $\rho_{\Q_p[x_1]\times\Q_p[x_2]}^{(\infty)}(x_1,x_2)$ over $\rho_{\Q_p[x_1]}^{(\infty)}(x_1)\rho_{\Q_p[x_2]}^{(\infty)}(x_2)$ is approximately $||\Res(Z_1,Z_2)||$.

\begin{thm}\label{thm: estimate of two point correlation}
We have
\begin{equation}\label{eq: lower and upper bound of two point correlation}
||\Res(Z_1,Z_2)||(p^{-d_i};p^{-d_i})_\infty<\frac{\rho_{\Q_p[x_1]\times\Q_p[x_2]}^{(\infty)}(x_1,x_2)}{\rho_{\Q_p[x_1]}^{(\infty)}(x_1)\rho_{\Q_p[x_2]}^{(\infty)}(x_2)}<\frac{||\Res(Z_1,Z_2)||}{(p^{-d_i};p^{-d_i})_\infty}.
\end{equation}
\end{thm}

Recall from \Cref{prop: resultant positive power} that $||\Res(Z_1,Z_2)||$ is a positive integer power of $p^{-d_i}$. It is easy to prove that the right-hand side of \eqref{eq: lower and upper bound of two point correlation} is less than $1$ when $||\Res(Z_1,Z_2)||=p^{-kd_i}$ with $k\ge 2$. However, when $||\Res(Z_1,Z_2)||=p^{-d_i}$, this upper bound might be still larger than $1$, because one can take $p=2,d_i=1,$ and $||\Res(Z_1,Z_2)||=1/2$. To truly confirm that the eigenvalues in $\mathcal{U}_i$ repel each other, we make a slight improvement of the above estimate for the particular case $||\Res(Z_1,Z_2)||=p^{-d_i}$.

\begin{thm}\label{thm: refinement of two point correlation}
Suppose $||\Res(Z_1,Z_2)||=p^{-d_i}$. Then we have
$$p^{-d_i}<\frac{\rho_{\Q_p[x_1]\times\Q_p[x_2]}^{(\infty)}(x_1,x_2)}{\rho_{\Q_p[x_1]}^{(\infty)}(x_1)\rho_{\Q_p[x_2]}^{(\infty)}(x_2)}\le p^{-d_i}+p^{-2d_i}.$$
\end{thm}

The rest of this section aims to prove \Cref{thm: estimate of two point correlation} and \Cref{thm: refinement of two point correlation}. The following lemma provides a useful elementary bound for expectations of determinants. For example, we may take $Z=Z_1,Z=Z_2$, or $Z=Z_1Z_2$.

\begin{lemma}\label{lem: elementry bound for expectation}
Let $Z\in\Z_p[x]$ be monic, the residue of which is a power of $F_i$. Then we have
$$(p^{-d_i};p^{-d_i})_\infty<\lim_{n\rightarrow\infty}\E||\det(Z(A))||<1.$$
\end{lemma}

\begin{proof}
The statement in \Cref{thm: limit distribution of cok} implies that
$$\lim_{n\rightarrow\infty}\mathbf{P}(||\det(Z(A))||=1)=\lim_{n\rightarrow\infty}\mathbf{P}(Z(A)\text{ invertible})=(p^{-d_i};p^{-d_i})_\infty.$$
On the one hand, when $Z(A)$ is not invertible, we always have $||\det(Z(A))||\le 1/p$. This proves the upper bound. On the other hand, there exists a nontrivial finite $\Z_p[x]/(Z)$-module, whose probability in \Cref{thm: limit distribution of cok} is strictly positive. This proves the lower bound.
\end{proof}

\begin{proof}[Proof of \Cref{thm: estimate of two point correlation}]
By \Cref{thm:limit_cor_fns_exist}, we have
$$\frac{\rho_{\Q_p[x_1]\times\Q_p[x_2]}^{(\infty)}(x_1,x_2)}{\rho_{\Q_p[x_1]}^{(\infty)}(x_1)\rho_{\Q_p[x_2]}^{(\infty)}(x_2)}=||\Res(Z_1,Z_2)||\cdot\lim_{n\rightarrow\infty}\frac{\E||\det(Z_1(A))\det(Z_2(A))||}{\E||\det(Z_1(A))||\cdot\E||\det(Z_2(A))||}.$$
Then we apply \Cref{lem: elementry bound for expectation}. On the one hand, 
$$\E||\det(Z_1(A))||>(p^{-d_i};p^{-d_i})_\infty$$
and
$$\E||\det(Z_1(A))\det(Z_2(A))||\le\E||\det(Z_2(A))||,$$ 
so we get the upper bound. On the other hand, $$\E||\det(Z_1(A))\det(Z_2(A))||>(p^{-d_i};p^{-d_i})_\infty$$ 
and 
$$\E||\det(Z_1(A))||,\E||\det(Z_2(A))||<1,$$ 
so we get the lower bound.
\end{proof}

The following lemma is crucial for our proof of \Cref{thm: refinement of two point correlation}.

\begin{lemma}\label{lem: DVR when large resultant}
If $||\Res(Z_1,Z_2)||=p^{-d_i}$, then $\Z_p[x_1]=\mathcal{O}_{K_1}$ and $\Z_p[x_2]=\mathcal{O}_{K_2}$.
\end{lemma}

\begin{proof}
We only prove the claim for $\Z_p[x_1]$, as it is symmetric. It suffices to show the implication 
\begin{equation}\label{eq:val_to_ring}
    \val f(x_1) \geq 0 \Rightarrow f(x) \in \Z_p[x]
\end{equation}
for all $f \in \Q_p[x]$ such that $\deg f\le \deg Z_1-1$. The reason is that \eqref{eq:val_to_ring} implies $\{f(x_1): f \in \Q_p[x], \val f(x_1)\ge 0\} \subset \Z_p[x_1]$, and the reverse inclusion is clear, while the left-hand side is the unit ball $\mc{O}_{K_1}$ of $\Q_p[x_1]$ by definition. 

Since the residue $\overline Z_1(x) \in \F_p[x]$ is a power of $F_i$, one of the following must hold:
\begin{enumerate}
\item $\deg Z_1=d_i$.
\item $Z_1=Y^k+pZ_3$, where $k\ge 2$, and $Y\in\Z_p[x]$ is monic and a lift of $F_i$.
\end{enumerate}
We will prove these two cases separately.

If $\deg Z_1=d_i$, then for any $f\in\Z_p[x]$ with non-zero residue and $\deg f\le d_i-1$, we have $||\Res(f,Z_1)||=1$ by \Cref{prop: resultant positive power}, which implies $\val f(x_1)=0$. Now suppose $f\in\Q_p[x],\deg f\le d_i-1$ has a coefficient that is not in $\Z_p$. Then we can find an integer $m\ge 1$ such that $p^mf\in\Z_p[x]$ has non-zero residue, and therefore 
$$\val f(x_1)=\val(p^mf(x_1))-m=-m<0.$$
This verifies the implication \eqref{eq:val_to_ring}.

If $Z_1=Y^k+pZ_3$, we claim that the residue of $Z_3$ is not divisible by $F_i$. Suppose for the sake of contradiction that this is not true. Then $\val Z_3(x_2)>0$. Combining this with the fact that $\val(Z_1(x_2))=\frac{1}{\deg Z_2}\val\Res(Z_1,Z_2)\le 1$, we have 
$$\val(Y^k(x_2))=\val(Z_1(x_2)-pZ_3(x_2))=\val(Z_1(x_2))$$
by the strong triangle inequality. Therefore, $||\Res(Y^k,Z_2)||=||\Res(Z_1,Z_2)||=p^{-d_i}$. However, $||\Res(Y^k,Z_2)||=||\Res(Y,Z_2)||^k$ and $||\Res(Y,Z_2)||$ must be a positive integer power of $p^{-d_i}$, which is a contradiction. 

Hence the residue of $Z_3$ is not divisible by $F_i$, and $\val Z_3(x_1)=0$. Also, we have $\val Y^k(x_1)=1+\val Z_3(x_1)=1$, which implies $\val Y(x_1)=1/k$. Now, for any $f\in\Z_p[x]$ with non-zero residue and $\deg f\le \deg Z_i-1=kd_i-1$, we can write the expression
$$f=Y^{k_0}Z_4+pZ_5,$$
where $0\le k_0<k$ is the number of times $F_i$ appears in the decomposition of the residue of $f$, $Z_4\in\Z_p[x]$ with residue not divisible by $F_i$, and $Z_5\in\Z_p[x]$. In this case, we have $\val(Y^{k_0}(x_1)Z_4(x_1))=k_0\val Y(x_1)+\val Z_4(x_1)=k_0/k<1$, and
$$\val f(x_1)=\val(Y^{k_0}(x_1)Z_4(x_1)+pZ_5(x_1))=k_0/k,$$
which is strictly less than $1$. Following the same proof by contradiction as in the $\deg Z_1=d_i$ case, we once again confirmed the implication \eqref{eq:val_to_ring}.
\end{proof}

\begin{proof}[Proof of \Cref{thm: refinement of two point correlation}]
By \Cref{lem: DVR when large resultant}, we have $\Z_p[x_1]=\mathcal{O}_{K_1}$, $\Z_p[x_2]=\mathcal{O}_{K_2}$ are discrete valuation rings with residue field size $p^{d_i}$. Therefore, applying \Cref{thm:compute_char_poly_det}, we have
$$\lim_{n\rightarrow\infty}\E||\det(Z_1(A))||=\lim_{n\rightarrow\infty}\E||\det(Z_2(A))||=1-p^{-d_i},$$ and
$$\frac{\rho_{\Q_p[x_1]\times\Q_p[x_2]}^{(\infty)}(x_1,x_2)}{\rho_{\Q_p[x_1]}^{(\infty)}(x_1)\rho_{\Q_p[x_2]}^{(\infty)}(x_2)}=p^{-d_i}\frac{\lim_{n\rightarrow\infty}\E||\det(Z_1(A)Z_2(A))||}{(1-p^{-d_i})^2}.$$
On the one hand, applying the Cauchy-Schwarz inequality, we obtain the upper bound by showing
\begin{align}
\begin{split}
\lim_{n\rightarrow\infty}\E||\det(Z_1(A)Z_2(A))||&\le\lim_{n\rightarrow\infty}\left(\E||\det(Z_1(A))||^2\right)^{1/2}\cdot\lim_{n\rightarrow\infty}\left(\E||\det(Z_2(A))||^2\right)^{1/2}\\
&=(1-p^{-d_i})(1-p^{-2d_i}).
\end{split}
\end{align}
Here, the second line comes from \Cref{prop: limit expectation over generator}. On the other hand, by \Cref{lem: elementry bound for expectation},
\begin{align}
\begin{split}
\lim_{n\rightarrow\infty}\E||\det(Z_1(A)Z_2(A))||&>(p^{-d_i};p^{-d_i})_\infty\\
&>(1-p^{-d_i})\cdot\left(1-\sum_{k=2}^\infty p^{-kd_i}\right)\\
&\ge (1-p^{-d_i})^2.
\end{split}
\end{align}
Here the second line follows by the generalized Bernoulli's inequality
\begin{equation}
    \prod_{i=1}^k (1+x_i) > 1 + \sum_{i=1}^k x_i \quad \quad \quad \quad \text{ for all $x_1,\ldots,x_k \in (-1,0)$}.
\end{equation}
This gives the lower bound.
\end{proof}

%% file: average_number_high_degree.tex
\section{Average numbers of eigenvalues in extensions of degree {$r > 2$}}\label{sec:average_number_of_eigenvalues_in_high_degree_extensions}

Let $K/\Q_p$ be a finite extension. In this section, we study the orders of magnitude of the expectation of eigenvalues that fall inside $K$ but not any proper subfield of $K$. We have already done this exactly for $K$ of degree $1$ and $2$, but in general we are still able to compute accurate estimates when the residual degree is high. We apply the following notations from \cite[Section 4]{caruso2022zeroes}. Let

\begin{enumerate}
\item $K/\Q_p$ be a finite extension of degree $r$.
\item $\mathcal{O}_K$ be the ring of integers of $K$, equipped with the probability Haar measure.
\item $e$ be the ramification degree of $K/\Q_p$.
\item $f$ be the residue degree of $K/\Q_p$, so $r=ef$, and the residue field of $\mathcal{O}_K$ is isomorphic to $\F_{p^f}$.
\item For all $d\ge 1$, $\mathcal{G}_d=\{\alpha\in\F_{p^d}\mid \F_{p^d}=\F_p[\alpha]\}$ be the set of generators of $\F_{p^d}$ over $\F_p$. Denote by $G_d=\#\mathcal{G}_d$ the cardinality of $\mathcal{G}_d$.
\item $\mu:\Z_{>0}\rightarrow\{0,1,-1\}$ be the M\"obius function.
\item $\mathcal{G}_K=\{x\in\mathcal{O}_K\mid \mathcal{O}_K=\Z_p[x]\}$ be the set of generators of $\mathcal{O}_K$.
\item For all $\alpha\in \F_{p^f}$, $U_\alpha$ be the open subset of $\mathcal{O}_K$ whose image in the residue field is $\alpha$.
\item $V=\mathcal{O}_K\backslash\mathcal{G}_K$ be the complementary set of $\mathcal{G}_K$ in $\mathcal{O}_K$, and $V_\alpha=V\cap U_\alpha$ for all $\alpha\in\F_{p^f}$.
\end{enumerate}

The goal of this section is to prove the following theorem. 

\begin{thm}\label{thm: expectation of higher degree extension}
We have 
$$-p^{-f}\le\frac{\lim_{n\rightarrow\infty}\E[Z_{\mathcal{O}_K^{\new},n}]}{||\Disc_{K/\Q_p}||}-\sum_{d\mid f}\mu(\frac{f}{d})p^{d-f}< \frac{1+\tau(f)}{1-p^{-f}}p^{-f},$$
where $\tau(f)=\#\{d\in\Z_{>0}\mid d\text{ divides }f\}$ is the number of positive divisors of $f$.
\end{thm}

\begin{lemma}\label{lem: Mobius inverse formula}
$$G_f=\sum_{d\mid f}\mu(\frac{f}{d})p^d.$$
\end{lemma}

\begin{proof}
This is the M\"obius inversion of the relation $\sum_{m\mid f}G_m=p^f$, which is clear.
\end{proof}

\begin{prop}\label{prop: integral of density_generators}
We have 
$$-p^{-f}\le\frac{\int_{\mathcal{G}_K}\rho_K^{(\infty)}(x)dx}{||\Disc_{K/\Q_p}||}-\sum_{d\mid f}\mu(\frac{f}{d})p^{d-f}\le 0.$$
\end{prop}

\begin{proof}
Following \Cref{thm: limit density over generator}, we have
$$\frac{\int_{\mathcal{G}_K}\rho_K^{(\infty)}(x)dx}{||\Disc_{K/\Q_p}||}=\int_{\mathcal{G}_K}dx=\sum_{\alpha\in\mathcal{G}_f}\int_{\mathcal{G}_K\cap U_\alpha}dx.$$
By \cite[Lemma 4.4]{caruso2022zeroes}\footnote{In the notation of \cite[Lemma 4.4]{caruso2022zeroes}, we are taking $q=p$, and $\lambda_K=\mu_K$ denotes the Haar probability measure over $\mathcal{O}_K$.}, we have for all $\alpha\in \mc{G}_f$,
$$-p^{-2f}\le\int_{\mathcal{G}_K\cap U_\alpha}dx-p^{-f}\le 0.$$
Summing over all $\alpha\in \mathcal{G}_f$ and applying \Cref{lem: Mobius inverse formula}, we have
$$-G_fp^{-2f}\le \sum_{\alpha\in\mathcal{G}_f}\int_{\mathcal{G}_K\cap U_\alpha}dx-\sum_{d\mid f}\mu(\frac{f}{d})p^{d-f}\le 0.$$
Then the inequality is clear by the trivial bound $G_f\le p^f$.
\end{proof}

\begin{cor}\label{cor: estimate of limit density}
Suppose $x\in\mathcal{O}_K^{\new}$ satisfies $\F_p[x\pmod{p}]=\F_{p^d}$, where $d\mid f$. When $x\notin\mathcal{G}_K$, we have
$\rho_K^{(\infty)}(x)\le\frac{p^{-d}+2p^{-2d}}{1-p^{-f}}||\Disc_{K/\Q_p}||$.
\end{cor}

\begin{proof}
Applying the estimate in \Cref{lem: estimate of orbital integral} to the expression of $\rho_K^{(\infty)}(x)$, we have
\begin{align}
\begin{split}
\rho_K^{(\infty)}(x)
&\le\frac{\#(\Lambda\backslash\Mod_{\Z_p[x]})}{1-p^{-f}}||\Delta_\sigma(x)||^2\\
&=\frac{||\Disc_{K/\Q_p}||}{1-p^{-f}}\cdot\frac{\#(\Lambda\backslash\Mod_{\Z_p[x]})}{\#(\mathcal{O}_K/\Z_p[x])^2}\\
&\le\frac{p^{-d}+2p^{-2d}}{1-p^{-f}}||\Disc_{K/\Q_p}||.
\end{split}
\end{align}
Here the second line comes from the fact that $||\Delta_\sigma(x)||^2=\frac{||\Disc_{K/\Q_p}||}{\#(\mathcal{O}_K/\Z_p[x])^2}$, which is proved on page 10 of \cite{caruso2022zeroes}.
\end{proof}

\begin{prop}\label{prop: integral of density_Gf complement}
We have
$$\sum_{\alpha\in\mathcal{G}_f}\frac{\int_{U_\alpha\backslash\mathcal{G}_K}\rho_K^{(\infty)}(x)dx}{||\Disc_{K/\Q_p}||}\le \frac{p^{-2f}+2p^{-3f}}{1-p^{-f}}.$$
\end{prop}

\begin{proof}
When $x\in U_\alpha$ for some $\alpha\in\mathcal{G}_f$, one always has $\F_p[ x \pmod{p}]=\F_{p^f}$. Applying the upper bound of the density function in \Cref{cor: estimate of limit density}, we have
\begin{align}
\begin{split}
\sum_{\alpha\in\mathcal{G}_f}\frac{\int_{U_\alpha\backslash\mathcal{G}_K}\rho_K^{(\infty)}(x)dx}{||\Disc_{K/\Q_p}||}&\le\frac{p^{-f}+2p^{-2f}}{1-p^{-f}}\cdot\sum_{\alpha\in\mathcal{G}_f}\int_{U_\alpha\backslash\mathcal{G}_K}dx\\
&\le\frac{p^{-f}+2p^{-2f}}{1-p^{-f}}\cdot\sum_{\alpha\in\mathcal{G}_f}p^{-2f}\\
&\le\frac{p^{-2f}+2p^{-3f}}{1-p^{-f}}.
\end{split}
\end{align}
Here the second line comes from \cite[Lemma 4.4]{caruso2022zeroes}, and the third line holds because $G_f\le p^f$.
\end{proof}

\begin{prop}\label{prop: integral of density_Gd}
Let $1\le d<f$ be a divisor of $f$. Then we have
$$\sum_{\alpha\in\mathcal{G}_d}\frac{\int_{U_\alpha}\rho_K^{(\infty)}(x)dx}{||\Disc_{K/\Q_p}||}\le\frac{1+2p^{-d}}{1-p^{-f}}p^{-f}.$$
\end{prop}

\begin{proof}
Applying the upper bound of the density function in \Cref{cor: estimate of limit density}, we have
\begin{align}
\begin{split}
\sum_{\alpha\in\mathcal{G}_d}\frac{\int_{U_\alpha}\rho_K^{(\infty)}(x)dx}{||\Disc_{K/\Q_p}||}&\le\frac{p^{-d}+2p^{-2d}}{1-p^{-f}}\cdot\sum_{\alpha\in\mathcal{G}_d}\int_{U_\alpha}dx\\
&=\frac{p^{-d}+2p^{-2d}}{1-p^{-f}}\cdot\sum_{\alpha\in\mathcal{G}_d}p^{-f}\\
&\le\frac{1+2p^{-d}}{1-p^{-f}}p^{-f}.
\end{split}
\end{align}
Here the third line holds because $G_d\le p^d$.
\end{proof}

Combining the estimate from \Cref{prop: integral of density_generators}, \Cref{prop: integral of density_Gf complement} and \Cref{prop: integral of density_Gd}, we can turn back to our proof of \Cref{thm: expectation of higher degree extension}.

\begin{proof}[Proof of \Cref{thm: expectation of higher degree extension}]
The lower bound is a direct corollary of \Cref{prop: integral of density_generators}. For the upper bound, we have 
\begin{align}\label{eq: upper bound of three parts}
\begin{split}
\frac{\lim_{n\rightarrow\infty}\E[Z_{\mathcal{O}_K^{\new},n}]}{||\Disc_{K/\Q_p}||}&=\frac{\int_{\mathcal{G}_K}\rho_K^{(\infty)}(x)dx}{||\Disc_{K/\Q_p}||}+\sum_{\alpha\in\mathcal{G}_f}\frac{\int_{U_\alpha\backslash\mathcal{G}_K}\rho_K^{(\infty)}(x)dx}{||\Disc_{K/\Q_p}||}+\sum_{\substack{1\le d<f \\ d|f}}\sum_{\alpha\in\mathcal{G}_d}\frac{\int_{U_\alpha}\rho_K^{(\infty)}(x)dx}{||\Disc_{K/\Q_p}||}\\
&\le\sum_{d\mid f}\mu(\frac{f}{d})p^{d-f}+\frac{p^{-2f}+2p^{-3f}}{1-p^{-f}}+\sum_{\substack{1\le d<f \\ d|f}}\frac{1+2p^{-d}}{1-p^{-f}}p^{-f},\\
\end{split}
\end{align}
which is the combination of the three upper bounds in \Cref{prop: integral of density_generators}, \Cref{prop: integral of density_Gf complement} and \Cref{prop: integral of density_Gd} in sequence. Notice that
\begin{align}\label{eq: remain of upper bound}
\begin{split}
\frac{p^{-2f}+2p^{-3f}}{1-p^{-f}}+\sum_{\substack{1\le d<f \\ d|f}}\frac{1+2p^{-d}}{1-p^{-f}}p^{-f}&=\frac{p^{-f}}{1-p^{-f}}\left(\sum_{\substack{1\le d<f \\ d|f}}1+\sum_{\substack{1\le d<f \\ d|f}}2p^{-d}+p^{-f}+2p^{-2f}\right)\\
&\le\frac{p^{-f}}{1-p^{-f}}\left(\tau(f)-1+\sum_{1\le d\le 2f}2p^{-d}\right)\\
&<\frac{1+\tau(f)}{1-p^{-f}}p^{-f}.
\end{split}
\end{align}
The estimates in \eqref{eq: upper bound of three parts} and \eqref{eq: remain of upper bound} together give the proof. 
\end{proof}

\begin{proof}[Proof of \Cref{thm:roots_in_extensions_intro}]
For the case when $K/\Q_p$ is quadratic, the limit $\lim_{n\rightarrow\infty}\E[\Z_{\mathcal{O}_K^{\new},n}]$ is given in \Cref{thm:number_of_roots_in_quadratic_extension}. For the more general case that $K/\Q_p$ is any finite extension, the estimate is given in \Cref{thm: expectation of higher degree extension}.
\end{proof}

%% file: all_Zp.tex
\section{The probability that all eigenvalues are in $\Z_p$}\label{sec:all_zp}

Denote by $\mathscr{E}_n$ (resp. $\mathscr{E}_n^{\GL}$) the event that all eigenvalues of $A$ lie in $\Z_p$, where $A\in\Mat_n(\Z_p)$ (resp. $A\in\GL_n(\Z_p)$) is Haar-distributed. The main result of this section is the following theorem, which is just \Cref{thm: all eigenvalues in Z_p intro} together with the analogous result for $\GL_n(\Z_p)$.

\begin{thm}\label{thm: all eigenvalues in Z_p}
We have 
$$\log_p\mathbf{P}(\mathscr{E}_n)=-\frac{n^2}{2(p-1)}-\frac{1}{2}n\log_p n+O(n),$$
and 
$$\log_p\mathbf{P}(\mathscr{E}_n^{\GL})=-\frac{pn^2}{2(p-1)^2}-\frac{1}{2}n\log_p n+O(n).$$
Here, the implicit constants depend on $p$. We abuse notation and use $\mathbf{P}$ for both the additive Haar measure on $\Mat_n(\Z_p)$ and the Haar measure on $\GL_n(\Z_p)$.
\end{thm}

The main challenge to prove \Cref{thm: all eigenvalues in Z_p} is to deal with the integration of the correlation function, which by \Cref{thm:Zp_coulomb_intro} has the form of Coulomb gas. Luckily, Caruso proved in \cite[Theorem 5.8 (5)]{caruso2022zeroes} that the same Coulomb gas structure also appears when considering random polynomials. For the purpose of connecting these two random settings, in this section, we will always let $Y_n\in\Z_p[x]$ be the random polynomial
$$Y_n(x)=x^n+a_{n-1}x^{n-1}+\cdots+a_0,$$
where $a_0,\ldots,a_{n-1}$ are i.i.d. Haar-distributed in $\Z_p$. Based on the result of Caruso, the asymptotic of $\mathbf{P}(\mathscr{E}_n)$ should proceed from $\mathbf{P}(\mathscr{E}_n^{\poly})$, where $\mathscr{E}_n^{\poly}$ is the event that all the roots of $Y_n$ lie in $\Z_p$.  Since the asymptotic of $\mathbf{P}(\mathscr{E}_n^{\poly})$ is already worked out by Buhler-Goldstein-Moews-Rosenberg in \cite[Theorem 5.1]{buhler2006probability} (in their paper, they write $r_n$ for $\mathbf{P}(\mathscr{E}_n^{\poly})$), we are able to deduce the same form of asymptotic for $\mathbf{P}(\mathscr{E}_n)$. In the end, we will derive the asymptotic of $\mathbf{P}(\mathscr{E}_n^{\GL})$.

The following definition is the ``monic Haar polynomial" analog of the correlation function $\rho_{\Q_p^n}^{(n)}(x_1,\ldots,x_n)$ we studied in earlier sections.

\begin{defi}
Denote by $\rho_{\Q_p^n}^{(n),\poly}(x_1,\ldots,x_n)$ the correlation function of roots of $Y_n(x)$ over the étale algebra $\Q_p^n$.
\end{defi}

The following proposition is a slightly different version of \cite[Theorem 5.8 (5)]{caruso2022zeroes} (in the original setting, all coefficients of the random polynomial are Haar-distributed, including the highest-degree one).

\begin{prop}\label{prop: density function of poly}
For all distinct $x_1,\ldots,x_n\in\Z_p$, we have
$$\rho_{\Q_p^n}^{(n),\poly}(x_1,\ldots,x_n)=\prod_{1\le i<j\le n}||x_i-x_j||.$$
\end{prop}

The proof of \Cref{prop: density function of poly} is based on the following lemma, which is an analog of \Cref{thm: points on variety_main theorem}.
\begin{lemma}\label{lem: points on variety_poly}
Let $x_1,\ldots,x_n\in\Z_p$ be distinct. For all tuples $s=(s_1,\ldots,s_n)\in \N^n$ such that $$s_i>n\val\Disc(x_1,\ldots,x_n)=n\sum_{1\le j<k\le n}\val(x_j-x_k)$$ 
for all $1\le i\le n$, we have
$$p^{\sum_{i=1}^ns_i}\cdot\mathbf{P}(\val(Y_n(x_i))\ge s_i,\forall 1\le i\le n)=\prod_{1\le i<j\le n}||x_i-x_j||^{-1}.$$
\end{lemma}

\begin{proof}[Proof of \Cref{prop: density function of poly}, based on \Cref{lem: points on variety_poly}]
The proof is the same argument as was used to derive \Cref{thm: joint distribution} from \Cref{thm: points on variety_main theorem} in \Cref{sec: joint distribution}. We omit the details.
\end{proof}

\begin{proof}[Proof of \Cref{lem: points on variety_poly}]
Denote by $\Omega_{n-1}$ the set of polynomials of degree $\leq n-1$ with coefficients in $\Z_p$ (the polynomials need not be monic), equipped with the Haar probability measure. Then the polynomial 
$$H_{n-1}(x):=Y_n(x)-(x-x_1)(x-x_2)\cdots(x-x_n)$$
is distributed by the Haar measure on $\Omega_{n-1}$. Denote by
$$\Omega_{n-1}^s:=\{H_{n-1}\in\Omega_{n-1}:\val(H_{n-1}(x_i))\ge s_i,\forall 1\le i\le n\},$$
which is a $\Z_p$-submodule of $\Omega_{n-1}$.
It suffices to show
\begin{equation}\label{eq: measure of H in Omega}
\mathbf{P}(H_{n-1}\in\Omega_{n-1}^s)=p^{-\sum_{i=1}^ns_i}\prod_{1\le i<j\le n}||x_i-x_j||^{-1}.
\end{equation}
Consider the $\Z_p$-module homomorphism
\begin{align}
\begin{split}
\phi:\Omega_{n-1}&\rightarrow \Z_p^n\\
P(x)&\mapsto (P(x_1),\ldots,P(x_n)).
\end{split}
\end{align}
Then $\phi(\Omega_{n-1})=\Z_p[(x_1,\ldots,x_n)]$. Also, for all $1\le i\le n$, we have $s_i>\sum_{j\ne i}\val(x_j-x_i)$, thus the polynomials
$$\frac{p^{s_i}}{\prod_{j\ne i}(x_i-x_j)} \prod_{j \neq i} (x-x_j)\in\Z_p[x], \quad \quad 1 \leq i \leq n$$
form a $\Z_p$-basis of $\Omega_{n-1}^s$, and $\phi(\Omega_{n-1}^s)=\prod_{i=1}^n p^{s_i}\Z_p$. Therefore, we have
\begin{align}
\begin{split}
\text{LHS}\eqref{eq: measure of H in Omega}&=\frac{1}{\#(\Omega_{n-1}:\Omega_{n-1}^s)}\\
&=\frac{1}{\#(\phi(\Omega_{n-1}):\phi(\Omega_{n-1}^s))}\\
&=\frac{1}{\#(\Z_p[(x_1,\ldots,x_n)]:\prod_{i=1}^n p^{s_i}\Z_p)}\\
&=p^{-\sum_{i=1}^n s_i}\cdot\#(\Z_p^n:\Z_p[(x_1,\ldots,x_n)]).
\end{split}
\end{align}
Finally, notice that a $\Z_p$-basis of $\Z_p[(x_1,\ldots,x_n)]$ has the form $(x_1^i,\ldots,x_n^i)$, where $0\le i\le n-1$. Therefore, 
$$\#(\Z_p^n:\Z_p[(x_1,\ldots,x_n)])=||\det
\begin{pmatrix}1 & \cdots & x_1^{n-1}\\
\vdots & \ddots & \vdots\\
1 & \cdots & x_n^{n-1}\\
\end{pmatrix}||^{-1}=\prod_{1\le i<j\le n}||x_i-x_j||^{-1}.$$
This completes the proof.
\end{proof}

Based on the expression of correlation function in \Cref{prop: density function of poly}, let
$$Z_{\Z_p^n,n}^{\poly}:=\#\{(x_1,\ldots,x_n)\in\Z_p^n: x_1,\ldots,x_n\text{ are distinct roots of }Y_n\}$$
denote the number of $n$-tuples of distinct roots of $Y_n$ that lie in $\Z_p$. The expectation of  $Z_{\Z_p^n,n}^{\poly}$ can be written as the integral of the correlation function
$$\E[Z_{\Z_p^n,n}^{\poly}]=\int_{\Z_p^n}\rho_{\Q_p^n}^{(n),\poly}(x_1,\ldots,x_n)dx_1\cdots dx_n=\int_{\Z_p^n}\prod_{1\le i<j\le n}||x_i-x_j||dx_1\cdots dx_n.$$
Furthermore, it is clear that 
$$\mathbf{P}(\mathscr{E}_n^{\poly})=\frac{1}{n!}\E[Z_{\Z_p^n,n}^{\poly}]=\frac{1}{n!}\int_{\Z_p^n}\prod_{1\le i<j\le n}||x_i-x_j||dx_1\cdots dx_n.$$ 
This reduces our problem to studying asymptotics of $\mathbf{P}(\mathscr{E}_n^{\poly})$, for which we may cite the following result.

\begin{thm}[{\cite[Theorem 5.1]{buhler2006probability}}]\label{thm: asymptotic of log_poly}
There exists constants $C_1=C_1(p),C_2=C_2(p)$ such that for all $n\ge 1$, we have
$$C_1n\le\log_p\mathbf{P}(\mathscr{E}_n^{\poly})+\frac{n^2}{2(p-1)}+\frac{1}{2}n\log_p n\le C_2n.$$
\end{thm}

Now, we are ready to prove the $\Mat_n(\Z_p)$ case of \Cref{thm: all eigenvalues in Z_p}.

\begin{proof}[Proof of \Cref{thm: all eigenvalues in Z_p}, $\Mat_n(\Z_p)$ case]
Recall from \Cref{thm: joint distribution} that
\begin{align}
\begin{split}
\rho_{\Q_p^n}^{(n)}(x_1,\ldots,x_n)&=\frac{(1-p^{-1})\cdots(1-p^{-n})}{(1-p^{-1})^n}\prod_{1 \leq i < j \leq n}||x_i-x_j||\\
&=\frac{(1-p^{-1})\cdots(1-p^{-n})}{(1-p^{-1})^n}\rho_{\Q_p^n}^{(n),\poly}(x_1,\ldots,x_n).
\end{split}
\end{align}
Therefore, we have
\begin{align}\label{eq: probability of En and Enpoly}
\begin{split}
\mathbf{P}(\mathscr{E}_n)&=\frac{1}{n!}\E[Z_{\Z_p^n,n}]\\
&=\frac{(1-p^{-1})\cdots(1-p^{-n})}{(1-p^{-1})^n}\cdot\frac{1}{n!}\int_{\Z_p^n}\prod_{1\le i<j\le n}||x_i-x_j||dx_1\cdots dx_n\\
&=\frac{(1-p^{-1})\cdots(1-p^{-n})}{(1-p^{-1})^n}\cdot\frac{1}{n!}\E[Z_{\Z_p^n,n}^{\poly}]\\
&=\frac{(1-p^{-1})\cdots(1-p^{-n})}{(1-p^{-1})^n}\cdot\mathbf{P}(\mathscr{E}_n^{\poly}).
\end{split}
\end{align}
Combining the above relation with \Cref{thm: asymptotic of log_poly} and transferring the term $\log_p\frac{(1-p^{-1})\cdots(1-p^{-n})}{(1-p^{-1})^n}$ into $O(n)$, we obtain the proof.
\end{proof}

Let us turn to the asymptotic of $\mathbf{P}(\mathscr{E}_n^{\GL})$. The following lemma expresses $\mathbf{P}(\mathscr{E}_n^{\GL})$ in the form of $\mathbf{P}(\mathscr{E}_{b_i}^{\poly}),1\le i\le p-1$, making the asymptotic easier to handle.

\begin{lemma}\label{lem: express GL by poly}
We have
\begin{equation}\label{eq: GL and poly}
\mathbf{P}(\mathscr{E}_n^{\GL})=\frac{1}{(1-p^{-1})^n}\sum_{b_1,\ldots,b_{p-1}}\prod_{i=1}^{p-1}\left(\mathbf{P}(\mathscr{E}_{b_i}^{\poly})\cdot p^{-b_i(b_i+1)/2}\right).
\end{equation}
Here, the sum ranges all $b_1,\ldots,b_{p-1}\in\Z_{\ge 0}$ such that $b_1+\cdots+b_{p-1}=n$. 
\end{lemma}

\begin{proof}
Recall from \Cref{thm: joint distribution_GL} that
$$\rho_{\Q_p^n}^{(n),\GL}(x_1,\ldots,x_n)=\frac{1}{(1-p^{-1})^n}\prod_{1\le i<j\le n}||x_i-x_j||.$$
Therefore, we have
\begin{align}\label{eq: partition the integral into cosets}
\begin{split}
\mathbf{P}(\mathscr{E}_n^{\GL})&=\frac{1}{(1-p^{-1})^n}\cdot\frac{1}{n!}\int_{(\Z_p^\times)^n}\prod_{1\le i<j\le n}||x_i-x_j||dx_1\cdots dx_n\\
&=\frac{1}{(1-p^{-1})^n}\cdot\frac{1}{n!}\sum_{1\le k_1,\ldots,k_n\le p-1}\int_{(k_1+p\Z_p)\times\cdots\times(k_n+p\Z_p)}\prod_{1\le i<j\le n}||x_i-x_j||dx_1\cdots dx_n.
\end{split}
\end{align}
When $k_i \neq k_j$ for a given $i,j$, we have $||x_i-x_j|| = 1$, so the Vandermonde reduces to $p-1$ Vandermondes corresponding to these classes. Given $b_1,\ldots,b_{p-1}\in\Z_{\ge 0}$ such that $b_1+\cdots+b_{p-1}=n$, there are
$$\binom{n}{b_1,\ldots,b_{p-1}}=\frac{n!}{b_1!\cdots b_{p-1}!}$$
choices of $k_1,\ldots,k_n\in\{1,2,\ldots,p-1\}$ such that for all $1\le i\le p-1$, we have $\#\{1\le j\le n:k_j=i\}=b_i$. Therefore, 
\begin{align}
\begin{split}
\text{RHS}\eqref{eq: partition the integral into cosets}&=\frac{1}{(1-p^{-1})^n}\sum_{b_1,\ldots,b_{p-1}}\prod_{i=1}^{p-1}\frac{1}{b_i!}\int_{(p\Z_p)^{b_i}}\prod_{1\le j< k\le b_i}||x_{j,i}-x_{k,i}||dx_{1,i}\cdots dx_{b_i,i}\\
&=\frac{1}{(1-p^{-1})^n}\sum_{b_1,\ldots,b_{p-1}}\prod_{i=1}^{p-1}\frac{p^{-b_i(b_i+1)/2}}{b_i!}\int_{\Z_p^{b_i}}\prod_{1\le j< k\le b_i}||x_{j,i}-x_{k,i}||dx_{1,i}\cdots dx_{b_i,i}\\
&=\frac{1}{(1-p^{-1})^n}\sum_{b_1,\ldots,b_{p-1}}\prod_{i=1}^{p-1}\left(\mathbf{P}(\mathscr{E}_{b_i}^{\poly})\cdot p^{-b_i(b_i+1)/2}\right).
\end{split}
\end{align}
Here as before, the sum ranges all $b_1,\ldots,b_{p-1}\in\Z_{\ge 0}$ such that $b_1+\cdots+b_{p-1}=n$. This completes the proof.
\end{proof}

The following lemma provides an estimate for the summation terms in the right-hand side of \eqref{eq: GL and poly}.

\begin{lemma}\label{lem: estimate for sum poly}
Let $C_2$ be the same as in \Cref{thm: asymptotic of log_poly}. Then for all $b_1,\ldots,b_{p-1}\in\Z_{\ge 0}$ with $b_1+\cdots+b_{p-1}=n$, we have
$$\sum_{i=1}^{p-1}\left(\log_p\mathbf{P}(\mathscr{E}_{b_i}^{\poly})-\frac{b_i(b_i+1)}{2}\right)\le-\frac{pn^2}{2(p-1)^2}-\frac{1}{2}n\log_p n+C_2n.$$
\end{lemma}

\begin{proof}
By applying \Cref{thm: asymptotic of log_poly} to each term in the sum, we have
\begin{equation}
    \label{eq:bound_all_b_i}
    \sum_{i=1}^{p-1}\left(\log_p\mathbf{P}(\mathscr{E}_{b_i}^{\poly})-\frac{b_i(b_i+1)}{2}\right)\le\sum_{i=1}^{p-1}\left(-\frac{pb_i^2}{2(p-1)}-\frac{1}{2}b_i\log_p b_i+(C_2-\frac{1}{2})b_i\right).
\end{equation}
Here we slightly extend our convention: when $b_i=0$, we set $\mathbf{P}(\mathscr{E}_0^{\poly})=1$ so that $\log_p \mathbf{P}(\mathscr{E}_0^{\poly})=0$.
Moreover, by continuous extension we naturally interpret $b_i\log_p b_i=0$ at $b_i=0$. With these conventions, the conclusion of \Cref{thm: asymptotic of log_poly}
remains valid even when some $b_i$ equals zero. Now, observe that the function
\begin{equation}
    f(b) = -\frac{pb^2}{2(p-1)}-\frac{1}{2}b\log_p b+(C_2-\frac{1}{2})b
\end{equation}
is concave. Hence by Jensen's inequality,
\begin{align}
\begin{split}
\text{RHS\eqref{eq:bound_all_b_i}}& = \sum_{i=1}^{p-1} f(b_i) \\ 
& \le (p-1) f\left(\frac{\sum_{i=1}^{p-1}b_i}{p-1}\right) \\ 
&=(p-1)\cdot\left(-\frac{p\frac{n^2}{(p-1)^2}}{2(p-1)}-\frac{n}{2(p-1)}\log_p \frac{n}{p-1}+(C_2-\frac{1}{2})\frac{n}{p-1}\right)\\
&=-\frac{pn^2}{2(p-1)^2}-\frac{1}{2}n(\log_pn-\log_p(p-1))+(C_2-\frac{1}{2})n\\
&\le -\frac{pn^2}{2(p-1)^2}-\frac{1}{2}n\log_p n+C_2n.
\end{split}
\end{align}
\end{proof}

Now, we are ready for the $\GL_n(\Z_p)$ case of \Cref{thm: all eigenvalues in Z_p}.

\begin{proof}[Proof of \Cref{thm: all eigenvalues in Z_p}, $\GL_n(\Z_p)$ case]
Let $C_1,C_2$ be the same as in \Cref{thm: asymptotic of log_poly}. Our goal is to provide elementary estimates for the right-hand side of \eqref{eq: GL and poly}. On the one hand, the case 
$$b_1=\ldots=b_{n-(p-1)\lfloor \frac{n}{p-1}\rfloor}=\lfloor \frac{n}{p-1}\rfloor+1,\quad b_{n+1-(p-1)\lfloor \frac{n}{p-1}\rfloor}=\ldots=b_{p-1}=\lfloor \frac{n}{p-1}\rfloor$$
already gives the lower bound 
\begin{align}
\begin{split}
\text{RHS}\eqref{eq: GL and poly}
&\ge \left(n-(p-1)\lfloor \frac{n}{p-1}\rfloor\right)\left(\log_p\mathbf{P}(\mathscr{E}_{\lfloor\frac{n}{p-1}\rfloor+1}^{\poly})-\frac{(\lfloor\frac{n}{p-1}\rfloor+1)(\lfloor\frac{n}{p-1}\rfloor+2)}{2}\right)\\
&+\left((p-1)(\lfloor \frac{n}{p-1}\rfloor+1)- n\right)\left(\log_p\mathbf{P}(\mathscr{E}_{\lfloor \frac{n}{p-1}\rfloor}^{\poly})-\frac{\lfloor \frac{n}{p-1}\rfloor(\lfloor \frac{n}{p-1}\rfloor+1)}{2}\right)\\
&-n\log_p(1-p^{-1})\\
&\ge (p-1)\Biggl(-\frac{(\frac{n}{p-1}+1)^2}{2(p-1)}-\frac{1}{2}\left(\frac{n}{p-1}+1\right)\log_p\left(\frac{n}{p-1}+1\right)+C_1\left(\frac{n}{p-1}+1\right)\\
&-\frac{(\frac{n}{p-1}+1)(\frac{n}{p-1}+2)}{2}\Biggr)-n\log_p(1-p^{-1})\\
&\ge-\frac{pn^2}{2(p-1)^2}-\frac{1}{2}n\log_p n+(C_1-100p)n.
\end{split}
\end{align}
Here, we get the fourth line by replacing every $\lfloor \frac{n}{p-1}\rfloor$ and $\lfloor \frac{n}{p-1}\rfloor+1$ with $\frac{n}{p-1}+1$. On the other hand, there are only $\binom{n+p-2}{p-2}$ choices of $b_1,\ldots,b_{p-1}\in\Z_{\ge 0}$ with $b_1+\cdots+b_{p-1}=n$. Applying \Cref{lem: estimate for sum poly}, we have
\begin{align}
\begin{split}
\log_p\mathbf{P}(\mathscr{E}_n^{\GL})&\le-n\log_p(1-p^{-1})+\log_p\binom{n+p-2}{p-2}-\frac{pn^2}{2(p-1)^2}-\frac{1}{2}n\log_p n+C_2n\\
&\le-\frac{pn^2}{2(p-1)^2}-\frac{1}{2}n\log_p n+(C_2+100p)n.
\end{split}
\end{align}
This completes the proof.
\end{proof}